\documentclass[10pt]{amsart}
\usepackage{amscd}
\usepackage{amssymb}
\usepackage{latexsym}
\usepackage{mathrsfs}
\usepackage{bbding,array,tabu}
\usepackage{graphicx}
\usepackage{mathscinet}
\usepackage{mathdots}
\usepackage{fullpage}
\usepackage{color}
\usepackage[dvipsnames]{xcolor}
\usepackage{mathtools}
\newsavebox{\mybox}
\usepackage[T1]{fontenc}
\usepackage{lscape}
\usepackage{url}

\newcommand{\tabincell}[2]{\begin{tabular}{@{}#1@{}}#2\end{tabular}}
\newcommand{\nn}{\nonumber}

\newcommand{\fr}{ {\mathcal F}_R(\alpha)}
\usepackage{array}
\newcolumntype{L}[1]{>{\raggedright\let\newline\\\arraybackslash\hspace{0pt}}m{#1}}
\newcolumntype{C}[1]{>{\centering\let\newline\\\arraybackslash\hspace{0pt}}m{#1}}
\newcolumntype{R}[1]{>{\raggedleft\let\newline\\\arraybackslash\hspace{0pt}}m{#1}}

\makeatletter
\providecommand*{\boxast}{%
  \mathbin{% as \boxplus and \boxtimes
    \mathpalette\@boxit{*}%
  }%
}
\newcommand*{\@boxit}[2]{%
  \sbox0{$\m@th#1\Box$}%
  \ifx#1\displaystyle \ht0=\dimexpr\ht0+.05ex\relax \fi
  \ifx#1\textstyle \ht0=\dimexpr\ht0+.05ex\relax \fi
  \ifx#1\scriptstyle \ht0=\dimexpr\ht0+.05ex\relax \fi
  \ifx#1\scriptscriptstyle \ht0=\dimexpr\ht0+.05ex\relax \fi
  \sbox2{$#1\vcenter{}$}% \ht2 is positionn of math axis
  \rlap{%
    \hbox to \wd0{%
      \hfill
      \raisebox{%
        \dimexpr.25\dimexpr\ht0+\dp0\relax-\ht1\relax
      }{$\m@th#1#2$}%
      \hfill
    }%
  }%
  \Box
}
\makeatother

\newtheorem{thm}[equation]{Theorem}
\newtheorem{lem}[equation]{Lemma}
\newtheorem{dfn}[equation]{Definition}
\newtheorem{prop}[equation]{Proposition}

\newtheorem{cnjture}[equation]{Conjecture}
\theoremstyle{remark}

\newtheorem{rmrk}[equation]{\bf Remark}

\numberwithin{equation}{subsection}
\numberwithin{table}{subsection}

\newcommand{\gam}[1]{\Gamma(#1)}
\newcommand{\A}{\mathbb A}
\newcommand{\Z}{\mathbb Z}
\newcommand{\R}{\mathbb R}
\newcommand{\C}{\mathbb C}
\newcommand{\Q}{\mathbb Q}

\newcommand{\bs}{\backslash}

\newcommand{\abs}[1]{\lvert #1 \rvert}

\newcommand{\Whit}[1]{W_{#1}}
\newcommand{\Mellin}[1]{\widetilde{W}_{#1}}
\newcommand{\gamm}[1]{\Gamma({\textstyle\frac{#1}{2}})}
\newcommand{\gammdelt}[2]{\Gamma({\textstyle\frac{#1}{2}}-#2)}

\DeclareMathOperator*{\Res}{Res}
\DeclareMathOperator{\re}{Re}
\DeclareMathOperator{\im}{Im}

\DeclareMathOperator{\GL}{GL}
\DeclareMathOperator{\SL}{SL}
\DeclareMathOperator{\Ad}{Ad}
\DeclareMathOperator{\diag}{diag}

\renewcommand{\(}{\left(}
\renewcommand{\)}{\right)}
\def\rddots{\displaystyle\cdot^{\displaystyle\cdot^{\displaystyle\cdot}}}

\begin{document}

\title{An Orthogonality Relation for $\GL(4, \mathbb R)$}
\author{Dorian Goldfeld \and Eric Stade \and Michael Woodbury}
%\authorrunning{Goldfeld, Stade, Woodbury}

%\institute{
%Dorian Goldfeld \at Department of Mathemtics\\ Columbia University \\ New York, NY 10027\\ USA \\\email{goldfeld@columbia.edu}
%\and
%Eric Stade \at Department of Mathematics \\ University of Colorado \\ Bolder, Colorado 80309\\ USA
%\\\email{stade@colorado.edu}
%\and
%Michael Woodbury \at Department of Mathematics \\ Columbia University \\ 
%New York, NY 10027\\ USA
%\\\email{woodbury@math.columbia.edu}
%}

%\date{Received: date / Accepted: date}

\begin{abstract}
Orthogonality is a fundamental theme in representation theory and Fourier 
analysis. An orthogonality relation for characters of finite abelian groups (now recognized as an orthogonality relation on $\GL(1)$) was used by Dirichlet to prove infinitely many primes in arithmetic progressions. Orthogonality relations for $\GL(2)$ and $\GL(3)$ have been worked on by many researchers  with a broad range of applications to number theory. We present here, for the first time, very explicit orthogonality relations for 
the real group $\GL(4, \mathbb R)$ with a power savings error term. The proof requires novel techniques in the computation of the geometric side of the Kuznetsov trace formula.
\end{abstract}

\keywords{Kuznetsov trace formula \and automorphic forms \and  spectral decomposition \and Poincar\'e series \and  Whittaker functions \and Kloosterman sums \and Eisenstein series.}

\thanks{Dorian Goldfeld is partially supported by Simons Collaboration Grant Number 567168.}
\maketitle

\tableofcontents

\section{\bf Introduction and Main Theorem}

Let $q > 1$ be an integer and let $\chi\hskip-2pt:\hskip-2pt (\Z/q\Z)^\times\to\C^\times$ be a Dirichlet character (mod $q$).  The classical orthogonality relation for Dirichlet characters states that for integers $m, n$ coprime to $q$,
       $$\frac{1}{\phi(q)}\sum_{\chi \hskip-6pt\pmod{q}} \chi(m) \overline{\chi(n)} = \begin{cases} 1 & \text{if} \;  m\equiv n\hskip-6pt\pmod{q},\\
     0 &\text{otherwise.}\end{cases}$$
 This orthogonality relation is the basis for Dirichlet's proof that there are infinitely many primes $p\equiv a\hskip-3pt\pmod{q}$ if $(a,q)=1.$    It has played an essential role in the modern development of analytic number theory.
 
 When they are lifted to the adele ring $\A$ over $\Q$, Dirichlet characters can be realized as automorphic representations of $\GL(1)$ (see chapter 2 in \cite{GH2011}). It is then very natural to try to generalize the above orthogonality relation to representations of higher rank reductive groups. When trying to do this, however, there is an immediate obstacle. In the case of $\GL(1)$, there are only finitely many characters $\hskip-6pt\pmod{q}$ for any fixed $q > 1.$ In higher rank, on the other hand, there will be infinitely many automorphic representations. It then becomes necessary to introduce a test function with rapid decay and define the orthogonality relation as an absolutely convergent  integral over the automorphic representations.
   
  The first successful attempt at obtaining an orthogonality relation for 
$\GL(2)$ was made by R. Bruggeman in 1978 (see \cite{Bruggeman1978}) who considered the orthonormal basis $\{\phi_j\}_{j=1,2,\ldots}$ of Maass cusp forms
    for $\SL(2,\Z)$ where
    $$\phi_j(z) = \sum_{n\ne 0} a_j(n) \sqrt{2\pi y}\; K_{it_j}(2\pi|n| 
y) \cdot e^{2\pi inx}, \qquad (z=x+iy \in \text{upper-half plane}),$$
 and $K_{it_j}$ is the modified $K$-Bessel function of the second kind while $a_j(n) \in \C$ are the Fourier coefficients of $\phi_j$.
 The Maass cusp form $\phi_j$ has Laplace eigenvalue $\lambda_j =1/4+t_j^2$  and is a Hecke eigenform.
  Each such Maass cusp form is associated to a unique irreducible cuspidal automorphic representation of $\GL(2)$. Then Bruggeman proved the following orthogonality relation for non-zero integers $m,n$:
$$
 \lim_{T\to \,\infty} \frac{4\pi^2 }{T}  \,\sum_{j=1}^\infty \frac{ a_j(m) \; \overline{a_j(n)} }{\cosh(\pi t_j)} \cdot e^{-\lambda_j/T} \; = \begin{cases} 1 & \text{if}\; m = n,\\
 0 & \text{if} \; m \ne n.\end{cases}
$$ 
  Other versions of $\GL(2)$ type orthogonality relations  were later obtained by P. Sarnak  \cite{Sarnak1984}, and, for the case of holomorphic Hecke modular forms, by Conrey-Duke-Farmer \cite{CDF1997} and J.P. Serre \cite{Serre1997}.
  
 An orthogonality relation for Maass cusp forms on $\GL(3, \mathbb R)$ was first proved independently  by Goldfeld--Kontorovich \cite{GK2013}  and 
  Blomer \cite{Blomer_2013} in 2013. 
 Further results on orthogonality relations for $\GL(3, \mathbb R)$ were obtained by Blomer-Buttcane-Raulf \cite{BBr_2014} and Guerreiro \cite{Guerreiro2015}. 
 In his 2013 thesis  (see \cite{Zhou2013}, \cite{Zhou2014}) Fan Zhou conjectured a very general orthogonality relation for $\GL(n)$ for $n \ge 2.$ 
We now describe  Zhou's conjecture.

 Fix $n\ge 2.$ A Maass cusp form for $\SL(n,\mathbb Z)$ is a smooth function $\phi: \GL(n, \mathbb R) \to \mathbb C$ which satisfies
 $\phi(gkr) = \phi(g)$
 for all $g\in \GL(n,\mathbb R),$ $k\in K=\text{O}(n,\mathbb R)$, and $r \in\mathbb R^\times.$  In addition $\phi$ is square integrable and is an eigenfunction of the Laplacian. If $\lambda$ denotes the Laplace eigenvalue of $\phi$  then $\lambda$ can be expressed in terms of Langlands parameters $\alpha = \left(\alpha_1, \ldots, \alpha_n\right) \in \mathbb C^n$ of $\phi$, where $\alpha_1 + \alpha_2+\cdots + \alpha_n = 0.$ The precise relation is given (see Section~6 in \cite{Miller_2002}) by
   $$\lambda =  \left(\frac{n^3 - n}{24} -\frac{\alpha_1^2+\alpha_2^2+\cdots+\alpha_n^2}{2}\right).$$
The Maass cusp form $\phi$ is said to be tempered at $\infty$ if the Langlands parameters $\alpha_1,\ldots, \alpha_n$ are all pure imaginary.
 
Let $\{\phi_j\}_{j=1,2,\ldots}$ denote an orthogonal basis of  Maass cusp forms for $\SL(n, \Z)$ with associated Langlands parameters $\alpha^{(j)} = \big(\alpha^{(j)}_1, \ldots, \alpha^{(j)}_n\big)$ and
$M^{th}$ Fourier coefficient $A_j(M)$ where $M = (m_1, m_2, \ldots, m_{n-1})$ with $m_1m_2\cdots m_{n-1}\ne 0.$.   We assume each Maass cusp form $\phi_j$ is normalized so that its first Fourier coefficient $A_j(1,1,\ldots,1) = 1.$
 Let
$$
\mathcal L_j := \underset {s=1}{ \text{Res}}\; L(s, \phi_j\times \overline \phi_j)
$$
 be the residue, at the edge of the critical strip, of the Rankin-Selberg 
L-function attached to $\phi_j\times \overline \phi_j$ which is  the value at $s = 1$ of the adjoint L-function $L(s, \mathrm{Ad} \;\phi_j).$  
 \vskip 5pt
 
For $T\to \infty$,  and Langlands parameters $\alpha = (\alpha_1, \alpha_2, \ldots, \alpha_n) \in \C^n$  of $\phi$,  let $h_{T}(\alpha)$ denote a good test function with exponential decay as $\sum_{k=1}^n \lvert \alpha_k \rvert^2 \to \infty$.  Here ``good'' means that $h_{T}$ is smooth, holomorphic in a region $-\eta <\re(\alpha_i) <\eta$ for some $\eta > 0$, 
invariant under permutation of the Langlands parameters $\alpha=(\alpha_1,\ldots,\alpha_n)$, real valued and positive, and is essentially supported on the Laplace eigenvalues of $\phi$ which are less than $T^2$.

 \begin{cnjture} {\bf (Orthogonality relation for $\GL(n, \mathbb R)$)}  \label{OrthoConj1}
  Let $\{\phi_j\}_{j=1,2,\ldots}$ denote an orthogonal basis of  Maass cusp forms for $\SL(n, \Z)$ as above.  Set  $M = (m_1, \ldots, m_{n-1})$ and  $M' = (m'_1, \ldots, m'_{n-1}) \in \Z^{n-1}_{+}.$ Let
   $h_{T}$ denote a good test function as above.  Then 
   \begin{equation*}
   \lim_{T\to\infty} \;\frac{\sum\limits_{j=1}^\infty A_j(M) \overline{A_j(M')} \; \frac{h_{T}\left(\alpha^{(j)}  \right)}{\mathcal L_j}}{  \sum\limits_{j=1}^\infty  \frac{h_{T}\left(\alpha^{(j)}  \right)}{\mathcal L_j}} \; = \; \begin{cases} 1 & \text{if $M = M',$}
   \\0 & \text{otherwise}.\end{cases}
   \end{equation*}
 \end{cnjture}

For applications it is important to determine the rate of convergence as $T\to\infty$ in the above asymptotic relation. With this in mind, we  reformulate Conjecture \ref{OrthoConj1} with an error term.\footnote{We adopt the standard convention that the constant implied by $\mathcal O_{M, M'}$ depends at most on $M$ and $M'$.} In this case, the orthogonality relation is expected to take  the form:
\begin{cnjture} \label{OrthConj2} For some constant $\theta < 1,$ we have
$$\sum\limits_{j=1}^\infty A_j(M) \overline{A_j(M')} \; \frac{h_{T}\left(\alpha^{(j)}  \right)}{\mathcal L_j} \; =  \;\delta_{M,M'}\sum\limits_{j=1}^\infty  \frac{h_{T}\left(\alpha^{(j)}  \right)}{\mathcal L_j} \; 
+ \; \mathcal O_{M,M'}\left( \sum\limits_{j=1}^\infty  \frac{h_{T}\left(\alpha^{(j)}  \right)}{\mathcal L_j}  \right)^\theta.$$
Here $\delta_{M,M'}$ is 1 or 0 depending on whether $M = M'$ or not.
\end{cnjture}

  In the above, since $\theta < 1$, the error term gives a power savings in the main term.
In the case $n=3$, this conjecture was proved in \cite{GK2013}  with 
Langlands parameters $\alpha = (\alpha_1, \alpha_2, \alpha_3) \in \mathbb C^3$, and the following choice of test function:
  \[ h_{T,R}(\alpha)  := e^{  \frac{\alpha_1^2+\alpha_2^2+\alpha_3^2}{T^2}  }\cdot
\frac{  \prod\limits_{1\le j \ne k \le 3}  \Gamma\left( \frac{2+R+\alpha_j-\alpha_k}{4}  \right)^2 }{\prod\limits_{1\le j \ne k\le 3}   \Gamma\left(\frac{1+\alpha_j-\alpha_k}{2}  \right)   }.  %\qquad\quad \big(R\ge 10 \; \text{\rm fixed}\big).
  \]
More precisely, it was shown in \cite{GK2013} that
$\sum\limits_{j=1}^\infty  \frac{h_{T,R}\left(\alpha^j  \right)}{\mathcal L_j} \sim c \,T^{5+3R}$ and $\theta =  \frac{3+3R +\varepsilon}{5+3R},$
for some constant $c > 0$, 
 and any fixed $\varepsilon > 0$ as $T\to\infty.$ Similar results were independently obtained by Blomer \cite{Blomer_2013} and improved later in  \cite{BBr_2014}, and more recently in \cite{BZ2016} where an interesting technique is developed to remove the arithmetic weight $\mathcal L_j$. 

Conjecture \ref{OrthConj2} has many important applications to low lying zeros, Katz--Sarnak conjectures on symmetry types of families of automorphic L-functions, Sato--Tate conjectures, etc. Such applications, for the special case of  $\GL(3, \mathbb R)$,
 are a major main theme in \cite{Blomer_2013}, \cite{BBr_2014}, \cite{BZ2016}, \cite{GK2013}, 
\cite{Guerreiro2015}, \cite{Zhou2013}, \cite{Zhou2014}.  See also \cite{ST2016} where results for these problems are obtained for general families 
of cohomological automorphic representations of reductive groups over number fields.  In this paper we focus only on the orthogonality conjecture  
as the techniques to obtain the above type applications from Conjecture \ref{OrthConj2} are very well established.

Shin--Templier \cite{ST2016} obtain their results by an application of the Arthur-Selberg trace formula, with polynomial dependence on the Hecke eigenvalue and a power saving on the error term. Matz-Templier \cite{MT2015} establish the analogous results for the family $\{\phi_j\}$ of Maass cusp forms for $\SL(n,\Z)$, and this is strengthened and generalized in Finis-Matz \cite{FM2019}. A variant of Theorem \ref{MainTheorem} is obtained in \cite{MT2015}, \cite{FM2019}, without the arithmetic weight $\mathcal L_j$, without the polynomial weight of size $T^{8R},$  and with different test functions which are indicator functions of $\alpha^{(j)}\in T\Omega$, and where the error term would be $O(|\ell m|^\frac12\cdot T^8)$. For comparison, note that in Theorem \ref{MainTheorem} (if the polynomial weights are removed) we obtain an error term of the form $\mathcal O\left(|\ell m|^{\frac25+\varepsilon}\cdot T^{6+\varepsilon}  \right)$.  Also note that our main term is of a different form than that of \cite{MT2015}, \cite{FM2019}, in that ours entails some high-order asymptotics (the terms $ \mathfrak c_2 T^{\,8+8R}$ and $ \mathfrak c_3 T^{\,7+8R}$).

Shortly after our paper first appeared, Subhajit Jana \cite{Jana2020} obtained a proof of the conjecture for compactly supported functions (on the 
geometric side) for  automorphic forms for $\text{PGL}_r(\mathbb Z)$ with 
$r\ge 2.$ Although a power savings error term is not given in Jana's paper the author has informed us that the method can give a power savings error term which is very far from optimal (even for GL(2)).
\vskip 8pt 

\subsection{\bf Main Theorem}
Let $\alpha = (\alpha_1, \alpha_2, \alpha_3, \alpha_4)\in \mathbb C^4$ and let $S_4$ denote the symmetric group on a set of size four.  The main 
result of this paper is a proof of Conjecture \ref{OrthConj2} for $\GL(4, 
\mathbb R)$ for the test function $h_{T,R}(\alpha)$ given by
\begin{align*} & e^{\frac{\alpha_1^2+\alpha_2^2+\alpha_3^2+\alpha_4^2}{T^2}}\, \prod\limits_{1\leq\, j \ne k\, \leq 4}
      \frac{\left( 
      \Gamma\left(\textstyle{\frac{2+R+\alpha_j - \alpha_k}{4}}\right)  \right)^2 
      }{
      \Gamma\left( \frac{1+\alpha_j-\alpha_k}{2}  \right)}\prod\limits_{\sigma\in S_4}  \Big(1+\alpha_{\sigma(1)}-\alpha_{\sigma(2)}-\alpha_{\sigma(3)}+\alpha_{\sigma(4)}\Big)^{\frac{R}{12}},
\end{align*}
where $T$ is a large positive number and $R$ (sufficiently large) is a fixed positive integer.
      
 \vskip 15pt     
 
 \begin{thm} {\bf (Main Theorem)}\label{MainTheorem} Let $\{\phi_j\}_{j=1,2,\ldots}$ denote an orthogonal basis of even Maass cusp forms for $\SL(4, \Z)$ (assumed to be tempered at $\infty$)
  with associated Langlands parameters 
 $$\alpha^{(j)} = \big(\alpha^{(j)}_1, \alpha^{(j)}_2, \alpha^{(j)}_3,\alpha^{(j)}_4\big)\;\in 
  \;\left(i \mathbb R\right)^4$$ and
$L^{th}$ Fourier coefficient $A_j(L)$  (as in (\ref{FourierExpansion})) where
$L = (\ell_1, \ell_2, \ell_3)\in \mathbb Z^3$. Let $\mathcal L_j = L(1, \Ad\, \phi_j)$. We assume each Maass cusp form $\phi_j$ is normalized so that its first Fourier coefficient $A_j(1,1,1) = 1.$  Let $\ell, m \in \mathbb Z$  with $\ell m\ne 0.$  Then,  for $T\to\infty$,
\begin{multline*}
\sum\limits_{j=1}^\infty  A_j(\ell,1,1) \, \overline{A_j(m,1,1)} \,\frac{h_{T,R}\left(\alpha^{(j)}   \right)}{\mathcal L_j} =    \delta_{\ell,m}\cdot\Big(\mathfrak c_1T^{\,9+8R} + \mathfrak c_2 T^{\,8+8R} + \mathfrak c_3 T^{\,7+8R}\Big) \\
+ \mathcal{O}_{\varepsilon,R}\bigg( |\ell m|^{{\frac25}+\varepsilon}  \cdot T^{\,6+8R+\varepsilon}  + \lvert \ell m\rvert^{\frac{7}{32}+\varepsilon}\cdot T^{5+8R+\varepsilon} + \lvert \ell m\rvert^{\frac{15}{2}}\cdot T^{4+8R+\varepsilon} \bigg),
\end{multline*}
where $\delta_{\ell,m}$ is the Kronecker symbol and  $\mathfrak c_1,\mathfrak c_2,\mathfrak c_3>0$ are absolute constants which depend at most on $R$.
Note that $h_{T,R}$ is of size $T^{8R}$ on the
relevant support. \end{thm}

\begin{rmrk}
The polynomial $\prod\limits_{\sigma\in S_4}  \Big(1+\alpha_{\sigma(1)}-\alpha_{\sigma(2)}-\alpha_{\sigma(3)}+\alpha_{\sigma(4)}\Big)^{\frac{R}{12}}$ is a new feature.  It has not appeared in previous versions of this method for $\GL(2,\R)$ and $\GL(3,\R)$, but we have included it because its inclusion improves the error terms.
\end{rmrk}
 
 \vskip 12pt
 \begin{rmrk}
 For $s\in\mathbb C$ with $ \re(s) > 5/2$, the L-function associated to  $\phi_j$ is given by
 \begin{align*}
 L(s, \phi_j)  = \sum_{m=1}^\infty A_j(m,1,1) m^{-s} = \prod_p \Bigg(1 - \frac{A_j(p,1,1)}{p^{s}} 
   + \frac{A_j(1,p,1)}{p^{2s}}
    - \frac{A_j(1,1,p)}{p^{3s}}
     +\frac{1}{p^{4s}}\Bigg)^{-1}.
     \end{align*}
This shows that Theorem \ref{MainTheorem} gives the orthogonality relation on $\GL(4, \mathbb R)$ for coefficients of cuspidal $L$-functions.  It is possible, using the Hecke relations, to obtain a more general version of Theorem \ref{MainTheorem} involving $A_L, A_M$ for arbitrary $L,M$ where $\prod\limits_{i=1}^3 \ell_i m_i\ne 0,$ but the formulas get quite complex and messy, so are omitted. 
 \end{rmrk}
 
 \vskip 10pt
 \begin{rmrk}
 For a tempered Maass cusp form with Langlands parameters $\alpha\in (i\mathbb R)^4$ note that $$h_{T,R}\left(\alpha\right) > 0$$ and is essentially supported on Laplace eigenvalues $\lambda < T^2.$ It is not necessary to assume all Maass cusp forms for $\SL(4,\mathbb Z)$ are tempered. A weaker version of Theorem \ref{MainTheorem}  can be proved that assumes that 
almost all Maass cusp forms (except for a set of zero density) are tempered.
 \end{rmrk}

\vskip 10pt
 
   The proof of Conjecture \ref{OrthConj2} for $\GL(4, \mathbb R)$ (with a strong power savings error term) has resisted all attempts up to now. Theorem \ref{MainTheorem} is the first orthogonality relation for $\GL(4, \mathbb R)$ obtained which has a strong power savings error term. 
    Many of the techniques used in the proof of the $\GL(3, \mathbb R)$ conjecture do not generalize in an obvious way and new difficulties arise for the first time. We now point out the obstacles that we faced in the last seven years of work on this paper with some indications of how we overcame them.
\vskip 10pt

    $\bullet$ {\it In the methods developed in \cite{GK2013} the Whittaker transform of a test function is estimated by first taking the Mellin transform of the Whittaker function and then taking the inverse Mellin transform to go back. This leads to multiple integrals involving ratios of Gamma functions which can be estimated by Stirling's asymptotic formula. When moving to $\GL(4, \mathbb R)$, however, the Mellin transform of the Whittaker function is much more complex and does not satisfy a simple recurrence relation as on $\GL(3, \mathbb R).$ The polynomials which appear in 
the recurrence formula in \cite{FG1993} are of large degree, and it did not seem possible to get good estimates for Mellin transforms of Whittaker 
functions via recurrence relations. 
} 
\vskip 8pt
$\bullet$ {\it Recent results (see Section 5.2) give precise control of the polynomials that appear in the recurrence formulae for Mellin transforms of shifted Whittaker functions allowing us to overcome the problem discussed in the previous bullet.}

\vskip 8pt
$\bullet$ {\it The classical Perron's formula allows one to obtain asymptotic formulae for the sum of coefficients of an L-function by computing a 
certain integral transform of that L-function and then  evaluating the integral transform by shifting the contour of integration. An important tool in the proof of Theorem \ref{MainTheorem} is a novel higher dimensional 
version of Perron's formula that gives asymptotic formulae for sums of terms arising in the cuspidal contribution to the trace formula.    In the case of $\GL(4, \mathbb R)$, the Perron type formula we develop involves a triple integral which requires shifting contours in 3 directions. It was necessary to generalize the method of Goldfeld-Kontorovich for finding the ``exponential zero set'' which gets repeatedly used for each shifted term. We also introduce a very precise bound for elementary integrals (see Appendix A) which turns out to be critical for accurately estimating the integrals over the shifted contours. 
}

\vskip 8pt
$\bullet$ {\it
Another difficulty is that the Langlands spectral decomposition is much more complex on $\GL(4, \mathbb R)$ with many more types of Langlands L-functions involving twists by Maass cusp forms of lower rank in the Levi components of the relevant parabolic subgroups. In order to obtain precise power savings error terms in the contribution of the continuous spectrum to the trace formula, it is necessary to have very explicit forms of the Fourier coefficients of the Eisenstein series. Although the Fourier coefficients are known in great generality (see, for example, Shahidi's book  \cite{Shahidi2010}) the archimedean  factors do not seem to have been worked out explicitly in the published literature.  In Section~3.2 we review 
\cite{GMW2019} where Borel Eisenstein series are used as a template to explicitly determine the non-constant Fourier coefficients of general Langlands Eisenstein series on $\GL(4, \mathbb R)$.}

%\pagebreak
%\newpage

\subsection{\bf Roadmap for the proof of the Main Theorem} 
       
       \vskip 5pt
       The proof of Theorem  \ref{MainTheorem}  is based on the Kuznetsov 
trace formula for $\GL(4, \mathbb R)$ which is worked out in Section~3.  The trace formula is the identity
       $\mathcal C = \mathcal M + \mathcal K - \mathcal E$
       where
       $$\mathcal C = \sum\limits_{j=1}^\infty  A_j(\ell,1,1) \, \overline{A_j(m,1,1)} \,\frac{h_{T,R}\left(\alpha^{(j)}   \right)}{\mathcal L_j}$$ is the cuspidal contribution. The main term $\mathcal M$ is computed in Proposition \ref{prop:mainterm} and is given by
       $$\mathcal M =  \delta_{L,M}\cdot\Big(\mathfrak c_1T^{\,9+8R} + \mathfrak c_2 T^{\,8+8R} + \mathfrak c_3 T^{\,7+8R} + \mathcal O\left( T^{\,6+8R}  \right)\Big).$$
The bound for the Kloosterman contribution $
\mathcal K$ is worked in Proposition \ref{prop:Iwbounds} (with $r=4$), while the bound for the continuous spectrum $\mathcal E$ is given in Theorem \ref{thm:EisensteinBound}. Combining these bounds with the main term $\mathcal M$ completes the proof. \qed

\section{\bf Whittaker functions, Maass cusp forms, and Poincar\'e series 
for $\SL(4,\mathbb Z)$}
  We review basic notation and the definitions of Whittaker functions, Maass cusp forms, and Poincar\'e series  following \cite{Goldfeld2006}.

\subsection{\bf Iwasawa Decomposition} Fix $n\ge 2$
and let $g\in \GL(n,\mathbb R).$ We have the Iwasawa decomposition
\begin{equation} \label{IwasawaDecomposition}
g = utkr
\end{equation}
where $u\in U_n(\mathbb R)$ and $k\in K = \text{O}(n,\mathbb R)$ and $r\in \mathbb R^\times$ and $t\in T$, the subgroup of diagonal matrices with positive entries. Then $t = t(g)$ can be uniquely chosen to take the form
\begin{equation}
t = \left(\begin{smallmatrix} \label{ToricElement} y_1y_2\cdots y_{n-1} 
& & &&\\
& \ddots &&&\\
&&y_1 y_2 &&\\
& & & y_1 &\\
& & & &1 \end{smallmatrix}\right)
\end{equation}
 for some $y = (y_1, y_2,\ldots, y_{n-1})$ with $y_i>0   \ (i=1,2,\ldots,n-1)$.

\subsection{\bf Spectral and Langlands parameters} \label{SLpars} 
In the context of $\SL(n,\R)$, we associate to a vector of complex numbers
 $s = (s_1,s_2,\ldots,s_{n-1}) \;\in\; \mathbb C^{n-1},$ 
the spectral parameters
 \[ v :=(v_1,v_2,\ldots,v_{n-1}), \qquad v_j := s_j - \frac1n, \quad j=1,\ldots,n-1, \]
and Langlands parameters
 \[ \alpha_i := \begin{cases} B_{n-1}(v) & \text{if}\; i=1,\\
B_{n-i}(v) - B_{n-i+1}(v) & \text{if}\; 1<i<n,\\
-B_1(v) & \text{if} \; i = n,    \end{cases} \]
where $${B_j(v) =\sum\limits_{i=1}^{n-1} b_{i,j} v_i\qquad\hbox{with}\qquad
b_{i,j} = \begin{cases} i j & \text{if $i + j \le n,$}\\
(n-i)(n-j) & \text{if $i + j \ge n.$}\end{cases}}$$Note that
\begin{equation}\label{eq:alphasum}
 \alpha_1 + \alpha_2 + \cdots + \alpha_n = 0.
\end{equation}

On the other hand, given $\alpha=(\alpha_1,\ldots,\alpha_n)\in \C^n$ satisfying \eqref{eq:alphasum}, it is straightforward to see that the Langlands parameters associated to 
\begin{equation}\label{eq:salpha}
 s_i = \frac{\alpha_i-\alpha_{i+1}+1}{n},
\end{equation} 
are given by $\alpha$.

To be completely explicit, in the special case of $\SL(4, \mathbb Z)$, the Langlands parameters $(\alpha_1,\alpha_2,\alpha_3,\alpha_4)$ associated 
to  $s = \frac14 + \left(v_1,v_2, v_3\right)$ are given by
 \[ \alpha_1 = 3v_1+2v_2+v_3, \quad \alpha_2 = -v_1+2v_2+v_3, \quad \alpha_3 = -v_1-2v_2+v_3, \quad \alpha_4 = -v_1-2v_2-3v_3; \]
 \[ v_1 = \frac{\alpha_1 - \alpha_2}{4}, \quad v_2 = \frac{\alpha_2 - 
\alpha_3}{4}, \quad v_3 = \frac{\alpha_3 - \alpha_4}{4}. \]

\subsection{\bf The $I_s$-function} Let $g \in \GL(n, \mathbb R)$ (with toric element given by (\ref{ToricElement})).  We define a power function in terms of either the spectral or Langlands parameters associated to $s= 
(s_1,s_2,\ldots,s_{n-1})\in \mathbb C^{n-1}$ via
\begin{equation*}%\label{eq:powerfunction1}
 I_s(g) = I_{s}(utkr) := \prod_{i=1}^{n-1} \prod_{j=1}^{n-1}
y_i^{b_{i,j} s_j},
\end{equation*}
or, equivalently,
\begin{equation*}%\label{eq:powerfunction2}
 I_s(g) = I_{s}(utkr) := \prod_{j=1}^{n-1}
y_j^{\frac{j(n-j)}{2}+\alpha_1+\cdots+\alpha_{j-1}}.
\end{equation*}

It is easy to see that if $s$ and $\alpha$ are associated to each other as in Section~\ref{SLpars}, then these two definitions are equivalent.  For example, when $n = 4$ and $s = (1/4+v_1,1/4+v_2,1/4+v_3) \in \mathbb C^3$, we have 
\begin{align*}
 I_s(g) & = y_1^{s_1+2s_2+3s_3} y_2^{2s_1+4s_2+2s_3} y_3^{3s_1+2s_2+s_3} \\ & = y_1^{\frac32 + v_1+2v_2+3v_3}y_2^{2+2v_1+4v_2+2v_3}y_3^{\frac32+3v_1+2v_2+v_3} = y_1^{\frac32 +\alpha_1+\alpha_2+\alpha_3}y_2^{2+\alpha_1+\alpha_2}y_3^{\frac32+\alpha_1}.
\end{align*}

\vskip 15pt
\subsection{\bf Additive character of $U_n(\mathbb R)$} Assume $n\ge 2.$ Fix $M = (m_1,m_2, \;\ldots,\; m_{n-1})\in \mathbb Z^{n-1}.$ Let $g\in \GL(n,\mathbb R)$ with Iwasawa decomposition $g = utkr,$ where
$$u \; = \; \left(\begin{smallmatrix} 1 & u_{1,2}& u_{1,3}& \cdots  & 
u_{1,n}\\
 & 1& u_{2,3} &\cdots &  u_{2,n}\\
& &\ddots & &  \vdots\\
& & &  1& u_{n-1,n}\\
& & & &  1\end{smallmatrix}\right).$$
Then associated to the vector $M$ we have an additive character $\psi_M:U_n(\mathbb R) \to\mathbb C$ defined by
\begin{equation} \label{CharacterofU4}
\psi_M(g) := \psi_M(u) := e^{2\pi i \big(m_1u_{1,2} \,+\, m_2 u_{2,3}\, + \;\,\cdots \;\, + \, m_{n-1}u_{n-1,n}\big)}.
\end{equation}

\subsection{\bf Jacquet's  Whittaker Function} Assume $n\geq 2$. Given Langlands parameters $\alpha = (\alpha_1,\alpha_2, \;\ldots, \;\alpha_n)$, let $s=s(\alpha)$ be defined as in Section~\ref{SLpars}.
\vskip 3pt
For $\re(v_i)> 0$  $(i=1,\ldots,n-1)$ and $w_{\mathrm{long}} = \left(\begin{smallmatrix} &&1\\
&\rddots&\\ 1&&
  \end{smallmatrix}   \right)$,  we define the completed    Whittaker function $W^{\pm}_{\alpha}: \GL(n,\mathbb R)\big/\left(
  \text{O}(n,\mathbb R)\cdot \mathbb R^\times\right) \to \mathbb C$   by the absolutely convergent integral
$$W^{\pm}_{\alpha}(g) := \prod_{1\leq j< k \leq n} \frac{\Gamma\big(\frac{1 + \alpha_j - \alpha_k}{2}\big)}{\pi^{\frac{1+\alpha_j - \alpha_k}{2}} }\cdot \int\limits_{U_4(\mathbb R)} I_{s}(w_{\mathrm{long}} ug) \,\overline{\psi_{1, \ldots, 1,\pm 1}(u)} \, du, $$
where $du$ is the Haar measure on $U_n(\mathbb R).$ 
The product of Gamma factors  is added so that   ${W^{\pm}_{\alpha}}$ is invariant under all permutations of the Langlands parameters $\alpha_1,\alpha_2, \ldots, \alpha_n.$ 

\begin{rmrk}
If $g$ is a diagonal matrix in $\GL(n,\mathbb R)$ then the value of $W^{\pm}_{\alpha}(g)$ is independent of sign, so we drop the $\pm$. We also drop the $\pm$ if the sign is $+1$.
\end{rmrk}

Let
 $\mathcal D^{n}$ denote the algebra of $\GL(n,\mathbb R)$-invariant differential operators on 
 $$
\mathfrak h^n := \GL(n,\mathbb R)\big/\left(\text{O}(n,\mathbb R)\cdot \mathbb R^\times\right).
 $$
 
It is well known that $I_{s}(g)$ is an eigenfunction of all $\delta \in \mathcal D^n$.  In particular, for the case of $\delta=\Delta$, the Laplacian, then 
  $\Delta I_{s} = \lambda_{\Delta}(\alpha )\cdot  I_{s},$
  where
  $$\lambda_{\Delta}(\alpha) =  \left(\frac{n^3 - n}{24} -\frac{\alpha_1^2+\alpha_2^2+\cdots+\alpha_n^2}{2}\right).$$

Define $\lambda_{\delta}(\alpha)\in\C$ by the eigenfunction equation $$\delta I_{s}(g) =  \lambda_{\delta}(\alpha)\cdot  I_{s}(g), \;\;\qquad \left(\delta\in \mathcal D^n, \; g\in \GL(n, \mathbb R)\right).$$ 
 Jacquet's Whittaker function for  $\GL(n, \mathbb R)$ is characterized (up to scalars) by the following properties:
 
 \begin{align} \label{Eigenfunction}
 & \bullet \; \delta W^{\pm}_{\alpha} (g) = \lambda_{\delta}(\alpha)\cdot W^{\pm}_{\alpha}(g), \qquad (\text{for all $\delta\in \mathcal D^n$, \, $g \in \GL(n,\R)$}),\\
 &\bullet \; W^{\pm}_{\alpha}(ug) = \psi_{1,\ldots,1,\pm1}(u)  \cdot W^{\pm}_{\alpha}(g), \qquad (\text{for all $u\in U_n(\R), \; g\in \GL(n, \R)$}),\nonumber\\
 &\bullet \; W^{\pm}_{\alpha}(\mathrm{diag}\big( y_1y_2\cdots y_{n-1} , \; \ldots \;  ,y_1, 1)\big) \;\text{has exponential decay as $y_i\to\infty$, $(1 \le i \le n-1)$},\nonumber\\
 &\bullet \; {W^{\pm}_{\alpha}}\;\text{has holomorphic continuation to all $\alpha\in\C^n$, for all $g\in \GL(n,\R)$,}\nonumber\\
 &\bullet \;  {W^{\pm}_{\alpha}}  =  {W^{\pm}_{\alpha'}}  \;\text{where 
$\alpha'$ is any permutation of $\alpha = (\alpha_1, \ldots, \alpha_n)$.}\nonumber
 \end{align}

\subsection{\bf Whittaker Transform
}\label{sec:WhittTransform} Assume $n\ge 2$. Let $v = \frac1n + \left(v_1,v_2, \ldots, v_{n-1}\right) \in \mathbb C^{n-1}$ with the associated Langlands parameters $\alpha = (\alpha_1,\alpha_2, \ldots, \alpha_n).$
   Set
$$y := (y_1, y_2,\ldots y_{n-1}), \qquad t(y) := \left(\begin{smallmatrix} y_1y_2\cdots y_{n-1} &&&&\\
 &\ddots&&&  \\
 & & y_1y_2 &&\\
&& &y_1&\\
&&&&1  \end{smallmatrix}\right).$$
 Let $f:\mathbb R_+^{n-1} \to \mathbb C$ be an integrable function. Then we define the Whittaker transform 
$f^\#: \mathbb R_+^n \to \mathbb C$ by
 \begin{equation} \label{WhittakerTransform}
 f^\#(\alpha) := \int\limits_{y_1=0}^\infty \cdots  \int\limits_{y_{n-1}=0}^\infty   f(y)\,W_\alpha\big(t(y)\big) \prod_{k=1}^{n-1} \frac{dy_k}{y_k^{k(n-k)+1}  },
 \end{equation}
 provided the above integral converges absolutely and uniformly on compact subsets of $\mathbb R_+^{n-1}$. Assume that $\alpha$ is tempered, i.e., 
$v_1, v_2, \ldots, v_{n-1}$ are all pure imaginary. The Whittaker transform was studied in \cite{GK2012} and the following explicit inverse Whittaker transform was obtained:  $$f(y) = \frac {1}{\pi^{n-1}} \int\limits_{\re(v_1)=0} \cdots \int\limits_{\re(v_{n-1})=0} f^\#(\alpha) W_{-\alpha}\big(t(y)\big) \frac{dv_1 dv_2\cdots  dv_{n-1}}{\prod\limits_{1\le k\ne\ell\le n}\Gamma\left( \frac{\alpha_k-\alpha_\ell}{2}  \right)},$$ provided the above integral converges absolutely and uniformly on compact subsets of $(i\mathbb R)^n$.

\subsection{\bf The inner product of two Whittaker functions} Assume $n\ge 2.$  Suppose that $\alpha = (\alpha_1, \ldots,\alpha_n)$ and 
$\beta = (\beta_1, \ldots,\beta_n)$ are Langlands parameters for which $\re(\alpha_j) = \re(\beta_k) = 0$ $(1\le j,k\le n)$.  Then 
\begin{equation} \label{StadesInnerProdFormula}
\int\limits_{y_1=0}^\infty \cdots  \int\limits_{y_{n-1}=0}^\infty  W_\alpha(y) \,\overline{W_\beta(y)}\; \prod_{j=1}^{n-1} y_j^{(n-j)s} \frac{dy_k}{y_k^{k(n-k)+1}} \; = \; \frac{\prod\limits_{j=1}^n\prod\limits_{k=1}^n\Gamma\big(\frac{s+\alpha_j-\beta_k}{2}  \big)}{2\pi^{s\frac{n(n-1)}{2}} \Gamma\left( \frac{ns}{2}  \right)},
\end{equation}
where the left side converges absolutely for Re$(s)$ sufficently large.
This is given in \cite{Stade2002}.

\subsection{\bf Fourier-Whittaker expansion of Maass cusp forms} 
Assume $n\ge 2.$ Fix Langlands parameters $\alpha = (\alpha_1, \alpha_2, \ldots,\alpha_n)\in \mathbb C^n$. Let 
$
\phi: \mathfrak h^n\to\mathbb C
$
be a Maass cusp form for $\SL(n, \mathbb Z)$ which satisfies
$\delta \phi(g) = \lambda_\alpha(\delta)\cdot \phi(g)$
for all $\delta\in \mathcal D^n$ and $g \in \GL(n,\R)$, as in (\ref{Eigenfunction}). Then for $g\in \GL(n, \mathbb R)$,  we have the following Fourier-Whittaker expansion:

\begin{equation} \label{FourierExpansion}
\phi(g) = \sum_{\gamma\in{U}_{n-1}(\mathbb Z)\backslash \SL_{n-1}(\mathbb Z)} \;\sum_{m_1=1}^\infty\cdots  \sum_{m_{n-2}=1}^\infty\;\sum_{m_{n-1}\ne0} \, \frac{A_\phi(M)}
{\prod\limits_{k=1}^{n-1}  |m_k|^{\frac{k(n-k)}{2}}} \; W^{\text{\rm sgn}(m_{n-1})}_{\alpha}\left(t(M) \bigg(\begin{matrix} \gamma &0\\0&1\end{matrix}\bigg)g\right), \end{equation}
where $M = (m_1, m_2, \; \ldots, \; m_{N-1})$
and  $A_\phi(M)$ is  the $M^{th}$ Fourier coefficient of $\phi$. This is proved in Section~9.1 of \cite{Goldfeld2015}.

\subsection{\bf First Fourier-Whittaker coefficient of  a Maass cusp form} For $n\ge 2$,  consider a  Maass cusp form $\phi$ for $\SL(n,\mathbb Z)$ with Fourier Whittaker expansion given by \ref{FourierExpansion}. Assume $\phi$ is a Hecke eigenform.   Let $A_\phi(1) := A_\phi(1,1,\ldots,1)$ denote the first Fourier-Whittaker coefficient of $\phi.$ Then we have
$$A_\phi(M) = A_\phi(1) \cdot\lambda_\phi(M)$$
where $\lambda_\phi(M)$ is the Hecke eigenvalue (see Section~9.3 in \cite{Goldfeld2015}), and $\lambda_\phi(1) = 1$.  

Recall also the definition of the adjoint L-function:  $L(s, \Ad\;\phi) := L(s, \phi\times\overline{\phi})/\zeta(s)$ where $L(s, \phi\times\overline{\phi})$ is the Rankin-Selberg convolution L-function as in \S12.1 of 
\cite{Goldfeld2015}.

\begin{prop} \label{PropFirstCoeff} Assume $n\ge 2.$ Let $\phi$ be a Maass cusp form for $\SL(n,\mathbb Z)$ with Langlands parameters $\alpha=(\alpha_1,\ldots,\alpha_n)$. Then the first coefficient $A_\phi(1)$ is given by 
$$|A_\phi(1)|^2 =  \frac{\mathfrak c_n\cdot\big\langle \phi, \,\phi\big\rangle}{L(1, \Ad \; \phi) \prod\limits_{1\le j\ne k\le n}\Gamma\big(\frac{1+\alpha_j-\alpha_k}{2}  \big)}$$
where $\mathfrak c_n \ne 0$ is a constant depending on $n$ only.
\end{prop}

\begin{proof}
We follow the Rankin-Selberg computations in \S12.1 of \cite{Goldfeld2015}.
\begin{align*}
\big\langle\phi, \,\phi\big\rangle & = \int\limits_{\SL(n,\mathbb Z)\backslash\mathfrak h^n} |\phi(g)|^2\, d^*g\\
& = \text{vol}\left(\SL(n,\mathbb Z)  \backslash\mathfrak h^n \right)\cdot \underset{s=1}{\text{Res} } \int\limits_{\SL(n,\mathbb Z)\backslash\mathfrak h^n} \phi(g) \overline{\phi(g)} E(g,s)\, d^*g
 \end{align*}
 where $E(g,s)$ is the maximal parabolic Eisenstein series. After unfolding and replacing $\phi$ with its Fourier-Whittaker expansion, we obtain
 \begin{align*}
 \int\limits_{\SL(n,\mathbb Z)\backslash\mathfrak h^n} \phi(g) \overline{\phi(g)} E(g,s)\, d^*g
 & = \frac{|A_\phi(1)|^2 L(s, \phi\times\overline{\phi})}{\zeta(ns)} \int\limits_{\mathbb R_+^{n-1}} W_\alpha(y) \overline{W_\alpha(y)}\, \prod_{j=1}^{n-1} \frac{dy_k}{y_k^{(k-s)(n-k)+1}}.
 \end{align*}
The proposition follows immediately from the formula \eqref{StadesInnerProdFormula}, since $L(1, \Ad\;\phi) = \underset{s=1}{ \text{Res}}\; L(s, \phi\times\overline{\phi})$.
\end{proof}

\subsection{\bf Vector or matrix notation depending on context}
\label{VectorMatrix} Given a vector $a = (a_1,a_2,a_3)$ in $\mathbb R^3$, we shall  define the toric element $t(a) := \mathrm{diag}(a_1a_2a_3, 
a_1 a_2 ,a_1,1).$ 
  
  Given a function $f:\mathbb R^3 \to \mathbb C$ we define
  $f(a) :=  f(a_1,a_2,a_3).$  
  On the other hand, if $f:\GL(4, \mathbb R) \to\mathbb C$ is a function defined on the group then we let
   $f(a) :=  f(t(a)),$ and more generally, for any $g_1, g_2\in \GL(4,\mathbb R)$ we define $f(g_1ag_2) :=
   f(g_1t(a)g_2).$  In other words, we may consider $a$ as a vector or a diagonal matrix depending on the context.
   
\vskip 15pt
\subsection{\bf Poincar\'e series for $\SL(4, \mathbb Z)$} Let $H:\mathfrak h^4  \to\mathbb C$ be a smooth test function satisfying $H(utkr) = H(t)$ (see \eqref{IwasawaDecomposition}). We assume that $H$ has sufficient decay properties so that the series defining the Poincar\'e series (given below) converges absolutely.  For $g\in \GL(4, \mathbb R)$, $M = (m_1,m_2, m_3)\in\mathbb Z_+^{3}$,  and $s = (s_1, s_2,s_3) \in \mathbb C^{3}$, with Re$(s_j)$ sufficiently large, the $\SL(4,\mathbb Z)$ Poincar\'e series is defined by
\begin{equation}\label{PoincareSeries}
P^M\left(g, s \right) := \sum_{\gamma\in U_4(\mathbb Z)\backslash \SL(4,\mathbb Z)} \psi_M(\gamma g)\, H\big(M \gamma g\big)\, I_s(\gamma g).
\end{equation}
{\bf Remark:} Following Section~\ref{VectorMatrix}, for $\psi_M$, we take 
$M = (m_1,m_2,m_3)$. On the other hand for $H(M\gamma g)$,  we take $M$ 
to be the diagonal matrix  $t(M)$.

\subsection{\bf Inner product of the Poincar\'e series with a Maass cusp form} Let $\phi$ be a Maass cusp form for $\SL(4,\mathbb Z)$ with Fourier 
expansion (\ref{FourierExpansion}). Let $P^M$ denote the Poincar\'e series (\ref{PoincareSeries}). The inner product is defined by
\begin{align*}
\big\langle P^M(*, s), \; \phi  \big\rangle\;  & := \int\limits_{\SL(4, 
\mathbb Z) \backslash\mathfrak h^4} P^M(g,s)\, \overline{\phi(g)} \; dg.
\end{align*}

It follows that
\begin{align*}
\big\langle P^M(*, s), \; \phi  \big\rangle\;& = 
 \int\limits_{U_4(\mathbb Z) \backslash\mathfrak h^4} \hskip-3pt\psi_M(x) 
H\big(My\big) I_s(y) \overline{\phi(xy)} \; dx_{1,2} dx_{1,3} dx_{1,4} dx_{2,3} dx_{2,4} dx_{3,4} \frac{dy_1 dy_2 dy_3}{y_1^4 y_2^5 y_3^4}\\
 & = \frac{\overline{A_\phi(M)}}
{m_1^{\frac32}  m_2^2  m_3^{\frac32}}\int\limits_{y_1=0}^\infty  \int\limits_{y_2=0}^\infty
 \int\limits_{y_3=0}^\infty  H\big(My\big)\, I_s(y)\cdot \overline{W_\alpha(My)}\; \frac{dy_1 dy_2 dy_3}{y_1^4 y_2^5 y_3^4}.
\end{align*}

We see that the above inner product picks out the $M^{th}$ Fourier coefficient of $\phi$ multiplied by a certain Whittaker transform of $H\big(My\big)\cdot I_s(y)$. Letting $s\to 0$, it follows from  (\ref{WhittakerTransform})
 that
 \begin{align} \label{InnerProductPoincare}
\lim_{s\to 0}\, \big\langle P^M(*, s), \; \phi  \big\rangle \; & = \;\frac{\overline{A_\phi(M)}}
{m_1^{\frac32} m_2^2 m_3^{\frac32}}\int\limits_{y_1=0}^\infty  \int\limits_{y_2=0}^\infty
 \int\limits_{y_3=0}^\infty  H\big(My\big)\cdot \overline{W_\alpha\big(My\big)}\; \frac{dy_1 dy_2 dy_3}{y_1^4 y_2^5 y_3^4}\\
 & = \; m_1^{\frac32} m_2^2 m_3^{\frac32} \cdot \overline{A_\phi(M)}
  \cdot
H^\#(\overline{\alpha}).\nonumber\end{align}

\subsection{\bf Fourier-Whittaker expansion of the Poincar\'e series}
\vskip 5pt
Let  $W_4 \cong S_4$ denote the Weyl group of $\GL(4, \mathbb R).$ For $w\in W_4$, we define
$$\Gamma_w := \Big(w^{-1}\cdot\hskip-8pt \phantom{U}^t U_4(\mathbb Z)\cdot w\Big)\cap U_4(\mathbb Z),$$
where 
${}^t U_4$ denotes the transpose of $U_4$.

We have the Bruhat decomposition
$$\GL(4, \mathbb R) = \bigcup_{w\in W_4} G_w, \qquad\quad \Big(G_w = U_4(\mathbb R) \cdot w \cdot T_4(\mathbb R) U_4(\mathbb R)\Big),$$
where $T_4(\mathbb R)$ is the subgroup of diagonal matrices in $\GL(4, \mathbb R).$

\begin{dfn}[Twisted Character] Let 
$$V_4 := \left\{v = \left(\begin{smallmatrix}
v_1& & &\\
& v_2 & &\\
& & v_3 &\\
& & & v_4
\end{smallmatrix}\right) \big | \;v_1,v_2,v_3,v_4\in\{\pm1\},  \; v_1v_2v_3v_4 = 1\right\}.$$ Let $M = (m_1,m_2,m_3) \in \mathbb Z^3,$ and consider $\psi_M$ an additive character of $U_4$. Then for $v \in V_4,$ we define the twisted character $\psi_M^v : U_4(\mathbb R) \to \mathbb C$ by 
$$\psi_M^v(g) := \psi_M\left(v^{-1} g v\right).$$
\end{dfn}

\begin{dfn}[Kloosterman Sum] \label{KloostSum} Fix $L=(\ell_1,\ell_2,\ell_3),\, M=(m_1,m_2,m_3)  \in\mathbb Z^3.$ Let $\psi_L, \psi_M$ be characters of $U_4(\mathbb R).$ Let $w \in W_4$ where $W_4$ is the Weyl group 
of $\GL(4).$ Let
$c = \left(\begin{smallmatrix} 1/c_3 & & &\\
& c_3/c_2 & &\\
& & c_2/c_1 &\\
& & & c_1  \end{smallmatrix}\right)$ with $c_i\in \mathbb N.$
 Then the Kloosterman sum is defined as
$$S_w(\psi_L, \psi_M, c) := \underset {\gamma = \beta_1cw\beta_2}{\sum\limits_{\gamma = U_4(\mathbb Z)\backslash \Gamma \cap G_w/\Gamma_w}} \psi_L(\beta_1) \,\psi_M(\beta_2),$$
with notation as in Definition~11.2.2 of \cite{Goldfeld2006}. The Kloosterman sum $S_w(\psi, \psi', c)$ is well defined and not equal to zero if and only if it satisfies the compatibility condition $\psi(cwuw^{-1}) = \psi'(u).$
\end{dfn}

 It follows from Theorem~11.5.4 of \cite{Goldfeld2006} that the $M^{th}$ Fourier coefficient of the Poincar\'e series $P^L\left(g, \,s \right)$ is 
given by
 
\begin{align} \label{FourierExp}
\int\limits_{U_4(\mathbb Z)\backslash U_4(\mathbb R)} P^L\left(ug, \,s \right) \cdot \overline{\psi_M(u)} \; d^*u \; = \; \sum_{w\in W_4}\sum_{v\in V_4}\sum_{c_1=1}^\infty \; \sum_{c_2=1}^\infty\; \sum_{c_3=1}^\infty  \frac{S_w(\psi_L,\psi_M^v,c) J_w(g; s,\psi_L,\psi_M^v,c)}{c_1^{4s_3} c_2^{4s_2} c_3^{4s_1}   },
  \end{align}
where 
$$J_w(g; s, \psi_L,\psi_M^v,c) = \int\limits_{U_w(\mathbb Z)\backslash {U}_w(\mathbb R)}\;
\int\limits_{\overline{U}_w(\mathbb R)} \psi_L(w u g)\, H\big(L c w u g\big)\, I_s(w u  g) \;\overline{\psi_M^v(u)} \; d^* u,$$
$${U}_w(\mathbb R) = \Big(w^{-1} \cdot  U_4(\mathbb R)\cdot w\Big) \cap 
U_4(\mathbb R), \qquad \overline{U}_w(\mathbb R) = \Big(w^{-1} \cdot \hskip-10pt\phantom{U}^t U_4(\mathbb R)\cdot w\Big) \cap U_4(\mathbb R),$$
and \hskip-8pt$\phantom{m}^t m$ denotes the transpose of a matrix $m$.

\section{\bf Kuznetsov Trace Formula for $\SL(4,\mathbb Z)$}

\subsection{\bf Choice of test function}

Let $\alpha = (\alpha_1,\alpha_2,\alpha_3,\alpha_4)\in(i\R)^4$ with $\alpha_1+\alpha_2+\alpha_3+\alpha_4 = 0$.  Let  $T > 1$ with $T\to\infty$ 
and $R\geq 14$ with $R$ fixed. We consider the test function
\begin{equation}\label{testfunctionsharp2}
p_{T,R}^\sharp(\alpha) \; := \;e^{\frac{\alpha_1^2+\alpha_2^2+\alpha_3^2+\alpha_4^2}{2T^2}}\,\mathcal F_R(\alpha)\,
      \prod_{1\leq\, j \ne k\, \leq 4}
      \Gamma\left(\textstyle{\frac{2+R+\alpha_j - \alpha_k}{4}}\right),
\end{equation}where
\begin{align}\label{Fr}&
 \mathcal F_R(\alpha)  := \left(\prod_{\sigma\in S_4}  \Big(1+\alpha_{\sigma(1)}-\alpha_{\sigma(2)}-\alpha_{\sigma(3)}+\alpha_{\sigma(4)}\Big) \right)^{\frac{R}{24}}\\ \nn
 & = \Big(\big(1+|\alpha_1+\alpha_2 -\alpha_3-\alpha_4|^2\big) 
 \big(1+|\alpha_1+\alpha_3 -\alpha_2-\alpha_4|^2\big)
 \big(1+|\alpha_1+\alpha_4 -\alpha_2-\alpha_3|^2\big)\Big)^{\frac{R}{6}}.
\end{align}

The function $p_{T,R}^\sharp(\alpha)$ defined in \eqref{testfunctionsharp2} generalizes the similar function defined in \cite{GK2013}. As before, the choice is motivated by the fact that we need $p_{T,R}^\sharp$ to be invariant under the Weyl group, and have meromorphic continuation in $\alpha\in\C^4$, while also requiring it to have enough exponential decay to kill the exponential growth of certain Gamma factors appearing in the denominator of the Kuznetsov trace formula.  The new feature is the introduction of the polynomial $\mathcal F_R(\alpha)$.

\vskip 10pt
By the inverse Lebedev-Whittaker transform (see \cite{GK2012}), we see that $p_{T,R}$ is given by
\begin{equation}\label{eq:testfunction2}
p_{T,R}(y) = p_{T,R}(y_1,y_2,y_3) = \frac{1}{\pi^{3}}
  \underset{\re(\alpha_j)=0}{\iiint} p^{\#}_{T,R}(\alpha) \;W_{\alpha}(y)
     \frac{ d\alpha_1\, d\alpha_2\, d\alpha_3
         }{\prod\limits_{1\leq\, j \ne k\, \leq 4} \Gamma\left(\frac{\alpha_j-\alpha_k}{2}\right) }.
\end{equation}

\subsection{\bf Setting up the  trace formula}

Set $L= (\ell_1,\ell_2,\ell_3), \; M = (m_1,m_2,m_3) \in \mathbb Z^3$ 
where we assume $\prod_{i=1}^3 \ell_im_i \ne 0.$   Consider the Poincar\'e series $P^L, P^M$, as defined in (\ref{PoincareSeries}) with the test 
function $H = p_{T,R}$. 
 
\begin{dfn}[Normalization factor $C_{L,M}$] \label{NormFactor} Let $\mathfrak c_4$ be as Proposition \ref{PropFirstCoeff}. Then we define
$C_{L,M} := \mathfrak c_4\cdot (\ell_1 m_1)^3  (\ell_2 m_2)^4  (\ell_3 m_3)^3.$
\end{dfn} 
 
 With the normalization factor $C_{L,M}$ defined above, the Kuznetsov trace formula is obtained by evaluating the inner product
$$C_{L,M}^{-1}\cdot \lim\limits_{s\to 0} \; \big\langle P^L(*, s), \; P^M(*, s)  \big\rangle\;  = \; C_{L,M}^{-1}\cdot \lim\limits_{s\to 0} \int\limits_{\SL(4, \mathbb Z) \backslash\mathfrak h^4} P^L(g,s)\, \overline{P^M(g,s)} \; dg$$
in two different ways. The first approach is to use spectral theory while 
the second uses geometry. For the spectral theory approach, we will need the following definition.
 
\begin{dfn}[Generic and non-generic automorphic forms]   Let $n\ge 2.$  An automorphic form for $SL(n,\mathbb Z)$ is said to be {\it generic} if there exists some nondegenerate character of $U_n(\mathbb R)$ for which the form has a nonzero Fourier coefficient.  Otherwise, the automorphic form is said to be {\it non-generic}. Note that a character $\psi_{m_1,m_2,\ldots,m_{n-1}}$ of $U_n(\mathbb R)$ is nondegenerate if all $m_1,m_2, \ldots,m_{n-1}$ are all nonzero.
\end{dfn}

In order to prove the Kuznetsov trace formula for $\GL(4,\R)$, we require 
the spectral expansion of the Poincar\'e series $P^M$ and $P^L$.  It will 
be shown (see the proof of Theorem \ref{ThmSpectralDecomp}) that these Poincar\'e series are orthogonal to the non-generic spectrum.  Therefore, the following theorem suffices for our needs.

\begin{thm} {\bf (Langlands spectral decomposition for $\SL(4,\mathbb Z)$)}
 \label{SpectralDecomp} Let $\phi_1, \phi_2, \phi_3, \ldots$ denote an orthogonal   basis of Maass  cusp forms for $\SL(4,\mathbb Z)$.   Assume that  $F, G \in \mathcal L^2(\SL(4,\mathbb Z)\backslash\mathfrak h^4)$ are orthogonal to the non-generic spectrum. Then for
$g\in \GL(4,\mathbb R)$ we have 
\begin{align*}
& F(g) = \sum_{j=1}^\infty \langle F, \phi_j\rangle \frac{\phi_j(g)}{\langle\phi_j, \phi_j\rangle} \; + \; \sum_{\mathcal P} \sum_{\Phi} c_{\mathcal P} \hskip-10pt\int\limits_{\re(s_1)=0}  \hskip-5pt \cdots \hskip-5pt\int\limits_{\re(s_{r-1})=0} \hskip-14pt\Big\langle F, E_{\mathcal P,\Phi}(*\, ,s)\Big\rangle E_{\mathcal P,\Phi}(g\, ,s)
\  ds_1 \cdots ds_{r-1};\\
&\nonumber
\\
&
\langle F, G\rangle =  \sum_{j=1}^\infty  \frac{\langle F, \phi_j\rangle \, \langle\phi_j, G\rangle}{\langle\phi_j, \phi_j\rangle} \; + \; \sum_{\mathcal P} \sum_{\Phi} c_{\mathcal P} \hskip-7pt\int\limits_{\re(s_1)=0}  \hskip-5pt \cdots \hskip-5pt\int\limits_{\re(s_{r-1})=0} \hskip-10pt\Big\langle F, \, E_{\mathcal P,\Phi}(*\, ,s)\Big\rangle      \Big\langle E_{\mathcal P,\Phi}(*\, ,s), \,G \Big\rangle \ ds_1 \cdots ds_{r-1};
\nonumber\end{align*}
where the sum over $\mathcal P$ ranges over  proper parabolic subgroups associated to partitions $4 =  \sum\limits_{k=1}^r n_k$ ($r>1$), and the sum over $\Phi$ (see Definition \ref{CuspFormforEis}) ranges over an orthonormal basis of Maass cusp forms associated to $\mathcal P$. Here  $s 
= (s_1, \ldots,s_r)$ where $\sum\limits_{k=1}^r n_k s_k = 0$ for the partition  $4=\sum\limits_{k=1}^r n_k.$ Furthermore, $c_\mathcal P$ 
is a fixed non-zero constant for each parabolic subgroup $\mathcal P.$
\end{thm}

\begin{proof} 
 For proofs of the spectral expansion for arbitrary reductive groups, see 
\cite{Arthur_1979}, \cite{Langlands1976}, and \cite{MW1995}.  

\vskip 3pt

 It has long been known that the residual spectrum of $\GL(n, \mathbb R)$ 
is non-generic (see, for example, \cite{LM2018}).  According to the results of Kostant, Casselman-Zuckerman, and Vogan (mentioned at the end of Vogan's paper \cite{V78}), generic representations have Gelfand-Kirillov dimension as large as possible.  It follows that Eisenstein series induced from the residual spectrum can never be generic.  Therefore, the spectral 
expansion here only involves terms coming from cusp forms, as claimed.
\end{proof}

 In particular since $P^L, P^M \in \mathcal L^2\left(\SL(4, \mathbb Z) \backslash\mathfrak h^4 \right)$,  the inner product can be computed with the spectral expansion of the Poincar\'e series. The geometric approach utilizes the Fourier Whittaker expansion of the Poincar\'e series which involve Kloosterman sums.
\vskip 3pt

The trace formula takes the following form.
\begin{equation} \label{TraceFormula}
\boxed{\boxed{ \underset{\mbox{\scriptsize{spectral side}}}{\underbrace{\mathcal C \; + \; \mathcal E}} \;\; = \;\; \underset{\mbox{\scriptsize{geometric side}}}{\underbrace{\mathcal M \; + \; \mathcal K }}.}}
\end{equation}

Here $\mathcal C$ is the cuspidal contribution, $\mathcal M$ is the main term coming from the identity element. Further, $\mathcal E$ = Eisenstein contribution, $\mathcal K$ = Kloosterman sum contribution. These will be small with the special choice of the test function $ p_{T,R}$. In the subsections that follow, we explicitly evaluate $\mathcal C, \mathcal E, \mathcal M,$  and  $\mathcal K$.

\subsection{\bf Cuspidal contribution $\mathcal C$ to the Kuznetsov trace 
formula}

\begin{prop} {\bf (Cuspidal contribution to the trace formula)}
\label{CuspContr} Fix  $L= (\ell_1,\ell_2,\ell_3)$ and $M = (m_1,m_2,m_3)$, where $\ell_1,\ell_2,\ell_3,m_1,m_2,m_3$ are non-zero rational integers. Let $\phi_1, \phi_2, \ldots$ denote an orthogonal basis of Maass cusp forms for $\SL(4, \mathbb Z)$ with spectral parameters $\alpha^{(1)}, 
\alpha^{(2)}, \ldots$, respectively, ordered by Laplace eigenvalue.  Let $A_j(L)$ and  $A_j(M)$ denote the $L^{th}$ and $M^{th}$ Fourier coefficients of $\phi_j$ and assume that $A_j(1,1,1) = 1$ for all $j=1,2,\ldots$. Let $\mathcal L_j = L(1, \Ad\; \phi_j).$
Then the  cuspidal contribution to the trace formula  \eqref{TraceFormula} is given by
$$\mathcal C =\sum_{j=1}^\infty  \frac{A_j(L)\overline{A_j(M)}\cdot\left|  p_{T,R}^\#\left(\, \overline{\alpha^{(j)}}\,\right)\right|^2}{ \mathcal L_j \prod\limits_{1\le j \ne k\le 4} \Gamma\left( \frac{1+\alpha_j-\alpha_k}{2}  \right)}$$
where  $p_{T,R}^{\#}$ is given by \eqref{testfunctionsharp2}.
 \end{prop}
\vskip 10pt
\begin{proof}It follows from Theorem~\ref{SpectralDecomp} that for $g\in\GL(4,\R)$, the  cuspidal contribution to the trace formula \eqref{TraceFormula} is given by
$$ \sum_{j=1}^\infty \big\langle P^L(*,0), \;\phi_j\big\rangle\cdot \frac{\phi_j(g)}{\langle\phi_j, \phi_j\rangle}.$$
Now, since $A_j(1,1,1) = 1$, Proposition~\ref{PropFirstCoeff} implies that 
$$\big\langle\phi_j, \phi_j\big\rangle= \mathfrak c_4^{-1}\cdot \mathcal L_j\prod\limits_{1\le j\ne k\le 4}\Gamma\left(\frac{1+\alpha_j-\alpha_k}{2}  \right).$$
The cuspidal contribution to the trace formula is given by
$$\mathcal C  \; :=  \; C_{L,M}^{-1}\cdot \sum_{j=1}^\infty \frac{\big\langle P^L(*,0), \;\phi_j\big\rangle \cdot \overline{\big\langle P^M(*,0), \;\phi_j\big\rangle}}{ \langle\phi_j, \phi_j\rangle}.$$
The proposition immediately follows from the inner product formula (\ref{InnerProductPoincare}).\end{proof}

\subsection{\bf Geometric side of the Kuznetsov trace formula}
Next we consider the geometric side of the trace formula (\ref{TraceFormula}). This is computed with the Fourier-Whittaker expansion of the Poincar\'e series given in (\ref{FourierExp}).

\begin{prop}[Geometric side of the trace formula] \label{GeomSide}
Fix $L= (\ell_1,\ell_2,\ell_3)$, $M = (m_1,m_2,m_3)$ with $C_{L,M}\ne 
0$. Then
$$\lim\limits_{s\to 0} \; \big\langle P^L(*, s), \; P^M(*, s)  \big\rangle\;  = \; \sum\limits_{w\in W_4} \mathcal I _w$$
where we have
\begin{align} \label{eq:IwDef}
\mathcal I_w & :=  \;  \sum_{v\in V_4}\sum_{c_1=1}^\infty \; \sum_{c_2=1}^\infty\; \sum_{c_3=1}^\infty   S_w(\psi_L,\psi_M^v,c)\int\limits_{y_1=0}^\infty \; \int\limits_{y_2=0}^\infty\; \int\limits_{y_3=0}^\infty\;\;\int\limits_{ {U}_w(\mathbb Z)\backslash  {U}_w(\mathbb R)}\;\int\limits_{\overline{{U}}_w(\mathbb R)}\\ \nonumber
&
\hskip 160pt
\cdot
\psi_L(w u y)\,\overline{\psi_M^v(u)}\;
  p_{T,R}(Lc w u y)\,\overline{ p_{T,R}(My)}\;  d^* u
\,\frac{dy_1\,dy_2\,dy_3}{y_1^4 \,y_2^5 \,y_3^4}.
\end{align} 
\end{prop}

\begin{proof}
We  compute the  inner product
\begin{align*}&
\lim_{s\to0}\big\langle P^L\left(*, \,s \right), \, P^M\left(*, \,s \right) \big\rangle  = \lim_{s\to0} \int\limits_{\SL(4,\mathbb Z)\backslash\mathfrak h^4} P^L\left(g, s \right) \cdot  \overline{P^M\left(g, s \right)} \; dg\\
  = &\lim_{s\to0} \int\limits_{U_4(\mathbb Z)\backslash\mathfrak h^4} P^L\left(g, s \right) \cdot  \overline{\psi_M(g)\,  p_{T,R}(Mg)\, I_s(g)} \; dg\\
 =& \lim_{s\to0} \int\limits_{y_1=0}^\infty \; \int\limits_{y_2=0}^\infty\; \int\limits_{y_3=0}^\infty\;\left( \;\int\limits_{U_4(\mathbb Z)\backslash U_4(\mathbb R)} P^L\left(uy, s \right) \cdot \overline{\psi_M(u)} \;du\right) \cdot \overline{ p_{T,R}(My)\, I_s(y)}\;\,\frac{dy_1\,dy_2\,dy_3}{y_1^4 \,y_2^5 \,y_3^4}.
\end{align*}

It follows from (\ref{FourierExp}) that
\begin{align*}
\lim_{s\to0}\big\langle P^L\left(*, \,s \right), \, P^M\left(*, \,s \right) \big\rangle
&\\
& \hskip-100pt = \lim_{s\to0} \sum_{w\in W_4}\sum_{v\in V_4}\sum_{c_1=1}^\infty \; \sum_{c_2=1}^\infty\; \sum_{c_3=1}^\infty   \frac{S_w(\psi_L,\psi_M^v,c)}{c_1^{4s_1} c_2^{4s_2} c_3^{4s_3}   }  \int\limits_{y_1=0}^\infty \; \int\limits_{y_2=0}^\infty\; \int\limits_{y_3=0}^\infty\;\;\int\limits_{{U}_w(\mathbb Z)\backslash {{U}}_w(\mathbb R)}\;\int\limits_{\overline{{U}}_w(\mathbb R)}\\
&
\hskip-37pt
\cdot
\psi_L(w u y)\,\overline{\psi_M^v(u)}\;
 p_{T,R}(Lc w u y)\,\overline{ p_{T,R}(My)}\; I_s(w u  y) \;\overline{ I_s(y)}\;
  d^* u
\,\frac{dy_1\,dy_2\,dy_3}{y_1^4 \,y_2^5 \,y_3^4}\\
& \hskip-100pt =  \; \sum_{v\in V_4}\sum_{c_1=1}^\infty \; \sum_{c_2=1}^\infty\; \sum_{c_3=1}^\infty   S_w(\psi_L,\psi_M^v,c)\int\limits_{y_1=0}^\infty \; \int\limits_{y_2=0}^\infty\; \int\limits_{y_3=0}^\infty\;\;\int\limits_{ {U}_w(\mathbb Z)\backslash {{U}}_w(\mathbb R)}\;\int\limits_{\overline{{U}}_w(\mathbb R)}\\
&
\hskip 40pt
\cdot
\psi_L(w u y)\,\overline{\psi_M^v(u)}\;
  p_{T,R}(Lc w u y)\,\overline{ p_{T,R}(My)}\;  du
\,\frac{dy_1\,dy_2\,dy_3}{y_1^4 \,y_2^5 \,y_3^4}\\
&\hskip-100pt  = \sum_{w\in W_4} \mathcal I _w,
\end{align*}
as claimed.
\end{proof}

 \subsection{\bf Main term $\mathcal M$ in the Kuznetsov trace formula}

Let $w_1$ denote the $4\times 4$  identity matrix. The main term $\mathcal M = I_{w_1}$ in the trace formula (\ref{TraceFormula}) can now be easily computed.

 \begin{prop}{\bf(Main term in the trace formula)}\label{prop:mainterm} There exist  fixed constants $\mathfrak c_1,\mathfrak c_2, \mathfrak c_3>0$ (depending only on R) such that the main term $\mathcal M$ in the trace 
formula \eqref{TraceFormula}  is given by
 $$\mathcal M =  \delta_{L,M}\cdot\Big(\mathfrak c_1T^{9+8R} + \mathfrak c_2 T^{8+8R} + \mathfrak c_3 T^{7+8R} + \mathcal O\left( T^{6+8R}  \right)\Big).$$
 
  \end{prop}

\begin{proof} 
The Kloosterman sum in Definition~\ref{KloostSum}
for the special case of the trivial Weyl group element $w_1$ is identically zero unless  $c = (1,1,1)$ in which case $S_{w_1}\big(\psi_M, \psi_L^\nu, (1,1,1)\big) = 1$. It follows from (\ref{eq:IwDef}) and the normalization (by $C_{L,M}$) of the cuspidal contribution $\mathcal C$ that
\begin{align*}
\mathcal M & = C_{L,M}^{-1} \cdot \mathcal I_{w_1}\\
& = C_{L,M}^{-1}\cdot \int\limits_{y_1=0}^\infty \; \int\limits_{y_2=0}^\infty\; \int\limits_{y_3=0}^\infty\;\left(\;\int\limits_{{U}_w(\mathbb Z)\backslash {{U}}_w(\mathbb R)} \hskip-5pt\psi_{L}(u)\overline{\psi_M(u)}\; d^*u\right) p_{T,R}(Ly)\, \overline{p_{T,R}(My)} \;\frac{dy_1\,dy_2\,dy_3}{y_1^4 \,y_2^5 \,y_3^4}.
\end{align*}

Next
\begin{align*}
 \mathcal M & = \delta_{L,M}\cdot \mathfrak c_4\int\limits_{y_1=0}^\infty \; \int\limits_{y_2=0}^\infty\; \int\limits_{y_3=0}^\infty | p_{T,R}(y)|^2 \,\frac{dy_1\,dy_2\,dy_3}{y_1^4 \,y_2^5 \,y_3^4} = \big\langle  p_{T,R},  p_{T,R}\big\rangle \\
 & = \; \delta_{L,M}\cdot \mathfrak c_4\iiint\limits_{\re(\alpha_j)=0}  \frac{\big| p_{T,R}^\#(\alpha)\big|^2 \;d\alpha_1 \, d\alpha_2 \, d\alpha_3}{\prod\limits_{1\le j\ne k\le 4}\Gamma\left(\frac{\alpha_j-\alpha_k}{2}  \right) }                     \\
 & = \delta_{L,M}\, \mathfrak c_4\cdot \big\langle  p_{T,R}^\#, \; p_{T,R}^\#\big\rangle,
 \end{align*}
where the second representation of $\mathcal M$ in terms of the norm of $ 
p_{T,R}^\#$ follows from the Plancherel formula in Corollary 1.9 of \cite{GK2012}.

It now follows from \eqref{testfunctionsharp2} that
$$\mathcal M = \delta_{L,M}\cdot \mathfrak c_4\iiint\limits_{\re(\alpha_j)=0} \frac{ \Big|e^{\frac{\alpha_1^2+\alpha_2^2+\alpha_3^2+\alpha_4^2}{2T^2}}
   \mathcal{F}_R(\alpha)  \prod\limits_{1\leq\, j \ne k\, \leq 4}
      \Gamma\left(\textstyle{\frac{2+R+\alpha_j - \alpha_k}{4}}\right)\Big|^2}{   \prod\limits_{1\le j\ne k\le 4}\Gamma\left(\frac{\alpha_j-\alpha_k}{2}  \right)     } \; d\alpha_1 d\alpha_2 d\alpha_3.$$

Let $\alpha_j = i\tau_j$ with $\tau_j\in\mathbb R,$  where $\tau_4 = -\tau_1-\tau_2-\tau_3.$  We see that
\begin{align*}
 |\mathcal{F}_R(\alpha)|^2 & :=\Big(\big(1+|\tau_1+\tau_2 -\tau_3-\tau_4|^2\big) 
 \big(1+|\tau_1+\tau_3 -\tau_2-\tau_4|^2\big)
 \big(1+|\tau_1+\tau_4 -\tau_2-\tau_3|^2\big)\Big)^{\frac{R}{3}}.
\end{align*}
 It then follows from Stirling's asymptotic formula 
$|\Gamma(\sigma+it)|^2   \sim     2\pi \lvert t \rvert^{2\sigma-1}\, e^{-\pi\abs{t}}$ 
that
\begin{align*}\mathcal M & \sim  \delta_{L,M} \cdot \mathfrak c_4\\&\cdot 
\iiint\limits_{\mathbb{R}^3}    e^{\frac{-\tau_1^2-\tau_2^2-\tau_3^2-\tau_4^2}{2T^2}} 
 \Big(\big(1+ \big|\tau_1+\tau_2-\tau_3-\tau_4\big|^2\big)
   \big(1+ \big|\tau_1-\tau_2+\tau_3-\tau_4\big|^2\big)
   \big(1+ \big|\tau_1-\tau_2-\tau_3+\tau_4\big|^2\big)\Big)^{\frac{R}{3}}
   \\
   &
   \hskip 20pt
   \cdot
      \Big(\big(1 + |\tau_1 - \tau_2|\big) 
      \big(1 + |\tau_1 - \tau_3|\big) 
      \big(1 + |\tau_2 - \tau_3|\big) 
      \big(1 + |2 \tau_1 + \tau_2 + \tau_3|\big)\\
      &
      \hskip 150pt
      \cdot 
     \big (1 + |\tau_1 + 2 \tau_2 + \tau_3|\big) 
     \big (1 + |\tau_1 + \tau_2 + 2 \tau_3|\big)\Big)^{1+R}\; d\tau_1 d\tau_2 d\tau_3.
      \end{align*}
      Next, make the change of variables
      $$\tau_1\to \tau_1 T, \quad \tau_2\to \tau_2 T, \quad \tau_3\to \tau_3 T.$$
      It follows that as $T\to\infty$ we have
      $\mathcal M  \sim  \mathfrak c_1 \cdot  \delta_{L,M} \, T^{9+8R},$  
    where
      \begin{align*} \mathfrak c_1 & = \mathfrak c_4 \iiint\limits_{\mathbb{R}^3}    e^{\frac{-\tau_1^2-\tau_2^2-\tau_3^2-\tau_4^2}{2}} \Big(\big(1+ \big|\tau_1+\tau_2-\tau_3-\tau_4\big|^2\big)
   \big(1+ \big|\tau_1-\tau_2+\tau_3-\tau_4\big|^2\big)
   \big(1+ \big|\tau_1-\tau_2-\tau_3+\tau_4\big|^2\big)\Big)^{\frac{R}{3}}
   \\
   &
   \cdot
      \Big(\big(1 + |\tau_1 - \tau_2|\big) 
      \big(1 + |\tau_1 - \tau_3|\big) 
      \big(1 + |\tau_2 - \tau_3|\big) 
      \big(1 + |2 \tau_1 + \tau_2 + \tau_3|\big)\\
      &
      \hskip 150pt
      \cdot 
     \big (1 + |\tau_1 + 2 \tau_2 + \tau_3|\big) 
     \big (1 + |\tau_1 + \tau_2 + 2 \tau_3|\big)\Big)^{1+R}\; d\tau_1 d\tau_2 d\tau_3.
      \end{align*}
This method of proof can be extended by using additional terms in Stirling's asymptotic expansion for the Gamma function to obtain additional terms in the asymptotic expansion of $\mathcal M.$
\end{proof}

\subsection{\bf Kloosterman term $\mathcal K$ in the Kuznetsov trace formula}

It immediately follows from Proposition \ref{GeomSide} that
$$\mathcal K = C_{L,M}^{-1}\cdot \underset{w\ne w_1}{\sum\limits_{w\in W_4}} \mathcal I_w.$$

\subsection{\bf Eisenstein contribution $\mathcal E$ to the Kuznetsov trace formula}
\label{EisensteinTerm}
This section is based on \cite{GMW2019}. 
There are 4 standard non-associate parabolic subgroups  on $\GL(4)$ corresponding to the partitions
$$
4 \; = \; 3+1 \; =  \; 2+2 \;   = \; 2+1+1 \;  = \; 1+1+1+1.
$$

Consider the minimal parabolic subgroup  
$$\mathcal P_{\text{\rm Min}}:= \left\{\left(\begin{smallmatrix}  
*&* &* &*\\
&  *&*&*\\
& & *&* \\
& & & *\end{smallmatrix}\right) \subset \GL(4,\mathbb R)   \right\}$$
associated to the partion 4 = 1+1+1+1. 
 Let $s = \frac14 + (s_1, s_2,s_3)$ with $s\in\mathbb C^{3}$.  The minimal parabolic Eisenstein series  for $\Gamma = \SL(4,\mathbb Z)$ is defined by
$$E_{\mathcal P_{\text{\rm Min}}}(g, s) := \sum_{\gamma\, \in \,\left(\mathcal P_{\text{\rm Min}}\cap \Gamma    \right)\backslash\Gamma} I_s(\gamma g), \qquad (g\in \GL(4,\mathbb R), \;\; \re(s)\gg 1).$$ 

\vskip 8pt 
For the other 3 partitions $4  =  3+1  =   2+2    = 2+1+1$, consider the parabolic subgroups
$$\mathcal P_{3,1}:=\left\{\left(\begin{smallmatrix}  
*&* &* &*\\
*&  *&*&*\\
*&* & *&* \\
& & & *\end{smallmatrix}\right)  \right\}, \quad 
\mathcal P_{2,2}:=\left\{\left(\begin{smallmatrix}  
*&* &* &*\\
*&  *&*&*\\
& & *&* \\
& &* & *\end{smallmatrix}\right)  \right\},\quad
\mathcal P_{2,1,1}:=\left\{\left(\begin{smallmatrix}  
*&* &* &*\\
*&  *&*&*\\
& & *&* \\
& & & *\end{smallmatrix}\right)  \right\},
$$
respectively. For each of these parabolic subgroups (denoted $\mathcal P$), associated to a partition with $4 = n_1+n_2 + \cdots +n_r$,  we 
may define an infinite family of Eisenstein series $E_{\mathcal P, \Phi}$ 
(as in \cite{GMW2019}) given by
$$E_{\mathcal P, \Phi}(g,s) := \sum_{\gamma\,\in\, (\mathcal P \cap \Gamma)\backslash \Gamma}  \Phi(\gamma g)\cdot |\gamma g|^s_{_{\mathcal P}}, 
\qquad (g\in \GL(4,\mathbb R), \;\; \re(s)\gg 1).$$
Here $\Phi$ runs over cusp forms on the Levi components 
$m = \left(\begin{smallmatrix} m_1& 0&\cdots & 0\\
0 & m_2 &\cdots & 0\\
\vdots& \vdots &\ddots & \vdots\\
0 & 0& \cdots & m_r \end{smallmatrix}   \right)$ of $g$ (with $m_i\in \GL(n_i,\mathbb R)$)  and
 $|g|_{_{\mathcal P}}^s = \prod\limits_{j=1}^r |\text{\rm det}(m_j)|^{s_j}$  is a toric character as  in \cite{GMW2019}. 
\begin{dfn} \label{CuspFormforEis} The cusp forms $\Phi$ and complex $s$ values associated to Eisenstein series $E_{\mathcal P, \Phi}(g,s)$ for $\SL(4, \mathbb Z)$ are given as follows.
\vskip 4pt
\noindent
$\bullet$  Let $\mathcal P = \mathcal P_{3,1}$, then $\Phi$ runs over Maass cusp forms for $\SL(3,\mathbb Z)$ and $s=(s_1,s_2)$ with $3s_1 + s_2 = 0.$ 

\vskip 4pt
\noindent
$\bullet$ Let 
$\mathcal P = \mathcal P_{2,2}$, then $\Phi$ runs over pairs $\Phi = (\phi_1, \phi_2)$ where $\phi_1,\phi_2$ are Maass cusp forms for $\SL(2, \mathbb Z)$ and $s=(s_1,s_2)$ with $2s_1+2s_2 = 0.$

\vskip 4pt
\noindent
$\bullet$
Let $\mathcal P = \mathcal P_{2,1,1}$, then $\Phi$ runs over Maass cusp 
forms for $\SL(2,\mathbb Z)$ and $s = (s_1,s_2,s_3)$ with $2s_1+s_2+s_3=0.$
\end{dfn}

Now $E_{\mathcal P, \Phi}$  has a Fourier-Whittaker expansion similar to \eqref{FourierExpansion}. If $M = (m_1,m_2,m_3)$ with $m_1$, $m_2$ and $m_3$ positive integers, then the $M^{th}$ Fourier-Whittaker coefficient of $E_{\mathcal P, \Phi}(*, s)$ is given by
\begin{align*}
\int\limits_{U_4(\mathbb Z)\backslash U_4(\mathbb R)} E_{\mathcal P}(uy, s)\, \overline{\psi_M(u) }\; du \; = \; \frac{A_{E_{\mathcal P, \Phi}}(M, s)}
{|m_1|^{\frac32} |m_2|^2 |m_3|^{\frac32}} \; W_\alpha\big(M y\big),
\end{align*}
where $\alpha$ denotes the Langlands parameter of $E_{\mathcal P, \Phi}(g,s).$
Here $$A_{E_{\mathcal P, \Phi}}(M, s) = A_{E_{\mathcal P, \Phi}}\big((1,1,1),s\big) \cdot \lambda_{E_{\mathcal P, \Phi}}(M, s)$$
and  $\lambda_{E_{\mathcal P, \Phi}}(M,s)$
  is the $M^{th}$ Hecke eigenvalue of $E_{\mathcal P, \Phi}$.

\begin{prop}{\bf (Inner product of the Poincar\'e series  $P^M$ with $E_{\mathcal P, \Phi}$)} \label{InnerProductEisenstein} Consider an Eisenstein series $E_{\mathcal P, \Phi}$ for $\SL(4, \mathbb Z).$  Fix $M= (m_1,m_2,m_3)\in\mathbb Z_+^3$. Let $P^M$ be the Poincar\'e series defined in 
\eqref{PoincareSeries} with test function $ p_{T,R}:\mathfrak h^4\to\mathbb C$ (as in  \eqref{testfunctionsharp2}). Then 
$$\lim_{\delta\to 0}\Big\langle P^M(*, \, \delta),\; E_{\mathcal P, \Phi}(*, \; s)\Big\rangle = m_1^{\frac32} m_2^2 m_3^{\frac32} \cdot \overline{A_{E_{\mathcal P, \Phi}}(M,s)}\cdot
p_{T,R}^\#(\overline{\alpha}).$$
\end{prop}

\begin{proof} The proof is similar to the proof of (\ref{InnerProductPoincare}).  One begins by taking $\text{\rm Re}(s)$ very large so that there 
is no problem with convergence and the result follows by analytic continuations in $s$.\end{proof}

\begin{thm} {\bf (Spectral decomposition for the inner product of  Poincar\'e series)} \label{ThmSpectralDecomp} Fix $L=(\ell_1,\ell_2,\ell_3), \, M= (m_1,m_2,m_3) \in\mathbb Z^3$ . 

 Let $\phi_1,\, \phi_2, \,\phi_3, \; \ldots$ denote a basis of Maass cusp 
forms for $\SL(4, \mathbb Z)$ with spectral parameters $\alpha^{(1)}$, $\alpha^{(2)}$, $\alpha^{(3)}, \;\ldots$, respectively, ordered by Laplace eigenvalue and normalized so that the first Fourier coefficient $A_j(1,1,1) = 1$ for all $j=1,2,\ldots$.  Set $\mathcal L_j = L(1,\rm Ad \; \phi_j).$ 

Let  $\mathcal P$ range over  parabolics  associated to partitions $4 = 
n_1+\cdots +n_r,$ and  $\Phi$ range over an orthonormal basis of Maass cusp forms associated to $\mathcal P$. Let $E_{\mathcal P,\Phi}(s)$ denote the Langlands Eisenstein series for $\SL(4,\mathbb Z)$ with Langlands parameter $\alpha_{_{\mathcal P,\Phi}}(s)$ and $L^{th}$, $M^{th}$ Fourier coefficients $A_{E_{\mathcal P,\Phi}}(L, s)$, $A_{E_{\mathcal P,\Phi}}(M, s)$, respectively. Then 
\begin{align*}
&
C_{L,M}^{-1}\cdot \lim\limits_{\delta\to 0} \; \Big\langle P^L(*, \delta), \; P^M(*, \delta)  \Big\rangle  =  \;\sum_{j=1}^\infty \frac{A_j(L)\overline{A_j(M)}\cdot\left| p_{T,R}^\#\left(\, \alpha^{(j)}\,\right)\right|^2}{ \mathcal L_j \prod\limits_{1\le j \ne k\le 4} \Gamma\left( \frac{1+\alpha_j-\alpha_k}{2}  \right)}\\
&
\hskip 5pt 
+ \; \sum_{\mathcal P} \sum_{\Phi} c_{\mathcal P}\hskip-4pt \int\limits_{\re(s_1) = 0}\cdots \int\limits_{\re(s_{r-1}) = 0} \hskip-6pt 
A_{E_{\mathcal P,\Phi}}(L, s) \, \overline{A_{E_{\mathcal P,\Phi}}(M, s)}\cdot\Big| p_{T,R}^\#\big(\alpha_{_{(\mathcal P,\Phi)}}(s)\big)\Big |^2\; 
ds_1 \cdots ds_{r-1},\end{align*}
for constants   $c_{\mathcal P} >0$.
\end{thm}

\begin{proof}
This follows immediately from Proposition  \ref{CuspContr}, Theorem~\ref{SpectralDecomp}, Proposition \ref{InnerProductEisenstein}, and the rapid decay of the $p_{T,R}$ function.  Note that the Poincar{\'e} Series $P^L$, $P^M$ are orthogonal to  the non-generic spectrum.  This can be seen as 
follows:
\[ \int\limits_{\SL(4,\Z)\backslash  \mathfrak{h}^4} \hskip -10pt P^M(g) \overline{\phi(g)}\, dg \ = 
\int\limits_{U_4(\Z)\backslash \mathfrak{h}^4} \hskip -10pt \psi(g) H(Mg) 
I_s(g) \overline{\phi(g)}\, dg \ = \int\limits_{U_4(\R)\backslash \mathfrak{h}^4} \hskip -10pt H(Mg) I_s(g) \hskip -10pt \int\limits_{U_4(\Z)\backslash U_4(\R)} \hskip -10pt \psi(n)\, \overline{\phi(n g)}\, dn\, dg. \]
The last integral is zero if $\phi$ is non-generic.
\end{proof}

 \begin{prop} \label{EisTermKuz}
 {\bf (Eisenstein term $\mathcal E$  in the Kuznetsov trace formula)}  With the notation of Theorem \ref{ThmSpectralDecomp},   the Eisenstein term 
 $\mathcal E$ in the Kuznetsov trace formula is given by
 $$\mathcal E = \sum_{\mathcal P} \sum_{\Phi} c_{\mathcal P} \hskip-4pt\int\limits_{\re(s_1) = 0}\cdots \int\limits_{\re(s_{r-1}) = 0}  \hskip-6pt
A_{E_{\mathcal P,\Phi}}(L, s) \, \overline{A_{E_{\mathcal P,\Phi}}(M, s)}\cdot\Big| p_{T,R}^\#\big(\alpha_{_{(\mathcal P,\Phi)}}(s)\big)\Big |^2\; 
ds_1 \cdots ds_{r-1}.$$
\end{prop}

\begin{proof}
The proof follows from Proposition  \ref{ThmSpectralDecomp}. 
\end{proof}

\section{\bf Bounding the geometric side}

Recall from \eqref{TraceFormula} and \eqref{eq:IwDef} that the Kloosterman contribution to the Kuznetsov trace formula is given by
 \[\mathcal K \; = \; C_{L,M}^{-1} \sum\limits_{w\in W_4} \mathcal I_w \] where $C_{L,M}= \mathfrak c_4 \cdot (\ell_1 m_1)^3(\ell_2m_2)^4(\ell_3m_3)^3$, and
\begin{align}\label{eq:newIw} 
\mathcal I_w & =  \;  \sum_{v\in V_4}\sum_{c_1=1}^\infty \; \sum_{c_2=1}^\infty\; \sum_{c_3=1}^\infty   S_w(\psi_L,\psi_M^v,c)\int\limits_{y_1=0}^\infty \; \int\limits_{y_2=0}^\infty\; \int\limits_{y_3=0}^\infty\;\;\int\limits_{\rm U_w(\mathbb Z)\backslash {U}_w(\mathbb R)}\;\int\limits_{\overline{U}_w(\mathbb R)}
\\  \nonumber
&
\hskip 120pt
\cdot
\psi_L(w u y)\,\overline{\psi_M^v(u)}\;
  p_{T,R}(Lc w u y)\,\overline{ p_{T,R}(My)}\;  d^* u
\,\frac{dy_1\,dy_2\,dy_3}{y_1^4 \,y_2^5 \,y_3^4},
\end{align} 
\vskip -10pt\noindent
and
$c = \left(\begin{smallmatrix} 1/c_3 & & &\\
& c_3/c_2 & &\\
& & c_2/c_1 &\\
& & & c_1  \end{smallmatrix}\right).$

\vskip 4pt
The main term in the Kuznetsov trace formula is given in
 Proposition~\ref{prop:mainterm} and consists of the first term (corresponding to the identity Weyl element $w =w_1$) on the geometric side of the trace formula.
  In this section we bound $I_w$ in each of the remaining cases.

  \vskip 5pt

We remark first, by Friedberg \cite{friedberg1987poincare}, that $\mathcal{I}_w=0$ unless $w$ is \emph{relevant}.  That is $w$ must be of the form
 $w = \left( \begin{smallmatrix} & & I_{n_1} \\ & \iddots & \\ I_{n_k} & & \end{smallmatrix}\right),$
where $I_{n_i}$ is the identity matrix of size $n_i\times n_i$ and $4=\sum\limits_{i=1}^k n_i$ with $1\le k\le 4$ and $n_i\in \mathbb Z_{\ge 0}.$ The relevant Weyl group elements are therefore,
 \[ \begin{array}{lclclcl}
 w_1 = \left( \begin{smallmatrix} 1 & & & \\ & 1 & & \\ & & 1 & \\ & & & 1 \end{smallmatrix}\right), 
 & \phantom{X} & 
 w_2 = \left( \begin{smallmatrix} & & & 1 \\ 1 & & & \\ & 1 & & \\ & & 1 &  \end{smallmatrix}\right), 
 & \phantom{X} & 
 w_3 = \left( \begin{smallmatrix} & 1 & & \\ & & 1 & \\ & & & 1 \\ 1 & & &  \end{smallmatrix}\right), 
 & \phantom{X} & 
 w_4 = \left( \begin{smallmatrix} & & 1 & \\ & & & 1 \\ 1 & & & \\ & 1 & & \end{smallmatrix}\right), 
 \\ \\
 w_5 = \left( \begin{smallmatrix} & & & 1 \\ & 1 & & \\ & & 1 & \\ 1 & & & \end{smallmatrix}\right), 
 & \phantom{X} & 
 w_6 = \left( \begin{smallmatrix} & & 1 & \\ & & & 1 \\ & 1 & & \\ 1 & & & \end{smallmatrix}\right), 
 & \phantom{X} & 
 w_7 = \left( \begin{smallmatrix} & & & 1 \\ & & 1 & \\ 1 & & & \\ & 1 & &  \end{smallmatrix}\right), 
 & \phantom{X} & 
 w_8 = \left( \begin{smallmatrix} & & & 1 \\ & & 1 & \\ & 1 & & \\ 1 & & & \end{smallmatrix}\right). 
 \end{array}\]
 
To motivate the strategy for  bounding $\mathcal I_w,$  we consider \eqref{eq:newIw} (choosing $L=M=1$) and then take absolute values to obtain
 \begin{equation} \label{Example-Iw} | I_w| \ll 
 \sum_{c\in \Z_{+}^3}\; \int\limits_{\R_{+}^3} \left| p_{T,R}(y)\right| \frac{\delta_w(y)}{\delta(y)} \int\limits_{U_w(\mathbb R)} \big | p_{T,R}(cwyu)  \big | \;du\, d^\times y,  \qquad\;\; \Big(y = \text{\rm diag}(y_1y_2y_3, y_1y_2,y_1,1)\Big),
 \end{equation}
  for each Weyl group element $w$ where $\delta(y) = y_1^3y_2^4y_3^3$ is the modular character and $\delta_w(y)$ is the 
Jacobian of $U_w(\mathbb R) \ni u\mapsto yuy^{-1}$. It can be checked that $\frac{\delta_w(y)}{\delta(y)} = \delta^{-\frac12}(y)\cdot \delta^{-\frac12} (wyw^{-1}).$ Theorem \ref{th:pTRbound} is required to prove the absolute convergence of (\ref{Example-Iw}).
 The trivial bound of the Kloosterman sum is given by
$$S(1,1,c)\ll \delta^{-\frac12}(c) \ll c_1c_2c_3.$$  From this, one obtains the absolute 
convergence of the sum over $c$ provided that  $p_{T,R}(y) \ll \delta^{1+\varepsilon}(y).$ Fortunately, then one checks that the $u$-integral is also absolutely convergent and
$$\int\limits_{U_w(\mathbb R)} \big|  p_{T,R}(cwyu) \big|\; du \ll \delta^{1+\varepsilon}(c) \cdot \delta^{1+\varepsilon}(wyw^{-1}   ).$$
Then the remaining $y$ integral becomes
$$\int\limits_{\mathbb R_{+}^3} \big| p_{T,R}(y)  \big| \cdot \delta^{\frac12+\varepsilon}\big(wy w^{-1}  \big)\cdot \delta^{-\frac12}(y) \; d^\times y.$$
So, to make the $y$ integral convergent one also needs
$p_{T,R}(y) \ll \delta^{\frac12}(y)\cdot \delta^{-\frac12}\left(wyw^{-1}\right)\, \hbox{min}(y^{-\varepsilon},y^{\varepsilon}),$
for all Weyl elements $w$. Note that $\hbox{min}(y^{-\varepsilon},y^{\varepsilon})$ is needed as buffer as in the proof of Theorem~\ref{BufferTheorem} and that for $w = w_{\text{long}}$ this bound coincides with $\delta^{1+\varepsilon}(y).$ One also notes that $\delta^{\frac12}(y)$ is the trivial bound of $p_{T,R}(y)$.

The main technical tool for bounding the geometric side, the proof of which occupies all of Section~\ref{sec:boundsforpTR}, is the following theorem.

\begin{thm}\label{th:pTRbound}
Let $0<\varepsilon<1$ and $r\in \Z_{\ge0}$.  There exists $R$ sufficiently large such that the following is true.  Suppose $a_1,a_2,a_3$ satisfy the conditions
$$\boxed{ \varepsilon\le |2r_j-a_j |\le1-\varepsilon\quad (1\le j\le 3)}\mbox{\quad if\quad  } (r_1,r_2,r_3)\in \big\{ (0,r,0),(r,r,0),(0,r,r),(r,r,r) \big\},$$or  $$\boxed{\varepsilon\le |2r-a_1 |\le1-\varepsilon,\quad 
\varepsilon\le |a_2 |\le1-\varepsilon,\quad  -1+\varepsilon\le a_3<-\varepsilon} \quad\mbox{if}\quad (r_1,r_2,r_3)=(r,0,0),$$
or  $$\boxed{-1+\varepsilon\le a_1<-\varepsilon ,\quad \varepsilon\le |a_2 |\le1-\varepsilon,\quad \varepsilon\le |2r-a_3 |\le1-\varepsilon} \quad\mbox{if}\quad (r_1,r_2,r_3)=(0,0,r).$$Then we have the bound
 \[ \boxed{p_{T,R}(y) \ll 
 y_1^{\frac32+a_1}y_2^{2+a_2}y_3^{\frac32+a_3} T^{\varepsilon+4R+9+\sum\limits_{j=1}^3 (\delta_{0,r_j}-r_j)}}. \]
Here, $\delta_{0,r_j}$ is equal to $1$ if $r_j=0$ and is zero otherwise.  The implicit constant depends on $\varepsilon$ and $R$.
\end{thm}

\begin{proof}
We will show in Section~\ref{sec:boundsforpTR} that $p_{T,R}(y)$ can be expressed as a sum of a ``shifted $p_{T,R}$ term,'' ``single residue terms,'' ``double residue terms'' and ``triple residue terms'' (see \eqref{eq:pTRdecomposition}).  A bound for the shifted $p_{T,R}$ term is given in Proposition~\ref{prop:GL4pTRemptyset}.  {Note that in the proof of Proposition~\ref{prop:GL4pTRemptyset} the first set of inequalities for $a_j$ ($j=1,2,3$) is shown to hold for any $(r_1,r_2,r_3)$ with $r_j\geq 0$.}  Bounds for the single residue terms are given in Proposition~\ref{prop:singleresidue1} and Proposition~\ref{prop:singleresidue2}.  {See the proof of Lemma~\ref{lem:pTr_1alpha2shift} for details of why the further inequalities given above are required.}  Bounds for the double residue terms are given in Proposition~\ref{prop:doubleresidue12} and Proposition~\ref{prop:doubleresidue13}.  Finally, bounds for the triple residue terms are given in Proposition~\ref{prop:tripleresidue}.  Combining these bounds gives the desired result.
\end{proof}

\begin{prop}\label{prop:Iwbounds}
Let $\mathcal{I}_w$ be as above.  Let $M=(m_1,m_2,m_3), \;  L=(\ell_1,\ell_2,\ell_3) \in \mathbb Z^3$, where $C_{L,M}\ne 0$ is the normalization factor given in Definition \ref{NormFactor}. Let $r\geq 1$ be an integer.  Then for $R$ sufficiently large and any $\varepsilon>0$, we have 
 \[ C_{L,M}^{-1}\big| \mathcal{I}_{w_j} \big| \ll_{\varepsilon, R} \;(\ell_1m_1)^{2r-1/2}(\ell_2m_2)^{2r-1}(\ell_3m_3)^{2r-1/2} B_j(T), \]
where
\[ B_j(T) = 
 \begin{cases}
 T^{\varepsilon+8R+20-4r} & \mbox{ if }j=2,3,4, \\
 T^{\varepsilon+8R+19-5r} & \mbox{ if }j=6,7, \\
 T^{\varepsilon+8R+18-6r} & \mbox{ if }j=5,8.
 \end{cases}\]
\end{prop}

\begin{rmrk}
Given Proposition~\ref{prop:mainterm} (the asymptotic formula for the main term with asymptotic error terms), we need $20-4r<6$ in order for the contribution from $\mathcal{I}_{w}$ to be smaller than the main term $\mathcal{M}$.  In other words, we will need to take $r\geq 4$. In order to minimize the contribution of $L$ and $M$ in the Main Theorem \ref{MainTheorem}, we take $r=4.$

\end{rmrk}

\begin{proof} 
The main idea of the proof of Proposition \ref{prop:Iwbounds}  is to apply Theorem \ref{th:pTRbound} to each of the two instances of $p_{T,R}$ which appear on the right hand side of \eqref{eq:newIw}. 
 Before doing so, we make the change of variables $u\mapsto yuy^{-1}$.  Note that, by definition, $d(yuy^{-1})= \delta_w(y)\,du$ for any $u\in \overline{U}_w$.  %After doing this, we also make the change of variables $y\mapsto L^{-1}y$.  
We see that taking absolute values in \eqref{eq:newIw}, we obtain
\begin{align}\label{eq:newerIw} 
\lvert \mathcal I_w \rvert & \ll  \;  \sum_{v\in V_4}\sum_{c_1=1}^\infty \; \sum_{c_2=1}^\infty\; \sum_{c_3=1}^\infty   \lvert S_w(\psi_L,\psi_M^v,c) \rvert \int\limits_{y_1=0}^\infty \; \int\limits_{y_2=0}^\infty\; \int\limits_{y_3=0}^\infty\;\;\int\limits_{{U}_w(\mathbb Z)\backslash {U}_w(\mathbb R)}\;\int\limits_{\overline{U}_w(\mathbb R)}
\\  \nonumber
&
\hskip 80pt
\cdot
  \delta_w(y) \cdot
  \lvert p_{T,R}(Lc w y u) \rvert \, \cdot \, \lvert p_{T,R}(My) \rvert \;  d^* u
\,\frac{dy_1\,dy_2\,dy_3}{y_1^4 \,y_2^5 \,y_3^4}
\end{align} 

For the term $p_{T,R}(Lcwyu)$, we need the Iwasawa form:
 \[ Lcwyu =: u'tk \]
so that we can apply Theorem~\ref{th:pTRbound} for a particular choice of 
integers $r_1,r_2,r_3$ and parameters $2r_j-1<a_j<2r_j$ for each $j=1,2,3$. 

For notational purposes, we write
\begin{align*}
 t = t(Lcwuy) =: Y=(Y_1,Y_2,Y_3).
\end{align*}
It is easy to see that $Y_i$ factors as
 \[ Y_i =: Y_i(w,c)Y_i(w,L)Y_i(w,y)Y_i(w,u) \]
with each factor being an expression in the entries of $c$, $L$, $y$ or $u$ respectively.  (Note that we are following the notation of \eqref{ToricElement} for the element $Y$.)  By inspection, we see that \begin{equation}\label{eq:YLi}
 Y_1(w,L) = \ell_1^{\frac32+a_1}, \qquad Y_2(w,L) = \ell_2^{2+a_2}, \qquad Y_3(w,L) = \ell_3^{\frac32+a_3}
\end{equation}
for all $w$.  Moreover,
\begin{equation}\label{eq:YuXi}
 Y_{1}(w,u) = \frac{\sqrt{\xi_2(w,u)}}{\xi_1(w,u)}, \qquad 
 Y_{2}(w,u) = \frac{\sqrt{\xi_1(w,u)\cdot \xi_3(w,u)}}{\xi_2(w,u)}, \qquad 
 Y_{3}(w,u) = \frac{\sqrt{\xi_2(w,u)}}{\xi_3(w,u)},
\end{equation}
where each function $\xi_j(w,u)$ is strictly positive for any $u\in U(\R)$.

For example, in the case of $w=w_8$, the long element, we find that if
\[ u   = \left(\begin{matrix}1 & x_{12} & x_{13} & x_{14}\\
0 & 1 & x_{23} & x_{24}\\
0&0&1 & x_{34}\\
0 & 0 & 0 & 1
\end{matrix}\right), \quad  y = (y_1,y_2,y_3) =  \left(\begin{matrix} 
y_1 y_2 y_3 & 0 & 0 & 0\\
0 & y_1y_2 & 0 & 0\\
0 & 0 & y_1 & 0\\
0 & 0 & 0 & 1
 \end{matrix}\right), \]
then
\begin{equation}\label{eq:Yiw8}
 Y_{1}(w_8,y) = \frac{1}{y_3},\qquad Y_{2}(w_8,y) =  \frac{1}{y_2},\qquad Y_{3}(w_8,y) =  \frac{1}{y_1},
\end{equation}
and
\begin{align*}
 \xi_1(w_8,u) & = 1+x_{12}^2+x_{13}^2 + x_{14}^2,
 \\
 \xi_2(w_8,u) & = 1+x_{23}^2+x_{24}^2 + \big( x_{12} x_{24} - x_{14} \big)^2 
+\big( x_{12} x_{23} - x_{13} \big)^2
+\big( x_{13} x_{24} - x_{14} x_{23}   \big)^2,
 \\
 \xi_3(w_8,u) & = 1+x_{34}^2 + \big(x_{23} x_{34} - x_{24}\big )^2  + \big( x_{12} x_{23} x_{34} - x_{13} x_{34} - x_{12} x_{24} + x_{14}\big )^2.
\end{align*}

\vskip 12pt

Returning now to \eqref{eq:newerIw}, we apply Theorem~\ref{th:pTRbound} for $(r,r,r)$ and yet to be determined values $a_1,a_2,a_3$ within the permissible bounds ($\varepsilon\leq \lvert 2r-a_j\rvert \leq 1-\varepsilon$) to the term $p_{T,R}(u'tk)$.  To the other term we apply the theorem for a choice $(r_1',r_2',r_3')$ and values $b_1,b_2,b_3$.  As we will see, values of $(r_1',r_2',r_3')$ and $b_j'$ will be forced upon us in order to guarantee the convergence of the sum over $c$, the integral over $(y_1,y_2,y_3)$ and the integral over $u$.  In fact, it will become evident that these values are determined by $w$ and $a_1,a_2,a_3$.

Independent of the choice of $w$, we define $$I_0 := (0,1], \qquad I_1 = (1,\infty),$$
hence
$$\int\limits_{y_1=0}^\infty \; \int\limits_{y_2=0}^\infty\; \int\limits_{y_3=0}^\infty = \sum_{i,j,k\in\{0,1\}}
\int_{I_i} \int_{I_j}\int_{I_k}.$$
Using the bounds for $p_{T,R}$ (once with parameters $a_1,a_2,a_3$ and once with parameters $b_1,b_2,b_3$ as described above), it follows that
\begin{align}\label{eq:Iwsimplified}
 |\mathcal I_w| & \ll \; T^{\varepsilon + 8R+18-3r-r_1'-r_2'-r_3'+\sum\limits_{i=1}^{3}\delta_{r_i',0}}
 (\ell_1 m_1)^{\frac32} (\ell_2m_2)^{2} (\ell_3 m_3)^{\frac32}
 \ell_1^{a_1}\ell_2^{a_2}\ell_3^{a_3}m_1^{b_1}m_2^{b_2}m_3^{b_3}
\\ \nonumber & \hskip 24pt 
 \cdot\sum_{v\in V_4}\sum_{c_1=1}^\infty \; \sum_{c_2=1}^\infty\; \sum_{c_3=1}^\infty \frac{ \big |S_w(\psi_M,\psi_L^v,c)\big| }{  c_1^{1+2 a_1-a_2} c_2^{1-a_1+2 a_2-a_3} c_3^{1-a_2+2 a_3}  }
\\ \nonumber & \hskip -24pt
\cdot \sum_{i,j,k\in\{0,1\}}
\int_{I_i} \int_{I_j}\int_{I_k}
  Y_{1}(w, y)^{\frac32+a_1}  Y_{2}(w, y)^{2+a_2}   Y_{3}(w, y)^{\frac32 +a_3}
 \cdot   y_1^{\frac32+b_1} y_2^{2+b_2} y_3^{\frac32 +b_3} 
\cdot |\delta_w(y)|\, 
\\ \nonumber &
\hskip 12pt
\cdot
\,\frac{dy_1\,dy_2\,dy_3}{y_1^4 \,y_2^5 \,y_3^4}
\cdot \int\limits_{{U}_w(\mathbb Z)\backslash {U}_w(\mathbb R)}\;\int\limits_{\overline{U}_w(\mathbb R)}
\xi_1^{-\frac12-a_1+\frac{a_2}{2}}
 \cdot \xi_2^{-\frac12+\frac{a_1}{2}-a_2+\frac{a_3}{2}} \cdot \xi_3^{-\frac12+\frac{a_2}{2}-a_3}\;  |d^* u|.
\end{align} 

It is known (cf. \cite{Jacquet1967}) that the integral in $\xi$ will converge provided that each of the exponents of $\xi_j$ for $j=1,2,3$ is less than $-\frac12$.  Explicitly, we require that
\begin{equation}\label{eq:xiconditions}
 2a_3,2a_1 > a_2, \qquad \mbox{and} \qquad 2a_2 > a_1+a_3.
\end{equation}
Note that for any $r\geq 2$ the first of the conditions in \eqref{eq:xiconditions} is ensured for \emph{any} choice of $a_1,a_2,a_3$ such that
\begin{equation}\label{eq:Abounds}
 0< \lvert 2r-a_j \rvert < 1 \qquad \mbox{for each }j=1,2,3.
\end{equation}
In this range, the second condition (that $2a_2>a_1+a_3$) imposes an additional constraint, but we will see that it can be easily satisfied as well.

We next want to determine what additional restrictions on $a_1,a_2,a_3$ must be satisfied to guarantee the convergence of the sum over $c$.  We require the following Kloosterman bounds.

\begin{prop} \label{BufferTheorem} Let $w \ne w_1.$ The sum
 \[ \mathcal{K}(M,L;w,a) := \sum_{v\in V_4}\sum_{c_1=1}^\infty \sum_{c_2=1}^\infty \sum_{c_3=1}^\infty \frac{\lvert S_w(\psi_M,\psi_L^v,c) 
\rvert}{c_1^{1+2a_1-a_2}c_2^{1-a_1+2a_2-a_3}c_3^{1-a_2+2a_3}}, \]
which appears on the right hand side of \eqref{eq:Iwsimplified}, satisfies the following bounds (with all implied constants independent of $M$, $L$,  $c$, and $a$).

 For $2\le j\le 8$ and any $\varepsilon>0$, we have\begin{align*} & |\mathcal{K}(M,L;w_j,a)|\ll  \sum_{c_1=1}^\infty  \frac{1}{c_1^{ 2a_1-a_2} } 
\sum_{c_2=1}^\infty  \frac{1}{ c_2^{ -a_1+2a_2-a_3} }\sum_{c_3=1}^\infty \frac{1}{ c_3^{ -a_2+2a_3}}. \end{align*} \end{prop}
 
 \begin{proof} As noted in Appendix B, we have
 \begin{align}\lvert S_{w_j}(\psi_M,\psi_L^v,c) \rvert\ll   c_1c_2c_3\quad (2\le j\le8).\label{kloos-bound}\end{align}Applying this  to the definition of $\mathcal{K}(M,L,c;w,a)$  gives the desired result.  \end{proof}
 
 Appendix B   gives slightly better bounds for the Kloosterman sums in terms of $c_1,c_2,c_3$ at the expense of some powers of $M,L.$ The bounds in Appendix B will likely be relevant in other applications.  However, we thank the referee for pointing out that 
  the  proof of Proposition \ref{BufferTheorem} only requires the 
   trivial bound for Kloosterman sums proved in \cite[Theorem 0.3 (i)]{DR1998}.

Clearly, the series in the above proposition will converge as long as each $c_j$, for $1\le j\le 3$, occurs to an exponent less than $-1$.  This will be guaranteed by the following  additional set of restrictions: 
\begin{equation}\label{eq:KloosConditions}
  -2a_1+a_2 < -1, \quad  a_1-2a_2+a_3< -1, \quad  a_2-2a_3<-1.
\end{equation}
Again, this is easy to arrange given that we may take any values of $a_j$ 
satisfying the bounds of \eqref{eq:Abounds}.  Note, moreover, that \eqref{eq:xiconditions} is a consequence of \eqref{eq:KloosConditions}.

The next step is to show that the integral over $y=(y_1,y_2,y_3)$ also converges.  An elementary calculation shows that  
 \begin{align*}& \frac{y_1^{\frac32+b_1}y_2^{2+b_2}y_3^{\frac32+b_3} Y_1(w,y)^{\frac32+a_1}Y_2(w,y)^{2+a_2}Y_3(w,y)^{\frac32+a_3}\delta_w(y)}{y_1^4y_2^5 y_3^4} dy_1dy_2dy_3 \\&= y_1^{b_1-e_1(w;a)} y_2^{b_2-e_2(w;a)} y_3^{b_3-e_3(w;a)} \frac{dy_1}{y_1} \frac{dy_2}{y_2} \frac{dy_3}{y_3},\end{align*}
where $e_j(w;a)$ is given by
\begin{equation}\label{eq:ej-table}
 \begin{array}{|c||c|c|c||c|} \hline
 w & e_1(w;a) & e_2(w;a) & e_3(w;a) & (r_1',r_2',r_3') \\ \hline
 w_2 & a_3 & -a_1+a_3 & -a_2+a_3 & (r,0,0) \\
 w_3 & a_1-a_2 & a_1-a_3 & a_1 & (0,0,r) \\
 w_4 & a_2-a_3 & a_2 & -a_1+a_2 & (0,r,0) \\
 w_5 & a_3 & a_1-a_2+a_3 & a_1 & (r,r,r) \\
 w_6 & a_2-a_3 & a_2 & a_1 & (0,r,r) \\
 w_7 & a_3 & a_2 & -a_1+a_2 & (r,r,0) \\
 w_8 & a_3 & a_2 & a_1 & (r,r,r) \\ \hline
 \end{array}
\end{equation}
The result for $w_8$, for example, is established using \eqref{eq:Yiw8} and the fact that $\delta_{w_8}(y) = y_1^3 y_2^4 y_3^3$.

This means that for each $j=1,2,3$ and $k=0,1$ we are searching for $a_1,a_2,a_3$ satisfying \eqref{eq:Abounds} and \eqref{eq:KloosConditions}.  Simultaneously, we search for $r_1',r_2',r_3'\in \Z_{\geq 0}$ and $b_1,b_2,b_3$ satisfying the necessary bounds given in Theorem~\ref{th:pTRbound} such that the integral
 \[ \int_{I_k} y_j^{{b_j-e_j(w;a)}} \frac{dy_j}{y_j} \]
converges.  Since $I_0=(0,1]$ and $I_1=(1,\infty)$, we require that $b_j-e_j(w;a)$ be positive if $k=0$ and negative if $k=1$.  With this in mind, we let $\varepsilon>0$ and choose $b_1,b_2,b_3$ such that
 \[ b_j - e_j(w;a) = \begin{cases}  +\varepsilon & \text{if}\; k=0\\
-\varepsilon &\text{if}\; k=1. \end{cases} \]
Note that having made this choice for $b_j$, the value of $r_j'$ is forced to be that given in the final column of the table above.  More precisely, given $(r_1',r_2',r_3')$ as indicated in the final column of the table, there exists some choice $a_1,a_2,a_3$ and $\varepsilon>0$ for which conditions \eqref{eq:Abounds} and \eqref{eq:KloosConditions} are satisfied, 
and for which the necessary bounds on $b_j=e_j(w;a)\pm \varepsilon$ are 
also satisfied.  

Except in the case that $w=w_2$ or $w=w_3$, the conditions on $b_j$ are that
 \[ 0 < \lvert 2r_j'(w)-e_j(w;a)\pm \varepsilon \rvert < 1. \]
It is not hard to see that in these cases there exists an appropriate choice of parameters $a_1,a_2,a_3$.

On the other hand, if $w=w_2$, for example, we see that $(r_1',r_2',r_3')=(r,0,0)$, and so the restriction on the quantities $b_2=a_3-a_1\pm\varepsilon$ and  $b_3=a_3-a_2\pm\varepsilon$ (which plays the role of $a_3$ in Theorem~\ref{th:pTRbound}) are
\begin{align}
 0 < \lvert b_2 \rvert = \lvert a_1-a_3 \pm \varepsilon \rvert < 1, \label{eq:b-conditions} \\
 -1+\varepsilon \leq b_3 = a_3-a_2\pm\varepsilon \leq -\varepsilon.
\end{align}
But we're assuming  that $a_1-2a_2+a_3<-1$ and $a_1-a_3<1-\varepsilon$ (by equations (\ref{eq:KloosConditions}) and (\ref{eq:b-conditions})), which together imply that $ a_3-a_2<-\varepsilon$ in any event.  A similar analysis applies to the above restriction on $b_1$ in the case that $w=w_3$ and $(r_1',r_2',r_3')=(0,0,r)$.

The values of $r_j'$ for $j=1,2,3$ determine the exponent of $T$ in \eqref{eq:Iwsimplified}, hence the values of $B_j(T)$ in the statement of the proposition. To complete the proof, we note that  $a_j<2r+1$ and $b_j<2r_j'(w)+1\leq 2r+1$, and therefore, the powers of $\ell_1,\ell_2,\ell_3,m_1,m_2,m_3$ in \eqref{eq:Iwsimplified} are as claimed.
\end{proof}

\section{\bf Mellin transform of the $\GL(4)$ Whittaker function}
Let $\alpha=(\alpha_1,\alpha_2,\alpha_3,\alpha_4) \in (i\R)^{4}$ with $\alpha_1+ \alpha_2+\alpha_3+\alpha_4 = 0$.  Let $y=(y_1,y_2,y_3)\in \R^3$ and $s=(s_1,s_2,s_3)\in \C^3$ with $\re(s_j)> \varepsilon$ for $j=1,2,3$ and $\varepsilon>0$. Define the Mellin transform (denoted $\Mellin{\alpha}(s)$) of the $\GL(4,\R)$ Whittaker function $\Whit{\alpha}$ defined in Section~\ref{sec:WhittTransform} by the absolutely convergent integral
\begin{equation}\label{eq:Mellin}
 \Mellin{\alpha}(s) := \int\limits_0^\infty \int\limits_0^\infty \int\limits_0^\infty y_1^{s_1-\frac32}y_2^{s_2-2}y_3^{s_3-\frac32} \Whit{\alpha}(y) \; \frac{dy_1 dy_2 dy_3}{y_1 y_2 y_3}. 
\end{equation}

\subsection{\bf Formula for Mellin transform of the Whittaker function}
Here we present an expression for  $\Mellin{\alpha}(s)$ as an integral, over an additional variable $t\in \C$, of a ratio of Gamma functions involving $s$ and $\alpha$.  For the purposes of this section, we will assume that $\alpha$ is in general position, meaning that $\alpha_j\neq \alpha_k$ for all $j\neq k$.  As explained in Remark~\ref{rmk:generalposition}, for our purposes this assumption is not prohibitive.

In the context of $\GL(4)$, this expression was first given in \cite{Stade1995}; that result was later generalized to $GL(n)$, in Theorem~3.1 of \cite{Stade2001}. The formula for $n=4$ takes the form $\Mellin{\alpha}(s)=2^{-3} \pi^{-s_1-s_2-s_3}\widehat{W}_{\frac{\alpha}{2}}\left(\frac{s}{2}\right)$, where, assuming that $\re(s_j)\geq \varepsilon>0$ for each $j=1,2,3$, 
\begin{align}\label{eq:StadeMellinGL4}&
  \widehat{W}_{\alpha}(s)  =  \Gamma(s_1+\alpha_1)\Gamma(s_1+\alpha_2) \Gamma(s_2-\alpha_1-\alpha_2) \Gamma(s_2+\alpha_1+\alpha_2)\Gamma(s_3-\alpha_1)\Gamma(s_3-\alpha_2) \cdot\frac{1}{2\pi i}   \\ \nonumber
 &  \cdot
 \int\limits_{\re(t)= -\varepsilon'}  \hskip -6pt\frac{\Gamma(-t+\alpha_3)\Gamma(-t+\alpha_4)\Gamma(t+s_1)\Gamma(t+s_2+\alpha_1)\Gamma(t+s_2+\alpha_2)\Gamma(t+s_3+\alpha_1+\alpha_2) }{\Gamma(t+s_1+s_2+\alpha_1+\alpha_2)\Gamma(t+s_2+s_3)} \; dt,
\end{align}
assuming that $0<\varepsilon'<\varepsilon$.

\subsection{\bf Poles and residues of $\Mellin{\alpha}$}\label{poles-res}
Following  \cite{Stade2001}, we introduce the sets
\begin{align}\label{Pjdef}
P_1 & := \{-\alpha_1,-\alpha_2,-\alpha_3,-\alpha_4\},\\
P_2 & := \{\pm(\alpha_1+\alpha_2),\pm(\alpha_1+\alpha_3),\pm(\alpha_2+\alpha_3)\},\nonumber\\
P_3 & := \{ \alpha_1,\alpha_2,\alpha_2,\alpha_4\}.\nonumber
\end{align}Note that if $p\in P_k$ then $-p\in P_{4-k}$.

It is shown in \cite[Propositions 12 and 14]{ST2018} that the sets $P_1,\;P_2$, and $P_3$ determine the poles of $\Mellin{\alpha}(s)$, in the following sense.
 \begin{prop} \label{prop:mellres} The Mellin transform $\Mellin{\alpha}(s)$ extends to a meromorphic function of the variable $s=(s_1,s_2,s_3)\in\C^3$ with the following properties.
 
\vskip 5pt
{\bf (a)} \;$\Mellin{\alpha}(s)$ has a pole at $s_j=p-2\delta$ for each 
$p\in P_j$ and $\delta\in \Z_{\ge0}$.   Moreover, $\Mellin{\alpha}(s)$ has no other poles or polar divisors in $\C^3$.

\vskip 5pt
{\bf (b)} The residue of $\Mellin{\alpha}(s)$ at any of the above poles is a polynomial times a ratio of Gamma functions.  More specifically, for each $\delta\in \Z_{\ge0}$, there exists a polynomial $Q_\delta(b,c,d;e,f,g)$, of degree at most $3\delta$, such that the following are true. 
\begin{align}\label{eq:s1res}&\Res_{s_1=-\alpha_1-2\delta }{\Mellin{\alpha}(s)} =
\frac{1}{\gam{\textstyle\frac{  s_2+s_3+\alpha_1}{2}+\delta  } }\left[\;\prod_{k=2 }^4
\gammdelt{\alpha_k-\alpha_1}{\delta }
\gamm{s_2+\alpha_1+\alpha_k}\gamm{ s_3-\alpha_k}{}
\right] 
\\&
\hskip 90pt
\cdot Q_{\delta }\left( \frac{s_3-\alpha_2}{2}, \frac{s_3-\alpha_3}{2}, \frac{s_3 - \alpha_4}{2}; \frac{1 + \alpha_1 -s_2+s_3}{2}, \frac{s_3 - \alpha_1 }{2}- \delta , \frac{s_2+s_3 + \alpha_1}{2}  \right)
,\nn\end{align} 
\begin{align}\label{eq:s2res}&\Res_{s_2=-\alpha_1-\alpha_4-2\delta }{\Mellin{\alpha}(s)}
 =\left[\;\prod_{j\in\{1,4\}} \gamm{s_1+\alpha_j} \right]
 \left[\prod_{k=2}^3   \gamm{s_3-\alpha_j} \prod_{ j\in\{1,4\}}  \gammdelt{\alpha_k-\alpha_j}{\delta }\right]  \\&
 \hskip 80pt
 \cdot Q_{\delta }\biggl(   \frac{s_1+\alpha_1}{2}, \frac{\alpha_3-\alpha_4}{2}-\delta ,\frac{s_3-\alpha_2}{2};\frac{1+\alpha_1-\alpha_2}{2},\frac{s_1+\alpha_3}{2}-\delta ,\frac{s_3-\alpha_4}{2}-\delta  \biggr),\nn\end{align} 
and
\begin{align}\label{eq:s3res}&\Res_{s_3=\alpha_1-2\delta }{\Mellin{\alpha}(s)}  =\frac{1}{\gam{\textstyle\frac{  s_1+s_2-\alpha_1}{2}+\delta  } }\left[\; \prod_{k=2}^4\gammdelt{\alpha_1-\alpha_k}{\delta }\gamm{s_1+\alpha_k} 
\gamm{s_2-\alpha_1-\alpha_k}\right]  
\\&
\hskip 80pt
\cdot Q_{\delta }\biggl( \frac{s_1+\alpha_2}{2}, \frac{s_1+\alpha_3}{2}, \frac{s_1+\alpha_4}{2}; \frac{1 - \alpha_1 -s_2+s_1}{2}, \frac{s_1+\alpha_1 }{2}- \delta , \frac{s_1+s_2- \alpha_1}{2}  \biggr) \nn.\end{align} 
 \end{prop}
The formulas for the remaining residues are found from each of the above by permuting $\alpha$ (i.e., applying Weyl group transformations).

It may be seen from \cite[Equation (43)]{ST2018} that $Q_0$ is the constant polynomial $Q_0\equiv 1$.  Thus the case $\delta=0$ of the above proposition yields the following special cases, also deduced in \cite[Theorem 3.2]{Stade2001}:
 \begin{equation}\label{eq:M4residue1}
 \Res_{s_1=-\alpha_1}{\Mellin{\alpha}(s)} =  \frac{\prod\limits_{k=2}^4 \Gamma(\frac{\alpha_k-\alpha_1}{2}) \Gamma(\frac{s_2+\alpha_1+\alpha_k}{2}) \Gamma(\frac{s_3-\alpha_k}{2})}
{\Gamma(\frac{s_2+s_3+\alpha_1}{2})},
\end{equation}
\begin{equation}\label{eq:M4residue2}
\Res_{s_2=-\alpha_1-\alpha_4}{\Mellin{\alpha}(s)}
 = \left[\;\prod_{j\in\{1,4\}} \gamm{s_1+\alpha_j} \right]
 \left[\prod_{k=2}^3   \gamm{s_3-\alpha_j} \prod_{ j\in\{1,4\}}  \gamm{\alpha_k-\alpha_j}\right] ,
\end{equation}
and
\begin{equation}\label{eq:M4residue3}
 \Res_{s_3=\alpha_1}{\Mellin{\alpha}(s)} =  \frac{\prod\limits_{k=2}^4 \gamm{\alpha_1-\alpha_k}\gamm{s_1+\alpha_k} 
\gamm{s_2-\alpha_1-\alpha_k}}
{\Gamma(\frac{s_1+s_2-\alpha_1}{2})}.
\end{equation}

In Section~\ref{sec:doubleresidues} below, we will also need to consider double residues of $\Mellin{\alpha}(s)$ --- that is, residues in any one of the three $s_j$'s of residues in either of the others.  To this end we 
have the following, which is proved in \cite[Proposition 15]{ST2018}.
  \begin{prop}\label{resprop}
  \begin{enumerate}
   \item[\rm(a)] For each $\delta_1,\delta_2 \in\Z_{\ge0}$, there is a polynomial $f_{\delta_1 ,\delta_2 }(s_3,\alpha)$, of degree at most $2\delta_1 +\delta_2 $, such that
    \begin{align}\label{eq:s2s1res} &
   \Res_{s_2=-\alpha_1-\alpha_4-2\delta_2 }\biggl(\Res_{s_1=-\alpha_1-2\delta_1 }{ {\Mellin{\alpha}(s)}} \biggr) \\&= {\textstyle\gam{\frac{\alpha_4-\alpha_1}{2}-\delta_1 }}\left[\;\prod_{k=2 }^3\textstyle\gam{\frac{\alpha_k-\alpha_1}{2}-\delta _1}\gam{\frac{\alpha_k-\alpha_4}{2}-\delta_2 }\gam{\frac{s_3-\alpha_k}{2}}
\right]  f_{\delta_1 ,\delta_2 }(s_3,\alpha).\nn \end{align}
  \item[\rm(b)] For each $\delta_1 ,\delta_3\in\Z_{\ge0}$, there is a polynomial $g_{\delta_1 ,\delta_3}(s_2,\alpha)$, of degree at most $2\delta_1 +\delta_3$, such that
    \begin{align}\label{eq:s3s1res} &
\Res_{s_3=\alpha_2-2\delta_3}\biggl(\Res_{s_1=-\alpha_1-2\delta_1}{ {\Mellin{\alpha}(s)}} \biggr) \\&= {\textstyle\gam{\frac{\alpha_2-\alpha_1}{2}-\delta_1 }}\left[\;\prod_{k=3 }^4\textstyle\gam{\frac{\alpha_k-\alpha_1}{2}-\delta_1 }\gam{\frac{\alpha_2-\alpha_k}{2}-\delta_3}\gam{\frac{s_2+\alpha_1+\alpha_k}{2}}
\right]  g_{\delta_1 ,\delta_3}(s_2,\alpha) \nn. \end{align} \end{enumerate} \end{prop}

\subsection{\bf Shift equations}
We would like to have an expression similar to \eqref{eq:StadeMellinGL4} for $\Mellin{\alpha}(s)$ in the region where $\re(s)<0$.  Although we cannot use the right hand side of \eqref{eq:StadeMellinGL4} directly when $\re(s)<0$, we can use the fact that $\Mellin{\alpha}(s)$ satisfies the following shift equation which (as will be shown) is a direct corollary of Propositions~7 and 9 from \cite{ST2018}.

\begin{prop}\label{prop:GL4shifteq}
Let $s=(s_1,s_2,s_3)\in\C^3$.  Let $\alpha=(\alpha_1,\alpha_2,\alpha_3,\alpha_4)\in \C^4$ with $\alpha_1+\alpha_2+\alpha_3+\alpha_4=0$.  Suppose that $r_1,r_2,r_3\geq 0$ are integers.  Then
 \[ \left| \Mellin{\alpha}(s_1,s_2,s_3) \right| \ll \sum_{\ell=0}^{r_1}\sum_{k=0}^{r_2}\sum_{j=0}^{r_3} \left| \frac{Q_{j,k,\ell}^{r_1,r_2,r_3}}{\mathcal{B}_1^{r_1}\mathcal{B}_2^{r_2}\mathcal{B}_3^{r_3}} \Mellin{\alpha}\big(s_1+2(r_1+j+k),s_2+2r_2,s_3+2(r_3+\ell)\big)\right|, \]
where
\begin{align*}
 \mathcal{B}_1 := \mathcal{B}_1(\alpha,s) & := \big(s_1+\alpha_1\big) 
\big(s_1+\alpha_2\big) \big(s_1+\alpha_3\big) \big(s_1+\alpha_4\big), 
 \\
 \mathcal{B}_2 := \mathcal{B}_2(\alpha,s) & := \big(s_2+\alpha_1+\alpha_2\big) \big(s_2+\alpha_1+\alpha_3\big) \big(s_2+\alpha_1+\alpha_4\big)
 \\ & \qquad \cdot
 \big(s_2+\alpha_2+\alpha_3\big) \big(s_2+\alpha_2+\alpha_4\big) \big(s_2+\alpha_3+\alpha_4\big), 
 \\
 \mathcal{B}_3 := \mathcal{B}_3(\alpha,s) & := \big(s_3-\alpha_1\big) 
\big(s_3-\alpha_2\big) \big(s_3-\alpha_3\big) \big(s_3-\alpha_4\big), 
\end{align*}
and $Q_{j,k,\ell}^{r_1,r_2,r_3}$ is a polynomial in $\alpha$ and $s$ with 
combined degree $2(r_1+2r_2+r_3-j-k-\ell)$.
\end{prop}

\begin{proof}
We begin with Proposition~7(a) of \cite{ST2018} which states that
 \[ \Mellin{\alpha}(s_1,s_2,s_3) = \frac{q_0(\alpha,s)}{\mathcal{B}_1(\alpha,s_1)} \Mellin{\alpha}(s_1+2,s_2,s_3) + \frac{q_1(\alpha,s)}{\mathcal{B}_1(\alpha,s_1)} \Mellin{\alpha}(s_1+2,s_2,s_3+2),
 \]
where $\deg(q_i)=2-2i$.  Iterating this formula $r_1$ times gives
\begin{align} \label{eq:ShiftEquations1}
 \left| \Mellin{\alpha}(s_1,s_2,s_3) \right| \ll \sum_{\ell=0}^{r_1} \left| \frac{\mathcal{Q}_{1,r_1}^{(\ell)}(\alpha,s)}{\mathcal{B}_1(\alpha,s_1)^{r_1}} \Mellin{\alpha}(s_1+2r_1,s_2,s_3+2\ell)\right|.
\end{align}
In order to have an equality here, we would need to keep track of various 
shifts of the polynomial $\mathcal{B}_1(\alpha,s_1)$, but since we are only interested in bounds, the version presented here suffices.  It is easy 
to show that
 \[ \deg(\mathcal{Q}_{1,r_1}^{(\ell)})=(r_1-\ell)\deg(q_0)+\ell\deg(q_1) = 2(r_1-\ell). \]

Similarly, Proposition~9 of \cite{ST2018} states that
 \[ \Mellin{\alpha}(s_1,s_2,s_3) = \frac{p_0(\alpha,s)}{\mathcal{B}_2(\alpha,s_2)} \Mellin{\alpha}(s_1,s_2+2,s_3) + \frac{p_1(\alpha,s)}{\mathcal{B}_2(\alpha,s_2)} \Mellin{\alpha}(s_1+2,s_2+2,s_3),
 \]
where $\deg(p_k)=4-2k$.  Iterating this formula $r_2$ times gives
\begin{align} \label{eq:ShiftEquations2}
 \left| \Mellin{\alpha}(s_1,s_2,s_3) \right| \ll \sum_{k=0}^{r_2} \left| \frac{\mathcal{Q}_{2,r_2}^{(k)}(\alpha,s)}{\mathcal{B}_2(\alpha,s_2)^{r_2}} \Mellin{\alpha}(s_1+2k,s_2+2r_2,s_3)\right|,
\end{align}
and
 \[ \deg(\mathcal{Q}_{2,r_2}^{(k)})=(r_2-k)\deg(p_0)+k\deg(p_1) = 4(r_2-k)+2k = 4r_2-2k. \]

Via the change of variables $(s_1,s_3,\alpha)\mapsto (s_3,s_1,-\alpha)$ which preserves $\Mellin{\alpha}$ applied to \eqref{eq:ShiftEquations1}, we have
\begin{align} \label{eq:ShiftEquations3}
 \left| \Mellin{\alpha}(s_1,s_2,s_3) \right| \ll \sum_{j=0}^{r_3} \left| \frac{\mathcal{Q}_{1,r_3}^{(j)}(-\alpha,\overline{s})}{\mathcal{B}_3(\alpha,s_3)^{r_3}} \Mellin{\alpha}(s_1+2j,s_2,s_3+2r_3)\right|.
\end{align}

Applying equations \eqref{eq:ShiftEquations1}, \eqref{eq:ShiftEquations2} 
and \eqref{eq:ShiftEquations3} in succession gives the desired result with
 \[ \deg(Q_{j,k,\ell}^{r_1,r_2,r_3}) = \deg(\mathcal{Q}_{1,r_1}^{(\ell)}) + \deg(\mathcal{Q}_{2,r_2}^{(k)}) + \deg(\mathcal{Q}_{1,r_3}^{(j)}) = 
2(r_1+2r_2+r_3-\ell-k-j), \]
as claimed.
\end{proof}

\subsection{\bf Expressing $W_\alpha$ as the inverse Mellin transform of $\Mellin{\alpha}$}
Given the equation \eqref{eq:Mellin} for the Mellin transform of $\Whit{\alpha}$, we find by Mellin inversion that
\begin{equation}\label{eq:WhittasInt}
 W_\alpha(y) =\frac{1}{(2\pi i)^3} \iiint \limits_{\re(s)=u} y_1^{\frac32-s_1}y_2^{2-s_2}y_3^{\frac32-s_3} \Mellin{\alpha}(s)\, ds,
\end{equation}
provided that $u=(u_1,u_2,u_3)$ satisfies $u_j >0$ for $j=1,2,3$. 

As a matter of notation, for $u=(u_1,u_2,u_3)\in \R^3$ define
 \[ W_\alpha(y;u) :=  \frac{1}{(2\pi i)^3}\iiint\limits_{\re(s_j)=u_j} y_1^{\frac32-s_1}y_2^{2-s_2}y_3^{\frac32-s_3}\Mellin{\alpha}(s) \; ds, \]
 and if $I$ is a nonempty subset of $\{1,2,3\}$ of cardinality $r$ and $p_I=(p_i)_{i\in I}\in \C^r$ define 
 \[ \mathcal{R}_\alpha^{p_I}(y;u) := \frac{1}{(2\pi i)^{3-r}} \int\limits_{\substack{\re(s_j)=u_j\\ j\notin I}} \Res_{s_I=p_I}\left( y_1^{\frac32-s_1}y_2^{2-s_2}y_3^{\frac32-s_3}\Mellin{\alpha}(s) \right) \prod_{j\notin I} ds_j, \]
where $\Res_{s_I=p_I}$ is the operator which evaluates the iterated residue at each of the points $s_i=p_i$ with $i\in I$, i.e., 
 \[ \Res_{s_I=p_I}:= \Res_{s_{i_1}=p_{i_1}}\circ \cdots \circ \Res_{s_{i_q}=p_{i_q}}, \qquad\qquad \mbox{($I=(i_1,\ldots,i_q)$)}. \]
We call $\mathcal{R}_\alpha^{p_I}$ a \emph{single residue term} if $q=1$, a \emph{double residue term} if $q=2$, and a \emph{triple residue term} if $q=3$.  For $a=(a_1,a_2,a_3)$ with $a_j>0$, we call $W_\alpha(y;-a)$ the \emph{shifted main term}.

Note that if $u_1, u_2, u_3 > 0$, then $W_\alpha(y)=W_\alpha(y;u)$.  We 
now shift the lines of integration in the $s$-variable to the left passing poles at $\re(s_j)=0,-2,-4,\ldots$.  By the Cauchy Residue Theorem, this allows us to write $W_\alpha(y)$ in terms of $\mathcal{R}_\alpha^*$.  
For example, if $a=(a_1,a_2,a_3)$ with $0<a_j<2$ for each $j=1,2,3$, then 
\begin{align*}
 W_\alpha(y) & = W_\alpha(y;-a) + \sum_{\substack{p\in P_j \\ 1\leq j \leq 3}} \mathcal{R}_\alpha^{(p)}(y;-a) + \sum_{\substack{p\in P_j \\ q\in 
P_k \\ 1\leq j< k \leq 3}} \mathcal{R}_\alpha^{(p,q)}(y;-a) + \sum_{\substack{p\in P_1 \\ q\in P_2 \\ r\in P_3}} \mathcal{R}_\alpha^{(p,q,r)}(y;-a).
\end{align*}

Letting $p_1:=-\alpha_1-2\delta_1$, $p_2:=-\alpha_1-\alpha_2-2\delta_2$ and $p_3:=-\alpha_1-\alpha_2-\alpha_3-2\delta_3$ for some $\delta_1,\delta_2,\delta_3\in \Z_{\geq 0}$, as given in Section~\ref{poles-res}, there are polynomials $f_{\delta_j}^{p_j}(s,\alpha)$ such that 
\begin{equation} \label{eq:Rp1}
 \mathcal{R}_\alpha^{(p_1)}(y;-a) = \iint\limits_{\substack{\re(s_j)=-a_j\\j=2,3}}y_1^{\frac32+p_1}y_2^{2-s_2}y_3^{\frac32-s_3}\; \frac{\prod\limits_{j=2}^4 \Gamma(\frac{\alpha_j-\alpha_1}{2}) \Gamma(\frac{s_2+\alpha_1+\alpha_j}{2}) \Gamma(\frac{s_3-\alpha_j}{2})}
{\Gamma(\frac{s_2+s_3+\alpha_1}{2})}\; f_{\delta_1}^{p_1}(s,\alpha)\; ds_2 ds_3, 
\end{equation}
and
\begin{align} \label{eq:Rp2}
 \mathcal{R}_\alpha^{(p_2)}(y;-a) = \iint\limits_{\substack{\re(s_j)=-a_j\\j=1,3}} &
 y_1^{\frac32-s_1}y_2^{2+p_2}y_3^{\frac32-s_3}\; 
 \Biggl(\prod_{j=1}^2\prod_{k=3}^4 \textstyle{ \Gamma\bigl(\frac{\alpha_k-\alpha_j}{2}\bigr)}\Biggr) \; f_{\delta_2}^{p_2}(s,\alpha)
 \\ \nonumber & \qquad \cdot
 \textstyle{ \Gamma\(\frac{s_1+\alpha_1}{2}\) \Gamma\(\frac{s_1+\alpha_2}{2}\) \Gamma\(\frac{s_3-\alpha_3}{2}\) \Gamma\(\frac{s_3-\alpha_4}{2}\)}\; ds_1ds_3.
\end{align}
All other single residue terms can be obtained from these by observing that, first, the action of the Weyl group, which acts by permuting the set $(\alpha_1,\alpha_2,\alpha_3,\alpha_4)$, acts transitively on each of the 
sets $P_1,P_2,P_3$ and permutes the residue terms accordingly.  Explicitly, if $\sigma=(12)$ represents the permutation interchanging $1$ and $2$, then $\mathcal{R}_\alpha^{{(\sigma(p_1)-\delta_1)}}(y;-a)=\mathcal{R}_{\sigma(\alpha)}^{{(p_1-\delta_1)}}(y;-a)$ has the same expression as $\mathcal{R}_\alpha^{{(p_1-\delta_1)}}(y;-a)$ except that each instance of 
$\alpha_1$ and $\alpha_2$ are interchanged.

Moreover, the involution $(\alpha,s_1,s_2,s_3)\mapsto (-\alpha,s_3,s_2,s_1)$ interchanges the set of single residues at, say, $s_1=-\alpha_1-2\delta$ with those at $s_3=\alpha_1-2\delta$ respecting the formulas of Proposition~\ref{resprop}.  So, in particular, we have $\mathcal{R}_{\alpha}^{{(\alpha_1-2\delta)}}(y_1,y_2,y_3;-a)= \mathcal{R}_{-\alpha}^{{(-\alpha_1-2\delta)}}(y_3,y_2,y_1;-a)$. 

It can be similarly shown that for appropriate polynomials $f_{*,*}^{(*,*)}$ and $f_\delta^{(p_1,p_2,p_3)}$, every double residue term is equivalent to either
\begin{align} \label{eq:Rp1p2}
 \mathcal{R}_\alpha^{(p_1,p_2)}(y;-a) = \int\limits_{\re(s_3)=-a_3} & 
y_1^{\frac32+p_1}y_2^{2+p_2}y_3^{\frac32-s_3}
f_{\delta_1,\delta_2}^{(p_1,p_2)}(s_3,\alpha) \\ \nonumber & \cdot
 \Biggl(\prod_{j=1}^2 \prod_{k=j+1}^4 \textstyle{\Gamma\bigl(\frac{\alpha_k-\alpha_j}{2} \bigr)} \Biggr) 
 \textstyle{\Gamma\(\frac{s_3-\alpha_4}{2}\) \Gamma\(\frac{s_3-\alpha_3}{2}\)}\;  ds_3, 
\end{align}
or
\begin{align} \label{eq:Rp1p3}
 \mathcal{R}_\alpha^{(p_1,p_3)}(y;-a) = \int\limits_{\re(s_2)=-a_3} & 
y_1^{\frac32+p_1}y_2^{2-s_2}y_3^{\frac32+p_3}\; f_{\delta_1,\delta_3}^{(p_1,p_3)}(s_2,\alpha) \; \textstyle{\Gamma\(\frac{s_2+\alpha_1+\alpha_2}{2}\) \Gamma\(\frac{s_2+\alpha_1+\alpha_3}{2}\)}
 \\ \nonumber & \qquad \cdot
 \textstyle{ \Gamma\(\frac{\alpha_4-\alpha_3}{2}\) \Gamma\(\frac{\alpha_4-\alpha_2}{2}\) \Gamma\(\frac{\alpha_4-\alpha_1}{2}\) \Gamma\(\frac{\alpha_3-\alpha_1}{2}\) \Gamma\(\frac{\alpha_2-\alpha_1}{2}\)} \; ds_2,
\end{align}
and every triple residue term is equivalent to
\begin{equation} \label{eq:Rp1p2p3}
 \mathcal{R}_\alpha^{(p_1,p_2,p_3)}(y) = y_1^{3/2+p_1}y_2^{2+p_2}y_3^{3/2+p_3} \; f_\delta^{(p_1,p_2,p_3)}(\alpha) \prod_{1\leq j<k \leq 4}  \textstyle{\Gamma\bigl(\frac{\alpha_k-\alpha_j}{2}\bigr)}.
\end{equation}

It turns out that the bounds obtained from applying our methods only to the shifted integral $W_\alpha(y;-a)$ when $0<a_j<2$ are not sufficient for our needs.  Instead, we will need to obtain bounds for various shifts of the form $2r_j-1<a_j<2r_j$ for $r_1,r_2,r_3\in \Z_{\geq 0}$.  For a given choice of $(r_1,r_2,r_3)\in \Z_{\geq 0}^3$, we extrapolate the procedure outlined above obtaining poles corresponding to each case of 
 \[ s_j=p_j - 2\delta_j, \quad \mbox{ for each } \quad \delta_j=0,\ldots,r_j-1. \]
Doing so, we arrive at the following result.

%\pagebreak

\begin{prop}\label{prop:WhittakerDecomposition} Let $(r_1,r_2,r_3)\in \Z_{\geq 0}^3$ be fixed.
Let $a = (a_1,a_2,a_3)$ satisfy $2r_j-2<a_j<2r_j$ be fixed ($a_j\neq a_k$ if $j\neq k$). 
Let 
 \[ P_1=\{-\alpha_j\mid j=1,\ldots,4\}, \quad P_2=\{-\alpha_j-\alpha_k\mid 1 \leq j\neq k \leq 4\}, \quad P_3=\{\alpha_j\mid j=1,\ldots,4\}, \] as in \eqref{Pjdef}.
 Then, after shifting contours of integration, the $\GL(4)$-Whittaker function  is given by 
\begin{align*}
&W_\alpha(y)=  W_\alpha(y;-a) 
  +\hskip-12pt \sum_{\substack{p_j\in P_j-2\delta_j\\j=1,2,3\\ \delta_j\in\{0,1, \ldots, r_j-1\}}}\hskip-5pt \mathcal{R}_\alpha^{(p_j)}(y; -a) 
 \;\; +\hskip-12pt \sum_{\substack{1\leq j \neq k\leq 4 \\ p_j\in P_j-2\delta_j\\p_k\in P_k-2\delta_k \\ \delta_j\in\{0,1,\ldots, r_j-1\} \\ \delta_k\in\{0,1,\ldots, r_k-1\}}}\hskip-10pt \mathcal{R}_\alpha^{(p_j,p_k)}(y;-a) 
 \;\;+ \hskip-17pt\sum_{\substack{p_1\in P_1-2\delta_1\\ p_2\in P_2-2\delta_2 \\ p_3\in P_3-2\delta_3 \\ \delta_\ell\in\{0,1,\ldots,r_\ell-1\} \\ \ell=1,2,3}}\hskip-15pt \mathcal{R}_\alpha^{(p_1,p_2,p_3)}(y). \end{align*}
\end{prop}

\section{\bf Bounds for the test function $p_{T,R}$}\label{sec:boundsforpTR}

\vskip 10pt
Recall (see (\ref{eq:testfunction2})) that $p_{T,R}$ is given by
\begin{equation} \label{eq:testfunction3}
p_{T,R}(y) = p_{T,R}(y_1,y_2,y_3) = \frac{1}{\pi^{3}}
  \underset{\re(\alpha_j)=0}{\iiint} p^{\#}_{T,R}(\alpha) \;W_{\alpha}(y)
     \frac{ d\alpha_1\, d\alpha_2\, d\alpha_3
         }{\prod\limits_{1\leq\, j \ne k\, \leq 4} \Gamma\bigl(\frac{\alpha_j-\alpha_k}{2}\bigr) }.
\end{equation}
Note that in \cite{GK2012} instead of $W_\alpha(y)$ one actually has its complex conjugate $\overline{W}_{\alpha}(y)$.  However, one arrives at the formula here by noting that  $\alpha$ is purely imaginary and, therefore,  $\overline{W}_{\alpha}(y)=W_{-\alpha}(y)$.  Hence,  the change of variables $\alpha\mapsto -\alpha$ which leaves $p_{T,R}^\sharp(\alpha)$ and the measure ${ d\alpha_1\, d\alpha_2\, d\alpha_{3}}\Big/{\prod\limits_{1\leq\, j \ne k\, \leq 4} \Gamma\bigl(\frac{\alpha_j-\alpha_k}{2}\bigr) }$ invariant, leads to the given formula.

\subsection{\bf Decomposition of $p_{T,R}$ in terms of poles and residues 
of $\Mellin{\alpha}$}
Define \begin{equation}\label{gammaRalpha}\Gamma_R(\alpha):=\prod_{1\leq \ell\neq m \leq 4}\frac{\Gamma(\frac{2+R+\alpha_\ell-\alpha_m}{4})}{\Gamma(\frac{\alpha_\ell-\alpha_m}{2})}.\end{equation} We now replace $W_\alpha(y)$, on the right side of \eqref{eq:testfunction3}, by the expression 
given in Proposition~\ref{prop:WhittakerDecomposition}.  It follows from the definition of the test function $p_{T,R}^\#$ given in \eqref{testfunctionsharp2}, that in doing so, we obtain a shifted $p_{T,R}$  term
\begin{align} \label{eq:pTRemptyset}
 p_{T,R}(y; -a):= \iiint\limits_{\re(\alpha_j)=0} e^{\frac{\alpha_1^2+\cdots+\alpha_4^2}{2T^2}}&\; W_\alpha(s;-a)\,\mathcal F_R(\alpha)\,
  \Gamma_R(\alpha)\, d\alpha,
\end{align}
single residue terms of the type
\begin{align*}
 p_{T,R}^{j,\delta}(y) := \iiint\limits_{\re(\alpha_j)=0} e^{\frac{\alpha_1^2+\cdots+\alpha_4^2}{2T^2} }\,\mathcal F_R(\alpha)\,{ \mathcal{R}_\alpha^{(p_j-2\delta)}(y;-a)}\, \Gamma_R(\alpha) \, d\alpha,
\end{align*}
double residue terms of the type
\begin{align*}
 p_{T,R}^{jk,\delta}(y) := \iiint\limits_{\re(\alpha_j)=0} e^{\frac{\alpha_1^2+\cdots+\alpha_4^2}{2T^2}} \, \mathcal{R}_\alpha^{(p_j-2\delta_j,p_k-2\delta_k)}(y;-a)\,\mathcal F_R(\alpha)\, \Gamma_R(\alpha) \, d\alpha,
\end{align*}
and triple residue terms of the type
\begin{align*}
 p_{T,R}^{123,\delta}(y) := \iiint\limits_{\re(\alpha_j)=0} e^{\frac{\alpha_1^2+\cdots+\alpha_4^2}{2T^2}}\,\mathcal{R}_\alpha^{(p_1-2\delta_1,p_2-2\delta_2,p_3-2\delta_3)}(y;-a)\,\,\mathcal F_R(\alpha)\,
 \Gamma_R(\alpha) \, d\alpha,
\end{align*}
where $p_1=-\alpha_1$, $p_2=-\alpha_1-\alpha_2$, $p_3=-\alpha_1-\alpha_2-\alpha_3$ and $\delta_j\in\{0,1,\ldots, r-1\}$.  Also, note that we 
use the notation $d\alpha:=d\alpha_1\,d\alpha_2\,d\alpha_3$.  In particular, by equations~\eqref{eq:Rp1}--\eqref{eq:Rp1p2p3}, we see that
\begin{align} \label{eq:pTr_1}
  p_{T,R}^{1,0}(y) = & \iiint\limits_{\re(\alpha_j)=0} e^{\frac{\alpha_1^2+\cdots+\alpha_4^2}{2T^2}}\; \iint\limits_{\substack{\re(s_j)=-a_j\\j=2,3}}y_1^{\frac32+p_1}y_2^{2-s_2}y_3^{\frac32-s_3}\,\mathcal F_R(\alpha) \, \Gamma_R(\alpha)
  \\ \nonumber & \hskip 48pt
  \cdot\frac{ \prod\limits_{j=2}^4 \Gamma(\frac{\alpha_j-\alpha_1}{2}) \Gamma(\frac{s_2+\alpha_1+\alpha_j}{2})  \Gamma(\frac{s_3-\alpha_j}{2}) }{ \Gamma(\frac{s_2+s_3+\alpha_1}{2}) } \; ds_2 ds_3\; d\alpha,
\end{align}
\begin{align} \label{eq:pTr_2}
 p_{T,R}^{2,0}(y) = & \iiint\limits_{\re(\alpha_j)=0} e^{\frac{\alpha_1^2+\cdots+\alpha_4^2}{2T^2}}\; \iint\limits_{\substack{\re(s_j)=-a_j\\j=1,3}} 
 y_1^{\frac32-s_1}y_2^{2+p_2}y_3^{\frac32-s_3}\; 
 \Biggl(\prod_{j=1}^2\prod_{k=3}^4 \textstyle{ \Gamma\bigl(\frac{\alpha_k-\alpha_j}{2}\bigr)}\Biggr)\mathcal F_R(\alpha)\, \Gamma_R(\alpha)
 \\ \nonumber & \hskip 96pt
 \cdot
 \textstyle{ \Gamma\(\frac{s_1+\alpha_1}{2}\) \Gamma\(\frac{s_1+\alpha_2}{2}\) \Gamma\(\frac{s_3-\alpha_3}{2}\) \Gamma\(\frac{s_3-\alpha_4}{2}\)}\; ds_1ds_3\; d\alpha,
\end{align}
\begin{align} \label{eq:pTr_12}
 p_{T,R}^{12,0}(y) = & \iiint\limits_{\re(\alpha_j)=0} e^{\frac{\alpha_1^2+\cdots+\alpha_4^2}{2T^2}}\; \int\limits_{\re(s_3)=-a_3} y_1^{\frac32+p_1}y_2^{2+p_2}y_3^{\frac32-s_3} \, \mathcal F_R(\alpha)\,\Gamma_R(\alpha)
 \\ \nonumber & \hskip 60pt
 \cdot \Biggl(\prod_{j=1}^2 \prod_{k=j+1}^4 \textstyle{\Gamma\bigl(\frac{\alpha_k-\alpha_j}{2} \bigr)} \Biggr) \textstyle{\Gamma\(\frac{s_3-\alpha_4}{2}\) \Gamma\(\frac{s_3-\alpha_3}{2}\)}\; ds_3 \; d\alpha,
\end{align}
\begin{align} \label{eq:pTr_13}
 p_{T,R}^{13,0}(y) = & \iiint\limits_{\re(\alpha_j)=0} e^{\frac{\alpha_1^2+\cdots+\alpha_4^2}{2T^2}}\; \int\limits_{\re(s_2)=-a_3} y_1^{\frac32+p_1}y_2^{2-s_2}y_3^{\frac32+p_3}\; \textstyle{\Gamma\(\frac{s_2+\alpha_1+\alpha_2}{2}\) \Gamma\(\frac{s_2+\alpha_1+\alpha_3}{2}\)}
 \\ \nonumber & \hskip 24pt \cdot
 \textstyle{ \Gamma\(\frac{\alpha_4-\alpha_3}{2}\) \Gamma\(\frac{\alpha_4-\alpha_2}{2}\) \Gamma\(\frac{\alpha_4-\alpha_1}{2}\) \Gamma\(\frac{\alpha_3-\alpha_1}{2}\) \Gamma\(\frac{\alpha_2-\alpha_1}{2}\)} \; ds_2\,\mathcal F_R(\alpha)\,\Gamma_R(\alpha) \; d\alpha,
\end{align}
\begin{align} \label{eq:pTr_123}
 p_{T,R}^{123,0}(y) = & \iiint\limits_{\re(\alpha_j)=0} e^{\frac{\alpha_1^2+\cdots+\alpha_4^2}{2T^2}}\; y_1^{3/2+p_1}y_2^{2+p_2}y_3^{3/2+p_3} \,
 \Biggl(\prod_{1\leq j<k \leq 4} \textstyle{\Gamma\bigl(\frac{\alpha_k-\alpha_j}{2}\bigr)}\Biggr)\,\mathcal F_R(\alpha)\, \Gamma_R(\alpha) \; d\alpha.
 \end{align}

Since we are integrating over $\alpha$ on the right hand side of \eqref{eq:testfunction2} and the integrand (and measure) is invariant under the action of the Weyl group, there are explicitly computable constants $c_1,c_{12},c_{123}$ such that

\pagebreak

\begin{align}\label{eq:pTRdecomposition}
p_{T,R}(y) &= p_{T,R}(y;-a) \;\, + \;\, c_1\hskip-15pt\sum_{\substack{1\leq j \leq 3\\\delta\in \{0,1,\ldots,r-1\}}} \hskip-5pt p_{T,R}^{j,\delta}(y;-a) \;\, + \;\, c_{12}\hskip-15pt\sum_{\substack{1\leq j< k \leq 3\\ 
\delta\in \{0,1,\ldots,r-1\}^2}} \hskip-5pt p_{T,R}^{jk,\delta}(y;-a)  \\ 
\nn&
\hskip 250pt
 +\;\, c_{123}\hskip-15pt\sum_{\substack{\delta\in \{0,1,\ldots,r-1\}^3}}\hskip-5pt p_{T,R}^{123,\delta}(y;-a).
\end{align}

\begin{rmrk}\label{rmk:generalposition}
There are no poles in the expressions of the residues given in (6.1.1) to 
(6.1.6) because the products of $\Gamma\big(\frac{\alpha_k-\alpha_j}{2}  \big)$ are cancelled by the Gamma functions in the denominator of $\Gamma_R(\alpha).$
\end{rmrk}

\vskip 10pt

\subsection{\bf Philosophy of the proof of Theorem \ref{th:pTRbound}}\label{sec:pTRBound}

Recall that Theorem \ref{th:pTRbound} gives sharp bounds for the $p_{T,R}$ function.Given \eqref{eq:pTRdecomposition}, this can be achieved by bounding each of the terms on the right hand side: $p_{T,R}(y;-a)$, $p_{T,R}^{1,\delta}(y;-a)$, $p_{T,R}^{2,\delta}(y;-a)$, $p_{T,R}^{12,\delta}(y;-a)$, $p_{T,R}^{23,\delta}(y;-a)$ and $p_{T,R}^{123,\delta}(y;-a)$ and each 
possible choice of $\delta$.  We obtain bounds for each of these terms in 
Sections~\ref{sec:shiftedmaintermbound}--\ref{sec:tripleresbound}, respectively. Theorem~\ref{th:pTRbound} then follows as an immediate consequence.

Before proceeding to those sections, we describe the idea in the special case that $a=(a_1,a_2,a_3)$ with $a_i<0$ for each $i=1,2,3$.  In this 
case there are no residue terms at all, i.e., we need only bound $p_{T,R}(y)=p_{T,R}(y;-a)$.  Then the basic idea is to insert the formula \eqref{eq:StadeMellinGL4} into \eqref{eq:WhittasInt} to get an expression for $W_\alpha(y)$ as an integral of the ratio of many Gamma functions.  We then insert this into \eqref{eq:testfunction2} and estimate each Gamma function using Stirling's approximation, which for fixed $\sigma$ and $\abs{t}\to\infty$ says that
\begin{equation} \label{eq:Stirling}
 \Gamma(\sigma+it)\sim \sqrt{2\pi} \lvert t \rvert^{\sigma-\frac12}\, e^{-\frac{\pi}{2}\abs{t}}.
\end{equation}
We call $\lvert t\rvert^{\sigma-\frac12}$ the \emph{polynomial factor} of 
$\Gamma(\sigma+it)$, and $e^{-\frac{\pi}{2}\abs{t}}$ is called the \emph{exponential factor}.  Through this process, we can replace the ratio of all of the Gamma factors by a rational function $\mathcal{P}(s,a)$   obtained as the product of the polynomial factors of the individual Gamma functions times $e^{-\mathcal{E}(s,\alpha)}$ which is the product of the exponential factors.

To be completely explicit, after making a simple change of variables in \eqref{eq:StadeMellinGL4}, it is easy to see that (up to a constant) for $\varepsilon'>0$ sufficiently small
\begin{equation}\label{eq:MellinGL4}
 \Mellin{\alpha}(s) = \int\limits_{\re(t)=-\varepsilon'}\frac{\Gamma_0(t,\alpha)\Gamma_1(t,s_1,\alpha) \Gamma_2(t,s_2,\alpha) \Gamma_3(t,s_3,\alpha)}{ \Gamma_{\mathrm den}(t,s,\alpha)}\, dt 
\end{equation}
where
\begin{align*}
 \Gamma_0(t,\alpha) & := \Gamma\Bigl( \frac{-t+\alpha_3}{2}\Bigr) \Gamma\Bigl( \frac{-t+\alpha_4}{2}\Bigr), \\
 \Gamma_1(t,s_1,\alpha) & := \Gamma\Bigl(\frac{s_1+\alpha_1}{2}\Bigr)  \Gamma\Bigl(\frac{s_1+\alpha_2}{2}\Bigr) \Gamma\Bigl(\frac{s_1+t}{2}\Bigr), \\
 \Gamma_2(t,s_2,\alpha) & := \Gamma\Bigl(\frac{s_2+\alpha_1+\alpha_2}{2}\Bigr)  \Gamma\Bigl(\frac{s_2+\alpha_1+t}{2}\Bigr) \Gamma\Bigl(\frac{s_2+\alpha_2+t}{2}\Bigr) \Gamma\Bigl(\frac{s_2+\alpha_3+\alpha_4}{2}\Bigr), \\
 \Gamma_3(t,s_3,\alpha) & := \Gamma\Bigl(\frac{s_3+\alpha_1+\alpha_2+t}{2}\Bigr)  \Gamma\Bigl(\frac{s_3+\alpha_1+\alpha_3+\alpha_4}{2}\Bigr) \Gamma\Bigl(\frac{s_3+\alpha_2+\alpha_3+\alpha_4}{2}\Bigr),\\
 \Gamma_{\mathrm den}(t,s,\alpha) & := \Gamma\Bigl(\frac{s_1+s_2+\alpha_1+\alpha_2+t}{2}\Bigr) \Gamma\Bigl(\frac{s_2+s_3+t}{2}\Bigr).
\end{align*}
Thus by combining \eqref{eq:WhittasInt} and \eqref{eq:MellinGL4} into \eqref{eq:testfunction2} as described above, we find that
\begin{align} \label{eq:pTRasGammas} 
 p_{T,R}(y) = & \iiint\limits_{\re(\alpha_j)=0} e^{\frac{\alpha_1^2+\alpha_2^2+\alpha_3^2+\alpha_4^2}{2T^2}} \iiint\limits_{\re(s)=\varepsilon} \int\limits_{\re(t)=-\varepsilon'} y_1^{\frac32-s_1}y_2^{2-s_2}y_3^{\frac32-s_3} \\ \nonumber
 & \hskip 36pt \cdot \frac{\Gamma_0(t,\alpha)\Gamma_1(t,s_1,\alpha)\Gamma_2(t,s_2,\alpha)\Gamma_3(t,s_3,\alpha)}{\Gamma_{\mathrm den}(t,s_1,s_2,s_3,\alpha)}\,\mathcal F_R(\alpha)\,\Gamma_R(\alpha)\;dt\; ds\, d\alpha.
\end{align}
Applying Stirling's bound to each of the Gamma functions in \eqref{eq:pTRasGammas} we see that (up to a constant factor depending at most on $R, \,\varepsilon$)
 \begin{align*}&
 p_{T,R}(y)   \ll y_1^{\frac32+a_1}y_2^{2+a_2}y_3^{\frac32+a_3} \iiint\limits_{\mathcal{R}_T(0)} \cdot \iiint\limits_{\re(s)=-a} \cdot \int\limits_{\re(t)=\varepsilon'} \mathcal{P}(s,\alpha)  \, \mathcal F_R(\alpha)\,\exp\left(-\frac{\pi}{4}\mathcal{E}(s,\alpha)\right) dt\; ds\; d\alpha, 
\end{align*}
where, for $\alpha_j=\kappa_j+i\tau_j$ ($j=1,2,3,4$),
\[
 \mathcal{R}_T(\kappa) := \left\{  (i\tau_1+\kappa_1,i\tau_2+\kappa_2,i\tau_3) \mid -\tau_1-\tau_2-\tau_3 \leq \tau_3\leq \tau_2\leq \tau_1 \leq 
T^{1+\varepsilon} \right\}.
\]
The Weyl group is isomorphic to $S_4$ and acts by permutations on the set 
$\{\alpha_j\}_{j=1}^{4}$ leaving the integrand for $p_{T,R}(y)$ invariant.  Hence it suffices to restrict the integration over $\alpha$ to the set $\mathcal{R}_T(0)$.  

We will prove below (see Lemma~\ref{lem:ExpZeroSet}) that the integration 
in $s$ and $t$ can also be restricted to a finite volume set $\mathcal{R}$ which we call the \emph{exponential zero set}.  As $s$ and $t$ vary within this set, most of the polynomial terms can be uniformly bounded by a power of something of the form $(1+\tau_k-\tau_j)$ where $j<k$.  We prove 
a very strong bound on the remaining terms (see Lemma~\ref{lemmaIntegralBound3}) which shows that it too is bounded by a product of similar factors.  This implies that
 \[ \iiiint\limits_{(s,t)\in \mathcal{R}} \mathcal{P}(s,\alpha)\; dt\, ds 
\ll \prod_{1\leq j<k\leq 4} \big( 1+\tau_j-\tau_k\big)^{b_{j,k}} \]
where each $b_{j,k}>0$.  Then the integration in $\alpha$ over the set $\mathcal{R}_T(0)$ can be estimated trivially.  The main difference between 
this outline and the actual proof is that instead of using \eqref{eq:MellinGL4} directly as above, we replace it by the expression on the right hand side of the result of Proposition~\ref{prop:GL4shifteq}.  Also, instead of dealing with $p_{T,R}$ itself as given in \eqref{eq:testfunction3}, we individually bound each term on the right hand side of \eqref{eq:pTRdecomposition}.

In order to describe the exponential zero set, note that since $\tau_1 \geq \tau_2 \geq \tau_3 \geq \tau_4$, if we let $\im(s_j)=:\xi_j$ and $\im(t_1)=:\rho$, the exponential term takes the simplified form
\begin{align*}
  \mathcal{E} = \mathcal{E}(s,\alpha) = & \; -6\tau_1-4\tau_2-2\tau_3
  %% Gamma_0 terms
  +\lvert \rho-\tau_4\rvert + \lvert \rho-\tau_3\rvert
  %% Gamma_1 terms
  +\lvert \xi_1+\tau_1\rvert +\lvert \xi_1+\tau_2\rvert +\lvert \xi_1+\rho\rvert\\ & \qquad
  %% Gamma_2 terms
  +\lvert \xi_2+\tau_1+\tau_2)\rvert +\lvert \xi_2+\tau_1+\rho)\rvert +\lvert \xi_2+\tau_2+\rho)\rvert +\lvert \xi_2+\tau_3+\tau_4)\rvert \\ & \qquad \qquad
  %% Gamma_3 terms
  +\lvert \xi_3+\tau_1+\tau_2+\rho\rvert + \lvert \xi_3+\tau_1+\tau_3+\tau_4\rvert+ \lvert \xi_3+\tau_2+\tau_3+\tau_4\rvert \\
   %% Gamma_den terms
& \qquad \qquad \qquad - \lvert \xi_1+\xi_2+\tau_1+\tau_2+\rho\rvert - \lvert \rho+\xi_2+\xi_3\rvert.
\end{align*}

An important observation is the fact that $p_{T,R}(y)$ has been defined in such a way that there is at worst polynomial growth in the integrand.  This means that the exponential factor is never negative, i.e., $\mathcal{E}(s,\alpha)\geq 0$ for all $s$.  Since there will be exponential decay for any choice of $s$ such that $\mathcal{E}(s,\alpha)>0$, for the purposes of bounding $p_{T,R}(y)$, we need only determine when $\mathcal{E} = 
0$.  Note that this set depends only on the imaginary parts of $s$ and $\alpha$.

The expression for $\mathcal{E}$ above involves 14 absolute value terms. We can remove each of the 14 absolute values by replacing $|x|$ with $\pm 
x$ depending on whether $x$ is positive or negative. This leads to an expression of the form
\begin{align} \label{eq:AbsValSum}
0 = & -6\tau_1-4\tau_2-2\tau_3
  %% Gamma_0 terms
  +\varepsilon_{t,1}(\rho-\tau_4) + \varepsilon_{t,2} (\rho-\tau_3)
  %% Gamma_1 terms
  +\varepsilon_{1,0} (\xi_1+\tau_1) +\varepsilon_{1,1} (\xi_1+\tau_2) +\varepsilon_{1,2} (\xi_1+\rho)\\ \nonumber & \qquad
  %% Gamma_2 terms
  +\varepsilon_{2,0} (\xi_2+\tau_1+\tau_2) +\varepsilon_{2,1} (\xi_2+\tau_1+\rho) +\varepsilon_{2,2} (\xi_2+\tau_2+\rho) +\varepsilon_{2,3} (\xi_2+\tau_3+\tau_4) \\ \nonumber & \qquad \qquad
  %% Gamma_3 terms
  +\varepsilon_{3,0} (\xi_3+\tau_1+\tau_2+\rho) + \varepsilon_{3,1} (\xi_3+\tau_1+\tau_3+\tau_4)+ \varepsilon_{3,2} (\xi_3+\tau_2+\tau_3+\tau_4) \\ \nonumber
   %% Gamma_den terms
& \qquad \qquad \qquad - \varepsilon_{2-\frac12} (\xi_1+\xi_2+\tau_1+\tau_2+\rho) - \varepsilon_{2+\frac12} (\rho+\xi_2+\xi_3),
\end{align}
where each of the 14 $\varepsilon$'s is equal to $\pm1$. For a particular 
choice of $\varepsilon$'s either the sum on the right hand side of \eqref{eq:AbsValSum} vanishes identically or not. If it does vanish, each $\varepsilon_*$ determines an inequality, and the set of $\tau_1,\tau_2,\tau_3,\rho,\xi_1,\xi_2,\xi_3$ which satisfy all of these inequalities simultaneously is contained in the exponential zero set.  The following lemma shows that there are three such choices of signs and each choice explicitly determines an exponential zero set.  

\begin{lem}\label{lem:ExpZeroSet}
Every solution $\varepsilon=(\varepsilon_{t,1},\varepsilon_{t,2},\ldots,\varepsilon_{2-\frac12},\varepsilon_{2+\frac12})\in (\pm 1)^{14}$ to \eqref{eq:AbsValSum} is of the form
 \begin{gather}
  \varepsilon_{t,1}=+1, \quad \varepsilon_{t,2}=+1, \label{eq:ExpZeroSetLemmat} \\
  \varepsilon_{1,0}=+1, \quad \varepsilon_{1,1}=\varepsilon_{2-\frac12}, \quad \varepsilon_{1,2}=-1, \label{eq:ExpZeroSetLemma1} \\
  \varepsilon_{2,0}=+1, \quad \varepsilon_{2,1}=\varepsilon_{2-\frac12}, \quad \varepsilon_{2,2}=\varepsilon_{2+\frac12}, \quad \varepsilon_{2,3}=+1,\label{eq:ExpZeroSetLemma2}\\ \varepsilon_{3,0}=+1, \quad\varepsilon_{3,1}=\varepsilon_{2+\frac12}, \quad \varepsilon_{3,2}=-1, \label{eq:ExpZeroSetLemma3} 
 \end{gather}
and $\varepsilon_{2-\frac12}\geq \varepsilon_{2+\frac12}$.

In particular, there are three possible exponential zero sets, which we denote as $\mathcal R_1, \mathcal R_2, \mathcal R_3$.  The first corresponds to the case of $(\varepsilon_{2-\frac12},\varepsilon_{2+\frac12})=(+1,+1)$:
\begin{align*}
\hskip -244pt  \mathcal R_1: \qquad\qquad\quad
\tau_4 \; \leq \; \rho & \; \leq \; \tau_3
\ \phantom{SOLUTION 1:}\\
 -\tau_2 \;\leq \; \xi_1 & \; \leq \;-\rho, \\
 -\tau_2 -\rho \;  \leq \; \xi_2 & \; \leq \; \tau_1+\tau_2  \\
 \tau_2 \; \leq \; \xi_3 & \; \leq \; \tau_1,
\end{align*}
the second corresponds to $(\varepsilon_{2-\frac12},\varepsilon_{2+\frac12})=(+1,-1)$:
\begin{align*}
\hskip -244pt \mathcal R_2:  \qquad\qquad\qquad
\tau_4 \; \leq \; \rho  & \; \leq \; \tau_3 
\ \phantom{SOLUTION 2:}\\
 -\tau_2 \;\leq \; \xi_1 & \; \leq \;-\rho, \\
 -\tau_1 -\rho \;  \leq \; \xi_2 & \; \leq \; -\tau_2 -\rho  \\
 -(\tau_1+\tau_2) -\rho \; \leq \; \xi_3 & \; \leq \; \tau_2,
\end{align*}
and the third corresponds to $(\varepsilon_{2-\frac12},\varepsilon_{2+\frac12})=(-1,-1)$:
\begin{align*}
\hskip -244pt \mathcal R_3: \qquad\qquad\qquad
\tau_4 \; \leq \; \rho  & \; \leq \; \tau_3 
\ \phantom{SOLUTION 3:}\\
 -\tau_1 \;\leq \; \xi_1 & \; \leq \; -\tau_2, \\
-(\tau_1+\tau_2) \;  \leq \; \xi_2 & \; \leq \; -\tau_1 -\rho  \\
 -(\tau_1+\tau_2) -\rho \; \leq \; \xi_3 & \; \leq \; \tau_2.
\end{align*}
\end{lem}

\begin{proof}
Suppose that $\varepsilon=(\varepsilon_{t,1},\varepsilon_{t,2},\ldots,\varepsilon_{2+\frac12})\in (\pm 1)^{14}$ is a solution to \eqref{eq:AbsValSum}.  If we replace every instance of $\tau_4$ in \eqref{eq:AbsValSum} with $-\tau_1-\tau_2-\tau_3$, notice that the coefficient of $\tau_3$ is $(\varepsilon_{t,1}-\varepsilon_{t,2}-2)$.  This immediately implies \eqref{eq:ExpZeroSetLemmat}, or equivalently, $\tau_4\leq \rho\leq \tau_3$.  

Recall that we are assuming that $\tau_1\geq \tau_2\geq \tau_3\geq \tau_4$.  This, together with the fact that $\tau_4\leq \rho\leq \tau_3$, implies that
 \[ \xi_1+\tau_1\geq \xi_1+\tau_2 \geq \xi_1+\rho. \]
Since $\varepsilon_{1,0}=-1$ implies that $(\xi_1+\tau_1)\leq 0$, it follows that this would also imply that $\varepsilon_{1,1}=\varepsilon_{1,2}=-1$.  But it can't be the case that all three $\varepsilon_{1,k}$ are $-1$ because if so, the coefficient of $\xi_1$ in \eqref{eq:AbsValSum} 
will not be zero.  The same argument implies that $\varepsilon_{1,3}=-1$, and that the same relations hold for $\varepsilon_{3,j}$.  Similarly, $+1=\varepsilon_{2,0}\geq\varepsilon_{2,1}\geq\varepsilon_{2,2}\geq \varepsilon_{2,3}=-1$.

Using this information, we now rewrite \eqref{eq:AbsValSum} as
\begin{align}\label{eq:AbsValSum2}
0 = & -5\tau_1-3\tau_2
  %% Gamma_1 terms
  + (\tau_1-\rho) +\varepsilon_{1,1} (\xi_1+\tau_2) 
  \\ \nonumber & \qquad
  %% Gamma_2 terms
  +2(\tau_1+\tau_2) +\varepsilon_{2,1} (\xi_2+\tau_1+\rho) +\varepsilon_{2,2} (\xi_2+\tau_2+\rho)
  \\ \nonumber & \qquad \qquad
  %% Gamma_3 terms
  +(2\tau_1+\tau_2+\rho) + \varepsilon_{3,1} (\xi_3+\tau_1+\tau_3+\tau_4) 
  \\ \nonumber & \qquad \qquad \qquad
   %% Gamma_den terms
  -\varepsilon_{2-\frac12} (\xi_1+\xi_2+\tau_1+\tau_2+\rho) + \varepsilon_{2+\frac12} (\rho+\xi_2+\xi_3)
  \\ \nonumber
 = &\ \varepsilon_{1,1}(\xi_1+\tau_2) + \varepsilon_{2,1}(\xi_2+\tau_1+\rho) + \varepsilon_{2,2}(\xi_2+\tau_2+\rho) + \varepsilon_{3,1}(\xi_3-\tau_2)\\ \nonumber & \qquad - \varepsilon_{2-\frac12}(\xi_1+\xi_2+\tau_1+\tau_2+\rho)-\varepsilon_{2+\frac12}(\xi_2+\xi_3+\rho).
\end{align}

Since the coefficient of $\xi_1$ is $(\varepsilon_{1,1}-\varepsilon_{2-\frac12})$, we see that \eqref{eq:ExpZeroSetLemma1} is satisfied.  By similarly looking at the coefficient of $\xi_3$, we see that \eqref{eq:ExpZeroSetLemma3} holds.  Using this, \eqref{eq:AbsValSum2} simplifies further to
 \[ 0 = (\varepsilon_{2,2}-\varepsilon_{2-\frac12})(\xi_2+\tau_1+\rho) + (\varepsilon_{2,3}-\varepsilon_{2+\frac12})(\xi_2+\tau_2+\rho),\]
which is obviously true if and only if $\varepsilon_{2,2}=\varepsilon_{2-\frac12}$ and $\varepsilon_{2,3}=\varepsilon_{2+\frac12}$.  This proves \eqref{eq:ExpZeroSetLemma2}.  Since it must be the case that $\varepsilon_{2,1}\geq \varepsilon_{2,2}$, it follows that $\varepsilon_{2-\frac12} \geq \varepsilon_{2+\frac12}$, as claimed.

Suppose that $\varepsilon$ corresponds to one of the three admissible solutions to \eqref{eq:AbsValSum2}.  Considering only the inequalities that are determined by the $\varepsilon_{j,k}$, the stated inequalities of the 
three solutions are immediate.  So, in order to complete the proof, we must show that the inequalities imposed by $\varepsilon_{2-\frac12}$ and $\varepsilon_{2+\frac12}$ are superfluous, i.e., they do not impose any further restriction on $\xi_1$, $\xi_2$ or $\xi_3$.  

We check this first in the case that $(\varepsilon_{2-\frac12},\varepsilon_{2+\frac12})=(+1,+1)$, for which 
\begin{align*}
 -\tau_2\leq &\ \ \xi_1\, \leq -\rho, \\
 -\tau_2-\rho\leq &\ \ \xi_2\, \leq \tau_1+\tau_2,\\
 -\tau_1-\tau_2-\rho\leq &\ \ \xi_3\, \leq \tau_2.
\end{align*}
Combining the first and second sets of inequalities, we see that
 \[ 0 \leq \tau_1-\tau_2 = (-\tau_2)+(-\tau_2-\rho)+\tau_1+\tau_2+\rho \leq \xi_1+\xi_2+\tau_1+\tau_2+\rho. \]
That is to say that $\varepsilon_{2-\frac12}$ must be $+1$.  In other words, the condition cut out by $\varepsilon_{2-\frac12}$ is already a consequence of the fact that $\varepsilon_{1,1}=\varepsilon_{2,1}=\varepsilon_{2,2}=+1$.  In like manner, we see that the inequality required by $\varepsilon_{2+\frac12}=+1$ is already true by combining the second and third inequalities above.

In each of the other two cases $(\varepsilon_{2-\frac12},\varepsilon_{2+\frac12})=(+1,-1)$ or $(-1,-1)$, one similarly shows that the inequalities imposed by $\varepsilon_{2\pm\frac12}$ are already satisfied given those imposed by $\varepsilon_{j,k}$.
\end{proof}

Lemma~\ref{lem:ExpZeroSet} allows us to restrict the integration in \eqref{eq:pTRasGammas} to the three possible bounded subsets in the $s$-variables,  and then the integral can be bounded.  However, the resulting bound 
is not strong enough for our application, thus requiring that we consider 
$a=(a_1,a_2,a_3)$ for which there exists at least one $i\in\{1,2,3\}$ with $a_i>0$.  In this case \eqref{eq:MellinGL4} is no longer valid.  Instead, we need to use Proposition~\ref{prop:GL4shifteq}.  Although this introduces some technical difficulties, the Gamma-functions which occur are the same as in \eqref{eq:MellinGL4}, and so Lemma~\ref{lem:ExpZeroSet} equally applies.  This is carried out in detail in Section~\ref{sec:shiftedmaintermbound}.

\subsection{\bf Bounds for the shifted $p_{T,R}$ term} \label{sec:shiftedmaintermbound}

$\phantom{xxx}$
\vskip 5pt

Recall that
\begin{equation}
p_{T,R}^\sharp(\alpha) \; := \;e^{\frac{\alpha_1^2+\alpha_2^2+\alpha_3^2+\alpha_4^2}{2T^2}}\,\mathcal F_R(\alpha)\,
      \prod_{1\leq\, j \ne k\, \leq n}
      \Gamma\left(\textstyle{\frac{2+R+\alpha_j - \alpha_k}{4}}\right).
\end{equation}

By the inverse Lebedev-Whittaker transform (see \cite{GK2012}), we see that $p_{T,R}$ is given by
\begin{equation}\label{pTRfunction}
p_{T,R}(y) = p_{T,R}(y_1,y_2,y_3) = \frac{1}{\pi^{3}}
  \underset{\re(\alpha_j)=0}{\iiint} p^{\#}_{T,R}(\alpha) \;W_{\alpha}(y)
     \frac{ d\alpha_1\, d\alpha_2\, d\alpha_3
         }{\prod\limits_{1\leq\, j \ne k\, \leq 4} \Gamma\left(\frac{\alpha_j-\alpha_k}{2}\right) }.
\end{equation}
To obtain a sharp bound for $p_{T,R}(y)$ we replace the Whittaker function $W_\alpha(y)$ on the right hand side of (\ref{pTRfunction}) by its inverse Mellin transform
$$W_\alpha(y) =  \frac{1}{(2\pi i)^3} \iiint\limits_{\text{Re}(s_j) = 
\varepsilon} y_1^{\frac32-s_1} y_2^{2-s_2} y_3^{\frac32-s_3} \, \widetilde W_{\alpha}(s) \, ds.$$
Following  (\ref{eq:pTRdecomposition}), we can then  shift the lines of integration $\re(s)=\varepsilon$  to the left to $\re(s) = -a =(-a_1,-a_2,-a_3)$ (with $a_1,a_2,a_3 > 0$)  and express  $p_{T,R}(y)$ as a sum 
of residues plus a shifted $p_{T,R}$ integral given by

\begin{equation}
p_{T,R}(y, -a) := \frac{1}{(2\pi^2 i)^3} \iiint\limits_{\text{Re}(s) = 
-a} y_1^{\frac32-s_1} y_2^{2-s_2} y_3^{\frac32-s_3} \underset{\re(\alpha_j)=0}{\iiint} p^{\#}_{T,R}(\alpha) \, \widetilde W_{\alpha}(s) \,\frac{ 
d\alpha_1\, d\alpha_2\, d\alpha_3
         }{\prod\limits_{1\leq\, j \ne k\, \leq 4} \Gamma\bigl(\frac{\alpha_j-\alpha_k}{2}\bigr) }\; ds.
         \end{equation}

\begin{prop}\label{prop:GL4pTRemptyset}
Let $\varepsilon>0$ and $R$ sufficiently large be fixed. Let $r_1,r_2,r_3$ be integers.  Given any choice of parameters $a_1,a_2,a_3$ for which $\varepsilon \leq \lvert 2r_j -a_j \rvert \leq 1-\varepsilon$ for each $j=1,2,3$, we have the bound
 \[ \big| p_{T,R}(y,-a) \big| \ll y_1^{\frac{3}{2} + a_1} y_2^{2 + a_2} y_3^{\frac{3}{2} + a_3}
\cdot  T^{\varepsilon+4R+ 9 + \sum\limits_{j=1}^3(\delta_{0,r_j} -r_j)},
 \]
where $\delta_{0,r_i}=1$ if $r_i=0$ and zero otherwise.  The implicit 
constant depends on $r, \varepsilon$ and $R$.
\end{prop}

\begin{proof}
In order to obtain good bounds for $p_{T,R}(y,-a)$ when 
\begin{equation}\label{eq:aibounds}
 a_i \in [2r_i-1+\varepsilon, 2r_i-\varepsilon ] \cup [2r_i+\varepsilon, 2r_i+1-\varepsilon],
\end{equation}
  it is necessary to have sharp bounds for the growth of $\widetilde W_{\alpha}(s)$ on the lines $\text{Re}(s_i)=-a_i.$ Such bounds are known when $\text{Re}(s_i)>0,$ and we can backtrack to this situation by use of the shift equation for $\widetilde W_{\alpha}(s)$ given in Proposition~\ref{prop:GL4shifteq}, which for any $t_1,t_2,t_3\geq 0$ and any $s_1,s_2,s_3\in\C$ states that
\begin{align}\label{eq:WhittakerShiftrrr}
\left| \Mellin{\alpha}(s_1,s_2,s_3) \right| \ll \sum_{\ell=0}^{t_1}\sum_{k=0}^{t_2}\sum_{j=0}^{t_3} \left| \frac{Q_{j,k,\ell}^{t_1,t_2,t_3}}{\mathcal{B}_1^{t_1}\mathcal{B}_2^{t_2}\mathcal{B}_3^{t_3}} \Mellin{\alpha}\big(s_1+2(t_1+j+k),s_2+2t_2,s_3+2(t_3+\ell)\big)\right|,
\end{align}
where $\deg(Q_{j,k,\ell}^{t_1,t_2,t_3})=2(t_1+2t_2+t_3-j-k-\ell)$.

Recall that
\begin{align}\label{eq:B1functiondef}
 \mathcal{B}_1 := \mathcal{B}_1(\alpha,s) & := \big(s_1+\alpha_1\big) 
\big(s_1+\alpha_2\big) \big(s_1+\alpha_3\big) \big(s_1+\alpha_4\big), 
 \\ \label{eq:B2functiondef}
 \mathcal{B}_2 := \mathcal{B}_2(\alpha,s) & := \big(s_2+\alpha_1+\alpha_2\big) \big(s_2+\alpha_1+\alpha_3\big) \big(s_2+\alpha_1+\alpha_4\big)
 \\ & \qquad \cdot \nonumber
 \big(s_2+\alpha_2+\alpha_3\big) \big(s_2+\alpha_2+\alpha_4\big) \big(s_2+\alpha_3+\alpha_4\big), 
 \\ \label{eq:B3functiondef}
 \mathcal{B}_3 := \mathcal{B}_3(\alpha,s) & := \big(s_3-\alpha_1\big) 
\big(s_3-\alpha_2\big) \big(s_3-\alpha_3\big) \big(s_3-\alpha_4\big). 
\end{align}

For each $i=1,2,3$, let 
 \[ t_i := \begin{cases} r_i & \mbox{ if }a_i<2r_i,\\ r_i+1 & \mbox{ if 
}a_i>2r_i. \end{cases}\]
We see that for $a=(a_1,a_2,a_3)$,
 \[ p_{T,R}(y,-a) = \sum_{j=0}^{t_3}\sum_{k=0}^{t_2}\sum_{\ell=0}^{t_1} p_{T,R}^{(j,k,\ell)}(y,-a), \]
where
\begin{align*} &
  \left| p_{T,R}^{(j,k,\ell)}(y, -a) \right|  \ll \;
  y_1^{\frac{3}{2} + a_1} y_2^{2 + a_2} y_3^{\frac{3}{2} + a_3}
  \underset{\re(\alpha_j)=0}{\iiint} e^{\frac{\alpha_1^2+\alpha_2^2+\alpha_3^2+\alpha_4^2}{2T^2}} \cdot \big| \mathcal F_R(\alpha)\big|\cdot \big| \Gamma_R(\alpha) \big| \;
  \underset {\re(s)=-a}  {\iiint}\;
   \left|\frac{Q_{j,k,\ell}^{t_1,t_2,t_3}}{\mathcal{B}_1^{t_1}\mathcal{B}_2^{t_2}\mathcal{B}_3^{t_3}} \right|
     \\ 
  & \cdot \int\limits_{\re(t)=-\varepsilon'}
  \left| \frac{ \Gamma_0(t,\alpha) \Gamma_1(t,s_1+2(t_1+j+k),\alpha) \Gamma_2(t,s_2+2t_2,\alpha) \Gamma_3(t,s_3+2(t_3+\ell),\alpha) }{ \Gamma_{\mathrm den}\big(t, (s_1+2(t_1+j+k),s_2+2t_2,s_3+2(t_3+\ell)),\alpha\big) }
  \right|\, dt\, ds\, d\alpha.
    \end{align*}
 
It follows that 
  \begin{align} \nonumber&
  p_{T,R}^{(j,k,\ell)}(y, -a)\\ &\nn \quad \ll \; y_1^{\frac{3}{2} + a_1} 
y_2^{2 + a_2} y_3^{\frac{3}{2} + a_3}
  \cdot T^{\varepsilon + 2t_1+2t_2+t_3-j-k-\ell)}
  \underset{\re(\alpha_j)=0}{\iiint} e^{\frac{\alpha_1^2+\alpha_2^2+\alpha_3^2+\alpha_4^2}{2T^2}} \;
  \underset {\re(s)=-a} {\iiint}\; 
  \left|\mathcal{B}_1^{{-t_1}} \mathcal{B}_2^{{-t_2}} \mathcal{B}_3^{{-t_3}} \right|
  \\&\nonumber 
  \cdot \int\limits_{\re(t)=-\varepsilon'}
 \left| \frac{ \Gamma\left( \frac{-t+\alpha_3}{2}\right) \Gamma\left( \frac{-t+\alpha_4}{2}\right) \Gamma\left(\frac{s_1+2t_1+2j+2k+\alpha_1}{2}\right)  \Gamma\left(\frac{s_1+2(t_1+j+k)+\alpha_2}{2}\right) \Gamma\left(\frac{s_1+2(t_1+j+k)+t}{2}\right)  }  
  {\Gamma\left(\frac{s_1+s_2+2({t_1+t_2}+j+k)+\alpha_1+\alpha_2+t}{2}\right) \Gamma\left(\frac{s_2+s_3+2({t_2+t_3}+\ell)+t}{2}\right) }\right|\nonumber
 \\   &
  \hskip 42pt
    \cdot 
    \textstyle{
    \left|\Gamma\left(\frac{s_2+2t_2+\alpha_1+\alpha_2}{2}\right)  \Gamma\left(\frac{s_2+2t_2+\alpha_1+t}{2}\right) \Gamma\left(\frac{s_2+2t_2+\alpha_2+t}{2}\right) \Gamma\left(\frac{s_2+2t_2+\alpha_3+\alpha_4}{2}\right)\right| 
    }
   \nonumber\\   &
   \hskip 24pt
   \cdot
    \textstyle{
    \left|\Gamma\left(\frac{s_3+2(t_3+\ell)+\alpha_1+\alpha_2+t}{2}\right)  \Gamma\left(\frac{s_3+2(t_3+\ell)+\alpha_1+\alpha_3+\alpha_4}{2}\right) \Gamma\left(\frac{s_3+2(t_3+\ell)+\alpha_2+\alpha_3+\alpha_4}{2}\right)\right|
    }
   \nonumber\\ 
   &
   \hskip 18pt  \cdot   \left| \prod_{\sigma\in S_4} 
  \big(1+\alpha_{\sigma(1)}+\alpha_{\sigma(2)}-\alpha_{\sigma(3)}-\alpha_{\sigma(4)}\big)
  \right|^{\frac{R}{24}}
    \prod\limits_{1\leq\, j < k\, \leq 4}\frac{ 
  \left| \Gamma\left(\textstyle{\frac{2+R+\alpha_j - \alpha_k}{4}}\right)\right|^2
     }{\left| \Gamma(\frac{\alpha_j-\alpha_k}{2})\right|^2  }\; dt\,ds \,d\alpha,
   \label{eq:pTRemptysetj}\end{align}
where  $-\varepsilon'=\re(t)$ is such that the real part of the arguments of all of the Gamma functions appearing here are positive.  Hence $0<\varepsilon'<\varepsilon$.

  Note that besides the presence of the additional polynomials $\mathcal{B}_i^{r_i}$, the Gamma factors occurring in \eqref{eq:pTRemptysetj} are the same as that of \eqref{eq:pTRasGammas}.  In any event, the exponential 
zero set is precisely the same as that which was determined in Lemma~\ref{lem:ExpZeroSet} and can be any one of 
the three exponential zero sets $\mathcal R_1, \mathcal R_2, \mathcal R_3$ as given in Lemma~\ref{lem:ExpZeroSet}.

Note that we can get from $\mathcal{R}_2$ to $\mathcal{R}_1$ by making the change of variables
 \[ (\xi_2,\xi_3)\mapsto (-\xi_3-\rho,\xi_2-\tau_1). \]
Similarly, to go from $\mathcal{R}_2$ to $\mathcal{R}_3$, one makes the change of variables
 \[ (\xi_1,\xi_2)\mapsto (\xi_2-\tau_1,-\xi_1+\rho). \]
Either of these transformations result in no change in \eqref{eq:WhittasInt} because it leaves both the measure $ds$ and the region over which we are integrating $\re(s)=-a$ invariant.  This implies any bound obtained 
for a given $\mathcal{R}_i$ holds for each of the other choices as well.

 Recall the exponential zero set  $\mathcal R_2$ given by
\begin{align*}
&\hskip 70pt \tau_4 \; \leq \; \rho   \; \leq \; \tau_3,\quad
 -\tau_2 \;\leq \; \xi_1  \; \leq \;-\rho,\\
& -\tau_1 -\rho \;  \leq \; \xi_2  \; \leq \; -\tau_2 -\rho,\quad
-(\tau_1+\tau_2) -\rho \; \leq \; \xi_3  \; \leq \; \tau_2.
 \end{align*}
 
 Recall, also,  that  $\text{\rm Im}(\alpha_j) = \tau_j$ (for $j=1,2,3,4$ with $\tau_4=-\tau_1-\tau_2-\tau_3$), $\text{\rm Im}(t) = \rho$, 
and $\text{\rm Im}(s_i) =\xi_i,$ (for $i=1,2,3$). Accordingly, we write $s = -a+i\xi,\; \alpha = i\tau$
with $\xi = (\xi_1,\xi_2,\xi_3)$ and  $\tau = (\tau_1, \tau_2,\tau_3).$

We now replace each Gamma factor in the integral (\ref{eq:pTRemptysetj}) using the Stirling bound~\eqref{eq:Stirling}.  A bound for the integral is given by integrating over $\mathcal R_2$ (the exponential zero set) and 
just using the polynomial bound coming from the
Stirling bound~\eqref{eq:Stirling}.
 It follows (up to a constant dependant on $a_k$, $R$ and $\varepsilon>0$) that $p_{T,R}^{(j,k,\ell)}(y, -a)$ is bounded by
\begin{align*}&
  y_1^{\frac{3}{2} + a_1} y_2^{2 + a_2} y_3^{\frac{3}{2} + a_3}
  \cdot T^{\varepsilon + 2(t_1+2t_2+t_3-j-k-\ell)}\\&\cdot
 \underset{|\tau_1|,|\tau_2|,|\tau_3|,|\tau_4| \ll T^{1+\varepsilon}} { \iiint\limits_{\tau_1\geq \tau_2\geq \tau_3\geq \tau_4}}\;\;\,
  \int\limits_{\rho=\tau_4}^{\tau_3}\;
  \int\limits_{\xi_1=-\tau_2}^{-\rho}\;
   \int\limits_{\xi_2=-\tau_1-\rho}^{-\tau_2-\rho} \;
  \int\limits_{\xi_3=-\tau_1-\tau_2 -\rho}^{\tau_2}  \big|\mathcal{B}_1^{-t_1} \mathcal{B}_2^{-t_2} \mathcal{B}_3^{-t_3}\big|
  \\ 
  & 
  \cdot 
  \frac{
   (1+|\tau_3-\rho|)^{-\frac12}\,
      \big(1+|\tau_4-\rho|\big )^{-\frac12}\,
      \,
        }  
  { \left(1+|\xi_1+\xi_2+\tau_1+\tau_2+\rho|   \right)^{-\frac12+j+k+t_1-\frac{a_1}{2}+t_2-\frac{a_2}{2} } 
\, ( 1+|\xi_2+\xi_3+\rho|   )^{-\frac12+\ell+t_3-\frac{a_3}{2}+t_2-\frac{a_2}{2} } }
 \\ &
  \cdot
   (1+|\xi_1+\tau_1|)^{j+k+t_1-\frac{a_1+1}{2}}\,  
        (1+|\xi_1+\tau_2|)^{j+k+t_1-\frac{a_1+1}{2}}\,
         \big(1+|\xi_1+\rho|\big)^{j+k+t_1-\frac{a_1+1}{2}}
         \\
         &
    \cdot 
    (1+|\xi_2+\tau_1+\tau_2|)^{t_2-\frac{a_2+1}{2}}\,
     (1+|\xi_2+\tau_1+\rho|)^{t_2-\frac{a_2+1}{2}}\,
      \big(1+|\xi_2+\tau_2+\rho|\big)^{t_2-\frac{a_2+1}{2}}\,
  \\
   &
   \cdot \big(1+|\xi_2+\tau_3+\tau_4|\big)^{t_2-\frac{a_2+1}{2}}
  \big(1+|\xi_3+\tau_1+\tau_2+\rho|\big)^{\ell+t_3-\frac{a_3+1}{2}}\,
  \big(1+|\xi_3+\tau_1+\tau_3+\tau_4|\big)^{\ell+t_3-\frac{a_3+1}{2}}
     \\&
   \cdot    \big (1+|\xi_3+\tau_2+\tau_3+\tau_4|\big)^{\ell+t_3-\frac{a_3+1}{2}} 
   \, K_R(\tau) \,d\xi_3\,d\xi_2\,d\xi_1\; d\rho\; d\tau_3\,d\tau_2\, d\tau_1,
    \nonumber\end{align*}
    where
\begin{align}\label{K_R(tau)}
  K_R(\tau) & :=
 \Big(1+\lvert \tau_1+\tau_2-\tau_3-\tau_4\rvert\Big)^{\frac{R}{3}}
 \Big(1+\lvert \tau_1+\tau_3-\tau_2-\tau_4\rvert\Big)^{\frac{R}{3}}
 \\ & \qquad \nonumber
 \Big(1+\lvert \tau_1+\tau_4-\tau_2-\tau_3\rvert\Big)^{\frac{R}{3}} 
 \cdot \prod_{1\leq j<k\leq 4}\big(1+\tau_j-\tau_k\big)^{1+\frac{R}{2}}.
\end{align}

\begin{rmrk} The exponents of $(1+|\tau_j-\rho|)$ etc. should be modified 
by $\pm\varepsilon$ as Re$(t)=-\varepsilon$ but we have absorbed these into the term $T^\varepsilon$.  To simplify the exposition we drop these $\varepsilon$'s in what follows since they do not affect the conclusion.
\end{rmrk}

\begin{lem}\label{lem:r2primecancel}
Suppose that $a_i = 2r_i + \gamma_i a_i'$ where 
 \[ \gamma_i = \begin{cases} -1 & \mbox{ if }a_i<2r_i \\ 1 & \mbox{ if } a_i>2r_i, \end{cases} \]
and $0<a_i'<1$ for each $i=1,2,3$.  Throughout the domain of integration, in the integrand for $p_{T,R}^{(j,k,\ell)}(y,-a)$ above, the expression
\begin{align*}
& \frac{
        \big(1+|\xi_2+\tau_1+\tau_2|\big)^{t_2-\frac{a_2+1}{2}}\,
        \big(1+|\xi_2+\tau_1+\rho|\big)^{t_2-\frac{a_2+1}{2}}\,
        \big(1+|\xi_2+\tau_2+\rho|\big)^{t_2-\frac{a_2+1}{2}}\,
        \big(1+|\xi_2+\tau_3+\tau_4|\big)^{t_2-\frac{a_2+1}{2}}
        }{
        \big(1+|\xi_1+\xi_2+\tau_1+\tau_2+\rho|   \big)^{t_2-\frac{a_2+1}{2}} \,
       \big( 1+|\xi_2+\xi_3+\rho|   \big)^{t_2-\frac{a_2+1}{2} }
        }
\end{align*}
is bounded by a constant multiple of $T^{\varepsilon+\frac{1+\gamma_2}{2}}$.  Similarly,
\begin{align*}
& \frac{
        \big(1+|\xi_1+\tau_1|\big)^{j+k+t_1-\frac{a_1+1}{2}}\,  
        \big(1+|\xi_1+\tau_2|\big)^{j+k+t_1-\frac{a_1+1}{2}}\,
        \big(1+|\xi_1+\rho|\big)^{j+k+t_1-\frac{a_1+1}{2}}
        }{
        \big(1+|\xi_1+\xi_2+\tau_1+\tau_2+\rho|   \big)^{j+k+t_1-\frac{a_1}{2}}
        } \ll  T^{\varepsilon+2(j+k)+\frac{1+\gamma_1}{2}},
\end{align*}
and
\begin{align*}
& \frac{
        \big(1+|\xi_3+\tau_1+\tau_2+\rho|\big)^{\ell+t_3-\frac{a_3+1}{2}}\,
        \big(1+|\xi_3+\tau_1+\tau_3+\tau_4|\big)^{\ell+t_3-\frac{a_3+1}{2}}
        \big(1+|\xi_3+\tau_2+\tau_3+\tau_4|\big)^{\ell+t_3-\frac{a_3+1}{2}} 
        }{
        \big( 1+|\xi_2+\xi_3+\rho|   \big)^{\ell+t_3-\frac{a_3}{2}}
        } 
\end{align*}
is bounded by a multiple of $T^{\varepsilon+2\ell+\frac{1+\gamma_3}{2}}$.
\end{lem}

\begin{proof}
Note that $t_i=r_i+\frac{1+\gamma_i}{2}$.  This allows us to simplify many of the exponents.  For example,
\begin{align} \label{eq:simpleboundlast}
 \big(1+\lvert \xi_i+f(\tau,\rho)\rvert\big)^{t_i-\frac{a_i+1}{2}}
 =
 \big(1+\lvert \xi_i+f(\tau,\rho)\rvert\big)^{\gamma_i\frac{1-a_i'}{2}}
 & \ll \big(1+\lvert \xi_i+f(\tau,\rho)\rvert\big)^{\frac{1+\gamma_i}{4}}
\end{align}
for each of the cases $i=1,2,3$ and $f(\tau,\rho)$ appearing above.

In a similar manner, one sees that
\begin{align} \label{eq:simpleboundfirst}
 \big(1+\lvert \xi_1+\xi_2+\tau_1+\tau_2\rvert \big)^{\frac12-j-k-t_1+\frac{a_2}{2}} & \ll \big(1+\lvert \xi_1+\xi_2+\tau_1+\tau_2\rvert \big)^{\frac12-j-k-\frac{1+\gamma_1}{4}}
 \\ \nonumber & \ll \frac{T^{\varepsilon+\frac12}}{\big(1+\lvert \xi_1+\xi_2+\tau_1+\tau_2\rvert \big)^{j+k+\frac{1+\gamma_1}{4}}}, \\ \nonumber
 \\ 
 \big(1+\lvert \xi_2+\xi_3+\rho\rvert \big)^{\frac12-\ell-t_3+\frac{a_3}{2}} & \ll
 \big(1+\lvert \xi_2+\xi_3+\rho\rvert \big)^{\frac12-\ell-\frac{1+\gamma_3}{4}} 
 \\ \nonumber & \ll \frac{T^{\varepsilon+\frac12}}{\big(1+\lvert \xi_2+\xi_3+\rho\rvert \big)^{\ell+\frac{1+\gamma_3}{4}}}.
\end{align}
Successively making the change of variables
 \[ \xi_1\mapsto \xi_1-\tau_2, \;\;\quad 
  \xi_2\mapsto \xi_2-\tau_1-\rho,\;\;\quad
  \xi_3\mapsto \xi_3 -\tau_1-\tau_2-\rho,\;\;\quad
  \rho\mapsto \rho+\tau_4, \]
and the substitutions
 \[ T_1 = \tau_1-\tau_2, \quad T_2 =\tau_2-\tau_3, \quad T_3 = \tau_3-\tau_4 = \tau_1+\tau_2+2\tau_3,\]
the bounds of integration in $p_{T,R}^{(j,k,\ell)}(y,-a)$ become
\[ 0\leq T_1,T_2,T_3 \leq T^{1+\varepsilon},\quad 0\leq \rho \leq T_3,\quad 0 \leq \xi_1\leq T_2+T_3-\rho, \quad 0\leq \xi_2 \leq T_1, \quad 0\leq 
\xi_3 \leq T_2+\rho, \]
and the terms involving $\xi_1$ are
\begin{align*}
\frac{
    \big(1+T_1+\xi_1\big)^{j+k+\frac{1+\gamma_1}{4}}\,  
        \big(1+\xi_1\big)^{j+k+\frac{1+\gamma_1}{4}}\,
        \big(1+T_2+T_3-\rho-\xi_1 \big)^{j+k+\frac{1+\gamma_1}{4}}  }
         { 
        \big(1+\xi_1+\xi_2 \big)^{j+k+\frac{1+\gamma_1}{4}} }, 
\end{align*}
which, since $0\leq j,k$ and $0\leq T_i\leq T^{1+\varepsilon}$ is bounded 
by $T^{\varepsilon+\frac12+2(j+k)+\frac{1+\gamma_1}{2}}$, as claimed.

The bound for $\xi_3$ follows in precisely the same fashion, and that for 
$\xi_2$ is similar.
\end{proof}

\vskip 12pt
Putting this together, we obtain the bound 
\begin{align*}&
\left| p_{T,R}^{(j,k,\ell)}(y, -a) \right| \; \ll \;
  y_1^{\frac{3}{2} + a_1} y_2^{2 + a_2} y_3^{\frac{3}{2} + a_3} \, T^{\varepsilon + 2(t_1+2t_2+t_3)+\sum\limits_{i=1}^3 \frac{1+\gamma_i}{2}}\\&\cdot
 \underset{|\tau_1|,|\tau_2|,|\tau_3|,|\tau_4| \ll T^{1+\varepsilon}} { \iiint\limits_{\tau_1\geq \tau_2\geq \tau_3\geq \tau_4}}\;\;\,
  \int\limits_{\rho=\tau_4}^{\tau_3}\;
  \int\limits_{\xi_1=-\tau_2}^{-\rho}\;
   \int\limits_{\xi_2=-\tau_1-\rho}^{-\tau_2-\rho} \;
  \int\limits_{\xi_3=-\tau_1-\tau_2 -\rho}^{\tau_2} \, (1+|\tau_3-\rho|)^{-\frac12}\,
     \big(1+|\tau_4-\rho|\big)^{-\frac12}
     \\ & 
    \cdot
        \Big[(1+|\xi_1+\tau_1|)  \big(1+|\xi_1+\tau_2|) (1+|\xi_1+\tau_3|\big) \big(1+|\xi_1+\tau_4|\big)\Big]^{-t_1}
         \\
         &
    \cdot 
        \Big[(1+|\xi_2+\tau_1+\tau_2|)  (1+|\xi_2+\tau_1+\tau_3|) (1+|\xi_2+\tau_1+\tau_4|)  \big(1+|\xi_2+\tau_2+\tau_3|\big)\Big]^{-t_2}
  \\
  &
  \cdot 
    \Big[\big(1+|\xi_2+\tau_2+\tau_4|\big) \big(1+|\xi_2+\tau_3+\tau_4|\big)\Big]^{-t_2}
  \\
   &
   \cdot
    \Big[\big(1+|\xi_3-\tau_1|\big) \big (1+|\xi_3-\tau_2|\big) \big(1+|\xi_3-\tau_3|\big) \big(1+|\xi_3-\tau_4|\big)\Big]^{-t_3} 
     \\
   &
 \cdot K_R(\tau)\, d\xi_3\,d\xi_2\,d\xi_1\; d\rho\; d\tau_3\,d\tau_2\, d\tau_1.
\nonumber\end{align*}

Next, we successively make the change of variables
  $$\xi_1\mapsto \xi_1-\tau_2, \;\;\quad 
  \xi_2\mapsto \xi_2-\tau_1-\rho,\;\;\quad
  \xi_3\mapsto \xi_3 -\tau_1-\tau_2-\rho,\;\;\quad
  \rho\mapsto \rho+\tau_4, $$
and the substitutions
   $$T_1 = \tau_1-\tau_2, \quad T_2 =\tau_2-\tau_3, \quad T_3 = \tau_3-\tau_4 = \tau_1+\tau_2+2\tau_3.$$
  Then by abuse of notation we may replace $K_R(\tau)$ by $K_R(T_1,T_2,T_3)$ where
\begin{align} \label{K_R(T)}
 K_R(T_1,T_2,T_3) & := 
 \Big(1+T_1+2T_2+T_3 \Big)^{\frac{R}{3}}
 \Big(1+T_1+T_3 \Big)^{\frac{R}{3}}
 \Big(1+\lvert T_1-T_3\rvert\Big)^{\frac{R}{3}}
 \\ \nonumber & \hskip 24pt
 \cdot
 \Big(1+T_1\Big)^{1+\frac{R}{2}}
 \Big(1+T_2\Big)^{1+\frac{R}{2}}
 \Big(1+T_3\Big)^{1+\frac{R}{2}}
 \\ \nonumber & \hskip 24pt
 \cdot
 \Big(1+T_1+T_2\Big)^{1+\frac{R}{2}}
 \Big(1+T_2+T_3\Big)^{1+\frac{R}{2}}
 \Big(1+T_1+T_2+T_3\Big)^{1+\frac{R}{2}}
 \\ \nonumber & \hskip -12pt
 \ll T^{\varepsilon+3+\frac{13R}{6}}\cdot 
 \Big(1+\lvert T_1-T_3\rvert\Big)^{\frac{R}{3}}
 \Big(1+T_1\Big)^{1+\frac{R}{2}}
 \Big(1+T_2\Big)^{1+\frac{R}{2}}
 \Big(1+T_3\Big)^{1+\frac{R}{2}}.
\end{align}    
It follows that
\begin{align*} %\label{p_{T,R}-Bound}
\left| p_{T,R}^{(j,k,\ell)}(y, -a) \right| & \:\ll \;
  y_1^{\frac{3}{2} + a_1} y_2^{2 + a_2} y_3^{\frac{3}{2} + a_3}
\, T^{\varepsilon+ \frac{13R}{6}+3} + 2(t_1+2t_2+t_3)+\sum\limits_{i=1}^3\frac{1+\gamma_i}{2}\\&\cdot
\iiint\limits_{0\leq T_1,T_2,T_3 \leq T^{1+\varepsilon}}\;\;\,
  \int\limits_{\rho=0}^{T_3}\;
  \int\limits_{\xi_1=0}^{T_2+T_3-\rho}\;
   \int\limits_{\xi_2=0}^{T_1} \;
  \int\limits_{\xi_3=0}^{T_2+\rho}  (1+T_3-\rho)^{-\frac12}\,
     \big(1+\rho\big)^{-\frac12}
 \\ & 
    \cdot
        \Big[(1+\xi_1+T_1)  \big(1+\xi_1) (1+|\xi_1-T_2|\big) \big(1+T_2+T_3-\xi_1\big)\Big]^{-t_1}
         \\
         &
    \cdot 
        \Big[\big(1+\xi_2+T_2+T_3-\rho\big)  \big(1+\xi_2+T_3-\rho\big) (1+|\xi_2-\rho|)  \Big]^{-t_2}
  \\
  &
  \cdot 
    \Big[\big(1+|\xi_2-T_1+T_3-\rho|\big) \big(1+T_1+\rho-\xi_2\big) \big(1+T_1-\xi_2+T_2+\rho\big)\Big]^{-t_2}
  \\
   & \hskip -24pt
   \cdot
    \Big[\big(1+T_1+T_2+\rho-\xi_3\big) \big (1+T_2+\rho-\xi_3\big) \big(1+|\xi_3-\rho|\big) \big(1+\xi_3+T_3-\rho\big)\Big]^{-t_3} 
    \; d\xi_3\,d\xi_2\,d\xi_1  \; d\rho
     \\
   & \hskip 48pt \cdot
    \big(1+\lvert T_1-T_3\rvert\big)^{\frac{R}{3}}
   \big(1+T_1\big)^{1+\frac{R}{2}}
   \big(1+T_2\big)^{1+\frac{R}{2}}
   \big(1+T_3\big)^{1+\frac{R}{2}}
\; dT_3\,dT_2\, dT_1.
\nonumber\end{align*}

Now for $0\leq \xi_1 \leq T_2+T_3-\rho$, $0\leq \xi_2\leq T_1$, $0\leq \xi_3 \leq T_2+\rho$ and $0\leq \rho\leq T_3$,
\begin{align*}
1+\xi_1+T_1 \; & \geq\; 1 + T_1 \\
1+T_2+T_3-\rho+\xi_2 \; & \geq\; 1 + T_2 +T_3-\rho \; \geq \; 1+T_2 \\
1+T_1+T_2+\rho-\xi_2 \; & \geq\; 1+ T_2+\rho \; \geq \; 1+T_2 \\
1+T_1+T_2+\rho-\xi_3\; & \geq\; 1 + T_1.
\end{align*}
Hence, we have the bounds
\begin{align*}
\big(1+\xi_1+T_1\big)^{-t_1} \; & \leq\; \big(1 + T_1\big)^{-t_1} \\
\big(1+T_2+T_3-\rho+\xi_2\big)^{-t_2} \; & \leq\; \big(1+T_2\big)^{-t_2} \\
\big(1+T_1+T_2+\rho-\xi_2\big)^{-t_2} \; & \leq\; \big(1+T_2\big)^{-t_2} \\
\big(1+T_1+T_2+\rho-\xi_3\big)^{-t_3}\; & \leq\; \big(1 + T_1\big)^{-t_3}.
\end{align*}

Inserting these bounds, we see that
\begin{align} \label{pTR-Bound}
\left| p_{T,R}^{(j,k,\ell)}(y, -a) \right| &\; \ll \;
  y_1^{\frac{3}{2} + a_1} y_2^{2 + a_2} y_3^{\frac{3}{2} + a_3}
  \ T^{\varepsilon + \frac{13R}{6} + 3 + 2(t_1+2t_2+t_3)+\sum\limits_{i=1}^3\frac{1+\gamma_i}{2}}\\\nn&\cdot
\iiint\limits_{0\leq T_1,T_2,T_3 \leq T^{1+\varepsilon}}\;
    \big(1+T_1\big)^{-t_1-t_3}
    \big(1+T_2\big)^{-2t_2}  \int\limits_{\rho=0}^{T_3} \big(1+T_3-\rho\big)^{-\frac12}
    \big(1+\rho\big)^{-\frac12}
    \\  \nonumber
  & 
  \cdot 
  \int\limits_{\xi_1=0}^{T_2+T_3-\rho}\;
        \big(1+\xi_1\big)^{-t_1} \big(1+|\xi_1-T_2|\big)^{-t_1} \big(1+T_2+T_3-\xi_1\big)^{-t_1} \; d\xi_1
         \\ \nonumber
         &
    \cdot 
   \int\limits_{\xi_2=0}^{T_1} \;
 \big(1+\xi_2+T_3-\rho\big)^{-t_2}
        \big(1+|\xi_2-\rho|\big)^{-t_2}
        \big(1+|\xi_2-T_1+T_3-\rho|\big)^{-t_2}
      \\\nn&\qquad\cdot  \big(1+T_1+\rho-\xi_2\big)^{-t_2}\; d\xi_2
  \\ \nonumber
   &
   \cdot
  \int\limits_{\xi_3=0}^{T_2+\rho}
    \big (1+T_2+\rho-\xi_3\big)^{-t_3}
    \big(1+|\xi_3-\rho|\big)^{-t_3}
    \big(1+\xi_3+T_3-\rho\big)^{-t_3}
    \; d\xi_3\, d\rho
     \\ \nonumber
   & \hskip 36pt
   \cdot   
       \big(1+\lvert T_1-T_3\rvert\big)^{\frac{R}{3}}
   \big(1+T_1\big)^{1+\frac{R}{2}}
   \big(1+T_2\big)^{1+\frac{R}{2}}
   \big(1+T_3\big)^{1+\frac{R}{2}}
  \; dT_3\,dT_2\, dT_1.
\nonumber\end{align}

Let us denote the integral in $\xi_i$   ($1\le i\le 3$), in \eqref{pTR-Bound}, by $\mathcal{I}_i(r_i)$.  We have
\begin{lem}\label{LemmaXi123rhoBound}
For each $i=1,2,3$, $\mathcal{I}_i(0) \ll T^{\varepsilon+1}$.  Otherwise, if $r\geq 1$, we have the following bounds:
\begin{align*}
\mathcal{I}_1(r), \mathcal{I}_3(r) & \ll \big(1+T_2\big)^{-r}\big(1+T_3\big)^{-r},
\\
\mathcal{I}_2(r) & \ll \Big(\big(1+T_3\big)^{-2r} +
    \big(1+T_1\big)^{-2r} \Big)
    \big(1+\lvert T_1-T_3\rvert\big)^{-r}.
\end{align*}
\end{lem}
\begin{proof}
The case of $r=0$ is clear given that $0\leq T_i \leq T^{1+\varepsilon}$.  Thus, we may assume now that $r\geq 1$.  By expanding the region of integration, since the integrand is positive, we have that
\begin{align*}
 \mathcal{I}_1 & \ll 
     \int_{\xi_1=0}^{T_2+T_3} 
     \big(1+\xi_1\big)^{-r} (1+|\xi_1-T_2|\big)^{-r} \big(1+T_2+T_3-\xi_1\big)^{-r}\; d\xi_1
 \\ & \qquad = 
     \int_{\xi_1=B_1}^{B_3} 
     \big(1+\lvert \xi_1-B_1\rvert \big)^{-r}
     \big(1+\lvert \xi_1-B_2\rvert \big)^{-r}
     \big(1+\lvert \xi_1-B_3\rvert \big)^{-r} \; d\xi_1,
\end{align*}
where
 \[ B_1 = 0, \qquad\qquad B_2=T_2, \qquad\qquad B_3=T_2+T_3.\]
Applying Lemma~\ref{lemmaIntegralBound3}, it follows that
 \[ \mathcal{I}_1 \ll \big(1+T_1\big)^{-r}\big(1+T_2\big)^{-r}\big(1+T_3\big)^{-r}. \]
Replacing $\rho$ by $T_3-\rho$ in $\mathcal{I}_1$, one finds that $\mathcal{I}_3=\mathcal{I}_1$, and hence the same bound holds for $i=3$.

To bound $\mathcal{I}_2$, we also expand the region of integration, whence 
\begin{align*}
 \mathcal{I}_2(r_2) & \ll
 \int\limits_{\xi_2=\rho-T_3}^{T_1+\rho} \big(1+\xi_2+T_3-\rho\big)^{-r}
 \big(1+|\xi_2-\rho|\big)^{-r}
 \\ & \hskip 60pt \cdot 
 \big(1+|\xi_2-T_1+T_3-\rho|\big)^{-r}
 \big(1+T_1+\rho-\xi_2\big)^{-r} \; d\xi_2
 \\ &  =
 \int\limits_{\xi_2=0}^{T_1} \prod_{\ell=1}^4 \big(1+\lvert \xi_2-B_\ell\rvert\big)^{-r}\; d\xi_2,
\end{align*}
where
 \[ B_1=\rho-T_3, \qquad B_2=\rho, \qquad B_3=\rho+T_1-T_3, \qquad B_4=\rho+T_1. \]
Note that 
 \[ B_2\leq B_3 \quad \Longleftrightarrow \quad T_3\leq T_1. \]
In either event, using Lemma~\ref{lemmaIntegralBound3} gives the claimed bound.
\end{proof}

Inserting these bounds into \eqref{pTR-Bound}, the integral in the $\rho$-variable is
 \[ \int\limits_{\rho=0}^{T_3} \big(1+\rho\big)^{-\frac12}\big(1+T_3-\rho\big)^{-\frac12}\; d\rho \ll 1, \]
where this bound is obtained via a direct application of Lemma~\ref{lemmaIntegralBound}.

\vskip 12pt
Plugging these bounds into \eqref{pTR-Bound}, it follows that $\big| p_{T,R}^{(j,k,\ell)}(y,-a)\big|$ is bounded by
\begin{align} \label{pTR-Bound2}
&
y_1^{\frac{3}{2} + a_1} y_2^{2 + a_2} y_3^{\frac{3}{2} + a_3}
\cdot  T^{\varepsilon+\frac{13R}{6}+ 3+  2(r_1+2r_2+r_3)+\delta_{0,t_1}+\delta_{0,t_2}+\delta_{0,t_3}+\sum\limits_{i=1}^3\frac{1+\gamma_i}{2}} 
 \\ \nonumber & \\ \nonumber
 & \hskip 24pt
    \cdot
\iiint\limits_{0\leq T_1,T_2,T_3 \leq T^{1+\varepsilon}}\;
    \big(1+T_1\big)^{1+\frac{R}{2}-t_1-t_3}
    \Big(\big(1+T_3\big)^{-2t_2} +
    \big(1+T_1\big)^{-2t_2} \Big)
     \\ \nonumber
   &
   \\ \nonumber
   &
   \hskip 84pt \cdot
   \big(1+T_3\big)^{1+\frac{R}{2}-t_1-t_3}       \big(1+\lvert T_1-T_3\rvert\big)^{\frac{R}{3}-t_2}
   \big(1+T_2\big)^{1+\frac{R}{2}-t_1-2t_2-t_3}
  \; dT_3\,dT_2\, dT_1.
\end{align} 

\vskip 12 pt
\begin{rmrk} \label{T-Remark}
As long as
 \[1+\frac{R}{2} \geq t_1+2t_2+t_3 \qquad \mbox{and} \qquad
    \frac{R}{3}\geq t_2, \]
    each of the terms in the integrand,  in \eqref{pTR-Bound2},  is of the form $\left(1+T_j\right)^{\alpha_j}$ for $j=1,2,3$ and some $\alpha_j\geq 0$.  Since the integral is over the domain $0\leq T_1,T_2,T_3\ll T^{1+\varepsilon}$, it follows that each of these terms is bounded by $T^{\alpha}$. 
\end{rmrk}

In light of the above remark, we obtain the bound
\begin{equation}\label{eq:riprimebound}
 \left | p_{T,R}^{(j,k,\ell)}(y, -a)\right | \; \ll \;  
 y_1^{\frac{3}{2} + a_1} y_2^{2 + a_2} y_3^{\frac{3}{2} + a_3}
\cdot  T^{\varepsilon+4R+ 9 +  \delta_{0,t_1}+\delta_{0,t_2}+\delta_{0,t_3} - (r_1+r_2+r_3)}. 
\end{equation}
Since $0\leq r_i\leq t_i$, it follows that $\delta_{0,t_i}\leq \delta_{0,r_i}$, and thus \eqref{eq:riprimebound} implies the desired result.
\end{proof}

%\newpage

\subsection{\bf Bounds for the single residue terms}There are two types of single residue terms that we need to deal with.  These are defined in equations \eqref{eq:pTr_1delta} and \eqref{eq:pTr_2delta}.  We show in Propositions~\ref{prop:singleresidue1} and \ref{prop:singleresidue2} respectively that the bounds from each of these is small.

The first single residue term that we need to consider is
\begin{align} \label{eq:pTr_1delta}
  & p_{T,R}^{1,\delta}(y;(-a_2,-a_3)) =  \iiint\limits_{\re(\alpha_j)=0} e^{\frac{\alpha_1^2+\cdots+\alpha_4^2}{2T^2}}\; \iint\limits_{\substack{\re(s_j)=-a_j\\j=2,3}}y_1^{\frac32-p_1}y_2^{2-s_2}y_3^{\frac32-s_3} 
\; 
 Q_{\delta}(s,\alpha)
  \\ \nonumber & \hskip 36pt
  \cdot
 \mathcal F_R(\alpha)\,\Gamma_R(\alpha)
  \cdot \frac{ \prod\limits_{j=2}^4 \Gamma(\frac{\alpha_j-\alpha_1}{2}-\delta) \Gamma(\frac{s_2+\alpha_1+\alpha_j}{2})  \Gamma(\frac{s_3-\alpha_j}{2}) }{ \Gamma(\frac{s_2+s_3+\alpha_1}{2}+\delta) } \; ds_2 ds_3\; d\alpha,
 \end{align}
where $Q_{\delta}$ is a polynomial (see Section \ref{poles-res}) of degree $\le3\delta$,  {$ \mathcal F_R(\alpha)$ is as in \eqref{Fr},  $\Gamma_R(\alpha)$ is as in \eqref{gammaRalpha},}  and $p_1=-\alpha_1-2\delta$ for $\delta=0,1,\ldots,r_1-1$.

\begin{prop}\label{prop:singleresidue1}
Let $r_1\geq 1$, $r_2,r_3\geq0$ be  integers, and $0<\varepsilon<1$. Suppose $a_1,a_2,a_3$ satisfy the hypotheses of Theorem \ref{th:pTRbound}.  If $0\le \delta\le r_1-1$, then
\begin{equation}\label{Prop-s1-Residue}
 \left| p_{T,r}^{1,\delta}(y;(-a_2,-a_3)) \right| \ll y_1^{\frac32+a_1}y_2^{2+a_2}y_3^{\frac32+a_3}  T^{\varepsilon+4R+9+\delta_{0,r_2}+\delta_{0,r_3}-(r_1+r_2+r_3)}.
\end{equation}
\end{prop}

\begin{proof}
In order to bound  $p_{T,R}^{1,\delta}(y;(-a_2,-a_3))$, we will need to shift the lines of integration in the $\alpha_1$ variable.  In doing so, we will pick up residues.  In other words, we may write
 \[ p_{T,R}^{1,\delta}(y;(-a_2,-a_3)) = p_{T,R}^{1,\delta}(y;(-a_2,-a_3),\kappa) + \sum\mbox{Residues}, \]
where 
\begin{align}\label{eq:pTr_1deltashift}
 & p_{T,R}^{1,\delta}(y;(-a_2,-a_3),\kappa) := 
 \iiint\limits_{\re(\alpha)=\kappa} e^{\frac{\alpha_1^2+\cdots+\alpha_4^2}{2T^2}}\; \iint\limits_{\substack{\re(s_j)=-a_j\\j=2,3}}y_1^{\frac32+\alpha_1+2\delta}y_2^{2-s_2}y_3^{\frac32-s_3} \, \mathcal F_R(\alpha)\,
 Q_{\delta}(s,\alpha)
  \\ \nonumber & 
  \cdot
  \frac{\Gamma\(\frac{\alpha_2-\alpha_1}{2}-\delta\)
      }{\Gamma\(\frac{\alpha_2-\alpha_1}{2}\)}
  \cdot
  \frac{\Gamma\(\frac{2+R+\alpha_1-\alpha_2}{4}\)
        \Gamma\(\frac{2+R+\alpha_2-\alpha_1}{4}\)
      }{\Gamma\(\frac{\alpha_1-\alpha_2}{2}\)}
  \cdot
  \frac{\Gamma\(\frac{2+R+\alpha_2-\alpha_3}{4}\)
        \Gamma\(\frac{2+R+\alpha_3-\alpha_2}{4}\)
      }{\Gamma\(\frac{\alpha_2-\alpha_3}{2}\)
        \Gamma\(\frac{\alpha_3-\alpha_2}{2}\)}
  \\ \nonumber & 
  \cdot
  \frac{\Gamma\(\frac{\alpha_3-\alpha_1}{2}-\delta\)
      }{\Gamma\(\frac{\alpha_3-\alpha_1}{2}\)}
  \cdot
  \frac{\Gamma\(\frac{2+R+\alpha_1-\alpha_3}{4}\)
        \Gamma\(\frac{2+R+\alpha_3-\alpha_1}{4}\)
      }{\Gamma\(\frac{\alpha_1-\alpha_3}{2}\)}
  \cdot
  \frac{\Gamma\(\frac{2+R+\alpha_2-\alpha_4}{4}\)
        \Gamma\(\frac{2+R+\alpha_4-\alpha_2}{4}\)
      }{\Gamma\(\frac{\alpha_2-\alpha_4}{2}\)
        \Gamma\(\frac{\alpha_4-\alpha_2}{2}\)}
  \\ \nonumber &   \cdot
  \frac{\Gamma\(\frac{\alpha_4-\alpha_1}{2}-\delta\)
      }{\Gamma\(\frac{\alpha_4-\alpha_1}{2}\)}
  \cdot
  \frac{\Gamma\(\frac{2+R+\alpha_1-\alpha_4}{4}\)
        \Gamma\(\frac{2+R+\alpha_4-\alpha_1}{4}\)
      }{\Gamma\(\frac{\alpha_1-\alpha_4}{2}\)}
  \cdot
  \frac{\Gamma\(\frac{2+R+\alpha_4-\alpha_3}{4}\)
        \Gamma\(\frac{2+R+\alpha_3-\alpha_4}{4}\)
      }{\Gamma\(\frac{\alpha_4-\alpha_3}{2}\)
        \Gamma\(\frac{\alpha_3-\alpha_4}{2}\)}
  \\ \nonumber & \cdot
  \frac{\Gamma\(\frac{s_2+\alpha_1+\alpha_2}{2}\) 
        \Gamma\(\frac{s_2+\alpha_1+\alpha_3}{2}\) \Gamma\(\frac{s_2+\alpha_1+\alpha_4}{2}\) \Gamma\(\frac{s_3-\alpha_2}{2}\) \Gamma\(\frac{s_3-\alpha_3}{2}\) \Gamma\(\frac{s_3-\alpha_4}{2}\)
      }{\Gamma(\frac{s_2+s_3+\alpha_1}{2}+\delta) }\,ds_2\, ds_3\, d\alpha.
\end{align}
and the residues that appear depend on the particular choice of $\kappa=(\kappa_1,\kappa_2,\kappa_3)$.  For example, if $\kappa_j=0$ for $j=1,2,3$, then equations \eqref{eq:pTr_1delta} and \eqref{eq:pTr_1deltashift} are the same, meaning there are no residues.

Our goal will be for the given value of $2r_j-1+\varepsilon\leq a_j \leq 2r_j-\varepsilon$ to shift the lines of integration of the $\alpha$ variables from $\re(\alpha)=0$ to $\re(\alpha)=\kappa=(\kappa_1,\kappa_2,\kappa_3)$ with 
 \[ \re(\alpha_1+2\delta)=\kappa_1+2\delta=a_1, \qquad \kappa_2=0=\kappa_3. \]
 
To help clarify the structure of the proof, we now give a brief outline of what is to follow.  As described above, in order for the exponent of $y_1$ to be correct, we need to shift the line of integration in the variable $\alpha_1$ from $\re(\alpha_1)=0$ to $\re(\alpha_1)=\kappa_1=a_1-2\delta$.  In Lemma~\ref{lem:shiftalpha1poles}, we identify the poles that are passed in making this shift, of which there are three types.  

After establishing this lemma, we bound the shifted integral \eqref{eq:pTr_1deltashift}.  The residue at any one of the three types of poles is essentially the same as that at any of the other two up to a simple transformation which doesn't effect the rest of the argument.  Hence, it suffices to pick any one of the types of poles from Lemma~\ref{lem:shiftalpha1poles}.  Having made a choice, we then show that it is, similar to before, necessary to shift in the variable $\alpha_2$ in order for the exponent of $y_1$ to be correct.  Unlike before, however, we are able to show in Lemma~\ref{lem:pTr_1alpha2shift} that in shifting $\alpha_2$, no further poles are encountered.  Hence, it suffices to bound the shifted terms, which we then do.  Since both \eqref{eq:pTr_1deltashift} and the shifted residue term satisfy \eqref{Prop-s1-Residue}, taken together, this proves the 
proposition.
 
\begin{lem}\label{lem:shiftalpha1poles}
In shifting the line of integration in $\alpha_1$ from $\re(\alpha_1)=0$ to $\re(\alpha_1)=a_1-2\delta$, poles are crossed only at the points $\alpha_1=q$ for
\begin{equation}\label{eq:qpoles}
 q \in\left\{\left.\begin{array}{c} -s_2-\alpha_2-2(r_2-\delta_2),\\-s_2-\alpha_3-2(r_2-\delta_2),\\-s_3-\alpha_2-\alpha_3-2(r_3-\delta_3)\end{array}\right| \begin{array}{c} 0 \leq \delta_j \leq r_j  \\ \delta_j\leq r_1-\delta \end{array}\right\}.
\end{equation}
\end{lem}

\begin{proof}
We first consider the ratio
 \[ \frac{\Gamma(\frac{\alpha_j-\alpha_1}{2}-\delta)}{\Gamma(\frac{\alpha_j-\alpha_1}{2})\Gamma(\frac{\alpha_1-\alpha_j}{2})}, \]
for $j=2,3,$ or $4$.  The numerator has a simple pole when $\frac{\alpha_j-\alpha_1}{2}-\delta$ is a nonpositive integer, but then $\frac{\alpha_j-\alpha_1}{2}$ is also an integer, so that there is a pole in the denominator as well. In other words, this ratio is holomorphic in $\alpha_1$.

For each of the terms
 \[ \Gamma\bigg(\frac{2+R+\alpha_j-\alpha_k}{4}\bigg), \qquad \mbox{($1\leq j\neq k\leq 4$)} \]
if $R$ is sufficiently large\footnote{We will come back to this later to determine exactly how big $R$ must be.} again no poles are crossed.

Examining the term $\Gamma(\frac{s_2+\alpha_1+\alpha_j}{2})$ with $j=2$ 
or $j=3$, assuming that $\re(\alpha_j)=0$ and $\re(s_2)=-a_2$, we see that poles occur at $\re(\alpha_1)=\kappa_1$, where $\kappa_1=a_2-2(r_2-\delta_2)$ for  nonnegative integers $r_2-\delta_2$. (Note that we have chosen this notation because it enumerates the $\delta_2$-th pole that is crossed as one starts at $\re(\alpha_1)=0$ and moves to $\re(\alpha_1)=a_1-2\delta$.)   Hence 
 \[ 0 \leq \re(\alpha_1) = a_2-2(r_2-\delta_2)  <2r_2+1-2(r_2-\delta_2) 
=1+2\delta_2\Longrightarrow \delta_2>-1/2 \Longrightarrow\delta_2\ge0 \](since $\delta_2$ is an integer).  And of course $\delta_2\le r_2$, since we are assuming that $r_2-\delta_2$ is nonnegative.

However, since the shift in $\alpha_1$ is going from $\kappa_1=0$ to $\kappa_1=a_1-2\delta$,  the largest value of $\delta_2$ that can in fact 
yield a pole is the one for which
 \[  a_2-2(r_2-\delta_2)\le a_1-2\delta< a_2-2(r_2-(\delta_2+1)) .\]
It's readily checked that these inequalities imply $\delta_2\le r_1-\delta$.

The details for evaluating the poles of $\Gamma(\frac{s_3-\alpha_4}{2})$ are similar.  We leave the details to the reader.
\end{proof}

\vskip 12pt

\noindent
\underline{\bf{Step 1}: Bounding the shifted integral $p_{T,R}^{1,\delta}$}
\vskip 10pt

Before dealing with the residues, we first bound $p_{T,R}^{1,\delta}(y;(-a_2,-a_3),(\kappa_1,0,0))$.

Notice that we may interchange $\tau_1,\tau_2$ and $\tau_3$ without affecting the integrand in \eqref{eq:pTr_1deltashift}.  Therefore, we may assume
 \[ -\tau_1-\tau_2-\tau_3 = \tau_4 \leq \tau_3 \leq \tau_2.  \]
It follows that the exponential factor is
\begin{align*}
 \mathcal{E} & = -2\tau_1-4\tau_2-2\tau_3-\lvert \xi_2+\xi_3+\tau_1\vert 
 + \lvert \xi_2+\tau_1+\tau_2\rvert 
 + \lvert \xi_2+\tau_1+\tau_3\rvert 
 \\ & \quad
 + \lvert \xi_2+\tau_1+\tau_4\rvert 
 + \lvert \xi_3-\tau_2\rvert 
 + \lvert \xi_3-\tau_3\rvert 
 + \lvert \xi_3-\tau_4\rvert,
\end{align*}
and using the method of Lemma~\ref{lem:ExpZeroSet} it is easy to show that there are two possible exponential zero sets:
\begin{align*}
 \mathcal{R}_+ & :=  \left\{ (-a_2+i\xi_2,-a_3+i\xi_3)\in \C^2\left| \begin{array}{cc} -\tau_1-\tau_3\leq \xi_2\leq \tau_2+\tau_3 \\ \tau_3\leq \xi_3 \leq \tau_2 \end{array}\right.\right\} \\
 \mathcal{R}_- & :=  \left\{ (-a_2+i\xi_2,-a_3+i\xi_3)\in \C^2\left| \begin{array}{cc} -\tau_1-\tau_2\leq \xi_2\leq -\tau_1-\tau_3 \\ \tau_4\leq 
\xi_3 \leq \tau_3 \end{array}\right.\right\}
\end{align*}
The change of variables $(\xi_2,\xi_3,\tau_2,\tau_4)\mapsto(\xi_3,\xi_2,\tau_4,\tau_2)$ relates $\mathcal{R}_+$ and $\mathcal{R}_-$, so it suffices to consider just the case of $\mathcal{R}_+$.

We replace the Gamma factors with their corresponding polynomial terms to 
obtain 
\begin{align} \label{pTr_1deltaShifted}&
  \left| p_{T,R}^{1,\delta}(y;(-a_2,-a_3),\kappa) \right| \\\nn& \ll  
  y_1^{\frac32+a_1}y_2^{2+a_2}y_3^{\frac32+a_3}
  \cdot  T^{\varepsilon+R+3\delta} \,
  \hskip -8pt \iiint\limits_{\substack{ \abs{\tau_1},\abs{\tau_2},\abs{\tau_3} \leq T^{1+\varepsilon} \\ 0 \leq \tau_1+\tau_2+2\tau_3 \leq \tau_1 + 2\tau_2 +\tau_3 }} 
  \hskip -8pt
  \big(1+\lvert\tau_1-\tau_2\rvert\big)^{\frac{R+1}{2}-\frac{a_1}{2}}  \big(1+\lvert\tau_1-\tau_3\rvert\big)^{\frac{R+1}{2}-\frac{a_1}{2}}
  \\ \nonumber &  
  \cdot
  \big(1+\lvert\tau_1-\tau_4\rvert\big)^{\frac{R+1}{2}-a_1+\delta}
  \big( 1+\tau_2-\tau_3 \big)^{1+\frac{R}{2}}
  \big( 1+\tau_2-\tau_4 \big)^{1+\frac{R}{2}}
  \big( 1+\tau_3-\tau_4 \big)^{1+\frac{R}{2}}
  \\ \nonumber & 
  \cdot
 \iint\limits_{\substack{-\tau_1-\tau_3\le\xi_2\le \tau_2+\tau_3\\ \tau_3\le\xi_3\le \tau_2} } 
  \frac{\big(1+\lvert\xi_2+\tau_1+\tau_2\rvert\big)^{-\frac{1+a_2}{2}+\frac{a_1}{2}-\delta}
        \big(1+\lvert\xi_2+\tau_1+\tau_3\rvert \big)^{-\frac{1+a_2}{2}+\frac{a_1}{2}-\delta}
      }{\big(1+\lvert\xi_2+\xi_3+\tau_1\rvert\big)^{-\frac{1+a_2+a_3-a_1}{2}} }
  \\ \nonumber &  
  \cdot        \big(1+ \lvert \xi_2-\tau_2-\tau_3\rvert \big)^{-\frac{1+a_2}{2}}
        \big(1+\lvert\xi_3-\tau_2\rvert\big)^{-\frac{1+a_3}{2}}
        \big(1+\lvert\xi_3-\tau_3\rvert\big)^{-\frac{1+a_3}{2}}
        \big(1+\lvert\xi_3-\tau_4\rvert\big)^{-\frac{1+a_3}{2}+\frac{a_1}{2}-\delta}
        \\\nn&\cdot d\xi_2 d\xi_3\; d\tau.
\end{align}We make the change of variables
 $ \xi_2 \mapsto \xi_2-\tau_1-\tau_3, \  \xi_3\mapsto \xi_3+\tau_3,  \mbox{ and }   T_j=\tau_j-\tau_{j+1}\ (1\le j\le 3), $
to get \begin{align}\label{ITintegral}
\left| p_{T,R}^{1,\delta}(y;(-a_2,-a_3),\kappa)\right|& \ll y_1^{\frac32+a_1}y_2^{2+a_2}y_3^{\frac32+a_3}\,  T^{\varepsilon+R+3\delta} 
  \hskip -12pt \iiint\limits_{0\leq \abs{T_1},T_2,T_3 \leq T^{1+\varepsilon}} \hskip-15pt
  \big(1+\lvert T_1 \rvert\big)^{ \frac{1+R-a_1}{2} }\big( 1+T_2 \big)^{1+\frac{R}{2}}  \big( 1+T_3 \big)^{1+\frac{R}{2}}
  \\ \nonumber & \hskip -100pt
  \cdot
   \big(1+\lvert T_1+T_2\rvert\big)^{\frac{1+R-a_1}{2}}
  \big(1+\lvert T_1+T_2+T_3 \rvert\big)^{\frac{1+R}{2} -a_1 +\delta}
%  \\ \nonumber & \hskip 42pt
%  \cdot
  \big( 1+T_2+T_3 \big)^{1+\frac{R}{2}}
  \\ \nonumber & 
  \hskip -57pt
  \cdot
 \int\limits_{\xi_3=0}^{T_2} \int\limits_{\xi_2=0}^{T_3}
  \frac{\big(1+\xi_2+T_2\big)^{-\frac{1+a_2}{2}+\frac{a_1}{2}-\delta}
        \big(1+\xi_2\big)^{-\frac{1+a_2}{2}+\frac{a_1}{2}-\delta}
        \big(1+T_3-\xi_2\big)^{-\frac{1+a_2}{2}}
      }{\big(1+\xi_2+\xi_3\big)^{-\frac{1+a_2+a_3-a_1}{2}} }
  \\ \nonumber & \hskip -68pt
  \cdot
        \big(1+T_2-\xi_3\big)^{-\frac{1+a_3}{2}}
        \big(1+\xi_3 \big)^{-\frac{1+a_3}{2}}
        \big(1+\xi_3+T_3\big)^{-\frac{1+a_3}{2}+\frac{a_1}{2}-\delta}
        \, d\xi_2\, d\xi_3\, dT_1\, dT_2\, dT_3.
\end{align}

We first examine \eqref{ITintegral} in the case  $r_2=r_3=0$.  In this case, we have $$1<1 +\varepsilon\le 2(\delta+1)-1+\varepsilon\le2r_1-1+\varepsilon\le a_1\le 2r_1+1-\varepsilon,$$and $ -1<-1+\varepsilon\le  a_j\le  1-\varepsilon<1\quad(2\le j\le 3).$
If we now denote the integral in $\xi_2$ and $\xi_3$, in \eqref{ITintegral}, by $\mathcal{I}(\tau,a,\delta)$ (where $\tau=(\tau_1,\tau_2,\tau_3)$ and $a=(a_1,a_2,a_3)$), then the above estimates on $a_1,a_2,a_3$  imply that \begin{align*} &\mathcal{I}(\tau,a,\delta) \ll T \int_{\xi_2=0 
} ^{ T_3 }   \int_{ \xi_3=0}^{T_2}  
  \big(1+ \xi_2+T_2\big)^{ \frac{a_1}{2}-\delta}
        \big(1+ \xi_2 \big)^{ \frac{a_1}{2}-\delta}
        \big(1+ \xi_2+\xi_3 \big)^{\frac{1 -a_1}{2}}  \big(1+ \xi_3+T_3\big)^{ \frac{a_1}{2}-\delta}
\; d\xi_3\; d\xi_2.\end{align*} {We note that $a_1/2-\delta>0$ and  $1-a_1<0$.  Two applications of Lemma \ref{lemmaIntegralBound3} (to the integrals in $\xi_3$ and $\xi_2$ respectively) then yield
 \begin{align*} &\mathcal{I}(\tau,a,\delta) \ll        T  \,        \big(1+ T_2 \big)^{\varepsilon+1}  \big(1+ T_2+T_3\big)^{a_1-2\delta} \big(1+ T_3 \big)^{\varepsilon+ \frac{3}{2}-\delta}  .\end{align*}Putting this back into \eqref{ITintegral} yields
\begin{align}\label{ITintegral001}
&\left| p_{T,R}^{1,\delta}(y;(-a_2,-a_3),\kappa)\right|  \ll y_1^{\frac32+a_1}y_2^{2+a_2}y_3^{\frac32+a_3}\,  T^{\varepsilon+1+R+3\delta} 
  \hskip -12pt \iiint\limits_{0\leq \abs{T_1},T_2,T_3 \leq T^{1+\varepsilon}} \hskip-15pt
  \big(1+\lvert T_1 \rvert\big)^{ \frac{1+R-a_1}{2} }\big( 1+T_2 \big)^{\varepsilon+2+\frac{R}{2}}  
  \\ \nonumber & 
  \cdot \big( 1+T_3 \big)^{\varepsilon+\frac52+\frac{R}{2}-\delta}
   \big(1+\lvert T_1+T_2\rvert\big)^{\frac{1+R-a_1}{2}}
  \big(1+\lvert T_1+T_2+T_3 \rvert\big)^{\frac{1+R}{2} -a_1 +\delta}
%  \\ \nonumber & \hskip 42pt
%  \cdot
  \big( 1+T_2+T_3 \big)^{1+\frac{R}{2}+a_1-2\delta}
  \\ \nonumber & \cdot   \, dT_1\,dT_2\, dT_3.
\end{align}
Assuming that $R$ is large enough to insure that the exponents in the above integrand are positive, we find that
\begin{align*}
\left| p_{T,R}^{1,\delta}(y;(-a_2,-a_3),\kappa)\right| & \ll y_1^{\frac32+a_1}y_2^{2+a_2}y_3^{\frac32+a_3}\,  T^{\varepsilon+1+R+3\delta} 
  \cdot T^3\cdot T^{\varepsilon+3R+7-a_1-2\delta}
  \\\nonumber&=y_1^{\frac32+a_1}y_2^{2+a_2}y_3^{\frac32+a_3}\, T^{\varepsilon+4R+11+\delta-a_1} 
  \\ \nonumber& \ll y_1^{\frac32+a_1}y_2^{2+a_2}y_3^{\frac32+a_3}\, {T^{\varepsilon+4R+11-r_1} },
  \end{align*}the last step because $\delta\le r_1-1$ and $2r_1-1+\varepsilon\le a_1$ imply $\delta-a_1\le -r_1$.}   Again, we're presently in the 
case $r_2=r_3=0$, whence $$11-r_1=9+\delta_{0,r_2}+\delta_{0,r_3}-r_1-r_2-r_3.$$
 So $p_{T,R}^{1,\delta}(y;(-a_2,-a_3),\kappa)$ satisfies the bound on $p_{T,r}^{1,\delta}(y;(-a_2,-a_3))$ given by \eqref{Prop-s1-Residue}, in this case.
    
   We now demonstrate that this is true for $r_2+r_3\ge r_1$ as well.  Because $$1+a_i>2r_1\ge0  \quad (1\le i\le 3)\hbox{\quad and \quad}\frac{a_1}{2}-\delta>\frac{2r_1-1}{2}-(r_1-1)=1/2>0,$$ we have the bounds
\begin{align*}
\big(1+\xi_j\big)^{-\frac{1+a_j}{2}} & \leq \big(1+\xi_j\big)^{-r_j}\quad 
(2\le j\le 3) ,
 \\
 \big(1+\xi_j+T_j\big)^{-\frac{1+a_j}{2}+\frac{a_1}{2}-\delta} & \leq \big(1+T_j\big)^{- r_j}\big(1+T_2+T_3\big)^{\frac{a_1}{2}-\delta} 
\quad (2\le j\le 3).
\end{align*}Further,\begin{align*} \frac{\big(1+\xi_2\big)^{\frac{a_1}{2}-\delta}}{\big(1+\xi_2+\xi_3\big)^{\frac{a_1}{2}-\frac{1+a_2+a_3}{2}}}
 &=\biggl(\frac{1+\xi_2}{1+\xi_2+\xi_3} \biggr)^{\frac{a_1}{2}-\delta}\ 
\big(1+\xi_2+\xi_3\big)^{\frac{1+a_2+a_3}{2}-\delta} 
  \ll  T^{\varepsilon+\frac{1+a_2+a_3}{2}-\delta},\end{align*}the last step because $r_2+r_3\ge r_1\Rightarrow \frac{1+a_2+a_3}{2}-\delta>0$.
 So \eqref{ITintegral} yields
 \begin{align*} 
&\left| p_{T,R}^{1,\delta}(y,(-a_2,-a_3),\kappa)\right|  \ll y_1^{\frac32+a_1}y_2^{2+a_2}y_3^{\frac32+a_3}\, T^{\varepsilon+R+2\delta  +\frac{1+a_2+a_3}{2}}
  \hskip -12pt \iiint\limits_{0\leq \abs{T_1},T_2,T_3 \leq T^{1+\varepsilon}} \hskip-15pt
  \big(1+\lvert T_1 \rvert\big)^{ \frac{1+R-a_1}{2} }\big( 1+T_2 \big)^{1+\frac{R}{2}-r_2}  
  \\ \nonumber &    \cdot\big( 1+T_3 \big)^{1+\frac{R}{2}-r_3}
   \big(1+\lvert T_1+T_2\rvert\big)^{\frac{1+R-a_1}{2}}
  \big(1+\lvert T_1+T_2+T_3 \rvert\big)^{\frac{1+R}{2} -a_1 +\delta}
%  \\ \nonumber & \hskip 42pt
%  \cdot
  \big( 1+T_2+T_3 \big)^{1+\frac{R}{2}+a_1-2\delta}
  \\ \nonumber & 
  \cdot
 \int\limits_{\xi_3=0}^{T_2}             \big(1+\xi_3 \big)^{-r_3}   \big(1+T_2-\xi_3\big)^{-r_3}\, d\xi_3 \int\limits_{\xi_2=0}^{T_3}
         \big(1+\xi_2\big)^{-r_2}
        \big(1+T_3-\xi_2\big)^{-r_2} \,d\xi_2         \, dT_1\, dT_2\, dT_3
\end{align*}or, by Lemma \ref{lemmaIntegralBound} applied to the above integrals in $\xi_2$ and $\xi_3$,
  
 \begin{align}\label{ITintegral3}
\left| p_{T,R}^{1,\delta}(y,(-a_2,-a_3),\kappa)\right|  \ll &y_1^{\frac32+a_1}y_2^{2+a_2}y_3^{\frac32+a_3}\,  T^{\varepsilon+R+2\delta+\frac{1+a_2+a_3}{2}} \\\nonumber&\cdot
  \iiint\limits_{0\leq \abs{T_1},T_2,T_3 \leq T^{1+\varepsilon}} \hskip-15pt
  \big(1+\lvert T_1 \rvert\big)^{ \frac{1+R-a_1}{2} }\big( 1+T_2 \big)^{1+\delta_{r_2,0}+\frac{R}{2}-r_2-r_3}  \big( 1+T_3 \big)^{1+\delta_{r_3,0}+\frac{R}{2}-r_2-r_3}
  \\ \nonumber &    \cdot   \big(1+\lvert T_1+T_2\rvert\big)^{\frac{1+R-a_1}{2}}
  \big(1+\lvert T_1+T_2+T_3 \rvert\big)^{\frac{1+R}{2} -a_1 +\delta}
%  \\ \nonumber & \hskip 42pt
%  \cdot
  \big( 1+T_2+T_3 \big)^{1+\frac{R}{2}+a_1-2\delta} 
       \\\nonumber &\cdot dT_1\, dT_2\, dT_3.
\end{align} 
Assuming again that $R$ is sufficiently large that the exponents in the above integrals in $T_1$, $T_2$, and $T_3$ are positive, we find that
\begin{align*} \left| p_{T,R}^{1,\delta}(y,(-a_2,-a_3),\kappa)\right|    & \ll
 y_1^{\frac32+a_1}y_2^{2+a_2}y_3^{\frac32+a_3}T^{\varepsilon+4R+\delta+\frac{1+a_2+a_3}{2}}\cdot T^3\cdot T^{\varepsilon+3R+\frac{9}{2}-a_1+\delta_{0,r_2}+\delta_{0,r_3}-2r_2-2r_3} .\end{align*}Since $a_j\le 2r_j+1$ for 
$2\le j\le3$ and $\delta-a_1\le r_1-1-( 2r_1-1)=-r_1$, we then conclude 
that
\begin{align*} \left| p_{T,R}^{1,\delta}(y,(-a_2,-a_3),\kappa)\right|    & \ll
 y_1^{\frac32+a_1}y_2^{2+a_2}y_3^{\frac32+a_3}T^{\varepsilon+4R+ 9-r_1+\delta_{0,r_2}+\delta_{0,r_3}-r_2-r_3} .\end{align*}  So, in the case $r_2+r_3\ge r_1$, we again have a bound on $p_{T,R}^{1,\delta}(y,(-a_2,-a_3),\kappa)$ that is consistent with the statement of Proposition \ref{prop:singleresidue1}.

\vskip 10pt

\noindent
\underline{\bf{Step 2}: Bounding the residue term at $\alpha_1=-s_2-\alpha_2-2(r_2-\delta_2)$}
\vskip 10pt

The next step of the proof is to show that all of the residues at the poles $\alpha_1=-s_2-\alpha_2-2(r_2-\delta_2)$ contribute smaller bounds.  
The residue we must bound is:
\begin{align}\label{eq:alpha1residue}
 &  \iiint\limits_{\substack{\re(\alpha_2)=0\\ \re(\alpha_3)=0 \\ \re(s_2)=-a_2}} e^{\frac{(s_2+\alpha_2+2(r_2-\delta_2))^2+\alpha_2^2+\alpha_3^2+(s_2-\alpha_3+2(r_2-\delta_2))^2}{2T^2}}\; \int\limits_{\re(s_3)=-a_3}y_1^{\frac32-s_2-\alpha_2-2(r_2-\delta_2-\delta)}y_2^{2-s_2}y_3^{\frac32-s_3} 
  \\ \nonumber   &
  \cdot
  \frac{\Gamma\(\frac{2+R+s_2-2\alpha_3+2(r_2-\delta_2)}{4}\)
        \Gamma\(\frac{2+R-s_2+2\alpha_3-2(r_2-\delta_2)}{4}\)
      }{\Gamma\(\frac{-s_2+2\alpha_3-2(r_2-\delta_2)}{2}\)
        \Gamma\(\frac{s_2-2\alpha_3+2(r_2-\delta_2)}{2}\)}
  \cdot
  \frac{\Gamma\(\frac{2+R+\alpha_2-\alpha_3}{4}\)
        \Gamma\(\frac{2+R+\alpha_3-\alpha_2}{4}\)
      }{\Gamma\(\frac{\alpha_2-\alpha_3}{2}\)}
  \\ \nonumber  &  \cdot
  \frac{\Gamma\(\frac{2+R-s_2+\alpha_2+\alpha_3-2(r_2-\delta_2)}{4}\)
        \Gamma\(\frac{2+R+s_2-\alpha_2-\alpha_3+2(r_2-\delta_2)}{4}\)
      }{\Gamma\(\frac{-s_2+\alpha_2+\alpha_3-2(r_2-\delta_2)}{2}\)}
  \cdot
  \frac{\Gamma\(\frac{s_2+2\alpha_2}{2}+r_2-\delta_2-\delta\)
      }{\Gamma\(\frac{s_2+2\alpha_2}{2}+r_2-\delta_2\)}
  \\ \nonumber & \cdot
  \frac{\Gamma\(\frac{2+R-s_2-2\alpha_2-2(r_2-\delta_2)}{4}\)
        \Gamma\(\frac{2+R+s_2+2\alpha_2+2(r_2-\delta_2)}{4}\)
      }{\Gamma\(\frac{-s_2-2\alpha_2-2(r_2-\delta_2)}{2}\)}
  \cdot
  \frac{\Gamma\(\frac{s_2+\alpha_2+\alpha_3}{2}+r_2-\delta_2-\delta\)
      }{\Gamma\(\frac{s_2+\alpha_2+\alpha_3}{2}+r_2-\delta_2\)}
  \\ \nonumber &   \cdot 
  \frac{\Gamma\(\frac{2+R-s_2-\alpha_2-\alpha_3-2(r_2-\delta_2)}{4}\)
        \Gamma\(\frac{2+R+s_2+\alpha_2+\alpha_3+2(r_2-\delta_2))}{4}\)
      }{\Gamma\(\frac{-s_2-\alpha_2-\alpha_3-2(r_2-\delta_2)}{2}\)}
  \cdot
  \frac{\Gamma\(\frac{2s_2+\alpha_2-\alpha_3+4(r_2-\delta_2)}{2}-\delta\)
      }{\Gamma\(\frac{2s_2+\alpha_2-\alpha_3+4(r_2-\delta_2)}{2}\)}
  \\ \nonumber &   \cdot
  \frac{\Gamma\(\frac{2+R-2s_2-\alpha_2+\alpha_3-4(r_2-\delta_2)}{4}\)
        \Gamma\(\frac{2+R+2s_2+\alpha_2-\alpha_3+4(r_2-\delta_2)}{4}\)
      }{\Gamma\(\frac{-2s_2-\alpha_2+\alpha_3-4(r_2-\delta_2)}{2}\)}
  \cdot
  \frac{\Gamma\(\frac{\alpha_3-\alpha_2}{2}-r_2+\delta_2\)
      }{\Gamma\(\frac{\alpha_3-\alpha_2}{2}\) }
  \\ \nonumber &  \cdot  
  \frac{\Gamma\(\frac{s_2-\alpha_2-\alpha_3}{2}\) 
      }{\Gamma\(\frac{s_2-\alpha_2-\alpha_3}{2}+r_2-\delta_2\) }
 \cdot 
  \frac{\Gamma\(\frac{s_3-\alpha_2}{2}\)
        \Gamma\(\frac{s_3-\alpha_3}{2}\) 
        \Gamma\(\frac{s_3-s_2+\alpha_3}{2}-r_2+\delta_2\)
      }{\Gamma(\frac{s_3-\alpha_2}{2}-r_2+\delta_2+\delta) } \bigl(\fr\, Q_\delta(s,\alpha)\bigr)\big|_{\alpha_1=-s_2-\alpha_2-2(r_2-\delta_2)}
      \\ \nonumber&\cdot ds_2 \,ds_3\,d\alpha.\nn
\end{align}
We need to shift the line of integration in $\alpha_2$ so that the real part of the exponent of $y_1$ is $a_1$.  That is, we require that
 \[ 2r_1-a_1'=a_1 = a_2-\re(\alpha_2)-2(r_2-\delta_2-\delta) = a_2'-\re(\alpha_2)+2(\delta_2+\delta), \]
where $0<a_1',a_2'<1$.  In other words, given the bounds from Lemma~\ref{lem:shiftalpha1poles}, we shift the line in $\alpha_2$ to
\begin{equation}\label{eq:alpha2shift}
 \re(\alpha_2) = \kappa_2 := a_2'-a_1'-2(r_1-\delta-\delta_2).
\end{equation}
 
\begin{lem}\label{lem:pTr_1alpha2shift}
In shifting the line of integration in \eqref{eq:alpha1residue} in the variable $\alpha_2$ from $\re(\alpha_2)=0$ to $\re(\alpha_2)=\kappa_2$ as in \eqref{eq:alpha2shift}, no poles are crossed.  \end{lem}

\begin{proof}
As in the proof of Lemma~\ref{lem:shiftalpha1poles}, we assume that $R$ is sufficiently large.  Then for each of the Gamma factors involving $R$, no poles are crossed.  Further, for reasons also described in that proof, 
each of the factors of the form
$$  \frac{ \Gamma\(-\frac{\alpha_2}{2}+z-\delta\) 
      }{\Gamma\(-\frac{\alpha_2}{2}+z\) \Gamma\(\frac{\alpha_2}{2}-z\) } $$is holomorphic in $\alpha_2$.  Thus, in moving the line of integration in $\alpha_2$, only the term $\Gamma\(\frac{s_3-\alpha_2}{2}\)$ might contribute poles.
 
Regarding this term, we consider two cases.  The first is when $r_2\ge r_1$.  In this case, since $\delta_2\le r_1-\delta$ by Lemma \ref{lem:shiftalpha1poles}, we have $\delta_2\le r_2-\delta$, so the pole in the denominator factor $\Gamma(\frac{s_3-\alpha_2}{2}-r_2+\delta+\delta_2)$      cancels the pole in $\Gamma\(\frac{s_3-\alpha_2}{2}\)$.

Finally, then, we need to consider the case $r_2<r_1$.  Of the triples $(r_1,r_2,r_3)$ under consideration, the only one satisfying this criterion 
is the triple $(r,0,0)$.  Recall that, in the case of this triple, we are 
assuming that $\re(a_3)<0.$
Also, note that our integral in $s_3$, in (\ref{eq:alpha1residue}), is over the line $\re(s_3)=-a_3$.  So on our original line of integration in 
$\alpha_2$, namely $\re(\alpha_2)=0$,  we have that $\re(s_3-\alpha_2)>0$.  But note that we are moving this line to the left,  since (\ref{eq:alpha2shift}) and the fact that $\delta_2=0$ (since $r_2=0$) together imply that the new line of integration is$$ \re(\alpha_2)  = a_2'-a_1'-2(r_1-\delta)\le a_2'-a_1'-2<0.$$Therefore, in moving this line we are only increasing the real part of $s_3-\alpha_2$, and consequently are not passing any poles of $\Gamma(\frac{s_3-\alpha_2}{2}).$
\end{proof}

In order to bound \eqref{eq:alpha1residue}, we first remark that the exponential factor is
\begin{align*}
 \mathcal{E} = - \lvert \xi_2-2\tau_3\rvert + \lvert \xi_3-\tau_3 \rvert + \lvert \xi_3-\xi_2+\tau_3 \rvert,
\end{align*}
and it follows that the exponential zero set is
 \[  \xi_2-\tau_3 \leq \xi_3\leq \tau_3\qquad \mbox{or} \qquad  \tau_3 \leq \xi_3 \leq \xi_2-\tau_3. \]
Using the first of these\footnote{For the other exponential zero set the answer is virtually identical.}, and assuming as before that $\tau_2\ge \tau_3$, it is not hard to see that the shifted version of \eqref{eq:alpha1residue}, to $\re(\alpha_2)=\kappa_2:=a_2-a_1-2(r_2-\delta-\delta_2)$, is bounded by
\begin{align}\label{new-a}
& 
y_1^{\frac32+a_1}y_2^{2+a_2}y_3^{\frac32+a_3}
\, T^{\varepsilon+R+3\delta}
  \iint\limits_{-\infty<\tau_3\le \tau_2<\infty}
    \big(1+ \tau_2-\tau_3\big)^{\frac{1+R+a_1-a_2}{2} -\delta} 
    \int\limits_{\xi_2=-\infty}^{2\tau_3} e^{-\frac{(\xi_2+\tau_2)^2+\tau_2^2+\tau_3^2+(\xi_2-\tau_3)^2}{2T^2}}\;
 \\ \nn& \hskip 12pt
 \cdot
    \big(1+ 2\tau_3-\xi_2\big)^{1+\frac{R}{2}}
    \big(1+\abs{\xi_2+2\tau_2}\big)^{ \frac{1+R-2a_1+a_2}{2}-r_2+\delta_2+\delta}
    \big(1+\tau_2+\tau_3-\xi_2 \big)^{ \frac{1+R+a_1-2a_2}{2}+r_2-\delta-\delta_2}
 \\ \nn& \hskip 84pt
 \cdot
    \big(1+\abs{2\xi_2+\tau_2-\tau_3}\big)^{ \frac{1+R-a_1-a_2}{2} + r_2-\delta_2 }
    \big(1+\abs{\xi_2+\tau_2+\tau_3}\big)^{\frac{1+R-a_1}{2} }
 \\ \nn& \hskip 30pt
 \cdot
    \int\limits_{\xi_3=\xi_2-\tau_3}^{ \tau_3}
    \big(1+\tau_2- \xi_3\big)^{ r_2-\delta_2-\delta}
    \big(1+\tau_3-\xi_3\big)^{-\frac{1+a_3}{2}}
    \big(1+ \xi_3-(\xi_2-\tau_3)\big)^{\frac{a_2-a_3-1}{2}-r_2+\delta_2} d\xi_3\; d\xi_2\; d\tau_2\; d\tau_3.
\end{align}
Now $$\big(1+\tau_2-\xi_3)^{-\delta}\le(1+\tau_2-\tau_3\big)^{-\delta} \qquad \mbox{and} \qquad   \big(1+ \xi_3-(\xi_2-\tau_3)\big)^{\delta_2+ 1/2} \le  \big(1+2\tau_3-\xi_2\big)^{\delta_2+1/2} ,$$so we find that \eqref{new-a} is 
\begin{align}\label{new-b}
& \ll
y_1^{\frac32+a_1}y_2^{2+a_2}y_3^{\frac32+a_3}
\,  T^{\varepsilon+R+3\delta}
  \iint\limits_{-\infty<\tau_3\le \tau_2<\infty}
    \big(1+ \tau_2-\tau_3\big)^{\frac{1+R+a_1-a_2}{2} -2\delta} 
    \int\limits_{\xi_2=-\infty}^{2\tau_3} e^{-\frac{(\xi_2+\tau_2)^2+\tau_2^2+\tau_3^2+(\xi_2-\tau_3)^2}{2T^2}}\;
 \\ \nn& \hskip 12pt
 \cdot
    \big(1+ 2\tau_3-\xi_2\big)^{ \frac{3+R}{2}+\delta_2}
    \big(1+\abs{\xi_2+2\tau_2}\big)^{ \frac{1+R-2a_1+a_2}{2}-r_2+\delta_2+\delta}
    \big(1+\tau_2+\tau_3-\xi_2 \big)^{ \frac{1+R+a_1-2a_2}{2}+r_2-\delta-\delta_2}
 \\ \nn& \hskip 84pt
 \cdot
    \big(1+\abs{2\xi_2+\tau_2-\tau_3}\big)^{ \frac{1+R-a_1-a_2}{2} + r_2-\delta_2 }
    \big(1+\abs{\xi_2+\tau_2+\tau_3}\big)^{\frac{1+R-a_1}{2} }
 \\ \nn& \hskip 30pt
 \cdot
    \int\limits_{\xi_3=\xi_2-\tau_3}^{ \tau_3}
    \big(1+\tau_2- \xi_3\big)^{ r_2-\delta_2}
    \big(1+\tau_3-\xi_3\big)^{-\frac{1+a_3}{2}}
    \big(1+ \xi_3-(\xi_2-\tau_3)\big)^{\frac{a_2-a_3-1}{2}-r_2} d\xi_3\; d\xi_2\; d\tau_2\; d\tau_3.
\end{align}
Since $r_2-\delta_2 \ge0$, $\frac{1+a_3}{2}>0$, and $\frac{a_2-a_3-2}{2}-r_2 <0$  we find, using Lemma A.0.3, that the integral in $\xi_3$, in \eqref{new-b}, is 
\begin{align*}
 &\ll
      \big(1+ 2\tau_3-\xi_2\big)^{\delta_{0,r_3}-\frac{1+a_3}{2}}\cdot 
    \big(1+{\tau_2+\tau_3-\xi_2}\big)^{r_2-\delta_2}.
\end{align*}
But then \eqref{new-b} is 
\begin{align}\label{new-c}
& \ll
y_1^{\frac32+a_1}y_2^{2+a_2}y_3^{\frac32+a_3}
\, T^{\varepsilon+R+3\delta}
  \iint\limits_{-\infty<\tau_3\le \tau_2<\infty}
    \big(1+ \tau_2-\tau_3\big)^{\frac{1+R+a_1-a_2}{2} -2\delta} 
    \int\limits_{\xi_2=-\infty}^{2\tau_3} e^{-\frac{(\xi_2+\tau_2)^2+\tau_2^2+\tau_3^2+(\xi_2-\tau_3)^2}{2T^2}}\;
 \\ \nn&  
 \cdot
    \big(1+ 2\tau_3-\xi_2\big)^{ \frac{2+R-a_3}{2}+\delta_2+\delta_{r_3,0}}
    \big(1+\abs{\xi_2+2\tau_2}\big)^{ \frac{1+R-2a_1+a_2}{2}-r_2+\delta_2+\delta}
    \big(1+\tau_2+\tau_3-\xi_2 \big)^{ \frac{1+R+a_1-2a_2}{2}+2r_2-\delta-2\delta_2}
 \\ \nn& 
 \cdot
    \big(1+\abs{2\xi_2+\tau_2-\tau_3}\big)^{ \frac{1+R-a_1-a_2}{2} + r_2-\delta_2 }
    \big(1+\abs{\xi_2+\tau_2+\tau_3}\big)^{\frac{1+R-a_1}{2} }
\, d\xi_2\, d\tau_2\, d\tau_3.
\end{align}Substituting $\tau_j\to \tau_jT$ ($1\le j\le 2$) and $\xi_3\to\xi_3T$, we find that \eqref{new-c} is 
\begin{align*} 
& \ll
y_1^{\frac32+a_1}y_2^{2+a_2}y_3^{\frac32+a_3}
\,T^{\varepsilon+R+3\delta+3+\frac{7}{2}+3R-\frac{2a_1+3a_2+a_3}{2}+2r_2-2\delta-\delta_2+\delta_{r_3,0}},  
 \end{align*}which, because $-a_j<1-2r_j$ for $1\le j\le 3$, is 
 \begin{align*} 
& \ll
y_1^{\frac32+a_1}y_2^{2+a_2}y_3^{\frac32+a_3}
\, T^{\varepsilon+4R+\delta+\frac{19}{2}-2r_1-r_2-r_3-\delta_2+\delta_{r_3,0}}\ll
y_1^{\frac32+a_1}y_2^{2+a_2}y_3^{\frac32+a_3}
\,  T^{\varepsilon +4R+\frac{17}{2}-r_1-r_2-r_3+\delta_{r_3,0}},  
 \end{align*}
where in the final step we have used the facts that $r_1-\delta\ge1$ and $\delta_2\ge0$.    So the residue term (\ref{eq:alpha1residue}) also satisfies the bound given in Proposition \ref{prop:singleresidue1}.  The proof of that proposition is therefore complete.  
\end{proof}

\vskip 12pt

The other type of single residue term that we have to consider is 
\begin{align} \label{eq:pTr_2delta}
 p_{T,R}^{2,\delta}(y;(-a_1,-a_3))& =  \iiint\limits_{\re(\alpha_j)=0}  e^{\frac{\alpha_1^2+\cdots+\alpha_4^2}{2T^2}}\; \iint\limits_{\substack{\re(s_j)=-a_j\\j=1,3}} 
 y_1^{\frac32-s_1}y_2^{2-p_2}y_3^{\frac32-s_3}\; 
 \(\prod_{j=1}^2\prod_{k=3}^4 \textstyle{ \Gamma\left(\frac{\alpha_k-\alpha_j}{2}-\delta\right)}\) 
 \\ \nonumber & \cdot   \mathcal F_R(\alpha)\,
 \Gamma_R(\alpha) \, Q_\delta(s,\alpha) \,
 \textstyle{ \Gamma\(\frac{s_1+\alpha_1}{2}\) \Gamma\(\frac{s_1+\alpha_2}{2}\) \Gamma\(\frac{s_3-\alpha_3}{2}\) \Gamma\(\frac{s_3-\alpha_4}{2}\)}\; ds_1ds_3
\; d\alpha
\end{align}
where $Q_{\delta}$ is a polynomial (see Section \ref{poles-res}) of degree $\le3\delta$ and $p_2=-\alpha_1-\alpha_2-2\delta$.

\begin{prop}\label{prop:singleresidue2}
Let $r_2\geq 1$, $r_1,r_3\geq0$ be integers, and $0<\varepsilon<1$. Suppose $a_1,a_2,a_3$ satisfy the hypotheses of Theorem \ref{th:pTRbound}.  If 
$0\le \delta\le r_2-1$, then
\begin{equation}\label{eq:pTr_2bound}
 \left| p_{T,R}^{2,\delta}(y;(-a_1,-a_3)) \right| \ll y_1^{\frac32+a_1} y_2^{2+a_2} y_3^{\frac32+a_3} \, T^{\varepsilon+4R+9+\delta_{0,r_1}+\delta_{0,r_3}-(r_1+r_2+r_3)}.
\end{equation}
\end{prop}

\begin{proof}
We first rearrange the terms on the right hand side of \eqref{eq:pTr_2delta} as follows.
\begin{align*} 
 p_{T,R}^{2,\delta}(y;(-a_1,-a_3)) = & \iiint\limits_{\re(\alpha_j)=0}  e^{\frac{\alpha_1^2+\cdots+\alpha_4^2}{2T^2}}\; \iint\limits_{\substack{\re(s_j)=-a_j\\j=1,3}} 
 y_1^{\frac32-s_1}y_2^{2+\alpha_1+\alpha_2+2\delta}y_3^{\frac32-s_3}\; 
 \mathcal F_R(\alpha)\, Q_\delta(s,\alpha) \,
 \\ &
 \cdot\(\prod_{j=1}^2\prod_{k=3}^4 
 \frac{ 
       \Gamma(\frac{\alpha_k-\alpha_j}{2}-\delta)
     }{
       \Gamma(\frac{\alpha_k-\alpha_j}{2})
     } 
 \frac{
       \Gamma(\frac{2+R+\alpha_j-\alpha_k}{4})
       \Gamma(\frac{2+R+\alpha_k-\alpha_j}{4})
     }{
       \Gamma(\frac{\alpha_j-\alpha_k}{2}) 
     }
  \)  
  \\ &  
  \cdot
  \frac{
        \Gamma(\frac{2+R+\alpha_1-\alpha_2}{4})
        \Gamma(\frac{2+R+\alpha_2-\alpha_1}{4})
      }{
        \Gamma(\frac{\alpha_1-\alpha_2}{2})
        \Gamma(\frac{\alpha_2-\alpha_1}{2})
      }
  \cdot
  \frac{
        \Gamma(\frac{2+R+\alpha_3-\alpha_4}{4})
        \Gamma(\frac{2+R+\alpha_4-\alpha_3}{4})
      }{
        \Gamma(\frac{\alpha_3-\alpha_4}{2})
        \Gamma(\frac{\alpha_4-\alpha_3}{2})
      }
  \\ & 
  \cdot
 \textstyle{ \Gamma\(\frac{s_1+\alpha_1}{2}\) \Gamma\(\frac{s_1+\alpha_2}{2}\) \Gamma\(\frac{s_3-\alpha_3}{2}\) \Gamma\(\frac{s_3-\alpha_4}{2}\)}\; ds_1\,ds_3
\, d\alpha.
\end{align*}

The proof follows the very same outline as that of Proposition~\ref{prop:singleresidue1}.  First, in order for the exponent of $y_2$ to match that 
in the statement of the proposition, we will shift the integration in the 
$\alpha$ variables from real part zero to $\re(\alpha_j) = \kappa_j$, such that
 \[ a_2 = \re(\alpha_1+\alpha_2+2\delta) = \kappa_1+\kappa_2+2\delta \]
lies in the interval $(2r_2-1,2r_2+1)$.  We do this by defining
 \[ \kappa_1:= a_2-2\delta \qquad\mbox{ and } \qquad \kappa_2:=0=:\kappa_3. \]
Let $\tau_j=\im(\alpha_j)$.  Note that since the above integral is invariant under each of the change of coordinates 
 \[ (\alpha_1,\alpha_2) \mapsto (\alpha_2,\alpha_1),\quad \mbox{and} \quad (\alpha_3,\alpha_4)\mapsto(\alpha_4,\alpha_3), \]
we may assume that $\tau_1\geq \tau_2$ and $\tau_3\geq \tau_4$.  

As before, we now use Stirling's formula to write the Gamma factors in the above integrand as a product of linear and exponential terms.  The exponential factor is $e^{\frac{\pi}{4}\mathcal{E}(\xi,\tau)}$ where
 \[ \mathcal{E} = 2\tau_1+2\tau_3 -\lvert \xi_1+\tau_2\rvert-\lvert \xi_1+\tau_1\rvert -\lvert \xi_3+\tau_1+\tau_2+\tau_3\rvert-\lvert \xi_3-\tau_3\rvert, \]
and from this we readily deduce that the exponential zero set is
 \[ \mathcal{R} = \big\{ (\xi_1,\xi_3)\mid -\tau_1\leq \xi_1 \leq -\tau_2,\, \tau_4\leq \xi_3 \leq \tau_3 \big\}. \]
The polynomial factor is
\begin{align*}
 \mathcal{P}(\xi,\tau) & = \big(1+\lvert \tau_1-\tau_2\rvert \big)^{1+\frac{R}{2}}
  \big(1+\lvert \tau_3-\tau_4\rvert \big)^{1+\frac{R}{2}} 
 \biggl(\prod_{j=1}^2\prod_{k=3}^4 
       \big( 1 + \lvert \tau_k-\tau_j\rvert \big)^{\frac{1+R}{2}-(\frac{\kappa_j-\kappa_k}{2})-\delta}
  \biggr)  
  \\ \nonumber & 
 \cdot \big(1+\lvert \xi_1+\tau_1\rvert \big)^{-\frac{1+a_1-\kappa_1}{2}}
\cdot \big(1+\lvert \xi_1+\tau_2\rvert \big)^{-\frac{1+a_1-\kappa_2}{2}}
  \cdot
 \big(1+\lvert \xi_3-\tau_3\rvert \big)^{-\frac{1+a_3+\kappa_3}{2}}
 \cdot
 \big(1+\lvert \xi_3-\tau_4\rvert \big)^{-\frac{1+a_3+\kappa_4}{2}}.
\end{align*}

Plugging in the values of $\kappa_j$ given above, and bounding the resulting terms as before we find that
\begin{align} \label{eq:beforechoosingkappa} 
& \left| p_{T,R}^{2,\delta}(y;(-a_1,-a_3),\kappa) \right| \ll 
 y_1^{\frac32+a_1}y_2^{2+a_2}y_3^{\frac32+a_3}
 T^{\varepsilon+R+3\delta} \hskip-12pt
 \iint\limits_{\substack{ -T^{1+\varepsilon}\leq \tau_4\leq \tau_3\leq T^{1+\varepsilon} \\ -T^{1+\varepsilon} \leq \tau_2\leq \tau_1\leq T^{1+\varepsilon}}} \;
    \big(1+\abs{\tau_2-\tau_3}\big)^{\frac{1+R}{2}-\delta}
 \\ \nonumber &
 \hskip 60pt 
  \cdot
    \big(1+\tau_1-\tau_2)^{1+\frac{R}{2}}
    \big(1+\tau_3-\tau_4)^{1+\frac{R}{2}}
    \big(1+\abs{\tau_1-\tau_4}\big)^{\frac{3+R}{2}-2r_2+\delta}
 \\ \nonumber & 
  \cdot
    \big(1+\abs{\tau_1-\tau_3}\big)^{\frac{2+R}{2}-r_2}
    \big(1+\abs{\tau_2-\tau_4}\big)^{\frac{2+R}{2}-r_2}
    \hskip-5pt\int\limits_{\xi_1=-\tau_1}^{-\tau_2}
    \big(1+\abs{\xi_1+\tau_1}\big)^{\frac{1}{2}-r_1+r_2-\delta}
    \big(1+\abs{\xi_1+\tau_2}\big)^{-r_1} \; d\xi_1
 \\ \nonumber & 
 \hskip 60pt
  \cdot
    \int\limits_{\xi_3=\tau_4}^{\tau_3}
    \big(1+\abs{\xi_3-\tau_3}\big)^{-r_3}
    \big(1+\abs{\xi_3-\tau_4}\big)^{\frac{1}{2}-r_3+r_2-\delta}
 \; d\xi_3
\; d\tau_2
\; d\tau_3.
\end{align}

Note that the integral in $\xi_3$ is bounded by
 $T^{\varepsilon+\frac{1}{2}+r_2-\delta}\big(1+\tau_3-\tau_4\big)^{\delta_{0,r_3}-r_3},$
and the integral in $\xi_1$ is bounded by
$T^{\varepsilon+\frac12+r_2-\delta}\big(1+\tau_1-\tau_2\big)^{\delta_{0,r_1}-r_1}.$ Plugging these bounds in and simplifying, we find that this shifted term is bounded as follows: 
\begin{align*}&
\left| p_{T,R}^{2,\delta}(y;(-a_1,-a_3),\kappa) \right| \ll  \
 y_1^{\frac32+a_1}y_2^{2+a_2}y_3^{\frac32+a_3}\,
 T^{\varepsilon+R+1+2r_2+\delta}\\&\cdot
 \iint\limits_{\substack{ -T^{1+\varepsilon}\leq \tau_4\leq \tau_3\leq T^{1+\varepsilon} \\ -T^{1+\varepsilon} \leq \tau_2\leq \tau_1\leq T^{1+\varepsilon}}} \;
    \big(1+\abs{\tau_2-\tau_3}\big)^{\frac{1+R}{2}-\delta} \big(1+\tau_1-\tau_2)^{1+\frac{R}{2}+\delta_{0,r_1}-r_1}
    \big(1+\tau_3-\tau_4)^{1+\frac{R}{2}+\delta_{0,r_3}-r_3}
 \\ \nonumber & 
  \cdot
    \big(1+\abs{\tau_1-\tau_4}\big)^{\frac{3+R}{2}-2r_2+\delta}
    \big(1+\abs{\tau_1-\tau_3}\big)^{\frac{2+R}{2}-r_2}
    \big(1+\abs{\tau_2-\tau_4}\big)^{\frac{2+R}{2}-r_2}
 \; d\tau_2
 \; d\tau_3
 \\ & \ll
 y_1^{\frac32+a_1}y_2^{2+a_2}y_3^{\frac32+a_3}\,
T^{\varepsilon+4R+10+\delta_{0,r_1}+\delta_{0,r_3}-r_1-r_2-r_3-(r_2-\delta)} \\ & \ll
 y_1^{\frac32+a_1}y_2^{2+a_2}y_3^{\frac32+a_3}
 \,T^{\varepsilon+4R+9+\delta_{0,r_1}+\delta_{0,r_3}-r_1-r_2-r_3},
\end{align*}
since $r_2-\delta\geq 1$.

In analogy to Lemma~\ref{lem:shiftalpha1poles}, it is easy to show that in shifting $\re(\alpha_1)$ from zero to $\kappa_1$, poles are crossed at $\alpha_1=q$ for
 \[ q\in \left\{\left.\begin{array}{c} -s_1-2(r_1-\delta_1) \\ -s_3-\alpha_2-\alpha_3 - 2(r_3-\delta_3) \end{array}\right| \begin{array}{c} 0\leq \delta_j \leq r_j  \\ \delta_j \leq r_2-\delta \end{array}\right\}. \]

In what follows, the method holds equally well for either of these two types of residues.  For concreteness, we consider the residue at $\alpha_1=-s_1-2(r_1-\delta_1)$.  This gives
\begin{align*} 
 & \iint\limits_{\substack{\re(\alpha_j)=\kappa_j \\ j=2,3}}  e^{\frac{(s_1+2(r_1-\delta_1))^2+\alpha_2^2+\alpha_3^2+(s_1+\alpha_2+\alpha_3+2(r_1-\delta_1))^2}{2T^2}}\;
 \iint\limits_{\substack{\re(s_j)=-a_j\\j=1,3}} 
 y_1^{\frac32-s_1}y_2^{2-s_1+\alpha_2-2(r_1-\delta_1-\delta)}y_3^{\frac32-s_3}\; 
 \\ & \qquad
 \cdot
    \left.\Big(  \mathcal F_R(\alpha)\,Q_\delta(s,\alpha) \Big)\right|_{\alpha_1=-s_1-2(r_1-\delta_1)} 
 \frac{ 
       \Gamma\left(\frac{\alpha_3-\alpha_2}{2}-\delta\right)
     }{
       \Gamma\left(\frac{\alpha_3-\alpha_2}{2}\right)
     } 
 \frac{
       \Gamma\(\frac{2+R+\alpha_2-\alpha_3}{4}\)
       \Gamma\(\frac{2+R+\alpha_3-\alpha_2}{4}\)
     }{
       \Gamma\(\frac{\alpha_2-\alpha_3}{2}\) 
     }
 \\ & \quad
 \cdot
 \frac{ 
       \Gamma\left(\frac{\alpha_3+s_1}{2}+r_1-\delta_1-\delta\right)
     }{
       \Gamma\left(\frac{\alpha_3+s_1}{2}+r_1-\delta_1\right)
     }
 \frac{
       \Gamma\(\frac{2+R-s_1-\alpha_3-2(r_1-\delta_1)}{4}\)
       \Gamma\(\frac{2+R+\alpha_3+s_1+2(r_1-\delta_1)}{4}\)
     }{
       \Gamma\(\frac{-s_1-\alpha_3}{2}-r_1+\delta_1\) 
     } 
  \\ & \quad 
  \cdot
 \frac{ 
       \Gamma\left(\frac{2s_1-\alpha_2-\alpha_3}{2}+2(r_1-\delta_1)-\delta\right)
     }{
       \Gamma\left(\frac{2s_1-\alpha_2-\alpha_3}{2}+2(r_1-\delta_1)\right)
     } 
 \frac{
       \Gamma\(\frac{2+R-2s_1+\alpha_2+\alpha_3}{4}-r_1+\delta_1\)
       \Gamma\(\frac{2+R+2s_1-\alpha_2-\alpha_3}{4}+r_1-\delta_1\)
     }{
       \Gamma\(\frac{-2s_1+\alpha_2+\alpha_3}{2}-2(r_1-\delta_1)\) 
     }
  \\ & \quad 
  \cdot
 \frac{ 
       \Gamma\left(\frac{s_1-2\alpha_2-\alpha_3}{2}+r_1-\delta_1-\delta\right)
     }{
       \Gamma\left(\frac{s_1-2\alpha_2-\alpha_3}{2}+r_1-\delta_1\right)
     } 
 \frac{
       \Gamma\(\frac{2+R-s_1-2(r_1-\delta_1)+2\alpha_2+\alpha_3}{4}\)
       \Gamma\(\frac{2+R+s_1+2(r_1-\delta_1)-2\alpha_2-\alpha_3}{4}\)
     }{
       \Gamma\(\frac{-s_1+2\alpha_2+\alpha_3}{2}-r_1+\delta_1\) 
     }
  \\ & \quad 
  \cdot
  \frac{
        \Gamma(\frac{2+R-s_1-2(r_1-\delta_1)-\alpha_2}{4})
        \Gamma(\frac{2+R+s_1+2(r_1-\delta_1)+\alpha_2}{4})
        \Gamma(\frac{2+R-s_1-2(r_1-\delta_1)+\alpha_2+2\alpha_3}{4})
      }{
        \Gamma(\frac{-s_1-2(r_1-\delta_1)-\alpha_2}{2})
        \Gamma(\frac{ s_1+2(r_1-\delta_1)+\alpha_2}{2})
        \Gamma(\frac{-s_1+\alpha_2+2\alpha_3}{2}-r_1+\delta_1)
      }
  \\ & \quad
  \cdot
  \frac{
        \Gamma(\frac{2+R+s_1+2(r_1-\delta_1)-\alpha_2-2\alpha_3}{4})
      }{
        \Gamma(\frac{ s_1-\alpha_2-2\alpha_3}{2}+r_1-\delta_1)
      }%\\&\quad
  \cdot
 \textstyle{ \Gamma\(\frac{s_1+\alpha_2}{2}\) \Gamma\(\frac{s_3-\alpha_3}{2}\) \Gamma\(\frac{s_3-s_1+\alpha_2+\alpha_3}{2}-r_1+\delta_1\)}\; ds_1ds_3
\; d\alpha.
\end{align*}

Strictly speaking this is what we want to bound in the case that $\kappa_2=\kappa_3=0$, but as before we need to shift the line of integration 
in the $\alpha_2$ variable to $\re(\alpha_2)=\kappa_2$ such that the real part of the exponent of $y_2$ is $2+a_2$.  This means that
 \[ a_2 = \re(-s_1+\alpha_2-2(r_1-\delta_1-\delta)) = a_1 + \kappa_2 - 2(r_1-\delta_1-\delta), \]
or in other words,
\begin{align}\label{eq:a2secondrequirement}
 \kappa_2 =  2(r_2-\delta-\delta_1)+a_2'-a_1'\quad\mbox{where}\quad -1<a_j'=a_j-2r_j<1 (j=1,2).
\end{align}
This implies that $\re(\alpha_2)$ gets shifted to the right.  Just as in the case of Lemma~\ref{lem:pTr_1alpha2shift}, one can show in this case that no poles are crossed in moving $\re(\alpha_2)$. So it suffices to bound the above for these values of $\kappa_2$ and $\kappa_3$.

The exponential factor is $e^{-\frac{\pi}{4}\mathcal{E}}$ where
\begin{align*}
 \mathcal{E} & = \lvert \xi_1-\tau_2-2\tau_3\rvert - \lvert \xi_3-\tau_3\rvert - \lvert \xi_3-\xi_1+\tau_2+\tau_3 \rvert,
\end{align*}
which leads to two exponential zero sets, the first of which is
 \[\mathcal{R}:\qquad \tau_3 \leq \xi_3 \leq \xi_1-\tau_2-\tau_3, \]
and the second is similar but the inequalities are reversed.

The polynomial factor (coming from the Gamma factors specifically) is
\begin{align*}
 \mathcal{P}  = &
    \big(1+\abs{\xi_1+\tau_2}\big)^{\frac{1+R+2 a_1' - a_2'}{2} -r_1+r_2- 
\delta-\delta_1}
    \big(1+\abs{\tau_2-\tau_3}\big)^{\frac{1+R- a_1' + a_2'}{2} -r_2+\delta_1}
\\\cdot & 
    \big(1+\abs{\xi_1+\tau_3}\big)^{\frac{1+R+a_1'}{2}-\delta-\delta_1}
    \big(1+\abs{-2 \xi_1+\tau_2+\tau_3}\big)^{\frac{1+R+a_1'+a_2'}{2}-r_2- \delta_1}
\\\cdot & 
    \big(1+\abs{-\xi_1+2\tau_2+\tau_3}\big)^{\frac{1+R-a_1'+2a_2'}{2} - 2 
r_2 + \delta + \delta_1 }
    \big(1+\abs{-\xi_1+\tau_2+2\tau_3}\big)^{1+\frac{R}{2}}
\\\cdot &
    \big(1+\abs{\xi_3-\tau_3}\big)^{\frac{-1+a_3'}{2}-r_3}
    \big(1+\abs{-\xi_1+\xi_3+\tau_2+\tau_3}\big)^{\frac{-1- a_2' + a_3'}{2} + r_2 - r_3 - \delta}.
\end{align*}

Note that the the presence of the exponential term means that we can restrict the integral to the set of $\lvert \tau_2\rvert,\lvert \tau_3\rvert,\lvert \xi_1\rvert \leq T^{1+\varepsilon}$.  Then  the integral that we seek to bound is  
\begin{align*}\ll
 &y_1^{\frac32+a_1}y_2^{2+a_2}y_3^{\frac32+a_3} \,T^{\varepsilon+R+3\delta+\frac{1}{2}} \hskip-10pt\iiint\limits_{\substack{\lvert \tau_2\rvert,\lvert \tau_3\rvert,\lvert \xi_1\rvert \leq T^{1+\varepsilon}\\ \tau_2+2\tau_3\leq \xi_1}} 
  \big(1+\abs{\xi_1+\tau_2}\big)^{\frac{1+R+2 a_1' - a_2'}{2} -r_1+r_2- \delta-\delta_1}
\\\cdot &     \big(1+\abs{\tau_2-\tau_3}\big)^{\frac{1+R- a_1' + a_2'}{2} 
-r_2+\delta_1}
    \big(1+\abs{\xi_1+\tau_3}\big)^{\frac{1+R+a_1'}{2}-\delta-\delta_1}
    \big(1+\abs{-2 \xi_1+\tau_2+\tau_3}\big)^{\frac{1+R+a_1'+a_2'}{2}-r_2- \delta_1}
\\\cdot & 
    \big(1+\abs{-\xi_1+2\tau_2+\tau_3}\big)^{\frac{1+R-a_1'+2a_2'}{2} - 2 
r_2 + \delta + \delta_1 }
    \big(1+\abs{-\xi_1+\tau_2+2\tau_3}\big)^{1+\frac{R}{2}}
     \\ & \quad \cdot
 \int\limits_{\xi_3=\tau_3}^{\xi_1-\tau_2-\tau_3} 
    \big(1+\abs{\xi_3-\tau_3}\big)^{-r_3}
    \big(1+\abs{-\xi_1+\xi_3+\tau_2+\tau_3}\big)^{r_2 - r_3 - \delta}\, d\xi_3\; d\xi_1\, d\alpha_2\, d\alpha_3
 \\ &  \qquad \qquad \ll \
 y_1^{\frac32+a_1}y_2^{2+a_2}y_3^{\frac32+a_3} \, T^{\varepsilon+4R+3+\frac{8}{2}+a_1'+\frac{3}{2}a_2'+\delta_{0,r_3}-r_1-r_2-r_3-(r_2-\delta)-\delta_1}
 \\ &  \qquad \qquad \ll \
 y_1^{\frac32+a_1}y_2^{2+a_2}y_3^{\frac32+a_3} \,T^{\varepsilon+4R+ \frac{17}2+\delta_{0,r_3}-r_1-r_2-r_3},
\end{align*}
since $r_2-\delta\ge1$, $\delta_1\geq 0$.
\end{proof}

\subsection{\bf Bounds for the double residue terms}
\label{sec:doubleresidues}
There are two types of double residue terms that we need to consider, namely $  p_{T,R}^{12,(\delta_1,\delta_2)}(y;-a_3) $ and $  p_{T,R}^{13,(\delta_1,\delta_3)}(y;-a_2) $.    They are obtained by taking residues, at $s_2=-\alpha_1-\alpha_4-2\delta_2$ and $s_3=\alpha_2-2\delta_3$ respectively, of the single residue term $  p_{T,R}^{1,(\delta_1)}(y) $ defined 
by (\ref{eq:pTr_1delta}).  

Specifically,  write
\begin{equation}p_1=-\alpha_1-2\delta_1,\quad  p_2=-\alpha_1-\alpha_4-2\delta_2,\quad p_3=\alpha_2-2\delta_3,\label{eq:pj}\end{equation}where  $0\le \delta_j\le r_j-1$ for $1\le j\le 3$.  Then we find from Proposition \ref{resprop}  that

\begin{align} \label{eq:pTr_12delta}
  p_{T,R}^{12,(\delta_1,\delta_2)}(y;-a_3) = & \iiint\limits_{\re(\alpha)=0} e^{\frac{\alpha_1^2+\cdots+\alpha_4^2}{2T^2}}\; \int\limits_{\re(s_3)=-a_3}y_1^{\frac32-p_1}y_2^{2-p_2}y_3^{\frac32-s_3} \; 
  \mathcal F_R(\alpha)\, \Gamma_R(\alpha) 
  \\ \nonumber & \hskip 48pt
  \cdot f_{\delta_1,\delta_2}(s_3,\alpha)\textstyle{
   \Gamma(\frac{\alpha_2-\alpha_1}{2}-\delta_1)  \Gamma(\frac{\alpha_3-\alpha_1}{2}-\delta_1) \Gamma(\frac{\alpha_4-\alpha_1}{2}-\delta_1) }
  \\ \nonumber & \hskip 48pt
  \cdot 
  \textstyle{  \Gamma(\frac{\alpha_2-\alpha_4}{2}-\delta_2) \Gamma(\frac{\alpha_3-\alpha_4}{2}-\delta_2) } \Gamma(\frac{s_3-\alpha_2}{2})\Gamma(\frac{s_3-\alpha_3}{2})  \; ds_3\; d\alpha,
 \end{align} and
\begin{align} \label{eq:pTr_13delta}
  p_{T,R}^{13,(\delta_1,\delta_3)}(y;-a_2) = & \iiint\limits_{\re(\alpha)=0} e^{\frac{\alpha_1^2+\cdots+\alpha_4^2}{2T^2}}\; \int\limits_{\re(s_2)=-a_2}y_1^{\frac32-p_1}y_2^{2-s_2}y_3^{\frac32-p_3} \,
  \mathcal F_R(\alpha)\,  \Gamma_R(\alpha) 
  \\ \nonumber & \hskip 48pt
  \cdot g_{\delta_1,\delta_3}( s_2,\alpha)\textstyle{
   \Gamma(\frac{\alpha_2-\alpha_1}{2}-\delta_1)  \Gamma(\frac{\alpha_3-\alpha_1}{2}-\delta_1) \Gamma(\frac{\alpha_4-\alpha_1}{2}-\delta_1) }
  \\ \nonumber & \hskip 48pt
  \cdot
  \textstyle{  \Gamma(\frac{\alpha_2-\alpha_3}{2}-\delta_3) \Gamma(\frac{\alpha_2-\alpha_4}{2}-\delta_3) }  \Gamma(\frac{s_2+a_1+\alpha_3}{2})\Gamma(\frac{s_2+a_1+\alpha_4}{2})   \; ds_2\; d\alpha,
 \end{align} where $\deg{f_{\delta_1,\delta_2}}\le 2\delta_1+\delta_2$  and $\deg{g_{\delta_1,\delta_3}}\le2\delta_1+\delta_3.$   (See Section \ref{poles-res}.)

In what follows, we show that the bounds on \eqref{eq:pTr_12delta} and \eqref{eq:pTr_13delta}
are ``small.''

We begin with \eqref{eq:pTr_12delta}.  We have:

\begin{prop}\label{prop:doubleresidue12} Let $r_1,r_2\geq 1$, $r_3\geq0$ be integers, and $0<\varepsilon<1$. Suppose $a_1,a_2,a_3$ satisfy the hypotheses of Theorem \ref{th:pTRbound}.  If $0\le \delta_j\le r_j-1$ for $1\le j\le2$, then
\begin{equation}\label{eq:Prop-s1s2-Residue}
 \left|   p_{T,R}^{12,(\delta_1,\delta_2)}(y;-a_3) \right| \ll y_1^{\frac32+a_1}y_2^{2+a_2}y_3^{\frac32+a_3} \,T^{\varepsilon+4R+\frac{13}{2}+\delta_{r_3,0}-r_1-r_2-r_3}.
\end{equation}\end{prop}

\begin{proof} The proof is similar, in spirit and in many of the details, 
to that of Proposition \ref{prop:singleresidue1}.

More specifically: to obtain the desired bound on  $p_{T,R}^{12,(\delta_1,\delta_2)}(y;-a_3)$, we will need to shift the lines of integration in both the $\alpha_1$ and $\alpha_2$ variables, so that the resulting exponents of $y_1$ and $y_2$ have real parts as stated in the  proposition.  In 
doing so, we will pick up residues.  That is, we will have
\begin{equation}\label{eq:KappaPlusResDouble}p_{T,R}^{12,(\delta_1,\delta_2)}(y;-a_3)  =  p_{T,R}^{12,(\delta_1,\delta_2)}(y;-a_3,\kappa) + \sum\mbox{Residues}, \end{equation}
where 
\begin{align}\label{eq:pTr_12deltashift}
 & p_{T,R}^{12,(\delta_1,\delta_2)}(y;-a_3,\kappa) := 
 \iiint\limits_{\re(\alpha)=\kappa} e^{\frac{\alpha_1^2+\cdots+\alpha_4^2}{2T^2}} \int\limits_{\re(s_3)=-a_3}y_1^{\frac32+\alpha_1+2\delta_1}y_2^{2+\alpha_1+\alpha_4+2\delta_2}y_3^{\frac32-s_3} 
  \mathcal{F}_R(\alpha) 
    f_{\delta_1,\delta_2}(s_3,\alpha)
  \\ \nonumber & 
  \cdot
  \frac{\Gamma(\frac{2+R+\alpha_1-\alpha_2}{4})
        \Gamma(\frac{2+R+\alpha_2-\alpha_1}{4})\Gamma(\frac{\alpha_2-\alpha_1}{2}-\delta_1)
      }{\Gamma(\frac{\alpha_1-\alpha_2}{2})\Gamma(\frac{\alpha_2-\alpha_1}{2})}
  \cdot
  \frac{\Gamma(\frac{2+R+\alpha_2-\alpha_3}{4})
        \Gamma(\frac{2+R+\alpha_3-\alpha_2}{4}) \Gamma(\frac{s_3-\alpha_2}{2})\Gamma(\frac{s_3-\alpha_3}{2}) 
      }{\Gamma(\frac{\alpha_2-\alpha_3}{2})
        \Gamma(\frac{\alpha_3-\alpha_2}{2})}
  \\ \nonumber & 
  \cdot
  \frac{
      }{}
  \frac{\Gamma(\frac{2+R+\alpha_1-\alpha_3}{4})
        \Gamma(\frac{2+R+\alpha_3-\alpha_1}{4})\Gamma(\frac{\alpha_3-\alpha_1}{2}-\delta_1)
      }{\Gamma(\frac{\alpha_1-\alpha_3}{2})\Gamma(\frac{\alpha_3-\alpha_1}{2})}
  \cdot
  \frac{\Gamma(\frac{2+R+\alpha_2-\alpha_4}{4})
        \Gamma(\frac{2+R+\alpha_4-\alpha_2}{4})  \Gamma(\frac{\alpha_2-\alpha_4}{2}-\delta_2) 
      }{\Gamma(\frac{\alpha_2-\alpha_4}{2})
        \Gamma(\frac{\alpha_4-\alpha_2}{2})}
  \\ \nonumber & 
  \cdot
  \frac{\Gamma(\frac{2+R+\alpha_1-\alpha_4}{4})
        \Gamma(\frac{2+R+\alpha_4-\alpha_1}{4})\Gamma(\frac{\alpha_4-\alpha_1}{2}-\delta_1)
      }{\Gamma(\frac{\alpha_1-\alpha_4}{2})\Gamma(\frac{\alpha_4-\alpha_1}{2})}
  \cdot 
  \frac{
        \Gamma(\frac{2+R+\alpha_3-\alpha_4}{4})\Gamma(\frac{2+R+\alpha_4-\alpha_3}{4})\Gamma(\frac{\alpha_3-\alpha_4}{2}-\delta_2)
      }{
        \Gamma(\frac{\alpha_3-\alpha_4}{2})\Gamma(\frac{\alpha_4-\alpha_3}{2})} \, ds_3\, d\alpha,
\end{align}
and the residues that appear depend on the particular choice of $\kappa=(\kappa_1,\kappa_2,\kappa_3)$.   We've grouped the Gamma factors, above, in a manner that will be convenient for what follows.

Because we want the exponents of $y_1$ and $y_2$, in \eqref{eq:pTr_12deltashift}, to have real parts $3/2+a_1$ and $2+a_2$ respectively, we will choose $\kappa=(\kappa_1,\kappa_2,\kappa_3)\in\R^3$  such that  
 \begin{equation*} \re(\alpha_1+2\delta_1)=\kappa_1+2\delta_1=a_1, \qquad \re(\alpha_1+\alpha_4+2\delta_2)=-\kappa_2-\kappa_3+2\delta_2=a_2 .  \end{equation*}Specifically, we will define
 \begin{equation}\kappa=(\kappa_1,\kappa_2,\kappa_3)=(a_1-2\delta_1,-a_2 +2\delta_2 ,0) \quad  (\hbox{and }\kappa_4=-\kappa_1-\kappa_2-\kappa_3). \label{newkappadouble}\end{equation}
For this value of $\kappa$, we will obtain an estimate of the desired magnitude for $p_{T,R}^{12,(\delta_1,\delta_2)}(y;-a_3,\kappa) .$  
 
 It will remain to estimate the residues that appear in \eqref{eq:KappaPlusResDouble}.  To do so we first identify the poles, cf. Lemma \ref{lem:shiftalpha1alpha2poles} below.  We then show that, for  these residue terms,  the desired exponents on the $y_j$'s can be obtained by shifting lines of integration, without passing additional poles.  Finally, using this information, we show that the residue terms are small.

%  \pagebreak
  
\vskip 12pt

\noindent
\underline{\bf{Step 1}: Bounding the shifted integral $p_{T,R}^{12,(\delta_1,\delta_2)}$}
\vskip 10pt

Before estimating the  residue terms in \eqref{eq:KappaPlusResDouble}, we 
obtain a bound of the desired magnitude on the shifted integral $ p_{T,R}^{12,(\delta_1,\delta_2)}(y;-a_3,\kappa)$, with $\kappa=(\kappa_1,\kappa_2,\kappa_3)$ as in \eqref{newkappadouble}.

Note that, from any of the grouped combinations of Gamma functions  in \eqref{eq:pTr_12deltashift} {\it except for the second one}, the contribution to the exponential factor in Stirling's formula is zero.  This is because absolute values of imaginary parts from the numerator of any such combination cancel those from the denominator. So, again,  the only one of these terms that contributes to the exponential factor is the second one, which contributes a factor of 
 \begin{equation*}  e^{-\frac{\pi}{2}(|\tau_2-\tau_3|/4+|\tau_3-\tau_2|/4+|\xi_3-\tau_2|/2+|\xi_3-\tau_3|/2-|\tau_2-\tau_3|/2-|\tau_3-\tau_2|/2)}=e^{-\frac{\pi}{4}( |\xi_3-\tau_2|+|\xi_3-\tau_3|-|\tau_2-\tau_3|)}.\end{equation*}
But the integrand in \eqref{eq:pTr_12deltashift}  is invariant under $\alpha_2\leftrightarrow \alpha_3$, so we may assume that
 $ \tau_2 \ge \tau_3$, whence the exponential factor in question is simply
 $$e^{-\frac{\pi}{4}( |\xi_3-\tau_2|+|\xi_3-\tau_3|-\tau_2+\tau_3)}.$$
 It is then easily seen that there is just one exponential zero set, namely
\begin{align*}
 \mathcal{R} & :=  \left\{ (-a_3+i\xi_3)\in \C |   \tau_3\leq \xi_3\leq 
 \tau_2 \right\} .\end{align*}

Replacing the Gamma factors with their corresponding polynomial terms, in 
\eqref{eq:pTr_12deltashift}, then gives 
\begin{align} \label{eq:pTr_12Shifted}&
  \left| p_{T,R}^{12,(\delta_1,\delta_2)}(y;-a_3,\kappa)   \right|  \ll  
  y_1^{\frac32+a_1}y_2^{2+a_2}y_3^{\frac32+a_3}
 \, T^{\varepsilon+R+2\delta_1+\delta_2}  \\&
 \cdot \iiint\limits_{\substack{ \abs{\tau_1},\abs{\tau_2},\abs{\tau_3} \leq T^{1+\varepsilon} \\ 0 \leq  \tau_2-\tau_3  }} 
  \hskip -8pt
  \big(1+\lvert\tau_1-\tau_2\rvert\big)^{\frac{R+1+\kappa_2-\kappa_1}{2}-\delta_1}\nn
  \big(1+\lvert\tau_1-\tau_3\rvert\big)^{\frac{R+1+\kappa_3-\kappa_1}{2} -\delta_1}
  \big(1+\lvert\tau_1-\tau_4\rvert\big)^{\frac{R+1+\kappa_4-\kappa_1}{2} -\delta_1}
  \\ \nonumber & 
  \cdot
  \big( 1+ \tau_2-\tau_3  \big)^{ \frac{2+R}{2} }
  \big( 1+\lvert\tau_2-\tau_4 \rvert\big)^{ \frac{R+1+\kappa_2-\kappa_4}{2}-\delta_2}
  \big( 1+|\tau_3-\tau_4 |\big)^{ \frac{R+1+\kappa_3-\kappa_4}{2}-\delta_2} 
  \\ \nonumber & 
  \cdot
 \int _{\xi_3=\tau_3}^{\tau_2}
   \big(1+\tau_2 - \xi_3\big)^{\frac{-1-a_3-\kappa_2}{2}}
        \big(1+\xi_3-\tau_3\big)^{\frac{-1-a_3-\kappa_3}{2}}
        \;  d\xi_3\; d\tau.
\end{align}The change of variables
 \[  \quad \xi_3\mapsto \xi_3+\tau_3, \quad \quad T_j=\tau_j-\tau_{j+1} 
\quad (1\le j\le3) \]
applied to \eqref{eq:pTr_12Shifted} then gives\begin{align} \label{pTr_12NewShifted}
  \left| p_{T,R}^{12,(\delta_1,\delta_2)}(y;-a_3,\kappa)   \right| & \ll  
  y_1^{\frac32+a_1}y_2^{2+a_2}y_3^{\frac32+a_3}
  \,T^{\varepsilon+R+2\delta_1+\delta_2}  \\&
 \cdot \iiint\limits_{ \substack{0\le \abs{T_1},T_2, \abs{T_3} \leq T^{1+\varepsilon}   }} 
  \hskip -8pt
  \big(1+\lvert T_1\rvert\big)^{\frac{R+1+\kappa_2-\kappa_1}{2}-\delta_1}\nn
  \big(1+\lvert T_1+T_2\rvert\big)^{\frac{R+1+\kappa_3-\kappa_1}{2} -\delta_1}
  \\ \nonumber & 
  \cdot   \big(1+\lvert T_1+T_2+T_3\rvert\big)^{\frac{R+1+\kappa_4-\kappa_1}{2} -\delta_1}
  \big( 1+ T_2 \big)^{ \frac{2+R}{2} }\\\nn&\cdot
  \big( 1+\lvert T_2+T_3 \rvert\big)^{ \frac{R+1+\kappa_2-\kappa_4}{2}-\delta_2}
  \big( 1+|T_3 |\big)^{ \frac{R+1+\kappa_3-\kappa_4}{2}-\delta_2} 
  \\ \nonumber & 
  \cdot
 \int _{\xi_3=0}^{T_2}
     \big(1+T_2-\xi_3\big)^{\frac{-1-a_3-\kappa_2}{2}}
        \big(1+\xi_3 \big)^{\frac{-1-a_3-\kappa_3}{2}}
        \;  d\xi_3\; dT_1\;dT_2\;dT_3. 
\end{align}Now Lemma \ref{lemmaIntegralBound} and the fact that $\kappa_3=0$ (cf. (\ref{newkappadouble})) tell us that
\begin{align}&\int _{\xi_3=0}^{T_2}
         \big(1+T_2-\xi_3\big)^{\frac{-1-a_3-\kappa_2}{2}}
        \big(1+\xi_3 \big)^{\frac{-1-a_3-\kappa_3}{2}}
        \;  d\xi_3 \nn\\& \ll 
 (1+T_2)^{-\min\big\{\frac{1+a_3+\kappa_2}{2},\,\frac{1+a_3+\kappa_3}{2},\,a_3+\frac{ \kappa_2+\kappa_3}{2}\big\}+\varepsilon}\nn
 \\&= \;
 (1+T_2)^{-\frac{1+a_3+\kappa_2}{2}+\frac{1}{2}\max\{0,\,\kappa_2,\,1-a_3 
\}+\varepsilon} .\label{T2int}\end{align}But$$\kappa_2=-a_2+2\delta_2\le 1-2r_2-\varepsilon+2(r_2-1)=-1-\varepsilon<0,$$and$$-2r_3+\varepsilon\le 1-a_3\le2-2r_3-\varepsilon;$$from this information, it follows that$$\max\{0,\,\kappa_2,\,1-a_3 \} \le 2\delta_{r_3,0}. $$So by \eqref{T2int},
 \begin{align}&\int _{\xi_3=0}^{T_2}
         \big(1+T_2-\xi_3\big)^{\frac{-1-a_3-\kappa_2}{2}}
        \big(1+\xi_3 \big)^{\frac{-1-a_3-\kappa_3}{2}}
        \;  d\xi_3  \ll   (1+T_2)^{-\frac{1+a_3+\kappa_2}{2}+\delta_{r_3,0}+\varepsilon} . \label{T2int-a}\end{align}Then \eqref{pTr_12NewShifted} 
gives
\begin{align}  &
  \left| p_{T,R}^{12,(\delta_1,\delta_2)}(y;-a_3,\kappa)   \right|  \ll  
  y_1^{\frac32+a_1}y_2^{2+a_2}y_3^{\frac32+a_3}
\, T^{\varepsilon+R+2\delta_1+\delta_2} \hskip-10pt\iiint\limits_{ \substack{0\le \abs{T_1},T_2, \abs{T_3} \leq T^{1+\varepsilon}   }} 
  \hskip -8pt
  \big(1+\lvert T_1\rvert\big)^{\frac{R+1+\kappa_2-\kappa_1}{2}-\delta_1} 
\\&
 \cdot \nn
  \big(1+\lvert T_1+T_2\rvert\big)^{\frac{R+1+\kappa_3-\kappa_1}{2} -\delta_1}
  \big(1+\lvert T_1+T_2+T_3\rvert\big)^{\frac{R+1+\kappa_4-\kappa_1}{2} -\delta_1}
  \\ \nonumber & 
  \cdot
  \big( 1+ T_2 \big)^{ \frac{ R+1-a_3-\kappa_2}{2} +\delta_{r_3,0}+\varepsilon}
  \big( 1+\lvert T_2+T_3 \rvert\big)^{ \frac{R+1+\kappa_2-\kappa_4}{2}-\delta_2}
  \big( 1+|T_3 |\big)^{ \frac{R+1+\kappa_3-\kappa_4}{2}-\delta_2} 
\; dT_1\;dT_2\;dT_3. \nn
\end{align}

It's straightforward to estimate the above integral, using the facts that, on the indicated domains of integration,$$ T_1, T_2, T_3  \ll T^{1+\varepsilon},$$and that the length of each domain of integration is also $\ll 
T^{1+\varepsilon}$.  We thereby find that
\begin{align*}&
  \left| p_{T,R}^{12,(\delta_1,\delta_2)}(y;-a_3,\kappa)   \right|  \ll  
  y_1^{\frac32+a_1}y_2^{2+a_2}y_3^{\frac32+a_3}
  \, T^{\varepsilon+4R+6+\delta_1+\delta_2+\delta_{r_3,0}-a_1-a_2-\frac{a_3}{2}}.\end{align*}
But we're assuming that $-a_3\le 1-2r_3$,  and  $$\delta_j-a_j\le \delta_j+1-2r_j=(\delta_j+1-r_j)-r_j\le -r_j \quad\mbox{for}\quad j=1,2.$$
So our above estimate reads
  \begin{align} \label{eq:pTr_12NewerShifted}   \left| p_{T,R}^{12,(\delta_1,\delta_2)}(y;-a_3,\kappa)   \right|  \ll  &
  y_1^{\frac32+a_1}y_2^{2+a_2}y_3^{\frac32+a_3}
\, T^{\varepsilon+4R+\frac{13}{2}+\delta_{r_3,0}-r_1-r_2-r_3}, \end{align}which gives us a bound of the desired magnitude on the  shifted integral 
 $p_{T,R}^{12,(\delta_1,\delta_2)}(y;-a_3,\kappa) $ in \eqref{eq:KappaPlusResDouble}.
\vskip 12pt
%\newpage
\noindent
\underline{\bf{Step 2}: Bounding the residue terms}
\vskip 10pt

 Our next step is to estimate the residues in \eqref{eq:KappaPlusResDouble}.  Recall:  these are the residues at the poles that one crosses in moving the lines of integration in \eqref{eq:pTr_12delta}, to transform it into  \eqref{eq:pTr_12deltashift}.
 
 We first locate these poles.

\begin{lem}\label{lem:shiftalpha1alpha2poles}Suppose the lines of integration, 
 in \eqref{eq:pTr_12delta}, are shifted  from $\re(\alpha_1,\alpha_2,\alpha_3)=(0,0,0)$ to $\re(\alpha_1,\alpha_2,\alpha_3)=( a_1-2\delta_1,-a_2+2\delta_2,0)$.  Then:
 
 \begin{enumerate}
 \item[\rm (a)] No poles are crossed in the $\alpha_1$ variable.
 
 \item[\rm (b)] For a fixed $s_3$, any pole crossed in the $\alpha_2$ variable belongs to the set \begin{equation}\label{eq:alpha12poles}
  \left\{ s_3+2\delta_3\,\vert \, \delta_3\in\Z_{\ge0},\,\max\{0,r_3-(r_2-\delta_2)\}\le \delta_3\le r_3  \right\}.
\end{equation}\end{enumerate}
\end{lem}

\begin{proof}  First we consider the factors
 \[\textstyle \Gamma\big(\frac{2+R+\alpha_j-\alpha_k}{4}\big)  \qquad \mbox{($1\leq j\neq k\leq 4$),} \]in the integrand of  \eqref{eq:pTr_12deltashift}.
If $R$ is sufficiently large, then no poles of these factors will be crossed in moving our lines of integration in $\alpha_1$ and $\alpha_2$.

Nor do any of the terms of the form
$$  \frac{ \Gamma(\frac{\alpha_k-\alpha_j}{2}-\delta_n)
      }{\Gamma(\frac{\alpha_j-\alpha_k}{2})\Gamma(\frac{\alpha_k-\alpha_j}{2})} $$give rise to any poles.  Indeed,   if the numerator of this  term  has a pole, then $\frac{\alpha_k-\alpha_j}{2}-\delta_n$ is a nonpositive integer, whence  $\frac{\alpha_k-\alpha_j}{2}\in\Z$, so this numerator 
pole will be cancelled by a pole from either $\Gamma(\frac{\alpha_j-\alpha_k}{2})$  or $\Gamma(\frac{\alpha_k-\alpha_j}{2})$ in the denominator.
      
The only factors remaining to consider are the factors $$\textstyle \Gamma(\frac{s_3-\alpha_2}{2})\quad\mbox{and}\quad \Gamma(\frac{s_3-\alpha_3}{2}),$$and since we are not shifting the line of integration in $\alpha_3$, we need only examine the former of these factors.  

In particular, part (a) of our lemma is proved.

Regarding $\Gamma(\frac{s_3-\alpha_2}{2})$:  for fixed $s_3$, this factor 
has poles, as a function of $\alpha_2$, whenever\begin{align}\alpha_2=s_3+2\delta_3 \quad(\delta_3\in\Z_{\ge0}).\label{eq:alpha2poles}\end{align}But  for such an $\alpha_2$ to lie between the initial line of integration $\re{\alpha_2}=0$ and the terminal line  $\re{\alpha_2}=-a_2+2\delta_2$, we must have
\begin{align}-a_2+2\delta_2\le\re{\alpha_2}=\re{(s_3+2\delta_3)}=-a_3+2\delta_3 \le0 .\label{eq:alpha2poles2}\end{align}

But$$-a_3+2\delta_3\ge -a_2+2\delta_2\implies \delta_3\ge \frac{a_3-a_2}{2}+\delta_2\ge r_3-r_2+\delta_2-1+\varepsilon,$$the last inequality because $2r_j-1+\varepsilon\le a_j\le 2r_j+1-\varepsilon$.  Since $\delta_3$ is an integer, this implies $\delta_3\ge r_3-r_2+\delta_2$.  On the other hand,\begin{equation}\label{delta3}-a_3+2\delta_3 \le0\implies \delta_3\le\frac{a_3}{2}\le r_3+\frac{1}{2}-\varepsilon,\end{equation}so that $\delta_3\le r_3$.  So part (b) of our lemma is proved.
\end{proof}

To complete our proof of Proposition \ref{prop:doubleresidue12}, then, we 
need only show that the residue at each of the above poles  in the variable $\alpha_2$ is sufficiently small.

For ease of notation, let us denote such a pole  by $\widehat{\alpha}_2$, 
for some fixed $\delta_3$ as described  in the above lemma. We also write 
$\widehat{\alpha}_4:=-\alpha_1-\widehat{\alpha}_2-\alpha_3$, and $\widehat{\alpha}:=(\alpha_1,\widehat{\alpha}_2,\alpha_3,\widehat{\alpha}_4)$.  Then by \eqref{eq:pTr_12delta} and \eqref{eq:pTr_12deltashift}, the residue at  $\widehat{\alpha}_2$ has the following form:
\begin{align}\label{eq:res-pTr_12deltashift}
 &\Res_{\alpha_2=\widehat{\alpha}_2}\biggl(p_{T,R}^{12,(\delta_1,\delta_2)}(y;-a_3,\kappa)\biggr)
 \\& \nn= \frac{(-1)^{\delta_3}}{\delta_3!}
 \iint\limits_{\substack{\re(\alpha_1)=a_1-2\delta_1\\\re(\alpha_3)=0}} e^{\frac{\alpha_1^2+\widehat{\alpha}_2^2+\alpha_3^2+\widehat{\alpha}_4^2}{2T^2}} \int\limits_{\re(s_3)=-a_3}y_1^{\frac32+\alpha_1+2\delta_1}y_2^{2-\widehat{\alpha}_2-\alpha_3+2\delta_2}y_3^{\frac32-s_3} 
  \mathcal{F}_R(\widehat{\alpha}) 
    f_{\delta_1,\delta_2}(s_3,\widehat{\alpha})
  \\ \nonumber & 
  \cdot
  \frac{\Gamma(\frac{2+R+\alpha_1-\widehat{\alpha}_2}{4})
        \Gamma(\frac{2+R+\widehat{\alpha}_2-\alpha_1}{4})\Gamma(\frac{\widehat{\alpha}_2-\alpha_1}{2}-\delta_1)
      }{\Gamma(\frac{\alpha_1-\widehat{\alpha}_2}{2})\Gamma(\frac{\widehat{\alpha}_2-\alpha_1}{2})}
  \cdot
  \frac{\Gamma(\frac{2+R+\widehat{\alpha}_2-\alpha_3}{4})
        \Gamma(\frac{2+R+\alpha_3-\widehat{\alpha}_2}{4})  \Gamma(\frac{s_3-\alpha_3}{2}) 
      }{\Gamma(\frac{\widehat{\alpha}_2-\alpha_3}{2})
        \Gamma(\frac{\alpha_3-\widehat{\alpha}_2}{2})}
  \\ \nonumber & 
  \cdot
  \frac{
      }{}
  \frac{\Gamma(\frac{2+R+\alpha_1-\alpha_3}{4})
        \Gamma(\frac{2+R+\alpha_3-\alpha_1}{4})\Gamma(\frac{\alpha_3-\alpha_1}{2}-\delta_1)
      }{\Gamma(\frac{\alpha_1-\alpha_3}{2})\Gamma(\frac{\alpha_3-\alpha_1}{2})}
  \cdot
  \frac{\Gamma(\frac{2+R+\widehat{\alpha}_2-\widehat{\alpha}_4}{4})
        \Gamma(\frac{2+R+\widehat{\alpha}_4-\widehat{\alpha}_2}{4})  \Gamma(\frac{\widehat{\alpha}_2-\widehat{\alpha}_4}{2}-\delta_2) 
      }{\Gamma(\frac{\widehat{\alpha}_2-\widehat{\alpha}_4}{2})
        \Gamma(\frac{\widehat{\alpha}_4-\widehat{\alpha}_2}{2})}
  \\ \nonumber & 
  \cdot
  \frac{\Gamma(\frac{2+R+\alpha_1-\widehat{\alpha}_4}{4})
        \Gamma(\frac{2+R+\widehat{\alpha}_4-\alpha_1}{4})\Gamma(\frac{\widehat{\alpha}_4-\alpha_1}{2}-\delta_1)
      }{\Gamma(\frac{\alpha_1-\widehat{\alpha}_4}{2})\Gamma(\frac{\widehat{\alpha}_4-\alpha_1}{2})}
  \cdot 
  \frac{
        \Gamma(\frac{2+R+\alpha_3-\widehat{\alpha}_4}{4})\Gamma(\frac{2+R+\widehat{\alpha}_4-\alpha_3}{4})\Gamma(\frac{\alpha_3-\widehat{\alpha}_4}{2}-\delta_2)
      }{
        \Gamma(\frac{\alpha_3-\widehat{\alpha}_4}{2})\Gamma(\frac{\widehat{\alpha}_4-\alpha_3}{2})} \\ \nonumber & 
\cdot ds_3\; d\alpha_3\;d\alpha_1.
\end{align}We want our bound on \eqref{eq:res-pTr_12deltashift} to contain the factor $y_1^{3/2+a_1}y_2^{2+a_2}y_3^{3/2+a_3}$, as usual.  To effect this, we will move the line of integration in $\alpha_3$, in \eqref{eq:res-pTr_12deltashift}, from $\re\alpha_3=0$ to $\re(2-\widehat{\alpha}_2-\alpha_3+2\delta_2)=2+a_2$, which is to say, to the line

\begin{equation}\re(\alpha_3)=-a_2-\widehat{\alpha}_2+2\delta_2=-a_2+a_3+2\delta_2-2\delta_3.\label{eq:realpha3}\end{equation}The crucial observation here is that, in moving this line, we do not cross any poles.  This is by arguments very similar to those employed in the proof of the above lemma.  The only additional argument we need to make here regards the term $\Gamma(\frac{s_3-\alpha_3}{2})$, in \eqref{eq:res-pTr_12deltashift}.  But if this factor contributes a pole, then $s_3-\alpha_3\in 2\Z$; since $s_3-\widehat{\alpha}_2=-2\delta_3$ is also in $2\Z$, we conclude that $\alpha_3-\widehat{\alpha}_2\in 2\Z$, whence the pole from either the term $\Gamma(\frac{\widehat{\alpha}_2-\alpha_3}{2})$ or the term $
        \Gamma(\frac{\alpha_3-\widehat{\alpha}_2}{2})$ in the denominator 
of \eqref{eq:res-pTr_12deltashift} cancels the pole from $\Gamma(\frac{s_3-\alpha_3}{2})$.

 So in estimating \eqref{eq:res-pTr_12deltashift}, we may replace the line of integration $\re(\alpha_3)=0$ with the line given by \eqref{eq:realpha3}.  The estimation is then similar to that of \eqref{eq:pTr_12deltashift}.
 
 Specifically:  as was the case with \eqref{eq:pTr_12deltashift}, the only grouped combination  of Gamma functions  in \eqref{eq:res-pTr_12deltashift}   that contributes to the exponential factor in Stirling's formula is the second one.  In the present case, since $\im{\widehat{\alpha}_2}=\im{s_3}=\xi_3$, these Gamma functions contribute a factor of 
  \begin{equation*}  e^{-\frac{\pi}{2}(|\xi_3-\tau_3|/2 -|\xi_3-\tau_3|/2)}=e^{0}.\end{equation*}
 In other words, our exponential zero set entails no restrictions on our integration in $\xi_3=\im(s_3)$.
 
We now write
 \begin{align}\label{lambdas}(\lambda_1,\lambda_2,\lambda_3,\lambda_4)&:=(\re{\alpha_1},\re{\widehat{\alpha}_2},\re{\alpha_3},\re{\widehat{\alpha}_4})\\&=(a_1-2\delta_1,-a_3+2\delta_3,-a_2+a_3+2\delta_2-2\delta_3,-a_1 
+ a_2 + 2\delta_1 - 2\delta_2).\nn\end{align}Then \eqref{eq:res-pTr_12deltashift} yields
 \begin{align}\label{eq:res-pTr_12deltashift-b}
 &\Res_{\alpha_2=\widehat{\alpha}_2}\biggl(p_{T,R}^{12,(\delta_1,\delta_2)}(y;-a_3,\kappa)\biggr)\ll y_1^{\frac32+a_1}y_2^{2+a_2}y_3^{\frac32+a_3} \, T^{\varepsilon+R+2\delta_1+\delta_2}
  \\&
 \cdot \iiint\limits_{\substack{  \tau_1, \tau_3, \xi_3\in\R }} e^{\frac{\alpha_1^2+\widehat{\alpha}_2^2+\alpha_3^2+\widehat{\alpha}_4^2}{2T^2}}   
\big(1+\lvert\tau_1-\xi_3\rvert\big)^{\frac{R+1+\lambda_2-\lambda_1}{2}-\delta_1}
  \big(1+\lvert\tau_1-\tau_3\rvert\big)^{\frac{R+1+\lambda_3-\lambda_1}{2} -\delta_1}
\nn  \\ \nonumber & 
  \cdot  \big(1+\lvert2\tau_1+\tau_3+\xi_3 \rvert\big)^{\frac{R+1+\lambda_4-\lambda_1}{2} -\delta_1}
  \big( 1+\vert \xi_3-\tau_3\vert  \big)^{ \frac{1+R-a_3-\lambda_3}{2} }
  \big( 1+\lvert2\xi_3+\tau_1+\tau_3 \rvert\big)^{ \frac{R+1+\lambda_2-\lambda_4}{2}-\delta_2}
  \\ \nonumber & 
  \cdot  \big( 1+|2\tau_3+\tau_1+\xi_3 |\big)^{ \frac{R+1+\lambda_3-\lambda_4}{2}-\delta_2} 
        \;  d\xi_3\; d\tau_3\;d\tau_1.
\end{align}The factor
$$e^{\frac{\alpha_1^2+\widehat{\alpha}_2^2+\alpha_3^2+\widehat{\alpha}_4^2}{2T^2}}$$ in \eqref{eq:res-pTr_12deltashift-b} is of exponential decay in $\tau_1$   if $|\tau_1| \gg T^{1+\varepsilon}$, and similarly for the variables $\tau_3$ and $\xi_3$.  So for our estimate, we may restrict attention to the domain where $|\tau_1|, |\tau_3|,|\xi_3| \ll T^{1+\varepsilon}$. On such a domain, each of the other factors in our integrand is $\ll T^{c+\varepsilon}$, where $c$ is the exponent on that factor.  So \eqref{eq:res-pTr_12deltashift-b} implies
\begin{align}\label{s2res}
 &\Res_{\alpha_2=\widehat{\alpha}_2}\biggl(p_{T,R}^{12,(\delta_1,\delta_2)}(y;-a_3,\kappa)\biggr)\ll y_1^{\frac32+a_1}y_2^{2+a_2}y_3^{\frac32+a_3} \, T^{\varepsilon+R+2\delta_1+\delta_2}
  \\&\nn
 \cdot \iiint\limits_{|\tau_1|, |\tau_3|,|\xi_3| \ll T^{1+\varepsilon}}T^{\varepsilon+4R+3+\frac{-a_3-3\lambda_1+2\lambda_2+\lambda_3-\lambda_4}{2}-3\delta_1-2\delta_2}  \;  d\xi_3\; d\tau_3\;d\tau_1  \\&\nn \ll y_1^{\frac32+a_1}y_2^{2+a_2}y_3^{\frac32+a_3} \,T^{\varepsilon+4R+6-a_1-a_2-a_3 +\delta_1+\delta_2+\delta_3}.
\end{align}By \eqref{delta3}, we have $\delta_3-a_3\le \frac{1}{2}-r_3$; also, $\delta_j-a_j\le r_j-1+1-2r_j=-r_j$ for $j=1,2$. Then \eqref{s2res} yields
\begin{align*}
 &\Res_{\alpha_2=\widehat{\alpha}_2}\biggl(p_{T,R}^{12,(\delta_1,\delta_2)}(y;-a_3,\kappa)\biggr)\ll y_1^{\frac32+a_1}y_2^{2+a_2}y_3^{\frac32+a_3} \,   T^{\varepsilon+4R+\frac{13}2 -r_1-r_2-r_3}.
\end{align*}
In other words, the sum of the residue terms in \eqref{eq:KappaPlusResDouble} also has a bound of the magnitude stipulated in Proposition \ref{prop:doubleresidue12}.  This completes the proof of that proposition.  \end{proof}

We now turn to our estimate of the term \eqref{eq:pTr_13delta}.  The analysis here is similar to that of \eqref{eq:pTr_12delta}, but different enough that some detail is merited.

We have:
\begin{prop}\label{prop:doubleresidue13} Let $r_1,r_3\geq 1$, $r_2\geq0$ be integers, and $0<\varepsilon<1$. Suppose $a_1,a_2,a_3$ satisfy the hypotheses of Theorem \ref{th:pTRbound}.  If $0\le \delta_j\le r_j-1$ for $j=1,3$, then
\begin{equation}\label{eq:Prop-s1s3-Residue}
 \left|   p_{T,R}^{13,(\delta_1,\delta_3)}(y;-a_2) \right| \ll    y_1^{\frac32+a_1}y_2^{2+a_2}y_3^{\frac32+a_3} \, T^{\varepsilon+4R+\frac{13}{2}+\delta_{r_2,0}-r_1-r_2-r_3}.
\end{equation}\end{prop}

\begin{proof} 

To obtain the desired bound on  $p_{T,R}^{13,(\delta_1,\delta_3)}(y;-a_3)$, we will need to shift the lines of integration in both the $\alpha_1$ and $\alpha_2$ variables, so that the exponents of $y_1$ and $y_3$  become $3/2+a_1$ and $3/2+a_3$ respectively.  In doing so, we will pick up residues, whence\begin{equation}\label{eq:KappaPlusResDouble2}p_{T,R}^{13,(\delta_1,\delta_3)}(y;-a_2)  =  p_{T,R}^{13,(\delta_1,\delta_3)}(y;-a_2,\kappa) + \sum\mbox{Residues}, \end{equation}
where 
\begin{align}\label{eq:pTr_13deltashift}
 & p_{T,R}^{13,(\delta_1,\delta_3)}(y;-a_2,\kappa) := 
 \iiint\limits_{\re(\alpha)=\kappa} e^{\frac{\alpha_1^2+\cdots+\alpha_4^2}{2T^2}} \int\limits_{\re(s_2)=-a_2}y_1^{\frac32+\alpha_1+2\delta_1}y_2^{2-s_2}y_3^{\frac32-\alpha_2+2\delta_3} 
  \cdot\fr\,
    g_{\delta_1,\delta_3}(s_2,\alpha)
  \\ \nonumber & 
  \cdot
  \frac{\Gamma(\frac{2+R+\alpha_1-\alpha_2}{4})
        \Gamma(\frac{2+R+\alpha_2-\alpha_1}{4})\Gamma(\frac{\alpha_2-\alpha_1}{2}-\delta_1)
      }{\Gamma(\frac{\alpha_1-\alpha_2}{2})\Gamma(\frac{\alpha_2-\alpha_1}{2})}
  \cdot
  \frac{\Gamma(\frac{2+R+\alpha_2-\alpha_3}{4})
        \Gamma(\frac{2+R+\alpha_3-\alpha_2}{4})  \Gamma(\frac{\alpha_2-\alpha_3}{2}-\delta_3)
      }{\Gamma(\frac{\alpha_2-\alpha_3}{2})
        \Gamma(\frac{\alpha_3-\alpha_2}{2})}
  \\ \nonumber & 
  \cdot
  \frac{
      }{}
  \frac{\Gamma(\frac{2+R+\alpha_1-\alpha_3}{4})
        \Gamma(\frac{2+R+\alpha_3-\alpha_1}{4})\Gamma(\frac{\alpha_3-\alpha_1}{2}-\delta_1)
      }{\Gamma(\frac{\alpha_1-\alpha_3}{2})\Gamma(\frac{\alpha_3-\alpha_1}{2})}
  \cdot
  \frac{\Gamma(\frac{2+R+\widehat{\alpha}_2-\alpha_4}{4})
        \Gamma(\frac{2+R+\alpha_4-\widehat{\alpha}_2}{4})  \Gamma(\frac{\widehat{\alpha}_2-\alpha_4}{2}-\delta_3) 
      }{\Gamma(\frac{\widehat{\alpha}_2-\alpha_4}{2})
        \Gamma(\frac{\alpha_4-\widehat{\alpha}_2}{2})}
  \\ \nonumber & 
  \cdot
  \frac{\Gamma(\frac{2+R+\alpha_1-\alpha_4}{4})
        \Gamma(\frac{2+R+\alpha_4-\alpha_1}{4})\Gamma(\frac{\alpha_4-\alpha_1}{2}-\delta_1)
      }{\Gamma(\frac{\alpha_1-\alpha_4}{2})\Gamma(\frac{\alpha_4-\alpha_1}{2})}\cdot \frac{\Gamma(\frac{2+R+\alpha_4-\alpha_3}{4})
        \Gamma(\frac{2+R+\alpha_3-\alpha_4}{4})\Gamma(\frac{s_2+\alpha_1+\alpha_3}{2})\Gamma(\frac{s_2+\alpha_1+\alpha_4}{2})
      }{\Gamma(\frac{\alpha_4-\alpha_3}{2})
        \Gamma(\frac{\alpha_3-\alpha_4}{2})} \\\nn&
\cdot \;
 ds_2\; d\alpha,
\end{align}
and the residues that appear depend on the particular choice of $\kappa=(\kappa_1,\kappa_2,\kappa_3)$.   As before, we've grouped the Gamma factors in an auspicious manner.

To obtain the desired exponents on $y_1$ and $y_3$, in \eqref{eq:pTr_13deltashift}, we will choose $\kappa=(\kappa_1,\kappa_2,\kappa_3)\in\R^3$  
such that  
 \begin{equation*} \re(\alpha_1+2\delta_1)=\kappa_1+2\delta_1=a_1, \qquad \re( -\alpha_2+2\delta_3)=-\kappa_2+2\delta_3=a_3,  \end{equation*}by putting\begin{equation}\kappa=(\kappa_1,\kappa_2,\kappa_3)=(a_1-2\delta_1,-a_3+2\delta_3 ,0) \quad(\hbox{and } \kappa_4=-\kappa_1-\kappa_2-\kappa_3). \label{newkappadouble2}\end{equation}For this value of $\kappa$, we will obtain an estimate of the desired magnitude for $p_{T,R}^{13,(\delta_1,\delta_3)}(y;-a_2,\kappa) .$
 Subsequently we will, as before, show that the residues  in \eqref{eq:KappaPlusResDouble2}  are small.

\vskip 12pt

\noindent
\underline{\bf{Step 1}: Bounding the shifted integral $p_{T,R}^{13,(\delta_1,\delta_2)}$}
\vskip 10pt

Here we estimate the term $ p_{T,R}^{13,(\delta_1,\delta_3)}(y;-a_2,\kappa)$, with $\kappa=(\kappa_1,\kappa_2,\kappa_3)$ as in \eqref{newkappadouble2}.

Only the last grouped combination of Gamma functions, in \eqref{eq:pTr_13deltashift}, contributes to the exponential factor in Stirling's formula, 
for the same reasons as we discussed in the proof of Proposition \ref{prop:doubleresidue12}.  In the present case, this last term contributes a factor of  \begin{equation*}  e^{-\frac{\pi}{2}(|\tau_4-\tau_3|/4+|\tau_3-\tau_4|/4+|\xi_2+\tau_1+\tau_3|/2+|\xi_2+\tau_1+\tau_4|/2-|\tau_4-\tau_3|/2-|\tau_3-\tau_4|/2)}=e^{-\frac{\pi}{4}( |\xi_2+\tau_1+\tau_3|+|\xi_2+\tau_1\tau_4|-|\tau_3-\tau_4|)}.\end{equation*}
As the integrand in \eqref{eq:pTr_13deltashift}  is invariant under $\alpha_3\leftrightarrow \alpha_4$, we may assume that
 $ \tau_3 \ge \tau_4$, so that the exponential factor in question equals
 $$e^{-\frac{\pi}{4}( |\xi_2+\tau_1+\tau_3|+|\xi_2+\tau_1+\tau_4|-\tau_3+\tau_4)}.$$
Then the corresponding  exponential zero set is seen to be 
\begin{align*}
 \mathcal{R} & :=  \left\{ (-a_2+i\xi_2)\in \C |   -\tau_1-\tau_3\leq \xi_2\leq   -\tau_1-\tau_4 \right\} .\end{align*}

We replace the Gamma factors in \eqref{eq:pTr_13deltashift} with their corresponding polynomial terms; we get
\begin{align} \label{eq:pTr_13Shifted}&
  \left| p_{T,R}^{13,(\delta_1,\delta_3)}(y;-a_2,\kappa)   \right|  \ll  
  y_1^{\frac32+a_1}y_2^{2+a_2}y_3^{\frac32+a_3}
 \, T^{\varepsilon+R+2\delta_1+\delta_3} \\&
 \cdot \iiint\limits_{\substack{ \abs{\tau_1},\abs{\tau_2},\abs{\tau_3} \leq T^{1+\varepsilon} \\ 0 \leq \tau_1+\tau_2+2\tau_3  }} 
  \hskip -8pt
  \big(1+\lvert\tau_1-\tau_2\rvert\big)^{\frac{R+1+\kappa_2-\kappa_1}{2}-\delta_1}\nn
  \big(1+\lvert\tau_1-\tau_3\rvert\big)^{\frac{R+1+\kappa_3-\kappa_1}{2}-\delta_1}
  \big(1+\lvert\tau_1-\tau_4\rvert\big)^{\frac{R+1+\kappa_4-\kappa_1}{2}-\delta_1}
  \\ \nonumber & 
  \cdot
  \big( 1+\lvert\tau_2-\tau_3 \rvert\big)^{ \frac{R+1+\kappa_2-\kappa_3}{2}-\delta_3}
  \big( 1+\lvert\tau_2-\tau_4 \rvert\big)^{ \frac{R+1+\kappa_2-\kappa_4}{2}-\delta_3}
  \big( 1+\tau_3-\tau_4 \big)^{  \frac{2+R}{2}} 
  \\ \nonumber & 
  \cdot
 \int _{\xi_2 =-\tau_1 -\tau_3}^{-\tau_1-\tau_4}
   \big(1+\xi_2+\tau_1+\tau_3\big)^{\frac{-1-a_2+\kappa_1+\kappa_3 }{2} }
        \big(1-(\xi_2+\tau_1+\tau_4)\big)^{\frac{-1-a_2+\kappa_1+\kappa_4}{2} }
        \;  d\xi_2\; d\tau.
\end{align}The change of variables
 \[  \quad \xi_2\mapsto \xi_2-\tau_1-\tau_3, \quad \quad T_j=\tau_j-\tau_{j+1} \quad (1\le j\le3) \]
applied to \eqref{eq:pTr_13Shifted} then gives\begin{align} \label{pTr_13NewShifted}
  \left| p_{T,R}^{13,(\delta_1,\delta_3)}(y;-a_2,\kappa)   \right| & \ll  
  y_1^{\frac32+a_1}y_2^{2+a_2}y_3^{\frac32+a_3}
 \, T^{\varepsilon+R+2\delta_1+\delta_3}  \\&
 \cdot \iiint\limits_{ \substack{0\le \abs{T_1},\abs{T_2}, T_3 \leq T^{1+\varepsilon}   }} 
  \hskip -8pt
  \big(1+\lvert T_1\rvert\big)^{\frac{R+1+\kappa_2-\kappa_1}{2}-\delta_1}\nn
  \big(1+\lvert T_1+T_2\rvert\big)^{\frac{R+1+\kappa_3-\kappa_1}{2}-\delta_1}
  \\ \nonumber & 
  \cdot  \big(1+\lvert T_1+T_2+T_3\rvert\big)^{\frac{R+1+\kappa_4-\kappa_1}{2}-\delta_1}
  \big( 1+\lvert T_2 \rvert\big)^{ \frac{R+1+\kappa_2-\kappa_3}{2}-\delta_3}\\\nn&\cdot
  \big( 1+\lvert T_2+T_3 \rvert\big)^{ \frac{R+1+\kappa_2-\kappa_4}{2}-\delta_3}
  \big( 1+T_3 \big)^{  \frac{2+R}{2}} 
  \\ \nonumber & 
  \cdot
 \int _{\xi_2=0}^{T_3}
         \big(1+\xi_2\big)^{\frac{-1-a_2+\kappa_1+\kappa_3 }{2} }
        \big(1+T_3-\xi_2 \big)^{\frac{-1-a_2+\kappa_1+\kappa_4}{2} }
        \;  d\xi_2\; dT_1\;dT_2\;dT_3.
\end{align}Now by Lemma \ref{lemmaIntegralBound}, we find that
\begin{align}\label{T3-int}& \int _{\xi_2=0}^{T_3}
         \big(1+\xi_2\big)^{\frac{-1-a_2+\kappa_1+\kappa_3 }{2} }
        \big(1+T_3-\xi_2 \big)^{\frac{-1-a_2+\kappa_1+\kappa_4}{2} }
        \;  d\xi_2\\&\nn \ll 
 (1+T_3)^{-\min\big\{\frac{1+a_2-\kappa_1-\kappa_3 }{2} ,\, \frac{1+a_2-\kappa_1-\kappa_4}{2} ,\,a_2-\frac{2\kappa_1+\kappa_3+\kappa_4 }{2}\big\}+\varepsilon}\\&=(1+T_3)^{\frac{-1-a_2+2\kappa_1+\kappa_3+\kappa_4 }{2}+\frac{1}{2}\max\{ -\kappa_1-\kappa_4 ,\, -\kappa_1-\kappa_3 ,\,1-a_2 \}+\varepsilon}.\nn\end{align}Further, it follows from  \eqref{newkappadouble2} that$$ -\kappa_1-\kappa_4 ,\, -\kappa_1-\kappa_3\le0,$$and that $1-a_2<0$ unless $r_2=0$, in which case $1-a_2<2$.  In either case, $1-a_2< 2\delta_{r_2,0}$.  So \eqref{pTr_13NewShifted} and \eqref{T3-int} give
 \begin{align} \label{pTr_13NewShifted-b}
  \left| p_{T,R}^{13,(\delta_1,\delta_3)}(y;-a_2,\kappa)   \right| & \ll  
  y_1^{\frac32+a_1}y_2^{2+a_2}y_3^{\frac32+a_3}
\, T^{\varepsilon+R+2\delta_1+\delta_3} \\\nn&\cdot\iiint\limits_{ \substack{0\le \abs{T_1},\abs{T_2}, T_3 \leq T^{1+\varepsilon}   }} 
  \hskip -8pt
  \big(1+\lvert T_1\rvert\big)^{\frac{R+1+\kappa_2-\kappa_1}{2}-\delta_1} 
\big(1+\lvert T_1+T_2\rvert\big)^{\frac{R+1+\kappa_3-\kappa_1}{2}-\delta_1} \\&
 \cdot \nn
  \big(1+\lvert T_1+T_2+T_3\rvert\big)^{\frac{R+1+\kappa_4-\kappa_1}{2}-\delta_1}  \big( 1+\lvert T_2 \rvert\big)^{ \frac{R+1+\kappa_2-\kappa_3}{2}-\delta_3}
  \\ \nonumber & 
  \cdot
  \big( 1+\lvert T_2+T_3 \rvert\big)^{ \frac{R+1+\kappa_2-\kappa_4}{2}-\delta_3}
  \big( 1+T_3 \big)^{  \frac{R+1-a_2+2\kappa_1+\kappa_3+\kappa_4 }{2}+\delta_{r_2,0}+\varepsilon } \; dT_1\;dT_2\;dT_3.
\end{align}
Then, because$$ T_1, T_2, T_3 \ll T^{1+\varepsilon}$$on the indicated domains of integration, and because the length of each domain of integration 
is also $\ll T^{1+\varepsilon}$, we see that
\begin{align*}&
  \left| p_{T,R}^{13,(\delta_1,\delta_3)}(y;-a_2,\kappa)   \right|  \ll  
  y_1^{\frac32+a_1}y_2^{2+a_2}y_3^{\frac32+a_3}
  \, T^{\varepsilon+4R+6+\delta_1+\delta_3+\delta_{r_2,0}-a_1-\frac{a_2}{2}-a_3}.\end{align*}
But our assumptions on the $a_j$'s and the $\delta_j$'s imply that $\delta_j-a_j\le -r_j$ for $j=1,3$ and $-\frac{a_2}{2}\le\frac{1}{2}-r_2$, so 
we conclude that\begin{align} \label{eq:pTr_13NewerShifted}   \left| p_{T,R}^{13,(\delta_1,\delta_3)}(y;-a_2,\kappa)   \right|  \ll  & \;
  y_1^{\frac32+a_1}y_2^{2+a_2}y_3^{\frac32+a_3}
 \, T^{\varepsilon+4R+\frac{13}{2}+\delta_{r_2,0}-r_1-r_2-r_3}, \end{align}which gives us a bound of the desired magnitude on the  shifted integral  $p_{T,R}^{13,(\delta_1,\delta_3)}(y;-a_2,\kappa) $ in \eqref{eq:KappaPlusResDouble2}.
\vskip 12pt
%\newpage
\noindent
\underline{\bf{Step 2}: Bounding the residue terms}
\vskip 10pt

Next, we estimate the residues in \eqref{eq:KappaPlusResDouble2}, which arise from moving the lines of integration in \eqref{eq:pTr_13delta}, to get  \eqref{eq:pTr_13deltashift}.
 
The locations of the poles in question are as follows.

\begin{lem}\label{lem:shiftalpha1alpha2poles-b}Suppose the lines of integration, 
 in \eqref{eq:pTr_13delta}, are shifted  from $\re(\alpha_1,\alpha_2,\alpha_3)=(0,0,0)$ to $\re(\alpha_1,\alpha_2,\alpha_3)=(a_1-2\delta_1,-a_3+2\delta_3 ,0)$.  Then:
 
 \begin{enumerate}
 \item[\rm (a)]For a fixed $s_2$ and $\alpha_3$, any pole crossed in the $\alpha_1$ variable belongs to the set \begin{equation}\label{eq:alpha13poles-a}
  \left\{ -s_2-\alpha_3-2\delta_2\,\vert \, \delta_2\in\Z_{\ge0},\,\max\{0,r_2-(r_1-\delta_1)\}\le \delta_2\le r_2  \right\}. \end{equation}
 \item[\rm (b)] For a fixed $s_2$ and $\alpha_3$, any pole crossed in the 
$\alpha_2$ variable belongs to the set \begin{equation}\label{eq:alpha13poles-b}
  \left\{ s_2-\alpha_3+2\delta_2\,\vert \, \delta_2\in\Z_{\ge0},\,\max\{0,r_2-(r_3-\delta_3)\}\le \delta_2\le r_2  \right\}.
\end{equation}\end{enumerate}
\end{lem}

\begin{proof}  As before, no poles will arise from the factors
 \[\textstyle \Gamma\left(\frac{2+R+\alpha_j-\alpha_k}{4}\right)  \qquad \mbox{($1\leq j\neq k\leq 4$)} \]in  \eqref{eq:pTr_13deltashift},
if $R$ is sufficiently large.

Nor will any of the terms of the form
$$  \frac{ \Gamma(\frac{\alpha_k-\alpha_j}{2}-\delta_n)
      }{\Gamma(\frac{\alpha_j-\alpha_k}{2})\Gamma(\frac{\alpha_k-\alpha_j}{2})} $$give rise to any poles, for the same reasons as before.      
The only terms remaining to consider are the factors $$\textstyle \Gamma(\frac{s_2+\alpha_1+\alpha_3}{2})\quad\mbox{and}\quad \Gamma(\frac{s_2+\alpha_1+\alpha_4}{2})=\Gamma(\frac{s_2-\alpha_2-\alpha_3}{2}).$$The former of these factors will give rise to poles when we shift the line of integration in $\alpha_1$; the latter will do so when we shift the line in $\alpha_2$.

Consider the first of these factors,  $\Gamma(\frac{s_2+\alpha_1+\alpha_3}{2})$.  For fixed $s_2$ and $\alpha_3$, this factor has poles, as a function of $\alpha_1$, whenever\begin{align}\alpha_1=-s_2-\alpha_3-2\delta_2 \quad(\delta_2\in\Z_{\ge0}).\label{eq:alpha1poles-b}\end{align}But  for such an $\alpha_1$ to lie between the initial line of integration $\re{\alpha_1}=0$ and the terminal line  $\re{\alpha_1}=a_1-2\delta_1$, we 
must have
\begin{align}0\le\re{\alpha_1}=\re{(-s_2-\alpha_3-2\delta_2)}=a_2-2\delta_2 \le a_1-2\delta_1 .\label{eq:alpha1poles2b}\end{align}But$$a_2-2\delta_2\le a_1-2\delta_1\implies \delta_2\ge \frac{a_2-a_1}{2}+\delta_1\ge 
r_2-r_1+\delta_1-1+\varepsilon.$$As $\delta_2$ is a nonnegative integer, we therefore have $\delta_2\ge \mbox{max}\{0, r_2-r_1+\delta_1\}$.  On the other hand,\begin{equation}\label{delta2}a_2-2\delta_2 \ge0\implies \delta_2\le\frac{a_2}{2}\le r_2+\frac{1}{2}-\varepsilon,\end{equation}so that $\delta_2\le r_2$.  So part (a) of our lemma is proved.

Part (b) is similar:  as a function of $\alpha_2$, $\Gamma(\frac{s_2-\alpha_2-\alpha_3}{2})$ has poles,  for fixed $s_2$ and $\alpha_3$, whenever\begin{align}\alpha_2=s_2-\alpha_3+2\delta_2 \quad(\delta_2\in\Z_{\ge0}).\label{eq:alpha2poles-b}\end{align}But  for such an $\alpha_2$ to lie between $\re{\alpha_2}=0$ and   $\re{\alpha_2}=-a_3+2\delta_3$, we must have
\begin{align}-a_3+2\delta_3\le\re{\alpha_2}=\re{(s_2-\alpha_3+2\delta_2)}=-a_2+2\delta_2 \le 0 .\label{alpha2poles2-c}\end{align}
We conclude from  \eqref{alpha2poles2-c} and the fact that $\delta_3$ is a nonnegative integer that
$$\max\{0,r_2-(r_3-\delta_3)\}\le \delta_2\le r_2,$$so part (b) of our lemma is proved.
\end{proof}

Therefore, to complete our proof of Proposition \ref{prop:doubleresidue13}, it will suffice to show that the residue at each of the above poles in 
$\alpha_1$ or $\alpha_2$ is sufficiently small.

We will consider the poles in $\alpha_1$ only; those in $\alpha_2$ may be 
treated in a very similar fashion.  
Let us, then, denote such a pole in $\alpha_1$ by $\widetilde{\alpha}_1$, 
for some fixed $\delta_2$ and $\alpha_3$ as described  in part (a) of the 
above lemma. We also write $\widetilde{\alpha}_4:=-\widetilde{\alpha}_1-\alpha_2-\alpha_3$, and $\widetilde{\alpha}:=(\widetilde{\alpha}_1,\alpha_2,\alpha_3,\widetilde{\alpha}_4)$.  Then by \eqref{eq:pTr_13delta} and \eqref{eq:pTr_13deltashift}, the residue at  $\widetilde{\alpha}_1$ has 
the following form:
\begin{align}\label{eq:res-pTr_13deltashift}
 &\Res_{\alpha_1=\widetilde{\alpha}_1}\biggl(p_{T,R}^{13,(\delta_1,\delta_3)}(y;-a_2,\kappa)\biggr)
 \\& \nn= \frac{(-1)^{\delta_2}}{\delta_2!}\hskip-5pt
 \iiint\limits_{\substack{\re(\alpha_2)=-a_3+2\delta_3\\ \re(\alpha_3)=0}} e^{\frac{\widetilde{\alpha}_1^2+\alpha_2^2+\alpha_3^2+\widetilde{\alpha}_4^2}{2T^2}} \int\limits_{\re(s_2)=-a_2}y_1^{\frac32 -s_2-\alpha_3-2\delta_2+2\delta_1}y_2^{2-s_2}y_3^{\frac32-\alpha_2+2\delta_3}  
 \mathcal F_R(\alpha)\,    g_{\delta_1,\delta_3}(s_2,\widetilde{\alpha})
  \\ \nonumber & 
  \cdot
  \frac{\Gamma(\frac{2+R+\widetilde{\alpha}_1-\alpha_2}{4})
        \Gamma(\frac{2+R+\alpha_2-\widetilde{\alpha}_1}{4})\Gamma(\frac{\alpha_2-\widetilde{\alpha}_1}{2}-\delta_1)
      }{\Gamma(\frac{\widetilde{\alpha}_1-\alpha_2}{2})\Gamma(\frac{\alpha_2-\widetilde{\alpha}_1}{2})}
  \cdot
  \frac{\Gamma(\frac{2+R+\alpha_2-\alpha_3}{4})
        \Gamma(\frac{2+R+\alpha_3-\alpha_2}{4})  \Gamma(\frac{\alpha_2-\alpha_3}{2}-\delta_3)
      }{\Gamma(\frac{\alpha_2-\alpha_3}{2})
        \Gamma(\frac{\alpha_3-\alpha_2}{2})}
  \\ \nonumber & 
  \cdot
  \frac{
      }{}
  \frac{\Gamma(\frac{2+R+\widetilde{\alpha}_1-\alpha_3}{4})
        \Gamma(\frac{2+R+\alpha_3-\widetilde{\alpha}_1}{4})\Gamma(\frac{\alpha_3-\widetilde{\alpha}_1}{2}-\delta_1)
      }{\Gamma(\frac{\widetilde{\alpha}_1-\alpha_3}{2})\Gamma(\frac{\alpha_3-\widetilde{\alpha}_1}{2})}
  \cdot
  \frac{\Gamma(\frac{2+R+\widehat{\alpha}_2-\widetilde{\alpha}_4}{4})
        \Gamma(\frac{2+R+\widetilde{\alpha}_4-\widehat{\alpha}_2}{4})  \Gamma(\frac{\widehat{\alpha}_2-\widetilde{\alpha}_4}{2}-\delta_3) 
      }{\Gamma(\frac{\widehat{\alpha}_2-\widetilde{\alpha}_4}{2})
        \Gamma(\frac{\widetilde{\alpha}_4-\widehat{\alpha}_2}{2})}
  \\ \nonumber & 
  \cdot
  \frac{\Gamma(\frac{2+R+\widetilde{\alpha}_1-\widetilde{\alpha}_4}{4})
        \Gamma(\frac{2+R+\widetilde{\alpha}_4-\widetilde{\alpha}_1}{4})\Gamma(\frac{\widetilde{\alpha}_4-\widetilde{\alpha}_1}{2}-\delta_1)
      }{\Gamma(\frac{\widetilde{\alpha}_1-\widetilde{\alpha}_4}{2})\Gamma(\frac{\widetilde{\alpha}_4-\widetilde{\alpha}_1}{2})}
 \cdot \frac{\Gamma(\frac{2+R+\widetilde{\alpha}_4-\alpha_3}{4})
        \Gamma(\frac{2+R+\alpha_3-\widetilde{\alpha}_4}{4}) \Gamma(\frac{s_2-\alpha_2-\alpha_3}{2})
      }{\Gamma(\frac{\widetilde{\alpha}_4-\alpha_3}{2})
        \Gamma(\frac{\alpha_3-\widetilde{\alpha}_4}{2})}
 \,   ds_2\; d\alpha_3\;d\alpha_2.
\end{align}
In order that our bound on \eqref{eq:res-pTr_13deltashift} contain the factor $y_1^{3/2+a_1}y_2^{2+a_2}y_3^{3/2+a_3}$,  we now move the line of integration in $\alpha_3$, in \eqref{eq:res-pTr_13deltashift}, from $\re\alpha_3=0$ to $\re( 3/2 -s_2-\alpha_3-2\delta_2+2\delta_1)=3/2+a_1$, or 
equivalently
\begin{equation}\re(\alpha_3)=-  a_1+a_2+2\delta_1-2\delta_2.\label{eq:realpha3-b}\end{equation}In moving this line, we do not cross any poles.  
This is by the same kinds of arguments as were used above.  In particular 
we note that, if the factor  $\Gamma(\frac{s_2-\alpha_2-\alpha_3}{2})$ has a pole, then $s_2-\alpha_2-\alpha_3\in 2\Z$; since $s_2+\widetilde{\alpha}_1+\alpha_3=s_2-\alpha_2-\widetilde{\alpha}_4=-2\delta_2$ is also in $2\Z$ (by assumption), we conclude that $\alpha_3-\widetilde{\alpha}_4\in 2\Z$, whence the pole from either the term $\Gamma(\frac{\widetilde{\alpha}_4-\alpha_3}{2})$ or the term $
        \Gamma(\frac{\alpha_3-\widetilde{\alpha}_4}{2})$ in the denominator of \eqref{eq:res-pTr_13deltashift} cancels the pole from $\Gamma(\frac{s_2-\alpha_2-\alpha_3}{2})$.

 So in estimating \eqref{eq:res-pTr_13deltashift}, we may replace the line of integration $\re(\alpha_3)=0$ with the line given by \eqref{eq:realpha3-b}.  The estimation then proceeds as follows.  First, the only grouped combination  of Gamma functions  in \eqref{eq:res-pTr_13deltashift}   
that contributes to the exponential factor in Stirling's formula is the last one.  Since $$\im{\widetilde{\alpha}_4}=\im{(-\widetilde{\alpha}_1-\alpha_2-\alpha_3)}=\im{(s_2-\alpha_2)}=\xi_2-\tau_2,$$ these Gamma functions contribute a factor of 
  \begin{equation*}  e^{-\frac{\pi}{2}(|\xi_2-\tau_2-\tau_3|/2 -|\xi_2-\tau_2-\tau_3|/2)}=e^{0}.\end{equation*}
So our exponential zero set here places no restrictions on our domain of integration in $\xi_2=\im(s_2)$.
 
Next, we write
 \begin{align}\label{lambdas-new}(\lambda_1,\lambda_2,\lambda_3,\lambda_4)&:=(\re{\widetilde{\alpha}_1},\re{{\alpha}_2},\re{\alpha_3},\re{\widetilde{\alpha}_4})\\&=(a_1-2\delta_1,-a_3+2\delta_3,-  a_1+a_2+2\delta_1-2\delta_2,-a_2 + a_3 +2\delta_2- 2\delta_3 ).\nn\end{align}Then \eqref{eq:res-pTr_13deltashift} yields
 \begin{align}\label{eq:res-pTr_13deltashift-b}
 &\Res_{\alpha_1=\widetilde{\alpha}_1}\biggl(p_{T,R}^{13,(\delta_1,\delta_3)}(y;-a_2,\kappa)\biggr)\ll y_1^{\frac32+a_1}y_2^{2+a_2}y_3^{\frac32+a_3} \, T^{\varepsilon+R+2\delta_1+\delta_3}
  \\&
 \cdot \iiint\limits_{\substack{  \tau_2, \tau_3, \xi_2\in\R }} e^{\frac{\widetilde{\alpha}_1^2+{\alpha}_2^2+\alpha_3^2+\widetilde{\alpha}_4^2}{2T^2}}   \big(1+\lvert\tau_2+\tau_3+\xi_2\rvert\big)^{\frac{R+1+\lambda_2-\lambda_1}{2}-\delta_1}
  \big(1+\lvert2\tau_3+\xi_2\rvert\big)^{\frac{R+1+\lambda_3-\lambda_1}{2} -\delta_1}
\nn  \\ \nonumber & 
  \cdot  \big(1+\lvert2\xi_2-\tau_2+\tau_3\rvert\big)^{\frac{R+1+\lambda_4-\lambda_1}{2} -\delta_1}
  \big( 1+\vert  \tau_3-\tau_2\vert  \big)^{ \frac{R+1+\lambda_2-\lambda_3}{2} -\delta_3}
  \big( 1+\lvert2\tau_2-\xi_2 \rvert\big)^{ \frac{R+1+\lambda_2-\lambda_4}{2}-\delta_3}
  \\ \nonumber & 
  \cdot  \big( 1+|\xi_2-\tau_2-\tau_3 |\big)^{ \frac{R+1-a_2-\lambda_2-\lambda_3}{2}} 
        \;  d\xi_2\; d\tau_3\;d\tau_2.
\end{align}The factor
$$e^{\frac{\widetilde{\alpha}_1^2+{\alpha}_2^2+\alpha_3^2+\widetilde{\alpha}_4^2}{2T^2}}$$ in \eqref{eq:res-pTr_13deltashift-b} is of exponential decay in $\tau_2$   if $|\tau_2| \gg T^{1+\varepsilon}$, and similarly for the variables $\tau_3$ and $\xi_2$.  So for our estimate, we need only consider the domain where $|\tau_2|, |\tau_3|,|\xi_2| \ll T^{1+\varepsilon}$. On this domain, each of the other factors in the integrand of \eqref{eq:res-pTr_13deltashift-b} is $\ll T^{c+\varepsilon}$, where $c$ is the exponent on that factor.  So \eqref{eq:res-pTr_13deltashift-b} implies
\begin{align*}
 &\Res_{\alpha_1=\widetilde{\alpha}_1}\biggl(p_{T,R}^{13,(\delta_1,\delta_3)}(y;-a_2,\kappa)\biggr)\ll y_1^{\frac32+a_1}y_2^{2+a_2}y_3^{\frac32+a_3} \, T^{\varepsilon+2\delta_1+\delta_3}
  \\&\nn
 \cdot \iiint\limits_{|\tau_2|, |\tau_3|,|\xi_2| \ll T^{1+\varepsilon}}   
T^{\varepsilon+3R+3+\frac{-a_2-3\lambda_1+2\lambda_2-\lambda_3}{2}-3\delta_1-2\delta_3}  \;  d\xi_2\; d\tau_3\;d\tau_2  \\&\nn \ll y_1^{\frac32+a_1}y_2^{2+a_2}y_3^{\frac32+a_3} \, T^{\varepsilon+4R+6-a_1-a_2-a_3 +\delta_1+\delta_2+\delta_3}.
\end{align*}So, by \eqref{delta2} and the assumptions $\delta_j-a_j\le r_j-1+1-2r_j\le -r_j$ for $j=1,3$, we have
\begin{align*}
 &\Res_{\alpha_1=\widetilde{\alpha}_1}\biggl(p_{T,R}^{13,(\delta_1,\delta_3)}(y;-a_2,\kappa)\biggr)\ll y_1^{\frac32+a_1}y_2^{2+a_2}y_3^{\frac32+a_3} \,  \, T^{\varepsilon+4R+\frac{13}{2} -r_1-r_2-r_3  }.
\end{align*}
So the sum of the residue terms in \eqref{eq:KappaPlusResDouble2} also has a bound of the magnitude described in Proposition \ref{prop:doubleresidue13}, and the proposition is proved.  \end{proof}

\subsection{\bf Bounds for the triple residue terms} \label{sec:tripleresbound}

There is only one type of triple residue term  to consider, namely $  p_{T,R}^{123,(\delta_1,\delta_2,\delta_3)}(y) $.    This term may be obtained by taking the residue, at $s_3= \alpha_2-2\delta_3$, of the double residue term $  p_{T,R}^{12,(\delta_1,\delta_2)}(y;-a_3) $ defined by \eqref{eq:pTr_12delta}.  Thus  

\begin{align} \label{eq:pTr_123delta}
  p_{T,R}^{123,(\delta_1,\delta_2,\delta_3)}(y ) = & \iiint\limits_{\re(\alpha)=0} e^{\frac{\alpha_1^2+\cdots+\alpha_4^2}{2T^2}}\;  y_1^{\frac32+\alpha_1+2\delta_1}y_2^{2+\alpha_1+\alpha_4+2\delta_2}y_3^{\frac32-\alpha_2+2\delta_3} \; 
 \mathcal F_R(\alpha)\,\Gamma_R(\alpha) 
  \\ \nonumber & \hskip 48pt
  \cdot h_{\delta_1,\delta_2,\delta_3}(\alpha)\textstyle{
   \Gamma(\frac{\alpha_2-\alpha_1}{2}-\delta_1)  \Gamma(\frac{\alpha_3-\alpha_1}{2}-\delta_1) \Gamma(\frac{\alpha_4-\alpha_1}{2}-\delta_1) }
  \\ \nonumber & \hskip 48pt
  \cdot 
  \textstyle{  \Gamma(\frac{\alpha_2-\alpha_4}{2}-\delta_2) \Gamma(\frac{\alpha_3-\alpha_4}{2}-\delta_2) }  \Gamma(\frac{\alpha_2-\alpha_3}{2}-\delta_3)  \;  d\alpha,
 \end{align}
 where $h_{\delta_1,\delta_2,\delta_3}$ is a polynomial of degree at most 
$2\delta_1+\delta_2$.
 
 We have
 
 \begin{prop}\label{prop:tripleresidue} Let $r_1,r_2,r_3$ be positive integers, and $0<\varepsilon<1$. Suppose $a_1,a_2,a_3$ satisfy the hypotheses of Theorem \ref{th:pTRbound}.  If $0\le \delta_j\le r_j-1$ for $1\le j\le3$, then
\begin{equation}\label{eq:Prop-s1s2s3-Residue}
 \left|   p_{T,R}^{123,(\delta_1,\delta_2,\delta_3)}(y) \right| \ll y_1^{\frac32+a_1}y_2^{2+a_2}y_3^{\frac32+a_3} \, T^{\varepsilon+4R+6-r_1-r_2-r_3}.
\end{equation}\end{prop}\begin{proof} To obtain the desired bound on  $p_{T,R}^{123,(\delta_1,\delta_2,\delta_3)}(y)$, we shift the lines of integration in  the $ \alpha_j$'s to $\re{(\alpha)}=\kappa$, where $\kappa=(\kappa_1,\kappa_2,\kappa_3)$ is such that the resulting exponents of $y_1,y_2,$ and $y_3$ have real parts as indicated in the  proposition.  We do so by choosing 
\begin{equation}\label{eq:kappatriple}\kappa=(a_1-2\delta_1,-a_3+2\delta_3,-a_2+a_3+2\delta_2-2\delta_3)\quad\hbox{(and $\kappa_4=-\kappa_1-\kappa_2-\kappa_3$)}.\end{equation}Then
\begin{align}\label{eq:pTr_123deltashift}&p_{T,R}^{123,(\delta_1,\delta_2,\delta_3)}(y)  =   
 \iiint\limits_{\re(\alpha)=\kappa} e^{\frac{\alpha_1^2+\cdots+\alpha_4^2}{2T^2}} y_1^{\frac32+\alpha_1+2\delta_1}y_2^{2+\alpha_1+\alpha_4+2\delta_2}y_3^{\frac32-\alpha_2+2\delta_3}  \mathcal F_R(\alpha)\,
    h_{\delta_1,\delta_2,\delta_3}(\alpha)
  \\  & \nn
  \cdot
  \frac{\Gamma(\frac{2+R+\alpha_1-\alpha_2}{4})
        \Gamma(\frac{2+R+\alpha_2-\alpha_1}{4})\Gamma(\frac{\alpha_2-\alpha_1}{2}-\delta_1)
      }{\Gamma(\frac{\alpha_1-\alpha_2}{2})\Gamma(\frac{\alpha_2-\alpha_1}{2})}
  \cdot
  \frac{\Gamma(\frac{2+R+\alpha_2-\alpha_3}{4})
        \Gamma(\frac{2+R+\alpha_3-\alpha_2}{4}) \Gamma(\frac{\alpha_2-\alpha_3}{2}-\delta_3) 
      }{\Gamma(\frac{\alpha_2-\alpha_3}{2})
        \Gamma(\frac{\alpha_3-\alpha_2}{2})}
  \\ \nonumber & 
  \cdot
  \frac{
      }{}
  \frac{\Gamma(\frac{2+R+\alpha_1-\alpha_3}{4})
        \Gamma(\frac{2+R+\alpha_3-\alpha_1}{4})\Gamma(\frac{\alpha_3-\alpha_1}{2}-\delta_1)
      }{\Gamma(\frac{\alpha_1-\alpha_3}{2})\Gamma(\frac{\alpha_3-\alpha_1}{2})}
  \cdot
  \frac{\Gamma(\frac{2+R+\alpha_2-\alpha_4}{4})
        \Gamma(\frac{2+R+\alpha_4-\alpha_2}{4})  \Gamma(\frac{\alpha_2-\alpha_4}{2}-\delta_2) 
      }{\Gamma(\frac{\alpha_2-\alpha_4}{2})
        \Gamma(\frac{\alpha_4-\alpha_2}{2})}
  \\ \nonumber & 
  \cdot
  \frac{\Gamma(\frac{2+R+\alpha_1-\alpha_4}{4})
        \Gamma(\frac{2+R+\alpha_4-\alpha_1}{4})\Gamma(\frac{\alpha_4-\alpha_1}{2}-\delta_1)
      }{\Gamma(\frac{\alpha_1-\alpha_4}{2})\Gamma(\frac{\alpha_4-\alpha_1}{2})}
  \cdot 
  \frac{
        \Gamma(\frac{2+R+\alpha_3-\alpha_4}{4})\Gamma(\frac{2+R+\alpha_4-\alpha_3}{4})\Gamma(\frac{\alpha_3-\alpha_4}{2}-\delta_2)
      }{
        \Gamma(\frac{\alpha_3-\alpha_4}{2})\Gamma(\frac{\alpha_4-\alpha_3}{2})} \; d\alpha.\nn
\end{align}
Notice that, in this case, {\it there are no poles crossed in moving the lines of integration}, and therefore no residue terms to consider.  This is for reasons encountered in prior situations:    if $R$ is large enough, then none of the terms
$$\textstyle  \Gamma(\frac{2+R+\alpha_k-\alpha_j}{4})\Gamma(\frac{2+R+\alpha_j-\alpha_k}{4})$$
in \eqref{eq:pTr_123deltashift} will give rise to any poles; moreover, any pole in the numerator of a factor $$  \frac{ \Gamma(\frac{\alpha_k-\alpha_j}{2}-\delta_n)
      }{\Gamma(\frac{\alpha_j-\alpha_k}{2})\Gamma(\frac{\alpha_k-\alpha_j}{2})}$$will be canceled by a pole in the denominator.
      
      So we need only bound the right-hand side of \eqref{eq:pTr_123deltashift}, and we do so in the usual way.  We get
      
      \begin{align} \label{eq:pTr_123Shifted}&
  \left| p_{T,R}^{123,(\delta_1,\delta_2,\delta_3)}(y)   \right|  \ll  
  y_1^{\frac32+a_1}y_2^{2+a_2}y_3^{\frac32+a_3}\,
T^{\varepsilon+R+2\delta_1+\delta_2}  \\&
 \cdot \iiint\limits_{\substack{ \abs{\tau_1},\abs{\tau_2},\abs{\tau_3} \leq T^{1+\varepsilon} }} 
  \hskip -8pt
  \big(1+\lvert\tau_1-\tau_2\rvert\big)^{\frac{R+1+\kappa_2-\kappa_1}{2}-\delta_1}\nn
  \big(1+\lvert\tau_1-\tau_3\rvert\big)^{\frac{R+1+\kappa_3-\kappa_1}{2} -\delta_1}
  \big(1+\lvert\tau_1-\tau_4\rvert\big)^{\frac{R+1+\kappa_4-\kappa_1}{2} -\delta_1}
  \\ \nonumber & 
  \cdot
  \big( 1+ \tau_2-\tau_3  \big)^{ \frac{R+1+\kappa_2-\kappa_3}{2}-\delta_3}
  \big( 1+\lvert\tau_2-\tau_4 \rvert\big)^{ \frac{R+1+\kappa_2-\kappa_4}{2}-\delta_2}
  \big( 1+|\tau_3-\tau_4 |\big)^{ \frac{R+1+\kappa_3-\kappa_4}{2}-\delta_2} 
\; d\tau\\\nn&
 \ll  
  y_1^{\frac32+a_1}y_2^{2+a_2}y_3^{\frac32+a_3}\,
 T^{\varepsilon+2\delta_1+\delta_2}   \iiint\limits_{\substack{ \abs{\tau_1},\abs{\tau_2},\abs{\tau_3} \leq T^{1+\varepsilon} }} 
  \hskip -8pt
 T^{ 3R+3 +\frac{-3\kappa_1+3\kappa_2+\kappa_3-\kappa_4}{2}-3\delta_1-2\delta_2-\delta_3} 
\; d\tau
\\\nn&
 \ll  
  y_1^{\frac32+a_1}y_2^{2+a_2}y_3^{\frac32+a_3}  \,T^{\varepsilon +4R+6 -a_1-a_2-a_3+\delta_1+\delta_2+\delta_3}.\end{align}
  But we're assuming $1\le \delta_j\le r_j-1$ and $-a_j\le 1-2r_j$, whence $\delta_j-a_j\le r_j$, for $1\le j\le 3$.  It follows that\begin{align*} \left| p_{T,R}^{123,(\delta_1,\delta_2,\delta_3)}(y)   \right|  \ll  &
  y_1^{\frac32+a_1}y_2^{2+a_2}y_3^{\frac32+a_3}
 \, T^{\varepsilon+4R+6-r_1-r_2-r_3}, \end{align*}which proves our proposition.\end{proof}

\section{\bf Bounding the contribution from the continuous spectrum}

Let $2\leq n\leq 4$.  Assuming that $\alpha=(\alpha_1,\ldots,\alpha_n)\in \C^n$, we define
\begin{equation}\label{eq:testfunctionsharp}
p_{T,R}^{\sharp,(n)}(\alpha) \; := \;
    \begin{cases}
      e^{\frac{\alpha_1^2+\cdots+\alpha_n^2}{2T^2}}
      \displaystyle{\prod_{1\leq\, j \ne k\, \leq n}}
      \Gamma\left(\textstyle{\frac{2+R+\alpha_j - \alpha_k}{4}}\right)
      & \mbox{ if }n=2,3,\\
      e^{\frac{\alpha_1^2+\cdots+\alpha_n^2}{2T^2}} \mathcal{F}_R(\alpha)
      \displaystyle{\prod_{1\leq\, j \ne k\, \leq n}}
      \Gamma\left(\textstyle{\frac{2+R+\alpha_j - \alpha_k}{4}}\right)
      & \mbox{ if }n=4,
      \end{cases}
\end{equation}
where
$\mathcal{F}_R(\alpha)$ is
as in (\ref{Fr}).  In the case of $n=4$ we will sometimes drop $n$ from 
the notation.  

Suppose that $\phi$ is a Maass cusp form for $\GL(n)$ with Langlands parameter $\alpha(\phi):=\alpha=(\alpha_1,\ldots,\alpha_n)\in \C^n$.  Then we define 
\begin{equation}\label{eq:hTRn}
 h_{T,R}^{(n)}(\phi) := \frac{ \big\lvert p_{T,R}^{\sharp,(n)}(\alpha) \big\rvert^2}{\prod\limits_{1\leq j\neq k\leq n} \Gamma\left(\frac{1+\alpha_j-\alpha_k}{2}\right)}.
\end{equation}

\begin{thm}[Weyl Law for $\GL(2)$ and $\GL(3)$]\label{th:WeylLaw}
Suppose that $n=2$ or $n=3$.  Let $\{\phi_1,\phi_2,\ldots\}$ be an orthogonal basis of Maass cusp forms for $\GL(n)$ ordered by eigenvalue.  Then there exists a constant $c_n$ such that
\begin{equation}\label{eq:WeylLaw}
 \sum_j \frac{h_{T,R}^{(n)}(\phi_j)}{L(1,\Ad{\phi_j})} \sim c_n\; T^{\frac{n(n-1)R}{2}+\frac{(n+2)(n-1)}{2}}.
\end{equation}
\end{thm}

\begin{rmrk}
In the case of $n=3$ this is Theorem~1.3 of \cite{GK2013}.  The case of 
$n=2$ is well known, but we remark that it can be proved by the same method as for $\GL(3)$ (see \cite{G2020}).  The point is that the main term 
for the left hand side of \eqref{eq:WeylLaw} is $\big\lvert p_{T,R}^{\sharp,(n)}\big\rvert^2$ which can be easily estimated using Stirling's estimate for the Gamma function.
\end{rmrk}

Suppose that $4=n_1+\cdots+n_r$ is a partition of $4$ and $\Phi=(\phi_1,\ldots,\phi_r)$ where, for $1\le j\le r$, $\phi_j$ is  a Maass cusp form for $\SL(n_j, \mathbb Z)$ if $n_j>1$, while $\phi_j$ is the constant function (properly normalized) if $n_j=1$.   Let $\mathcal{P}=\mathcal{P}_{n_1,\ldots,n_r}$.  Then we define
\begin{equation*}
\mathcal E_{\mathcal P,\Phi} :=  \int\limits_{\re(s_1)=0}\cdots \int\limits_{\re(s_{r-1})=0}  A_{E_{\mathcal P,\Phi}}(L, s)\cdot
 \overline{A_{E_{\mathcal P,\Phi}}(M, s)} \cdot\left| p_{T,R}^\#\big(\alpha_{_{\mathcal P, \Phi}}(s)\big)\right|^2 \; ds_1 \cdots ds_{r-1},
\end{equation*}
and
\begin{equation*}
\mathcal  E_{\mathcal P_{\rm\text{Min}}} :=
\int\limits_{\re(s_1)=0}\; \int\limits_{\re(s_2)=0}\; \int\limits_{\re(s_3)=0}   A_{E_{\mathcal P_{\rm \text{\rm Min}}}}(L, s)\, \overline{A_{E_{\mathcal P_{\rm \text{\rm Min}}}}(M, s)}\cdot
\left| p_{T,R}^\#\big(\alpha_{_{\mathcal P_{\rm Min}}}(s)\big)\right|^2 \; ds_1 \,ds_2 \,ds_3.
\end{equation*}

\begin{rmrk}
In the above integrals $\alpha_{_{\mathcal P, \Phi}}(s)$,  $\,\alpha_{_{\mathcal P_{\rm Min}}}(s)$ denote the Langlands parameters of the Eisenstein series $E_{\mathcal P,\Phi}(g,s)$,
$E_{\mathcal P_{\rm Min}}(g,s)$, respectively. Also, $A_{E_{\mathcal P, \Phi}}(L, s)$, $A_{E_{\mathcal P, \Phi}}(M, s)$ denote the $L^{th}$ and $M^{th}$ Fourier coefficient of  $E_{\mathcal P,\Phi}(g,s)$, and similarly for $E_{\mathcal P_{\rm Min}}(g,s)$.

\end{rmrk}
Thus, if we define
\[ \mathcal{E}_{\mathcal P} : = \sum_\Phi c_{L,M,\mathcal{P}}\cdot  \mathcal{E}_{\mathcal{P},\Phi}, \]
then the contribution to the Kuznetsov trace formula coming from the Eisenstein series (defined in Section~\ref{EisensteinTerm}) is given by
\[\mathcal{E}:= c_1 \mathcal{E}_{\mathcal P_{\rm\text{Min}}}+c_2\mathcal{E}_{\mathcal P_{2,1,1}}+c_3\mathcal{E}_{\mathcal P_{2,2}}+c_4\mathcal{E}_{\mathcal P_{3,1}},\]  
 for constants $c_1,c_2,c_3,c_4 > 0$.
\begin{thm}\label{thm:EisensteinBound}
Suppose the Ramanujan Conjecture (at $\infty$) for $\GL(n)$ with $n\leq 3$, i.e., the Langlands parameters are all purely imaginary.  Let $L=(\ell,1,1)$ and $M=(m,1,1)$.  Then 
\begin{align*} 
 \abs{\mathcal{E}_{\mathcal{P}_{\rm\text{Min}}}} & \ll_{\varepsilon} (\ell m)^\varepsilon\cdot T^{3+8R+\varepsilon}, \qquad
\abs{\mathcal{E}_{\mathcal P_{2,1,1}}}  \ll_{\varepsilon} (\ell m)^{{\frac{7}{64}}+\varepsilon}\cdot T^{2+8R+\varepsilon}, \\
\abs{\mathcal{E}_{\mathcal P_{2,2}}}  & \ll_{\varepsilon}  (\ell m)^{{\frac{7}{32}}+\varepsilon} \cdot  T^{5+8R+\varepsilon}, \qquad
 \abs{\mathcal{E}_{\mathcal P_{3,1}}}   \ll_{\varepsilon} (\ell m)^{{\frac25}+\varepsilon} \cdot  T^{6+8R+\varepsilon} ,  
\end{align*}
as $T\to\infty$ for any fixed $\varepsilon>0$.
\end{thm}

\begin{proof}

We shall require the following standard notation for completed L-functions. Let
$$\zeta^*(w) = \pi^{-\frac{w}{2}} \Gamma\left( \frac{w}{2} \right) \zeta(w) = \zeta^*(1-w), \qquad (w\in\mathbb C).$$
For a Maass cusp form $\phi$ on $\GL(2)$  with spectral parameter $\frac12 + v$, define the completed L-function $L^*(s, \phi)$  associated to $\phi$ by
$$L^*(w, \phi) :=  \pi^{-w}\Gamma\left( \frac{w+v}{2}  \right)  \Gamma\left( \frac{w-v}{2}  \right) L(w,\phi) = L^*(1-w,\phi), \qquad (w\in\mathbb C).$$
If $\phi_1, \phi_2$ are two Maass cusp forms on $\GL(2)$ with spectral parameters $\frac12+v, \frac12+v'$, respectively, then the completed L-function for the Rankin-Selberg convolution $L(w, \phi_1\times\phi_2)$ is given by
\begin{align} \label{RankinSelbergL}
L^*(w, \phi_1\times\phi_2) & =   \pi^{-2w}\Gamma\left( \frac{w+v+v'}{2} 
 \right)  \Gamma\left( \frac{w-v+v'}{2}  \right)\\
&
\hskip 90pt
\cdot
\Gamma\left( \frac{w+v-v'}{2}  \right)  \Gamma\left( \frac{w-v-v'}{2}  \right) L(w,\phi_1 \times \phi_2).
\nonumber
\end{align}
Finally, for a Maass cusp form $\phi$ on $\GL(3)$ with spectral parameter 
$\frac13 + (v, v')$ define the completed L-function 
$L^*(w, \phi)$  associated to $\phi$ by
$$L^*(w, \phi) :=  \pi^{-\frac{3w}{2}}\Gamma\left( \frac{w+v+2v'}{2}  \right)  \Gamma\left( \frac{w+v-v'}{2}  \right) \Gamma\left( \frac{w-2v-v'}{2}  \right) L(w,\phi) = L^*(1-w,\phi).$$
Recall the adjoint L-function of a Maass cusp form $\phi$  is defined by $L(w,\Ad \; \phi) :=\frac{ L(w,\phi\times \overline{\phi})}{\zeta(w)}.$
\vskip 4pt
The following table (see \cite{GMW2019}) lists, for each partition, the Maass form $\Phi$ (with its associated spectral parameters), the values of 
$s$-variables, Langlands parameters, and the Fourier-Whittaker coefficients of the $\SL(4, \mathbb Z)$ Eisenstein series $E_{\mathcal P,\Phi}$.Note that $E_{\mathcal P_{1,1,1,1}, \Phi} := E_{\mathcal P_{\rm Min}}.$

 {\extrarowsep=0.5mm
\begin{small}
\begin{table}[ht]
\begin{center}
\textbf{FOURIER COEFFICIENTS OF $\SL(4,\mathbb Z)$ LANGLANDS EISENSTEIN SERIES}
\vskip 8pt
\begin{tabular}{|c|c| C{5.45cm} |C{3.3cm} |} 
  \hline
 $\begin{matrix} \text{\rm Partition}\\ \text{Maass form $\Phi$}\\
 \text{Spectral pars.} \end{matrix}$& $\begin{matrix} $\text{\rm $s$ variables of $E_{\mathcal P, \Phi}$}$\\
 \alpha = \text{\rm Langlands pars.  }\end{matrix}$ & First Coefficient 
$A_{E_{\mathcal P, \Phi}}\big((1,1,1),s\big)$ (up to a constant factor)& $m^{th}$ Hecke eigenvalue $\lambda_{E_{\mathcal P, \Phi}}\big((m,1,1), s\big)$  \\ \hline
4 = 1+1+1+1  
  &
 $\begin{matrix} \\s=\frac14+(s_1,s_2,s_3,s_4)\\ \alpha_1 = 3s_1+2s_2+s_3\\ \alpha_2 = -s_1+2s_2+s_3\\ \alpha_3 = -s_1-2s_2+s_3\\ \alpha_4 
= -s_1-2s_2-3s_3
 \\
 \phantom{.}\end{matrix}$
%     & \hskip-80pt
%    $ \begin{matrix}& \Big(\zeta^*(1+4s_1) \, \zeta^*(1+4s_2) \\ 
%    & \hskip 50pt \cdot\zeta^*(1+4s_3)\, \zeta^*(1+4s_1+4s_2) \\
%    & \hskip 25pt \cdot \zeta^*(1+4s_2+4s_3)  \\
%     &
%     \hskip 85pt
%    \cdot \zeta^*(1+4s_1+4s_2+4s_3)\Big)^{-1}\end{matrix}$  
     &
     $\bigg(\prod\limits_{1\leq j<k\leq 4} \zeta^*(1+\alpha_j-\alpha_k)\bigg)^{-1}$
     &
     $\sum\limits_{c_1c_2c_3c_4=m} \hskip-15pt c_1^{\alpha_1} c_2^{\alpha_2} c_3^{\alpha_3} c_4^{\alpha_4}$
     
     \\ \hline
  $\begin{matrix} 4 = 2+1+1\\
  \phi \; \text{on $\GL(2)$}\\
  \frac12+v  \end{matrix}$
  &
  \hskip -2pt
  $\begin{matrix} \\
s= \big(1+s_1, -\frac12+s_2,s_3\big )\\
\alpha_1 = s_1+v\\
\alpha_2 = s_1-v\\
\alpha_3 = s_2\\
\alpha_4 = -2s_1-s_2\\
 \phantom{.}
\end{matrix}$
    &
    \hskip-12pt
    $\begin{matrix}
&
 \bigg(L(1, \Ad{\phi})^{\frac12}\;\left| \Gamma\left(\frac12+v   \right)\right|\\
 & \cdot \zeta^*(1+2s_1+2s_2) \\
&\hskip 20pt \cdot L^*(1+s_1+s_2, \phi)\\
&
\hskip 41pt\cdot L^*(1+3s_1+s_2, \phi)\bigg)^{-1}
% \bigg(L(1, \Ad\; \phi)^{\frac12}\;\left| \Gamma\left(\frac12+v   \right)\right|\\
% & \cdot \zeta^*(1+2s_1+2s_2)\\
%&\hskip 20pt \cdot L^*(1+s_1-s_2, \phi)\\
%&
%\hskip 41pt\cdot L^*(1+3s_1+s_2, \phi)\bigg)^{-1}
\end{matrix}$    
&
\hskip-10pt
$\begin{matrix}
\sum\limits_{c_1c_2c_3 =m} \lambda_\phi(c_1) c_1^{s_1}\\
\phantom{xxxxxx}\cdot\, c_2^{s_2}c_3^{-2s_1-s_2}
\end{matrix}$

    \\ \hline
 $\begin{matrix} \phantom{.}\\ 4= 2+2\\
 \Phi = (\phi_1, \, \phi_2)\\   \text{on}\\
  \GL(2)\times \GL(2)\\
 \frac12+v, \;\; \frac12+v' \\ 
 \phantom{.}  \end{matrix}$
  &
   $\begin{matrix}
s=\left(1+s_1,\, -1-s_1  \right)\\
\alpha_1 = s_1+v\\
\alpha_2 = s_1-v\\
\alpha_3 = -s_1+v'\\                                            
\alpha_4 = -s_1-v'
\end{matrix}$
   &
   \hskip-7pt
   $\begin{matrix}
   \bigg( L(1,\Ad\; \phi_1)^{\frac12} \, L(1,\Ad\; \phi_2)^{\frac12}\\
\cdot  \left| \Gamma\left(\frac12+v  \right)\; \Gamma\left(\frac12+v' \right)\right|
   \\
 \hskip 40pt  \cdot   L^*\big(1+2s_1,\phi_1\times\phi_2\big)\bigg)^{-1}
   \end{matrix}$
   &
   \hskip-20pt
   $\begin{matrix}
   \sum\limits_{c_1c_2=m}\hskip-5pt \lambda_{\phi_1}(c_1)\\
  \hskip 40pt \cdot \lambda_{\phi_2}(c_2) \left(\frac{c_1}{c_2}\right)^{s_1} 
   \end{matrix}$  
    \\ \hline
 $\begin{matrix} 4 = 3+1 \\
 \phi \;\text{on $\GL(3)$} \\
 \frac13 +(v,v') \end{matrix}$
 &
 $\begin{matrix} \phantom{.} \\
  s=\left(\frac12+s_1,-\frac32-3s_1\right)\\
 \alpha_1 = s_1+2v+v'\\
 \alpha_2 = s_1-v+v'\\
 \alpha_3 = s_1-v-2v'\\
 \alpha_4 = -3s_1\\
 \phantom{.}
 \end{matrix}$
  &
  $\begin{matrix}
  \hskip-25pt
  \bigg(L(1, \Ad \; \phi)^{\frac12}\cdot \left|\Gamma\left( \frac{1+3v}{2}  \right)\right| \\
  \cdot \left|\Gamma\left( \frac{1+3v'}{2}  \right) \Gamma\left( \frac{1+3v+3v'}{2}  \right)\right|
  \\
  \hskip 55pt
  \cdot L^*(1+4s_1, \phi)   \bigg)^{-1}
  \end{matrix}$
  &
  $\sum\limits_{c_1c_2 =m}\hskip-8pt \lambda_\phi(c_1, 1)\, c_1^{s_1} c_2^{-3s_1}$
   \\ \hline  
\end{tabular}

\end{center}
\end{table}
\end{small}
}

%\pagebreak

\begin{rmrk}
The formulas given here for the first coefficient are valid when the form 
$\phi$ or $\phi_j$ is a Maass cusp form.  \end{rmrk}

Following the above table of Fourier coefficients of $\SL(4,\mathbb Z)$ Langlands Eisenstein series, we now list the integrals arising in the contribution of the continuous spectrum  decomposition of the inner product of two Poincar\'e series given in Proposition \ref{EisTermKuz}.  For the rest of the proof, and for each partition of $4$, we will give the Langlands parameter $\alpha_{\mathcal P,\Phi}(s) := (\alpha_1,\alpha_2,\alpha_3,\alpha_4)$ for $(\mathcal P,\Phi)$ and then use Theorem~\ref{th:WeylLaw} to obtain the result.

In each case below we will use the fact that for any $1\leq j, k \leq 4$ and $\re(\alpha_j)=\re(\alpha_k)=0$,
\begin{equation}\label{eq:EisStirlingBound}
 \frac{\left\lvert \Gamma\left(\frac{2+R+\alpha_j-\alpha_k}{4}\right)\right\rvert^4 }{\left\lvert \Gamma\left(\frac{1+\alpha_j-\alpha_k}{2}\right)\right\rvert^2 } \; \sim \; c\; \big(1+\abs{\alpha_j-\alpha_k}\big)^R 
\end{equation}
for some constant $c$.  This follows trivially from Stirling's estimate.

We also will use the bound of Luo-Rudnick-Sarnak (see \cite{LRS1999}) for 
the $m^{th}$-Fourier coefficient of a  $\GL(n)$ ($n\geq2$) Maass cusp form $\phi$:
 \begin{equation} \label{LRSBound}
  \lambda_\phi(m,1,\ldots,1) \ll m^{\frac12-\frac{1}{n^2+1}+\varepsilon}. 
  \end{equation}
Note that this is proved in \cite{Goldfeld2015} as well.  In the special case of $\GL(2)$, Kim and Sarnak prove in the appendix of \cite{kim96} that
\begin{align}\label{KimBound}
 \lambda_\phi(m) \ll m^{\frac{7}{64} +\varepsilon}. 
\end{align}

\vskip 8pt\noindent
$\underline{\text{The   integral $\mathcal E_{\mathcal P_{\text{\rm Min}}}$:}}$
\vskip 5pt\noindent
The Langlands parameters  for $E_{\mathcal P_{Min}}(s)$ with $s=(s_1,s_2,s_3)$ these are given by: 
$$\alpha_{P_{\text{\rm Min}}}(s) = (\alpha_1,\alpha_2,\alpha_3,\alpha_4),$$
$$ \alpha_1 = 3s_1+2s_2+s_3,\;
 \alpha_2 = -s_1+2s_2+s_3,\;
 \alpha_3 = -s_1-2s_2+s_3,\;
 \alpha_4 = -s_1-2s_2-3s_3.$$

Therefore, using the above table   and the fact that  for any $\varepsilon>0$
 \begin{equation} \label{DivSumBound}
  \lambda_{E_{\mathcal{P}_{\rm \text{Min}}}}\big((m,1,1), s\big)= \sum\limits_{c_1c_2c_3c_4=m} \hskip-5pt c_1^{\alpha_1} c_2^{\alpha_2} c_3^{\alpha_3} c_4^{\alpha_4} \ll m^\varepsilon 
  \end{equation}
whenever $\re(\alpha_j)=0$ ($j=1,2,3,4$), we see that

\begin{align*}\mathcal  E_{\mathcal P_{\rm\text{Min}}} & =
\int\limits_{\re(s_1)=0}\; \int\limits_{\re(s_2)=0}\; \int\limits_{\re(s_3)=0}   A_{E_{\mathcal P_{\rm \text{\rm Min}}}}(L, s)\, \overline{A_{E_{\mathcal P_{\rm \text{\rm Min}}}}(M, s)}\cdot
\left| p_{T,R}^\#\big(\alpha_{\mathcal P_{\rm Min}}(s)\big)\right|^2 \; ds_1 ds_2 ds_3\\
 & \ll \int\limits_{\re(s_1)=0}\; \int\limits_{\re(s_2)=0}\; \int\limits_{\re(s_3)=0}
\frac{e^{\frac{\alpha_1^2+\alpha_2^2+\alpha_3^2+\alpha_4^2}{T^2}}
\cdot (\ell m)^\varepsilon }{ \prod\limits_{1\leq j<k\leq 4} \big|\zeta(1+\alpha_j-\alpha_k)\big|^2 }  \cdot 
 \prod_{1\leq j<k\leq 4} \frac{\left|\Gamma\left(\frac{2+R+\alpha_j-\alpha_k}{4}\right)\right|^4}{\left\lvert\Gamma\left(\frac{1+\alpha_j-\alpha_k}{2}\right)\right\rvert^2}
\\
& \hskip -10pt
 \cdot 
 \Big(\big(1+|\alpha_1+\alpha_2 -\alpha_3-\alpha_4|\big) 
 \big(1+|\alpha_1+\alpha_3 -\alpha_2-\alpha_4|\big)
 \big(1+|\alpha_1+\alpha_4 -\alpha_2-\alpha_3|\big)\Big)^{\frac{2R}{3}} \;  |ds_1 \,ds_2\, ds_3|.
\end{align*}

 If we make the change of variables
 $$ \alpha_1 = 3s_1+2s_2+s_3,\;
 \alpha_2 = -s_1+2s_2+s_3,\;
 \alpha_3 = -s_1-2s_2+s_3,\;
 \alpha_4 = -s_1-2s_2-3s_3$$
 in the above integral, it follows from the Jacobian transformation, that
 $ds_1\,ds_2\, ds_3$ maps to $\frac{d\alpha_1\,d\alpha_2\, d\alpha_3}{32}$.
 
Now, we have the Vinogradov (see \cite{Vinogradov1958}) bound
 \begin{equation} \label{ZetaInverseBound}
  \frac{1}{\lvert \zeta(1+it)\rvert} \ll \big(1+\abs{t}\big)^\varepsilon,\qquad (t\in\R),
  \end{equation}
which together with the above coordinate change imply that
\begin{align*}
& E_{\mathcal P_{\rm\text{Min}}} \ll  (\ell m)^\varepsilon \iiint\limits_{\substack{\re(\alpha)=0\\ \abs{\alpha}\leq T}}
\prod_{1\leq j<k\leq 4} \big(1+\abs{\alpha_j-\alpha_k}\big)^{R+\varepsilon} \\
&  \cdot 
 \Big(\big(1+|\alpha_1+\alpha_2 -\alpha_3-\alpha_4|\big) 
 \big(1+|\alpha_1+\alpha_3 -\alpha_2-\alpha_4|\big)
 \big(1+|\alpha_1+\alpha_4 -\alpha_2-\alpha_3|\big)\Big)^{\frac{2R}{3}} \, |d\alpha_1\, d\alpha_2\, d\alpha_3|.
\end{align*}

Next, make the change of variables
      $\alpha_1\to \alpha_1 T, \quad \alpha_2\to \alpha_2 T, \quad \alpha_3\to \alpha_3 T.$
      It easily follows that
     $$ E_{\mathcal P_{\rm\text{Min}}} \ll (\ell m)^\varepsilon\, T^{8R+3+\varepsilon}.$$

\vskip 8pt\noindent
$\underline{\text{The integral $\mathcal E_{\mathcal P_{2,1,1},\Phi}$:}}$
\vskip 5pt\noindent
We take $\phi$ to be a $\GL(2)$ Maass cusp form with spectral parameter $\frac12+v$ where $v\in\C$ is pure imaginary. 
The Langlands parameters $\alpha_{_{\mathcal P_{2,1,1},\Phi}}(s)$  for $E_{\mathcal P_{2,1,1},\Phi}(s)$ with $s=(s_1,s_2,s_3)$ are given by:
$ \alpha = (\alpha_1,\alpha_2,\alpha_3,\alpha_4)$  where
$$\alpha_1 = s_1+v,\quad
\alpha_2 = s_1-v,\quad
\alpha_3 = s_2,\quad
\alpha_4 = -2s_1-s_2.$$

It follows that
\begin{align*}
\mathcal E_{\mathcal P_{2,1,1},\Phi} & :=  \int\limits_{\re(s_1)=0}\; 
\int\limits_{\re(s_2)=0}  A_{E_{\mathcal P_{2,1,1},\Phi}}(L, s)\cdot
 \overline{A_{E_{\mathcal P_{2,1,1},\Phi}}(M, s)} \cdot\left| p_{T,R}^\#(\alpha)\right|^2 \; ds_1  ds_{2}\\
 &
 \\
 & = \int\limits_{\re(s_1)=0}\; \int\limits_{\re(s_2)=0} e^{\frac{(s_1+v)^2 + (s_1-v )^2 + s_2^2+ (2s_1+s_2 )^2}{T^2}}
\cdot \lambda_{E_{\mathcal P_{2,1,1},\Phi}}(L, s)   \,\overline{\lambda_{E_{\mathcal P_{2,1,1},\Phi}}(M, \;s)    } \\
&\hskip50pt \cdot\Big( \big(1+|4s_1|\big)\big(1+|2s_1+2s_2+2v| \big ) \big ( 1+|2s_1+2s_2-2v| \big ) \Big)^\frac{2R}{3}\\&
\hskip 10pt
\cdot
\frac{     \;\Big | 
      \Gamma\left(\textstyle{\frac{2+R+s_1-s_2-v}{4}}\right) \Gamma\left(\textstyle{\frac{2+R+s_1-s_2+v}{4}}\right) 
       \Gamma\left(\textstyle{\frac{2+R+3s_1+s_2-v}{4}}\right)
      \Gamma\left(\textstyle{\frac{2+R+3s_1+s_2+v}{4}}\right)\Big |^4
      }
      {  
\big |L^*(1+s_1-s_2, \phi)
L^*(1+3s_1+s_2, \phi)\big |^2}   \\
&
\hskip 80pt  \cdot \frac{ \Big |\Gamma\left(\textstyle{\frac{2+R+2v}{4}}\right) \Gamma\left(\textstyle{\frac{2+R+2s_1+2s_2}{4}}\right)\Big |^4   }{ L(1,\Ad \; \phi)  \left| \Gamma\left(\frac{1+2v}{2}\right)\right|^2
  \big|\zeta^*(1+2s_1+2s_2)\big |^2   }  \; ds_1 \,ds_2 
  \end{align*}
  from which we obtain the bound
 \begin{align*}
\mathcal E_{\mathcal P_{2,1,1},\Phi} & \ll \, \frac{ e^{\frac{2v^2}{T^2}} 
 \left|\Gamma\left(\textstyle{\frac{2+R+2v}{4}}\right) \right|^4 T^{2R}}{ 
 L(1, \Ad{\phi})  \left| \Gamma\left(\frac{1+2v}{2}\right)\right|^2
   } \int\limits_{\re(s_1)=0}\; \int\limits_{\re(s_2)=0}
   \hskip-10pt
    e^{\frac{2s_1^2 + s_2^2+ (2s_1+s_2 )^2}{T^2}}
 \lambda_{E_{\mathcal P_{2,1,1},\Phi}}(L, s)   \,\overline{\lambda_{E_{\mathcal P_{2,1,1},\Phi}}(M, \;s)    }\\
&
\hskip -50pt
\cdot
\frac{     \;\Big | 
      \Gamma\left(\textstyle{\frac{2+R+s_1-s_2-v}{4}}\right) \Gamma\left(\textstyle{\frac{2+R+s_1-s_2+v}{4}}\right) 
      \Gamma\left(\textstyle{\frac{2+R+3s_1+s_2-v}{4}}\right) 
      \Gamma\left(\textstyle{\frac{2+R+3s_1+s_2+v}{4}}\right)\Big |^4
      }
      {  
\big | \Gamma\left(\textstyle{\frac{1+s_1-s_2-v}{2}}\right) \Gamma\left(\textstyle{\frac{1+s_1-s_2+v}{2}}\right) 
      \Gamma\left(\textstyle{\frac{1+3s_1+s_2-v}{2}}\right) 
      \Gamma\left(\textstyle{\frac{1+3s_1+s_2+v}{2}}\right)L(1+s_1-s_2, \phi)
L(1+3s_1+s_2, \phi)\big |^2}\\
&
\hskip 115pt  \cdot \frac{ \Big |\Gamma\left(\textstyle{\frac{2+R+2s_1+2s_2}{4}}\right)\Big |^4   }{  \left|\Gamma\left( \frac{1+2s_2-2s_1}{2}  \right)\zeta(1+2s_1+2s_2)\right |^2   }  \; |ds_1 ds_2|. \end{align*}
Here 
$$\lambda_{E_{\mathcal P_{2,1,1},\Phi}}\Big((m,1,1), \;s\Big)  = \sum\limits_{c_1c_2c_3 =m} \lambda_\phi(c_1)\cdot c_1^{s_1} c_2^{s_2}c_3^{-2s_1-s_2} \; \ll \; m^{\frac{7}{64}+\varepsilon},
$$
by the bound  for $\lambda_\phi(c)$ given in \eqref{KimBound}.
It follows form \cite{HL1994}, \cite{HR1995}, \cite{GLS2004} that  $$L(1+it, \phi) \gg (1+|v|+|t|)^{-\varepsilon}.$$
It then follows from the above bound, Stirling's estimate for the Gamma function,   (\ref{ZetaInverseBound}), 
 that  
$$\mathcal E_{\mathcal P_{2,1,1},\Phi}\; \ll \;(\ell m)^{\frac{7}{64}+\varepsilon} \cdot T^{2+7R+\varepsilon} \cdot \frac{e^{\frac{2v^2}{T^2}}  \left|\Gamma\left(\textstyle{\frac{2+R+2v}{4}}\right) \right|^4 }{L(1,\Ad{\phi_j})  \left| \Gamma\left(\frac{1+2v}{2}\right)\right|^2} \; = \; (\ell m)^{\frac{7}{64}+\varepsilon} \cdot T^{2+8R+\varepsilon} \cdot   \frac{h_{T,R}^{(2)}(\phi_j)}{L(1,\Ad{\phi_j})}.$$  

   To bound $\mathcal E_{\mathcal P_{2,1,1}}$ we simply sum $\mathcal E_{\mathcal P_{2,1,1},\Phi}$ over all Maass cusp forms $\Phi$ for $\SL(2, \Z)$ using the Weyl law for $\GL(2)$ given in Theorem \ref{th:WeylLaw}. The 
stated result follows. 

\vskip 8pt\noindent
$\underline{\text{The integral $\mathcal E_{\mathcal P_{2,2},\Phi}$:}}$
\vskip 5pt\noindent
Here, we take  $\Phi = (\phi_1, \phi_2)$ to be Maass cusp forms for $\GL(2)$  with spectral parameters $\frac12+v, \; \frac12+v',$ respectively.
The Langlands parameters $\alpha_{_{\mathcal P_{2,2},\Phi }}$ for $E_{\mathcal P_{2,2},\Phi}(s)$ with $s=(s_1,s_2)$ are given by
$$\alpha_1 = s_1+v,\quad
\alpha_2 = s_1-v,\quad
\alpha_3 = -s_1+v',\quad                                         
\alpha_4 = -s_1-v'.$$

It follows that\begin{align*}
  p_{T,R}^{\sharp,(4)}(\alpha) & = e^{\frac{2s_1^2+v^2+{v'}^2}{T^2}} \cdot\Gamma\left( \frac{2+R+2v}{4} \right) \Gamma\left( \frac{2+R+2v'}{4} \right)\prod_{\delta,\delta'\in\{\pm1\}} \left\lvert \Gamma\Big(\frac{2+R+2s_1+\delta v+\delta' v'}{4}\Big)\right\rvert^2
  \\
  &
  \hskip 110pt
  \cdot \Big(\big(1+4|s_1|\big) \big(1 +2 |v+v'|\big) \big(1 +2 |v-v'|\big)\Big)^{\frac{R}{3}}.
  \end{align*}
Using this and the fact that
 \[ L^*(1+2s_1,\phi_1\times\phi_2) = \pi^{-2(1+2s_1)}L(1+2s_1,\phi_1\times\phi_2)\prod_{\delta,\delta'\in\{\pm1\}}\left\lvert\Gamma\Big(\frac{1+2s_1+\delta v +\delta'v'}{2}\Big)\right\rvert, \]
we see that
\begin{align*}
\mathcal E_{\mathcal P_{2,2},\Phi} & \ll \; T^{2R+\varepsilon}\frac{h_{T,R}^{(2)}(\phi_1)}{L(1,\Ad{\phi_1})} \frac{h_{T,R}^{(2)}(\phi_2)}{L(1,\Ad{\phi_2})} \int\limits_{\re(s_1)=0}  \frac{ \left| \lambda_{E_{\mathcal P_{2,2},\Phi}}(L, s)   \,\overline{\lambda_{E_{\mathcal P_{2,2},\Phi}}(M, 
\;s)} \right|}{|L(1+2s_1,\phi_1\times\phi_2)|^2} \\
 & \hskip 180pt 
 \cdot e^{\frac{4s_1^2}{T^2}} \prod_{\delta,\delta'\in\{\pm1\}} \frac{ \left\lvert\Gamma\left(\frac{2+R+2s_1+\delta v +\delta'v'}{4}\right)\right\rvert^4 }{ \left\lvert\Gamma\left(\frac{1+2s_1+\delta v +\delta'v'}{2}\right)\right\rvert^2 } \; |ds_1|
\end{align*}
where
$$\lambda_{E_{\mathcal P_{2,2},\Phi}}\Big((m,1,1), s\Big) =  \sum\limits_{c_1c_2=m} \lambda_{\phi_1}(c_1) \,\lambda_{\phi_2}(c_2)\cdot \left(\frac{c_1}{c_2}\right)^{s_1} \ll m^{\frac{7}{32}+\varepsilon}$$
by the bounds  for $\lambda_{\phi_1}(c), \lambda_{\phi_2}(c)$ given in \eqref{KimBound}.

It follows from  \cite{HR1995}, \cite{Moreno1985} that $$\frac{1}{L(1+2s_1,\phi_1\times\phi_2)} \ll  \big(1+\abs{s_1}+\abs{\nu}+\abs{\nu'}\big)^\varepsilon.$$
Note that the bound for the case that $|\text{\rm Im}(s_1)|$ is small involves the non-existence of Siegel zeros which is proved in \cite{HR1995} while the case when $|\text{\rm Im}(s_1)|$ is large was first proved in \cite{Moreno1985}. See also \cite{GLi2018}, \cite{HF2019}.

Applying these bounds, the Theorem~\ref{th:WeylLaw} bound, together with Stirling's bound to estimate the integral in $s_1$, we find
 \[ \lvert \mathcal E_{\mathcal P_{2,2}} \rvert \ll (\ell m)^{\frac{7}{32}+\varepsilon} \cdot T^{\varepsilon+6R+1} \sum_{(\phi_1,\phi_2)} \frac{h_{T,R}^{(2)}(\phi_1)}{L(1,\Ad{\phi_1})} \frac{h_{T,R}^{(2)}(\phi_2)}{L(1,\Ad{\phi_2})} \ll  (\ell m)^{\frac{7}{32}+\varepsilon} \cdot  T^{\varepsilon+6R+1}\, T^{R+2}\, T^{R+2} \]
as claimed.

\vskip 8pt\noindent
$\underline{\text{The integral $\mathcal E_{\mathcal P_{3,1},\Phi}$:}}$
\vskip 5pt\noindent
Let $\beta=(\beta_1,\beta_2,\beta_3)$  and $\frac13+(v, v')$ denote the 
Langlands and spectral parameters, respectively,  associated to a Maass cusp form $\phi$ on $\GL(3)$. Here $$\beta_1 = 2v+v', \quad \beta_2 = -v+v', \quad \beta_3 = -v-2v'.$$
The Langlands parameters $\alpha_{_{\mathcal P_{3,1},\Phi }}(s)$  for $E_{\mathcal P_{3,1},\Phi}(s)$ with $s=(s_1,-3s_1)$  are given by:
$ \alpha = (\alpha_1,\alpha_2,\alpha_3,\alpha_4)$  where
$$\alpha_1 = s_1+\beta_1,\quad
 \alpha_2 = s_1+\beta_2,\quad
 \alpha_3 = s_1+\beta_3,\quad
 \alpha_4 = -3s_1.$$
 
  Note that in this case, since
 \[ \sum_{j=1}^4 \alpha_j^2 = 9s_1^2 +\sum_{j=1}^3 (s_1+\beta_j)^2 = 12s_1^2+\sum_{j=1}^3 \beta_j^2\]
and $s_1,\beta_1,\beta_2,\beta_3$ are purely imaginary, we have
 \begin{align*}
  p_{T,R}^\sharp(\alpha) & = p_{T,R}^{\sharp,(3)}(\beta) \cdot e^{\frac{6s_1^2}{T^2}}\cdot \prod_{k=1}^3 \left\lvert \Gamma\left( \textstyle{\frac{2+R+4s_1-\beta_k}{4}} \right)\right\rvert^2
  \\
  & 
  \hskip-10pt
 \cdot \Big(\big(1+|\beta_1-\beta_2-\beta_3-4s_1|   \big)\big(1+|\beta_1+\beta_2-\beta_3+4s_1|   \big) \big(1+|\beta_1-\beta_2+\beta_3+4s_1|   \big)\Big)^{\frac{R}{3}}
  \end{align*}
where $p_{T,R}^{\sharp,(3)}(\beta)$ is defined by \eqref{eq:testfunctionsharp}.

This allows us to write
\begin{align*}
\mathcal E_{\mathcal P_{3,1},\Phi} & :=  \int\limits_{\re(s_1)=0}  A_{E_{\mathcal P_{3,1},\Phi}}(L, s)\cdot
 \overline{A_{E_{\mathcal P_{3,1},\Phi}}(M, s)} \cdot\left| p_{T,R}^\#(\alpha)\right|^2 \; ds_1 \\
 &
 \\
 & = \int\limits_{\re(s_1)=0}  \frac{\lambda_{E_{\mathcal P_{3,1},\Phi}}(L, s)\cdot
 \overline{\lambda_{E_{\mathcal P_{3,1},\Phi}}(M, s)}}{ L(1,\Ad{\phi}) \cdot |L^*(1+4s_1, \phi)|^2 } \cdot h_{T,R}^{(3)}(\beta) \cdot e^{\frac{12s_1}{T^2}} \prod_{k=1}^3 \Big\lvert \Gamma\Big(\frac{2+R+4s_1+\beta_k}{4}\Big)\Big\rvert^4
 \\
 &
 \hskip-27pt
 \cdot \Big(\big(1+|\beta_1-\beta_2-\beta_3-4s_1|   \big)\big(1+|\beta_1+\beta_2-\beta_3+4s_1|   \big) \big(1+|\beta_1-\beta_2+\beta_3+4s_1|   \big)\Big)^{\frac{2R}{3}}\, ds_1
\end{align*}
where, using \eqref{LRSBound},
\begin{equation} \label{GL3Bound}
\sum\limits_{c_1c_2 =m} {\lambda_\phi(c_1,1)}\cdot c_1^{s_1} c_2^{-3s_1} \ll m^{\frac{2}{5}+\varepsilon}.
\end{equation}

\vskip 5pt
In the above
\begin{align*}
L^*(1+4s_1, \phi) & =  \pi^{-\frac{3+12s_1}{2}} L(1+4s_1,\phi)\prod_{j=1}^3 \Gamma\Big(\frac{1+4s_1+\beta_j}{2}\Big).
\end{align*}
It follows from    
\cite{Moreno1985} and \cite{HR1995} that 
for every $\varepsilon > 0$
\begin{equation}\label{L1sBound}
L(1+4s_1,\phi) \; \gg_\varepsilon \;\frac{1}{\big(1+ |s_1|+|v|+|v'| \big)^{\varepsilon}}\end{equation}
where the implied constant in the $\gg_\varepsilon$ symbol is effective unless $\phi$ is a self-dual Maass cusp form that is not a symmetric square lift from $\GL(2).$
Note that the bound for the case that $|\text{\rm Im}(s_1)|$ is small involves the non-existence of Siegel zeros which is proved in \cite{HR1995} while the case when $|\text{\rm Im}(s_1)|$ is large was first proved in \cite{Moreno1985}. See also \cite{Sarnak2004}.

Let $\{\phi_1,\phi_2,\ldots\}$ be the Maass cusp forms for $\GL(3)$ ordered by eigenvalue, and set $\mathcal{L}_j:=L(1,\Ad{\phi_j})$.  It follows from (\ref{GL3Bound}), (\ref{L1sBound}) that
\begin{align*}
 \sum_j \mathcal{E}_{\mathcal{P}_{3,1},\phi_j} & =  \sum_j \frac{h_{T,R}^{\sharp,(3)}(\beta^{(j)})}{ \mathcal{L}_j } \int\limits_{\re(s_1)=0}  
\frac{\lambda_{E_{\mathcal P_{3,1},\phi_j}}(L, s)\cdot
 \overline{\lambda_{E_{\mathcal P_{3,1},\phi_j}}(M, s)}}{ |L(1+4s_1, \phi_j)|^2 } \; e^{\frac{12s_1^2}{T^2}} \prod_{k=1}^3 \frac{\left\lvert \Gamma\big(\frac{2+R+4s_1+\beta^{(j)}_k}{4}\big)\right\rvert^4}{\left|\Gamma\big(\frac{1+4s_1+\beta^{(j)}_k}{2}\big)\right|^2} \\
 &
 \hskip -40pt
\cdot \Big(\big(1+|\beta_1-\beta_2-\beta_3-4s_1|   \big)\big(1+|\beta_1+\beta_2-\beta_3+4s_1|   \big) \big(1+|\beta_1-\beta_2+\beta_3+4s_1|   \big)\Big)^{\frac{2R}{3}}\, ds_1
 \\
 & \ll \;(\ell m)^{\frac{2}{5}+\varepsilon}\cdot T^{2R}\; \sum_j \frac{h_{T,R}^{(3)}(\beta^{(j)})}{\mathcal{L}_j}\int\limits_{\re(s_1)=0} \frac{e^{\frac{12s_1^2}{T^2}}}{L(1+s_1,\phi_j)} \prod_{k=1}^3 \frac{\left\lvert \Gamma\big(\frac{2+R+4s_1+{\beta^{(j)}_k}}{4}\big)\right\rvert^4}{\left|\Gamma\big(\frac{1+4s_1+{\beta^{(j)}_k}}{2}\big)\right|^2}\; |ds_1|.
\end{align*}

Using Stirling's estimate for the Gamma functions here, it easy to see that 
 \[ \prod_{k=1}^3 \frac{\left\lvert \Gamma\big(\frac{2+R+4s_1+{\beta^{(j)}_k}}{4}\big)\right\rvert^4}{\left|\Gamma\big(\frac{1+4s_1+{\beta^{(j)}_k}}{2}\big)\right|^2} \ll \prod_{k=1}^3\big(1+|4s_1+2{\beta^{(j)}_k}|\big)^R. \]
Combining this with the previous bounds and Theorem~\ref{th:WeylLaw}, we find that 
\begin{equation}
\mathcal E_{\mathcal P_{3,1}}   \ll \sum_j \left\lvert \mathcal{E}_{\mathcal{P}_{3,1},\phi_j} \right\rvert \ll (\ell m)^{\frac25 +\varepsilon}\cdot T^{6+8R+\varepsilon}
\end{equation}
as claimed.
\end{proof}

\begin{rmrk}
As outlined in \cite{Blomer_2013}, it should be possible to obtain bounds 
for the more general case of $L=(\ell_1,\ell_2,\ell_3)$ and $M=(m_1,m_2,m_3)$ using the relations for the $\GL(4)$ Hecke operators.
\end{rmrk}

\vskip 25pt

\section*{\bf Acknowledgements}
 The authors would like to acknowledge several people whose contributions 
have enhanced this work.  First, we thank the referees for their careful reading of the paper and their many helpful suggestions.  Additionally, we thank Professor Stephen D. Miller of Rutgers University for a number of 
insightful and illuminating conversations.

\appendix
\section{\bf Integral bounds}\label{app:IntegralBound}

\begin{lem}\label{lemmaIntegralBound}
Suppose that $e,f$ are real numbers.  Then
 $$ \int\limits_{x=0}^T \frac{dx}{\big(1+T-x\big)^e \big(1+x\big)^{f}} \;  \ll \;
 (1+T)^{-\min\big\{e,\, f,\,e+f-1\big\}+\varepsilon} $$ for any $T\ge0$ and $\varepsilon>0$, and the implied constant does not depend on $T$.
\end{lem}
\begin{proof}
We consider the integrals
 \[ \int\limits_{x=0}^{T/2} (1+T-x)^{-e} (1+x)^{-f}\; dx \quad \mbox{and} \quad \int\limits_{x=T/2}^{T} (1+T-x)^{-e} (1+x)^{-f}\; dx \]
individually.  Since
 \[ \int\limits_{x=1}^{T/2}  x^{-f}\; dx \ll \begin{cases} T^{-f+1}+1 & 
\mbox{ if }f\neq 1, \\ \log{T}+1 &\mbox{ if }f=1,\end{cases} \] 
it follows in the case of $f\neq 1$ that 
 \[ \int\limits_{x=0}^{T/2} (1+T-x)^{-e} (1+x)^{-f}\; dx \ll (1+T)^{-e}\big(1+(1+T)^{-f+1} \big). \]
In like fashion, we find that if $e\neq 1$,
 \[ \int\limits_{x=T/2}^{T} (1+T-x)^{-e} (1+x)^{-f}\; dx \ll (1+T)^{-e-f+1}+(1+T)^{-f}. \]
Putting this together proves the result.  In the case that $e=1$ or $f=1$, the logarithm contributes $T^\varepsilon$ as claimed.
\end{proof}

\begin{lem}\label{lemmaIntegralBound2}
Assume $B_1\le B_2 \le \cdots \le B_k$.  Then for any $\varepsilon>0$
\begin{align*}
 \int\limits_{x=B_1}^{B_k}  \prod_{i=1}^k \big(1+|x-B_j | \big)^{-e_j}  \;\ll \; \big(1+B_k-B_1\big)^\varepsilon \;\sum_{j=1}^{k-1}  & \big(1+B_{j+1}-B_j\big)^{-\min\big\{e_j, e_{j+1},e_j+e_{j+1}-1\big\}  } \\
 & \quad \cdot \prod_{i\neq j,j+1}\big(1+|B_j^*(i)-B_i|\big)^{-e_i},
\end{align*}
where
 \[ B_j^*(i) := \begin{cases} B_j & \mbox{ if $i<j$ and $e_i>0$,} \\ B_{j+1} & \mbox{ if $i<j$ and $e_i<0$,} \\ B_{j+1} & \mbox{ if $i>j+1$ and $e_i>0$,} \\ B_j & \mbox{ if $i>j+1$ and $e_i<0$.} \end{cases}\]
\end{lem}
\begin{proof}
First, 
 \[ \int\limits_{x=B_1}^{B_k} = \sum_{j=1}^{k-1} \int\limits_{x=B_j}^{B_{j+1}}. \]
For every $j=1,2,\ldots,k-1$ we have
\begin{align*}
 \int\limits_{B_j}^{B_{j+1}} \prod_{i=1}^k \big(1+|x-B_i|\big)^{-e_i}dx 
&\ll  
   \int\limits_{B_j}^{B_{j+1}}  \big(1+x-B_j\big)^{-e_j} \big(1+B_{j+1}-x\big)^{-e_{j+1}} \\
   & \cdot \prod_{i=1}^{j-1} \big(1+{x-B_i}\big)^{-e_i} \prod_{\ell=j+2}^k \big(1+B_i-x\big)^{-e_\ell} \; dx.
\end{align*}
For each of the terms with $i< j$ and any $B_j\leq x\leq B_{j+1}$,
 \[ \big(1+x-B_i\big)^{-e_i} \ll \begin{cases} \big(1+B_j-B_i\big)^{-e_i} 
& \mbox{ if }e_i>0, \\ \big(1+B_{j+1}-B_i\big)^{-e_i} & \mbox{ otherwise.} \end{cases} \]
A similar bound holds for the terms with $\ell>j+1$.  So in order to complete the proof, we need the bound
 \[ \int\limits_{B_j}^{B_{j+1}} \big(1+x-B_j\big)^{-e_j} \big(1+B_{j+1}-x\big)^{-e_{j+1}}\; dx \ll \big(1+B_k-B_1\big)^\varepsilon \cdot \big(1+B_{j+1}-B_j\big)^{-\min\big\{e_j, e_{j+1},e_j+e_{j+1}-1\big\}}, \]
which follows from Lemma \ref{lemmaIntegralBound} by a simple change of variables.
\end{proof}

\begin{lem}\label{lemmaIntegralBound3}
Assume $B_1\leq B_2 \leq \cdots \leq B_k$.  Suppose that $1\leq j_{\mathrm{min}} < j_{\mathrm {max}} \leq k$.  Then for any $\varepsilon>0$, 
\begin{align*}
& \int\limits_{x=B_{j_{\mathrm min}}}^{B_{j_{\mathrm{max}}}}  \prod_{i=1}^k \big(1+|x-B_j | \big)^{-e_j} \;dx \ll \; \big(1+B_{j_{\mathrm{max}}}-B_{j_{\mathrm{ min}}}\big)^\varepsilon \\&\cdot\sum_{j=j_{\mathrm{min}}}^{j_{\mathrm{max}}-1}    \big(1+B_{j+1}-B_j\big)^{-\min\big\{e_j, e_{j+1},e_j+e_{j+1}-1\big\}  }
  \prod_{i\neq j,j+1}\big(1+|B_j^*(i)-B_i|\big)^{-e_i},
\end{align*}
where
 \[ B_j^*(i) := \begin{cases} B_j & \mbox{ if $i<j$ and $e_i>0$,} \\ B_{j+1} & \mbox{ if $i<j$ and $e_i<0$,} \\ B_{j+1} & \mbox{ if $i>j+1$ and $e_i>0$,} \\ B_j & \mbox{ if $i>j+1$ and $e_i<0$.} \end{cases}\]
\end{lem}
\begin{proof}
First, 
 \[ \int\limits_{x=B_{j_{\mathrm{min}}}}^{B_{j_{\mathrm{max}}}} = \sum_{j=j_{\mathrm{min}}}^{j_{\mathrm{max}}-1} \int\limits_{x=B_j}^{B_{j+1}}. \]
For every $j=1,2,\ldots,k-1$ we have
\begin{align*}
 \int\limits_{B_j}^{B_{j+1}} \prod_{i=1}^k \big(1+|x-B_i|\big)^{-e_i}dx 
&\ll  
   \int\limits_{B_j}^{B_{j+1}}  \big(1+x-B_j\big)^{-e_j} \big(1+B_{j+1}-x\big)^{-e_{j+1}} \\
   & \cdot \prod_{i=1}^{j-1} \big(1+{x-B_i}\big)^{-e_i} \cdot \prod_{\ell=j+2}^k \big(1+B_i-x\big)^{-e_\ell} \; dx.
\end{align*}
For each of the terms with $i< j$ and any $B_j\leq x\leq B_{j+1}$,
 \[ \big(1+x-B_i\big)^{-e_i} \ll \begin{cases} \big(1+B_j-B_i\big)^{-e_i} 
& \mbox{ if }e_i>0, \\ \big(1+B_{j+1}-B_i\big)^{-e_i} & \mbox{ otherwise.} \end{cases} \]
A similar bound holds for the terms with $\ell>j+1$.  So in order to complete the proof, we need the bound
 \[ \int\limits_{B_j}^{B_{j+1}} \big(1+x-B_j\big)^{-e_j} \big(1+B_{j+1}-x\big)^{-e_{j+1}}\; dx \ll \big(1+B_k-B_1\big)^\varepsilon \cdot \big(1+B_{j+1}-B_j\big)^{-\min\big\{e_j, e_{j+1},e_j+e_{j+1}-1\big\}}, \]
which follows from Lemma~\ref{lemmaIntegralBound} by a simple change of variables.
\end{proof}

\vskip 25pt

\section{\bf Kloosterman sums on $\GL(4)$\\ by Bingrong Huang}\label{app:bingrong}

\subsection{\bf Introduction} \label{sec: Introduction}

The classical Kloosterman sum is given by
\[
  S(m,n;c) = \underset{d\pmod{c}}{{\sum}^*} e\Big(\frac{md+n\bar{d}}{c}\Big),
\]
where $d\bar{d}\equiv1\pmod{c}$ and $e(x)=e^{2\pi i x}$,
which arises when one computes the Fourier expansion of the $\GL(2)$ Poincar\'{e} series.
Weil \cite{weil1948some} obtained the algebreo-geometric estimate
\[ |S(m,n;c)| \leq \gcd(m,n,c)^{1/2}  c^{1/2} \tau(c), \]
where $\tau(\cdot)$ is the divisor function.
%$S(m,n;c)$ arises when one computes the Fourier expansion of the $\GL(2)$ Poincar\'{e} series
%\[
%  P_m(z,s) = \sum_{\gamma\in \Gamma_\infty\backslash\Gamma} \Im(\gamma z)^s e(m \gamma z),
%\]
%where $m>0$, $z\in\mathbb{H}$ the upper half plane, $\Re(s)>1$, $\Gamma=\SL(2,\mathbb{Z})$, and $\Gamma_\infty=\{ \left(\begin{smallmatrix} 1 & n \\  0 & 1  \end{smallmatrix} \right):n\in\mathbb{Z}\}$, which plays
%a role in much work in number theory.
Bump, Friedberg and Goldfeld \cite{BFG1988} introduced Poincar\'{e} series for $\GL(n)$, $n\geq2$, and showed in the case $n=3$ that certain ``$\GL(3)$ Kloosterman sums'' arise in the Fourier expansion.
Friedberg \cite{friedberg1987poincare} and Stevens \cite{stevens1987poincare} extended this result to all $n$, studying $\GL(n)$ Poincar\'{e} series and their related $\GL(n)$ Kloosterman sums,
from the classical and adelic points of view respectively.
Friedberg, following the work of Larsen ($n=3$) \cite{BFG1988}, obtained upper bounds for $\GL(n)$ in certain cases.
Stevens \cite{stevens1987poincare} gave a nontrivial estimate for the $\GL(3)$ Kloosterman sum corresponding to the long element of the Weyl group.
By their results, we get nontrivial upper bounds for all $\GL(3)$ Kloosterman sums.
%In \cite{buttcane2013sums}, Buttcane worked out the dependence on the characters.

In this appendix, we consider all $\GL(4)$ Kloosterman sums.
We will write $\mathbb{Q}_p$ for the completion of $\mathbb{Q}$ at a place $p$ and write $\mathbb{A}$ for the adeles of $\mathbb{Q}$. Let $G=\GL(4)$.
Let $W$ be the Weyl group of $G$.
Let
$U=\left\{ \left( \begin{smallmatrix}
                        1 & * & * & * \\
                        0 & 1 & * & * \\
                        0 & 0 & 1 & * \\
                        0 & 0 & 0 & 1
                      \end{smallmatrix} \right) \right\}$
be the standard unipotent group and let
%$U=\left\{ \left( \begin{smallmatrix}
%                        1 & * & * & \cdots & * \\
%                        0 & 1 & * & \cdots & * \\
%                        0 & 0 & 1 & \cdots & * \\
%                        \vdots & \vdots & \vdots & & \vdots \\
%                        0 & 0 & 0 & \cdots & 1
%                      \end{smallmatrix} \right) \right\}$,
\begin{equation*}%\label{eqn: U_w}
  \text{ U}_w = (w^{-1}\cdot U\cdot w)\cap U, \quad
  \bar{U}_w = (w^{-1}\cdot{} ^t U \cdot w)\cap U, \quad w\in W.
%  \begin{split}
%       & U_w = (w^{-1}\cdot U\cdot w)\cap U, \\
%       & \bar{U}_w = (w^{-1}\cdot{} ^t U \cdot w)\cap U.
%  \end{split}
\end{equation*}
Let $c_1,\ldots,c_{3}$ be non-zero integers, and set
\[ c= \diag(1/c_{3},c_{3}/c_{2}, c_2/c_1, c_1). \]
Following Stevens \cite[\S2]{stevens1987poincare}, we define
\begin{equation*}%\label{eqn: C,X}
  C(cw)  := U(\mathbb{Q}_p)cwU(\mathbb{Q}_p)\cap G(\mathbb{Z}_p),  \quad
  X(cw)  := U(\mathbb{Z}_p)\bs C(cw)/\bar{U}_w(\mathbb{Z}_p).
%  \begin{split}
%    C(cw) & := U(\mathbb{Q}_p)cwU(\mathbb{Q}_p)\cap G(\mathbb{Z}_p),  \\
%    X(cw) & := U(\mathbb{Z}_p)\bs C(cw)/\bar{U}_w(\mathbb{Z}_p).
%  \end{split}
\end{equation*}
By the Bruhat decomposition we have natural maps
\begin{equation*}%\label{eqn: u,u'}
  u:  X(cw)\rightarrow U(\mathbb{Z}_p)\bs U(\mathbb{Q}_p),  \quad
  u':  X(cw)\rightarrow \bar{U}_w(\mathbb{Q}_p)/\bar{U}_w(\mathbb{Z}_p).
\end{equation*}
defined by the relation $x=u(x)cw u'(x)$ for $x\in X(cw)$.
% character
Let $\psi:U(\mathbb{A})/U(\mathbb{Q})\mapsto \mathbb{C}^*$ be a character 
of $U(\mathbb{A})$ which is trivial on $U(\mathbb{Q})$.
Every such character has the form $\psi=\psi_{N}$ for some $N=(n_1,n_2,n_3)\in\mathbb{Q}^{3}$ where $\psi_{N}$ is given by
\[
  \psi_{N}\left( \begin{smallmatrix}
                        1 & x_1 & * & * \\
                        0 & 1 & x_2 & * \\
                        0 & 0 & 1 & x_3 \\
                        0 & 0 & 0 & 1
                      \end{smallmatrix} \right)
  = \xi(n_1x_1+n_2x_2+n_3x_3)
\]
and $\xi:\mathbb{A}\rightarrow \mathbb{C}^*$ is the standard additive character.
We can write $\psi=\prod_p \psi_p$ where $\psi_p$ is a character of $U(\mathbb{Q}_p)$ which is trivial on $U(\mathbb{Z}_p)$.
The \emph{local Kloosterman sum} is defined by
\[
    Kl_p(\psi_p,\psi'_p;c,w) = \sum_{x\in X(cw)}\psi_p(u(x))\cdot \psi_p'(u'(x)).
\]
The \emph{global Kloosterman sum} is defined by
$
    Kl(\psi,\psi';c,w) = \prod_p Kl_p(\psi_p,\psi'_p,c,w).
$
Our main results for $Kl(\psi_{M},\psi_{N};c,w)$ are in the following table.

{\extrarowsep=0.5mm
%\begin{landscape}
\begin{small}\label{table-b11}
\begin{table}[h]
\begin{center}

\begin{tabular}{|c|C{4.2cm}| C{8cm} |} % p{7cm}<{\centering} % m{3cm}<{\centering} 垂直居中 加array宏包
  \hline
  % after \\: \hline or \cline{col1-col2} \cline{col3-col4} ...
  Weyl element & Compatibility conditions & Upper bounds of the Kloosterman sum \\ \hline
  $w_1=\left(\begin{smallmatrix}1&&&\\ &1&&\\ &&1&\\ &&&1\end{smallmatrix}\right)$  &
    \tabincell{c}{$m_1=n_1$,  $m_2=n_2$,\\ $m_3=n_3$; \\ $c_1=c_2=c_3=1$}
    &
    $\delta_{c_1,1}\delta_{c_2,1}\delta_{c_3,1}
    \delta_{m_1,n_1}\delta_{m_2,n_2}\delta_{m_3,n_3}$ \\ \hline
  $w_2=\left(\begin{smallmatrix}&&&-1\\ 1&&&\\ &1&&\\ &&1&\end{smallmatrix}\right)$ &
    \tabincell{c}{$n_1=\frac{c_1c_3m_2}{c_2^2}$, $n_2=\frac{c_2m_3}{c_1^2}$; \\
    $c_1|c_2|c_3$} &
    \tabincell{l}{
    $\bullet$ Friedberg \cite[Theorem C]{friedberg1987poincare}; \\
    $\bullet$ $\tau(c_1c_2c_3)^{\kappa_2}(m_1,c_3/c_2)^{1/4}(m_2,c_2/c_1)^{1/2}$\\
    \hskip 60pt $\cdot (m_3,n_3,c_1)^{3/4}(c_1c_2c_3)^{3/4}$} \\ \hline
  $w_3=\left(\begin{smallmatrix}&1&&\\ &&1&\\ &&&1\\ -1&&&\end{smallmatrix}\right)$  &
    \tabincell{c}{$n_3=\frac{c_1c_3m_2}{c_2^2}$, $n_2=\frac{c_2m_1}{c_3^2}$;\\
    $c_3|c_2|c_1$} &
    \tabincell{l}{
    $\bullet$ Similar to Friedberg \cite[Theorem C]{friedberg1987poincare}; \\
    $\bullet$ $\tau(c_1c_2c_3)^{\kappa_2}(m_3,c_1/c_2)^{1/4}(m_2,c_2/c_3)^{1/2}$\\
    \hskip 60pt $\cdot (m_1,n_1,c_3)^{3/4} (c_1c_2c_3)^{3/4}$} \\ \hline
  $w_4=\left(\begin{smallmatrix}&&1&\\ &&&1\\ 1&&&\\ &1&&\end{smallmatrix}\right)$  &
    \tabincell{c}{$n_1=\frac{m_3c_2}{c_1^2}$, $n_3=\frac{m_1c_2}{c_3^2}$;\\
    $c_1|c_2,\ c_3|c_2$}&
    $c_1c_2c_3$ \\ \hline
  $w_5=\left(\begin{smallmatrix}&&&-1\\ &1&&\\ &&1&\\ 1&&&\end{smallmatrix}\right)$  &
    $n_2=\frac{c_1c_3m_2}{c_2^2}$ &
    $c_1c_2c_3$ \\ \hline
  $w_6=\left(\begin{smallmatrix}&&1&\\ &&&1\\ &1&&\\ -1&&&\end{smallmatrix}\right)$  &
    \tabincell{c}{$n_3 = \frac{c_2 m_1}{c_3^2}$;\\ $c_3|c_2$} &
    $c_1c_2c_3$  \\ \hline
  $w_7=\left(\begin{smallmatrix}&&&-1\\ &&1&\\ 1&&&\\ &1&&\end{smallmatrix}\right)$  &
    \tabincell{c}{$n_1 =\frac{c_2 m_3}{c_1^2}$;\\ $c_1|c_2$} &
    $c_1c_2c_3$ \\ \hline  %\rule{0pt}{0.8cm}
  $w_8=\left(\begin{smallmatrix}&&&1\\ &&1&\\ &1&&\\ 1&&&\end{smallmatrix}\right)$  &
     &
    \tabincell{l}{$\bullet$ $\tau(c_1c_2c_3)^{\kappa_8}(m_1n_3,[c_1,c_2,c_3])^{1/2}$ \\ \ \ \ \
    $\cdot(m_2n_2,[c_1,c_2,c_3])^{1/2} (m_3n_1,[c_1,c_2,c_3])^{1/2}$\\ \ \ \ \
    $\cdot\min\left\{ [c_1,c_3]^{1/2}(c_1,c_3)c_2 \; , \; c_1c_3(c_1,c_3)c_2^{1/2} \right\}$ \\
    $\bullet$  $\tau(c_1c_2c_3)^{\kappa_8}(m_1n_3,[c_1,c_2,c_3])^{1/2}$ \\ \hskip 60pt
    $\cdot(m_2n_2,[c_1,c_2,c_3])^{1/2} $\\
    \hskip 60pt
    $\cdot (m_3n_1,[c_1,c_2,c_3])^{1/2} (c_1c_2c_3)^{9/10}$} \\ \hline
\end{tabular}
\vskip 8pt
\caption{Main results for $\GL(4)$ Kloosterman sums}
  \label{capt:kloos}
\end{center}
\end{table}
\end{small}
%\end{landscape}
}
%prove Conjecture 2 in Stevens \cite{stevens1987poincare} for the $GL(4)$ case ...?

It was shown in Friedberg \cite[\S1]{friedberg1987poincare} that the Kloosterman sums %$S(\psi_M,\psi_N;c,w)$
are non-zero only if $w\in W$ is of the form
$
  w = \left( \begin{smallmatrix}
          &   &   & I_{k_1} \\
          &   & I_{k_2} &  \\
          & \reflectbox{$\ddots$} &  &   \\
        I_{k_r} &  &  &
      \end{smallmatrix} \right),
$
where the $I_k$ are $k\times k$ identity matrices and $k_1+\cdots+k_r=n$ (may have some minus sign to make its determinant 1).

For the case $w=w_1$, we have
$Kl(\psi_{M},\psi_{N};c,w_1)= \delta_{c_1,1}\delta_{c_2,1}\delta_{c_3,1}
\delta_{m_1,n_1}\delta_{m_2,n_2}\delta_{m_3,n_3},$
where $\delta_{m,n}=1$ if $m=n$, and $\delta_{m,n}=0$ otherwise.

For the case $w=w_2$ or $w_3$, Friedberg \cite{friedberg1987poincare} gave some very nice bounds for
$\GL(n)$ Kloosterman sums attached to
$w = \left(\begin{smallmatrix}  & \pm 1 \\  I_{n-1} & \end{smallmatrix} 
\right)$.
For $n=3$, this is due to
Larsen, see \cite[Appendix]{BFG1988}.
Then Friedberg \cite[\S4]{friedberg1987poincare}
generalized it to all $n$.
In some applications, we may need to give a bound with power saving
in terms of all $c_1,c_2,c_3$.
One can modify Friedberg's proof to give such a bound in the case $n=4$. Note that the main situation is when $c_j=p^{ja},\ a\geq1,\ 1\leq j\leq 3$, in which case $(c_1c_2c_3)^{3/4}$ agrees with $c_j^{9/2j}$ in \cite{friedberg1987poincare}.
In the proof, we need Deligne's deep theorems from algebraic geometry \cite{Deligne1977}.
For the case $w=w_3$, one can use the involution operator $\iota: g \mapsto w_8 {}^t g^{-1} w_8$ to get the result.

By \cite[Theorem 0.3 (i)]{DR1998}, we have the ``trivial'' bound
\begin{equation}\label{eqn:trivialbound}
  Kl(\psi,\psi';c,w) \leq \# X(cw) \leq  c_1 c_2 c_3.
\end{equation}
We use this for $w=w_j$, with $4\leq j\leq 7$. In fact, this kind of bound holds for a general Kloosterman sum.

%%%%%%%%%%%%%%%%%%%%%%%%%%%%%%%%%%%%%%%%%%%%%%%%%%%%%%%%%%%%%%%%%%
%%                      Section                                 %%
%%%%%%%%%%%%%%%%%%%%%%%%%%%%%%%%%%%%%%%%%%%%%%%%%%%%%%%%%%%%%%%%%%
\subsection{\bf Stevens' approach}

In this section, we follow Stevens' approach \cite{stevens1987poincare} to bound the $\GL(4)$ long element Kloosterman sums.
For $w\in W$, we define $w(j)$, $j\in\{1,2,3,4\}$ by the formula
\[
  w\cdot e_j = \pm e_{w(j)},
\]
where $e_1,e_2,e_3,e_4$ is the standard basis of column vectors.
Let $\nu_1,\nu_2,\nu_3,\nu_1',\nu_2',\nu_3'\in\mathbb{Z}_p$
and define the characters $\psi,\psi'$ of $U(\mathbb{Q}_p)/U(\mathbb{Z}_p)$ by
\begin{equation}\label{eqn: psi}
  \begin{split}
     & \psi\left(\begin{smallmatrix} 1&u_1&*&*\\ &1&u_2&*\\ &&1&u_3\\ &&&1\end{smallmatrix}\right) = \xi(\nu_1u_1+\nu_2u_2+\nu_3u_3), \ \
     \psi'\left(\begin{smallmatrix} 1&u_1&*&*\\ &1&u_2&*\\ &&1&u_3\\ &&&1\end{smallmatrix}\right) = \xi(\nu_1'u_1+\nu_2'u_2+\nu_3'u_3).
  \end{split}
\end{equation}
Fix
\begin{equation}\label{eqn: c}
  c=\diag( p^{-t}, p^{t-r} ,  p^{r-s}, p^{s}).
%  c=\left(\begin{smallmatrix}
%         p^{-t} &  &  &  \\
%          & p^{t-r} &  &  \\
%          &  & p^{r-s} &  \\
%          &  &  & p^{s}
%       \end{smallmatrix}\right).
\end{equation}
We will use the same notation as in Stevens \cite[\S4]{stevens1987poincare}.
And we need Definition 4.9 and Theorem 4.10 in \cite{stevens1987poincare}.
Note that $n$ in \cite{stevens1987poincare} will be our $cw$.

\Big\}.

Our main result in this appendix is the following theorem.

\begin{thm}\label{thm: w_8}
  Let $Kl_p(\psi,\psi';c,w_8)$ be the local Kloosterman sum attached to the long element $w_8$.
  Let $\psi,\psi'$ be as in \eqref{eqn: psi}, $\ell=
  \max(r,s,t)$, $\varrho=\max(t,s)$, $\sigma=\min(t,s)$, and
  $$
    C_8 = 64   (|\nu_1\nu_3'|_p^{-1},p^\ell)^{1/2} (|\nu_2\nu_2'|_p^{-1},p^\ell)^{1/2} (|\nu_3\nu_1'|_p^{-1},p^\ell)^{1/2} (\varrho+1)(r+1)^2(\sigma+1)^2.
  $$
  Then
  \begin{equation}\label{eqn:Kl-w8}
    \begin{split}
      |Kl_p(\psi,\psi';c,w_8)|
      & \leq
        C_8 \min(p^{r+\sigma+\varrho/2},p^{\varrho+2\sigma+r/2}).
    \end{split}
  \end{equation}
  In particular, we have $|Kl_p(\psi,\psi';c,w_8)|\leq C_8 p^{9(t+r+s)/10}.$
\end{thm}

Suppose we are given $\alpha,\beta,\gamma\in\mathbb{Z}_p^\times$, and
nonnegative integers $a,b,c,d,e,f$
\footnote{Note that here we use $c$ in two meanings, one for a matrix, and another for a nonnegative integers.
However one can easily determine what does it mean in the context.}
satisfying
%\begin{equation}\label{eqn: a,b and alpha,beta}
\begin{gather}
  a\leq s, \quad d\leq s,\quad  e=s,  \quad f\leq r, \quad b+c\leq \max(t,f); \label{eqn: a,b,c,d,f}\\
  \left\{\begin{split}
       & x = -\gamma p^{r-f}\in \mathbb{Z}_p, \\
       & y=p^r(\beta p^{-s} - \gamma p^{-a-f}) \in \mathbb{Z}_p,  \\
       & z=p^t(\gamma p^{-f} - p^{-b-c})\in \mathbb{Z}_p;
  \end{split}\right. \label{eqn: y,z}\\
  \left\{\begin{split}
       & \lambda = p^r(\beta p^{-b-s} - \alpha\gamma p^{-d-f}) \in \mathbb{Z}_p^\times, \\
       & \mu = p^t(\gamma p^{-a-f} + \alpha p^{-d-c} - p^{-a-b-c} - \beta p^{-s})\in \mathbb{Z}_p^\times. \label{eqn: lambda,mu}
  \end{split}\right.
\end{gather}
Hence, by $\lambda\in\mathbb{Z}_p^\times$, we have
\begin{equation}\label{eqn: b}
  b \leq r, \quad  a+f\leq \max(r,s).
\end{equation}
Then there is an element $x_{a,b,c}^{d,f,\alpha,\beta,\gamma}\in X(cw_8)$ 
for which
\begin{equation}\label{eqn: u'(x)}
  u'(x_{a,b,c}^{d,f,\alpha,\beta,\gamma}) = \left(\begin{smallmatrix}
                                              1 & p^{-a} & \alpha p^{-d} & \beta p^{-s} \\
                                               & 1 & p^{-b} & \gamma p^{-f} \\
                                               &  & 1 & p^{-c} \\
                                               &  &  & 1
                                            \end{smallmatrix}\right) \pmod{U(\mathbb{Z}_p)}.
\end{equation}

Indeed, we have the matrix identity
\begin{equation}\label{eqn: MI for w_8}
  \left(\begin{smallmatrix}
    \mu^{-1} &&& \\
    z\lambda^{-1}& \mu\lambda^{-1} && \\
    x\beta^{-1} & y\beta^{-1} & \lambda\beta^{-1} & \\
    p^s & p^{s-a} & \alpha p^{s-d} & \beta
  \end{smallmatrix}\right) =
  \left(\begin{smallmatrix} 1&u_1&{U}_4&u_5 \\ &1&u_2&u_6 \\  &&1&u_3 \\ &&&1  \end{smallmatrix}\right) cw_8
  \left(\begin{smallmatrix} 1&p^{-a}&\alpha p^{-d}&\beta p^{-s} \\ &1&p^{-b}&\gamma p^{-f} \\  &&1& p^{-c} \\ &&&1  \end{smallmatrix}\right),
\end{equation}
where
\begin{equation}\label{eqn: m,n,l}
  \begin{array}{lll}
    u_1 = \mu^{-1}p^{r-t}(p^{-a-b}-\alpha p^{-d}),  & {U}_4 = -\mu^{-1} p^{s-r-a}, \\
    u_2 = \lambda^{-1}p^{t-r} (\alpha p^{s-c-d}-\beta),  & u_5 = \mu^{-1} p^{-s}, \\
    u_3 = -\beta^{-1}\gamma p^{r-s-f}, &  u_6 = \lambda^{-1}p^{t-s} (\gamma p^{-f} - p^{-b-c}).
  \end{array}
\end{equation}
Write \eqref{eqn: MI for w_8} as $g=ucw_8u'$. Then we have
\[
  g^\iota = w_8\!\ ^t g^{-1} w_8
  = \left(\begin{smallmatrix}
       \beta^{-1} &  &  &  \\
       -\alpha\lambda^{-1}p^{s-d} & \beta\lambda^{-1} &  &  \\
       (\alpha p^{r-d}-p^{r-a-b})\mu^{-1} & -y\mu^{-1} & \lambda\mu^{-1} &  \\
       p^t  & -p^{t-c} & -z & \mu
    \end{smallmatrix}\right),
\]
and its Bruhat decomposition is
%{\tiny
%\[
%  \left(\begin{smallmatrix}
%    1&-u_1&*&* \\ &1&-u_2&* \\  &&1&-u_3 \\ &&&1
%  \end{smallmatrix}\right)
%  \left(\begin{smallmatrix} p^{-s}&&&\\ &p^{s-r}&& \\ &&p^{r-t}& \\ &&&p^{t} \end{smallmatrix}\right) w
%  \left(\begin{smallmatrix}
%    1&-p^{-c}&p^{-b-c}-\gamma p^{-f}&\mu p^{-t} \\
%    &1&-p^{-b}&p^{-a-b}-\alpha p^{-d} \\  &&1&-p^{-a} \\ &&&1
%  \end{smallmatrix}\right).
%\]
%}
\[
  \left(\begin{smallmatrix}
    1&-u_1&*&* \\ &1&-u_2&* \\  &&1&-u_3 \\ &&&1
  \end{smallmatrix}\right)
  \left(\begin{smallmatrix} p^{-s}&&&\\ &p^{s-r}&& \\ &&p^{r-t}& \\ &&&p^{t} \end{smallmatrix}\right) w_8
  \left(\begin{smallmatrix}
    1&-p^{-c}&*&*\\ &1&-p^{-b}&* \\ &&1&-p^{-a} \\ &&&1
  \end{smallmatrix}\right).
\]
Since $g\in X(cw_8)$, we have $g^\iota\in X((cw_8)^\iota)\subseteq G(\mathbb{Z}_p)$, hence
\begin{equation}\label{eqn: a,b,c}
  c\leq t, \quad a+b\leq \max(r,d),\quad  \alpha p^{r-d}-p^{r-a-b}\in\mathbb{Z}_p.
\end{equation}

Let $\psi,\psi'$ be characters of $U(\mathbb{Q}_p)/U(\mathbb{Z}_p)$.
For $a,b,c$, and $d,f,\alpha,\beta,\gamma$ satisfying
\eqref{eqn: a,b,c,d,f}--\eqref{eqn: b} and \eqref{eqn: a,b,c}, let 
$
  X_{a,b,c}^{d,f,\alpha,\beta,\gamma}(cw_8)
              = T(\mathbb{Z}_p)*x_{a,b,c}^{d,f,\alpha,\beta,\gamma}
$
be the orbit through $x_{a,b,c}^{d,f,\alpha,\beta,\gamma}$, and let
$$
  S_{a,b,c}^{d,f,\alpha,\beta,\gamma}(\psi,\psi';c,w_8)
  = \sum_{x\in X_{a,b,c}^{d,f,\alpha,\beta,\gamma}(cw_8)}\psi(u(x))\psi'(u'(x))
$$
be the Kloosterman sum restricted to this orbit.
For $a,b,c$ satisfying \eqref{eqn: a,b,c,d,f}, \eqref{eqn: b} and \eqref{eqn: a,b,c}, let
$
  X_{a,b,c}(cw_8) = \bigcup\limits_{d,f,\alpha,\beta,\gamma} X_{a,b,c}^{d,f,\alpha,\beta,\gamma}(cw_8),
$
where $d,f$ run over nonnegative integers, and $\alpha,\beta,\gamma$ run over the elements of $\mathbb{Z}_p^\times$
satisfying \eqref{eqn: a,b,c,d,f}-\eqref{eqn: lambda,mu}. Let
$$
S_{a,b,c}(\psi,\psi';c,w_8) = \sum_{x\in X_{a,b,c}(cw_8)}\psi(u(x))\psi'(u'(x)).
$$

\begin{lem}\label{lemma: X(cw_8)}
  We have $X(cw_8) = \coprod_{a,b,c} X_{a,b,c}(cw_8)$, where $a,b,c\geq0$ run over integers
  satisfying \eqref{eqn: a,b,c,d,f}, \eqref{eqn: b}, and \eqref{eqn: a,b,c}.
\end{lem}

\begin{proof}
  See Stevens \cite[Lemmas 5.2 and 5.7]{stevens1987poincare}.
  % Note that $G(\mathbb{Z}_p)= U(\mathbb{Z}_p)\cdot B^-(\mathbb{Z}_p)\cdot U(\mathbb{Z}_p)$ is true for $G=\GL_n$.
\end{proof}

\begin{lem}\label{lemma:S<<w8}
  Let $\ell = \max(s,r,t)$, %$\sigma=\min(s,t)$
  and $a\leq s,\ b\leq r,\ c\leq t$ be nonnegative integers.
  Then
  \[
%    \begin{split}
%      |S_{a,b,c}(\psi,\psi';c,w)| \leq C & (|\nu_1\nu_3'|_p^{-1},p^\ell)^{1/2} (|\nu_2\nu_2'|_p^{-1},p^\ell)^{1/2} (|\nu_3\nu_1'|_p^{-1},p^\ell)^{1/2} \\
%            & \cdot (\sigma+1)(r+1) \cdot  p^{s+r+t-\ell+\frac{a+b+c}{2}}. %at least for the case $t\geq r$.
%    \end{split}
    \begin{split}
      & |S_{a,b,c}(\psi,\psi';c,w_8)|  \\
      & \hskip 10pt \leq
        64   (|\nu_1\nu_3'|_p^{-1},p^\ell)^{1/2} (|\nu_2\nu_2'|_p^{-1},p^\ell)^{1/2} (|\nu_3\nu_1'|_p^{-1},p^\ell)^{1/2}
          p^{-\frac{a+b+c}{2}}   \#(X_{a,b,c}(cw_8)).
    \end{split}
  \]
\end{lem}

\begin{proof}
  The involution $\iota$ sends $X_{a,b,c}(cw_8)$ to $X_{c,b,a}((cw_8)^\iota)$.
  Composing $\psi$ and $\psi'$ with $\iota$ has the effect of replacing
  $(\nu_1,\nu_2,\nu_3)$ by $(-\nu_3,-\nu_2,-\nu_1)$ and
  $(\nu_1',\nu_2',\nu_3')$ by $(-\nu_3',-\nu_2',-\nu_1')$.
  Applying $\iota$ to $cw_8$ reverses the roles of $t$ and $s$.
  Thus we may assume $t\geq s$ without loss of generality.

  Let $\ell=\max(r,t)$. The conditions \eqref{eqn: a,b,c,d,f}-\eqref{eqn: b}
  and \eqref{eqn: a,b,c} imply that the matrix entries of $u(x)$ and $u'(x)$
  lie in $p^{-\ell}\mathbb{Z}_p/\mathbb{Z}_p$ for every $x\in X(cw_8)$.
  Indeed, by Lemma \ref{lemma: X(cw_8)}, it is enough to verify this for $x=x_{a,b,c}^{d,f,\alpha,\beta,\gamma}$.
  Note that
  $\mu=p^{-s}  p^t(\alpha p^{s-d-c} - \beta) + p^{-a}  p^{t}(\gamma p^{-f} - p^{-b-c})
  =p^{-s}\lambda u_2 p^{r}+p^{-a}z \in \mathbb{Z}_p^\times$.
  We have $u_2\in p^{-r}\mathbb{Z}_p$. The claim is now easily verified.
  %by \eqref{eqn: a,b,c,d,f}-\eqref{eqn: lambda,mu} and \eqref{eqn: a,b,c}.

  Now let $\mathcal{S}$ be a finite subset of $\mathbb{Z}_{\geq0}^2\times 
(\mathbb{Z}_p^\times)^3$
  such that $X_{a,b,c}(cw_8)$ is the disjoint union of the
  $X_{a,b,c}^{d,f,\alpha,\beta,\gamma}(cw_8)$ with $(d,f,\alpha,\beta,\gamma)\in\mathcal{S}$.
  Then as in \cite[Th. 4.10]{stevens1987poincare} %Lemma \ref{lemma: stevens}
  we have
  \begin{equation}\label{eqn: S decomp}
    S_{a,b,c}(\psi,\psi';c,w_8)
    = p^{-3\ell}(1-p^{-1})^{-3} \sum_{(d,f,\alpha,\beta,\gamma)\in\mathcal{S}}
    \#(X_{a,b,c}^{d,f,\alpha,\beta,\gamma}(cw))
    S_{w_8}(\theta_{a,b,c}^{d,f,\alpha,\beta,\gamma};\ell),
  \end{equation}
  where $S_{w_8}$ is defined in \cite[Def. 4.9]{stevens1987poincare}, %in Definition \ref{defn: stevens},
  and $\theta_{a,b,c}^{d,f,\alpha,\beta,\gamma}: A_{w_8}(\ell)\rightarrow\mathbb{C}^\times$
  is the character given by
  \[
    \begin{split}
      \theta_{a,b,c}^{d,f,\alpha,\beta,\gamma}(\underline{\lambda}\times\underline{\lambda}')
      &  = e\left(\nu_1u_1\lambda_1+\nu_2u_2\lambda_2+\nu_3u_3\lambda_3+\nu_1'p^{-a}\lambda_1'+\nu_2'p^{-b}\lambda_2'+\nu_3'p^{-c}\lambda_3'\right) \\
      & = e\left(\frac{(-\nu_1\mu^{-1}p^{\ell+r-t}(p^{-a-b}-\alpha p^{-d}))\lambda_1 + (\nu_2\lambda^{-1}p^{\ell+t-r} (\alpha p^{s-c-d}-\beta))\lambda_2}{p^\ell} \right.\\
      & \hskip 60pt \left. +\frac{(\nu_3\beta^{-1}\gamma p^{\ell+r-s-f})\lambda_3 +\nu_1'p^{\ell-a}\lambda_1' +\nu_2'p^{\ell-b}\lambda_2' +\nu_3'p^{\ell-c}\lambda_3'}{p^\ell}\right).
    \end{split}
  \]
%  \[
%    \begin{split}
%     & \quad \theta_{a,b,c}^{d,f,\alpha,\beta,\gamma}(\underline{\lambda}\times\underline{\lambda}')  = e\left(\nu_1u_1\lambda_1+\nu_2u_2\lambda_2+\nu_3u_3\lambda_3+\nu_1'p^{-c}\lambda_1'+\nu_2'p^{-b}\lambda_2'+\nu_3'p^{-a}\lambda_3'\right) \\
%      & = e\left(\frac{(-\nu_1\beta^{-1}\gamma p^{t+r-s-f})\lambda_1 + (\nu_2\lambda^{-1}p^{2t-r} (\alpha p^{s-c-d}-\beta))\lambda_2}{p^t} \right.\\
%      & \qquad  \left. +\frac{(\nu_3\mu^{-1}p^{r}(p^{-a-b}-\alpha p^{-d}))\lambda_3 +\nu_1'p^{t-c}\lambda_1' +\nu_2'p^{t-b}\lambda_2' +\nu_3'p^{t-a}\lambda_3'}{p^t}\right),
%    \end{split}
%  \]
%  if $\ell=t$, and
%  \[
%    \begin{split}
%     & \quad \theta_{a,b,c}^{d,f,\alpha,\beta,\gamma}(\underline{\lambda}\times\underline{\lambda}') \\
%      & = e\left(\frac{(-\nu_1\beta^{-1}\gamma p^{2r-s-f})\lambda_1 + (\nu_2\lambda^{-1}p^{t} (\alpha p^{s-c-d}-\beta))\lambda_2}{p^r} \right.\\
%      & \qquad  \left. +\frac{(\nu_3\mu^{-1}p^{2r-t}(p^{-a-b}-\alpha p^{-d}))\lambda_3 +\nu_1'p^{r-c}\lambda_1' +\nu_2'p^{r-b}\lambda_2' +\nu_3'p^{r-a}\lambda_3'}{p^r}\right),
%    \end{split}
%  \]
%  if $\ell=r$.
  By Example 4.12 in Stevens \cite{stevens1987poincare}, we have
  \begin{equation}\label{eqn: S decomp to S_2}
    \begin{split}
      & S_{w_8}(\theta_{a,b,c}^{d,f,\alpha,\beta,\gamma};\ell) = S_2(\nu_1\mu^{-1}p^{\ell+r-t}(p^{-a-b}-\alpha p^{-d}),\nu_3'p^{\ell-c};p^\ell) \\
      & \hskip 40pt  \cdot S_2(\nu_2\lambda^{-1}p^{\ell+t-r} (\alpha p^{s-c-d}-\beta),\nu_2'p^{\ell-b};p^\ell)\cdot
                     S_2(-\nu_3\beta^{-1}\gamma p^{\ell+r-s-f},\nu_1'p^{\ell-a};p^\ell),
    \end{split}
  \end{equation}
  where $S_2$ is the classical $\GL(2)$-Kloosterman sum.
  By Weil \cite{weil1948some}, we have the inequality
  \begin{equation}\label{eqn: S_2}
    |S_2(m,n;p^{\ell})| \leq 2 (\gcd(|m|_p^{-1},|n|_p^{-1},p^\ell))^{1/2} 
p^{\ell/2},
  \end{equation}
  for $m,n\in\mathbb{Z}_p$. In order to apply this bound, we first note
  \[
    \begin{split}
       \gcd(|\nu_3 p^{\ell+r-s-f}|_p^{-1},|\nu_1'p^{\ell-a}|_p^{-1},p^\ell) & \leq \gcd(|\nu_3\nu_1'|_p^{-1},p^\ell)p^{\ell-a}, \\
       \gcd(|\nu_2 p^{\ell+t-r} (\alpha p^{s-c-d}-\beta)|_p^{-1},|\nu_2'p^{\ell-b}|_p^{-1},p^\ell) & \leq \gcd(|\nu_2\nu_2'|_p^{-1},p^\ell)p^{\ell-b}, \\
       \gcd(|\nu_1 p^{\ell+r-t}(p^{-a-b}-\alpha p^{-d})|_p^{-1},|\nu_3'p^{\ell-c}|_p^{-1},p^\ell) & \leq \gcd(|\nu_1\nu_3'|_p^{-1},p^\ell)p^{\ell-c}.
    \end{split}
  \]
  Hence we have
  \begin{equation*}\label{eqn: S_w_8 bound}
    |S_{w_8}(\theta_{a,b,c}^{d,f,\alpha,\beta,\gamma};\ell)|
    \leq 8 (|\nu_1\nu_3'|_p^{-1},p^\ell)^{1/2} (|\nu_2\nu_2'|_p^{-1},p^\ell)^{1/2}
            (|\nu_3\nu_1'|_p^{-1},p^\ell)^{1/2}   p^{3\ell-\frac{a+b+c}{2}}.
  \end{equation*}
  This inequality, together with \eqref{eqn: S decomp}, gives
  \begin{equation}\label{eqn: S_a,b,c}
    \begin{split}
      |S_{a,b,c}(\psi,\psi';c,w_8)| \leq 8 & (|\nu_1\nu_3'|_p^{-1},p^\ell)^{1/2}
      (|\nu_2\nu_2'|_p^{-1},p^\ell)^{1/2} (|\nu_3\nu_1'|_p^{-1},p^\ell)^{1/2} \\
      & \cdot (1-p^{-1})^{-3} p^{-\frac{a+b+c}{2}}  \sum_{(d,f,\alpha,\beta,\gamma)\in\mathcal{S}}  \#(X_{a,b,c}^{d,f,\alpha,\beta,\gamma}(cw_8)).
    \end{split}
  \end{equation}
  The sum appearing on the right hand side is equal to $\#(X_{a,b,c}(cw_8))$.
  Since $p\geq 2$ we have $(1-p^{-1})^{-3}\leq 8$, by \eqref{eqn: S_a,b,c}, we prove the lemma.
\end{proof}

\begin{proof}[Proof of Theorem \ref{thm: w_8}]
  By the involution $\iota$, we can assume $t\geq s$ without loss of generality.
  Let
  $$
    C=64   (|\nu_1\nu_3'|_p^{-1},p^\ell)^{1/2} (|\nu_2\nu_2'|_p^{-1},p^\ell)^{1/2} (|\nu_3\nu_1'|_p^{-1},p^\ell)^{1/2} (r+1)(s+1).
  $$
  At first, we deal with the case $t\geq r$.
  \begin{itemize}
  \item If $a+b+c\leq t$ and $d+f\leq r$, then $\#(d,f)\leq (s+1)(r+1)$, $\#(\alpha,\gamma,\beta)\leq p^{d+s+f}$, so
  \[
    \#(X_{a,b,c}(cw_8))\leq (r+1)(s+1)p^{a+b+c+d+f+s}\leq (r+1)(s+1)p^{r+s+a+b+c}.
  \]
  Hence by Lemma \ref{lemma:S<<w8}, we have
  $
     |S_{a,b,c}(\psi,\psi';c,w_8)|   \leq  C  p^{r+s+t/2}.
  $
  \item If $a+b+c\leq t$ and $d+f>r$, then we assume that $d+f=r+k$, $k\geq1$.
      Note that $d\leq s, f\leq r$, we have $k\leq s$.
      By \eqref{eqn: lambda,mu}, we have $b+s=d+f=r+k$.
      Since $\lambda\in\mathbb{Z}_p^\times$, we have $\#\{(\alpha,\gamma,\beta)\}\leq p^{d+f+(s-k)}=p^{r+s}$.
      Hence
      \[
        |S_{a,b,c}(\psi,\psi';c,w_8)|
        \leq  C  p^{r+s+\frac{a+b+c}{2}}
        \leq C   p^{r+s+t/2}.
      \]
  \item If $a+b+c>t$ and $d+f\leq r$, then by \eqref{eqn: lambda,mu} and a similar argument as above, we have $\#\{(\alpha,\gamma,\beta)\}\leq p^{(d-m)+f+s}$.
  Hence
  \[
    \begin{split}
      |S_{a,b,c}(\psi,\psi';c,w_8)|  & \leq  C  p^{d-m+f+s+\frac{a+b+c}{2}}
      \leq C   p^{r+s+t/2}.
    \end{split}
  \]
  \item If $a+b+c>t$ and $d+f>r$, then  we have $\#\{(\alpha,\gamma,\beta)\}\leq p^{(d-m)+f+(s-k)}$.
      Hence
      \[
        \begin{split}
       |S_{a,b,c}(\psi,\psi';c,w_8)| & \leq C  p^{d-m+f+s-k+\frac{a+b+c}{2}}
        \leq C  p^{r+s+t/2}.
        \end{split}
      \]
  \end{itemize}
  Note that in this case, we always have $r+s+t/2 \leq t+2s+r/2$.  Theorem \ref{thm: w_8} now follows from the equality
  $Kl_p(\psi,\psi';c,w_8)=\sum\limits_{a\leq s,b\leq r,c\leq t}S_{a,b,c}(\psi,\psi';c,w_8).$

  Now we handle the case $r>t$.
  By a similar argument as above, we obtain
  \[
    |S_{a,b,c}(\psi,\psi';c,w_8)| \leq C  p^{r+s+t/2}.
  \]
  If $t$ is small, this bound is not good enough.
  So we need to bound this in other way.
  \begin{itemize}
   \item If $f>t$, then by \eqref{eqn: y,z}, we have $b+c=f$, and $a+f\leq r$. By \eqref{eqn: lambda,mu},
       we have $\#(\alpha,\gamma)\leq p^{d+f-(a+f-t)}$.
       If $d+f\leq r$, then we have
       \[
        \begin{split}
         |S_{a,b,c}(\psi,\psi';c,w_8)| &
         \leq C p^{d+f-(a+f-t)+s+\frac{a+b+c}{2}}
          \leq C p^{t+s+d+\frac{b+c}{2}}
          \leq C p^{t+s+\frac{d}{2}+\frac{d+f}{2}}
         \leq C p^{t+3s/2+r/2}.
        \end{split}
        \]
    \item If $f>t$ and $d+f>r$, then by writing $d+f=r+k$, $1\leq k\leq 
s$, we have
        \[
            \begin{split}
            |S_{a,b,c}(\psi,\psi';c,w_8)| & \leq C  p^{d+f-(a+f-t)+s-k+\frac{a+b+c}{2}} \leq C  p^{t+s+\frac{d}{2}+\frac{d+f}{2}-k}
         \leq C  p^{t+3s/2+r/2}.
            \end{split}
        \]
  \item  If $f\leq t$ and $a+b+c>r$.   Since $\mu\in \mathbb{Z}_p^\times$, we have $\#(\alpha,\gamma)\leq p^{d+f-(a+b+c-t)}$.   Then by the same argument on the size of $d+f$, we have
     \[
         |S_{a,b,c}(\psi,\psi';c,w_8)|
         \leq C   p^{r+t+s-\frac{a+b+c}{2}}
         \leq C   p^{t+s+r/2}.
     \]
  \item If $f\leq t$, and $a+b+c\leq r$, then we have
    \[
    \begin{split}
      |S_{a,b,c}(\psi,\psi';c,w_8)| & \leq C  p^{d+f+s+\frac{a+b+c}{2}}  \leq C  p^{t+2s+r/2}.
    \end{split}
    \]
  \end{itemize}
  This proves \eqref{eqn:Kl-w8}.

  We now give a proof of the second claim.
  If $r+\sigma+\varrho/2\leq \varrho+2\sigma+r/2$, i.e., $r\leq \varrho+2\sigma$, then $\sigma+r\leq 4\varrho$, so $r+\sigma+\varrho/2 \leq 9(\varrho+r+\sigma)/10$. If $r+\sigma+\varrho/2 > \varrho+2\sigma+r/2$, i.e., $r>\varrho+2\sigma$, then $9\sigma< 3r$ and $\varrho+11\sigma < 4r$, so $\varrho+2\sigma+r/2 < 9(\varrho+r+\sigma)/10$. This proves $\min(p^{r+\sigma+\varrho/2},p^{\varrho+2\sigma+r/2})\leq p^{9(t+r+s)/10}$, as claimed, and hence Theorem \ref{thm: w_8}.
\end{proof}

\begin{rmrk}
  The result is not optimal.
  To improve the bound in some cases, one may use the stationary phase formulas as Dabrowski and Fisher did for $\GL(3)$, see \cite{dabrowski1997stationary}.
\end{rmrk}

\begin{rmrk}
  Stevens' method can be used to bound other Kloosterman sums as well. It's not too hard to prove bounds similar to \eqref{eqn:trivialbound}. But to improve these ``trivial'' bounds, one may need new ideas.
\end{rmrk}

{\sc Acknowledgements.} The author would like to thank Professor Dorian Goldfeld for suggesting this question, and for his valuable advice and constant encouragement.

\vskip 25pt
%He also wants to thank the referees and editors for their kind comments and valuable suggestions.

% \bibliographystyle{ieeetr}
\bibliographystyle{amsalpha}
%\nocite{*}
\bibliography{biblio}

\end{document}